\documentclass[a4paper]{article}
\title{The differential $\rm d_4(h_6^3)$ in the Adams spectral sequence for spheres}
\date{}
\author{Pomin Wu}

\usepackage[leqno]{amsmath}
\usepackage[lite]{amsrefs}
\usepackage{amssymb}
\usepackage{enumerate}    \renewcommand{\labelenumi}{(\arabic{enumi})}
                          \renewcommand{\labelenumii}{(\alph{enumii})}
                          \renewcommand{\labelenumiii}{(\roman{enumiii})}
\usepackage{hyperref}
\usepackage{sseq}         
\usepackage[all]{xy}
\usepackage{array}        
\usepackage{multicol}     
\usepackage{paralist}     
\usepackage{semantic}     

\usepackage[amsmath,amsthm,thmmarks,hyperref]{ntheorem}
\usepackage{changepage}   
\usepackage{aliascnt}     
\numberwithin{equation}{section}    

\theoremstyle{definition}
\newaliascnt{defn}{equation}
\newtheorem{defn}[defn]{Definition}
\aliascntresetthe{defn}

\theoremstyle{plain}

\newaliascnt{thm}{equation}
\newtheorem{thm}[thm]{Theorem}
\aliascntresetthe{thm}

\newaliascnt{lemm}{equation}
\newtheorem{lemm}[lemm]{Lemma}
\aliascntresetthe{lemm}

\newaliascnt{prop}{equation}
\newtheorem{prop}[prop]{Proposition}
\aliascntresetthe{prop}

\newaliascnt{coro}{equation}
\newtheorem{coro}[coro]{Corollary}
\aliascntresetthe{coro}

\newaliascnt{algo}{equation}

\aliascntresetthe{algo}

\newaliascnt{example}{equation}

\aliascntresetthe{example}

\newaliascnt{remark}{equation}
\newtheorem{remark}[remark]{Remark}
\aliascntresetthe{remark}

\makeatletter
\newtheoremstyle{marginonlynumber}%
{\item[\theorem@headerfont \llap{(##2)\theorem@separator}]}%
{\item[\theorem@headerfont \llap{##2} (##3)\theorem@separator]}%
\makeatother

\newlength{\notelabel}
\theoremstyle{marginonlynumber}
\theorembodyfont{}
\theoremseparator{}
\theoremprework{%
  \settowidth{\notelabel}{\textbf{(\thenote) }}
  \vskip\theorempreskipamount
  \begin{adjustwidth}{\notelabel}{}}
\theorempostwork{%
  \end{adjustwidth}%
  \vskip\theorempostskipamount}
\newaliascnt{note}{equation}
\newtheorem{note}[note]{Note}
\aliascntresetthe{note}

\newcommand{\U}{\Union}
\newcommand{\Union}{\cup}

\newcommand{\given}{\forall}

\renewcommand{\to}{\xrightarrow}

\renewcommand{\bar}{\overline}
\newcommand{\Sum}{\Summation}
\newcommand{\Summation}{\Sigma}

\newcommand{\R}{\Real}
\newcommand{\Real}{\mathbb{R}}

\newcommand{\Z}{\Integer}
\newcommand{\Integer}{\mathbb{Z}}

\newcommand{\N}{\NaturalNumber}
\newcommand{\NaturalNumber}{\mathbb{N}}

\newcommand{\dsum}{\DirectSum}
\newcommand{\DirectSum}{\oplus}

\newcommand{\tensor}{\TensorProduct}
\newcommand{\TensorProduct}{\otimes}

\newcommand{\Susp}{\Suspension}
\newcommand{\Suspension}{\Sigma}
\newcommand{\Loop}{\Omega}

\newcommand{\Smash}{\wedge}
\newcommand{\V}{\Wedge}
\newcommand{\Wedge}{\vee}

\newcommand{\homo}{\HomotopicTo}
\newcommand{\HomotopicTo}{\simeq}
\newcommand{\iso}{\IsomorphicTo}
\newcommand{\IsomorphicTo}{\cong}

\newcommand{\coker}{{\rm coker}}

\newcommand{\Ext}{{\rm Ext}}

\mathlig{|->}{\mapsto}
\mathlig{|-->}{\longmapsto}
\mathlig{<|}{\NormalSubgroup}

\newcommand{\NormalSubgroup}{\vartriangleleft}

\newcommand{\A}{\Steenrod}
\newcommand{\Steenrod}{\mathcal{A}}
\newcommand{\toda}[1]{\left<#1\right>}

\newcommand{\ssdropcirc}{\ssdrop{\circ}}
\newcommand{\ssdropast}{\ssdrop{\ast}}

\begin{document}

\thispagestyle{empty}
\begin{titlepage}
  \begin{center}
    \huge The differential $\rm d_4(h_6^3)$ in the Adams spectral sequence for spheres\\
    ~\\
    \Large By\\
    Pomin Wu\\
    ~\\
    Adviser\\
    Dr.\ W\^en-Hsiung Lin\\
    \vspace*{2in}
    \normalsize
    A Doctorial Dissertation\\
    Presented to\\
    Department of Mathematics\\
    at\\
    National Tsing-Hua University\\ 
    \vspace*{1in}
    In Partial Fulfillment\\
    of the Requirement for the Degree of\\
    Doctoral of Philosophy\\
    \vfill
    Sep, 2011
  \end{center}
\end{titlepage}

\pagenumbering{roman}

\begin{center}
  \textbf{Acknowledgements}
\end{center}
{\small
I would like to thank my advisor, Dr.\ W\^en-Hsiung Lin, for his infinite support and patience.
It may seem impossible for me to fully comprehend
every technique there is to learn from my teacher.
Yet it is among the non-technical things
that I believe I have learned the most
valuable lesson. The image
of a mathematician
which I have seen from Dr.\ Lin
will always
aspire me to pursuit a deeper understanding of the subject.
Therefore I feel I am most fortunate to have studied under his direction.
I would also like to thank Mr.\ Tai-Wei Chen for his help on some of the most complicated computations in this thesis.

I am also greatly in debt to many friends
who in the last few years
have supported me in one way or another.
You are my connections to this world, and
it is only by these connections that
a pursuit in abstract knowledge
acquires its full meaning.
Thank you for being here with me.

I must also thank my parents and my sister,
for all the support a family can give
to one of its members.
You are the closest thing
to a definition of me
in this world.
}
\newpage

\tableofcontents
\newpage


\pagenumbering{arabic}

\section{Introduction}\label{se:intro}
Let $\A$ denote the mod $2$ Steenrod algebra. Let $X$ be a stable complex
in the stable category (\cite{adams_non-existence_1960}).
Let $H_\ast(X)$ and $H^\ast(X)$ denote respectively the
mod $2$ reduced homology and cohomology of $X$.
The $\Ext$ groups $\Ext_\A^{\ast,\ast}( H^\ast(X),\Z/2)$,
which will be simply denoted by $\Ext_\A^{\ast,\ast}(X)$,
form the $E_2$ term of the mod $2$ Adams spectral sequence
(to be abbreviated as ASS) for
computing the $2$-primary stable homotopy groups $\pi_\ast^S(X)$ of $X$
(\cite{adams_structure_1958,adams_non-existence_1960}).
In the case $X=S^0$,
the sphere spectrum in stable dimension zero,
the $\Ext$ groups
$\Ext_\A^{\ast,\ast}( H^\ast(S^0),\Z/2)=\Ext_\A^{\ast,\ast}(\Z/2,\Z/2)$
will be simply denoted by $\Ext_\A^{\ast,\ast}$, and
the $2$-primary stable homotopy groups $\pi_\ast^S(S^0)$
by $\pi_\ast^S$.

For each $i\ge0$, let
$h_i\in \Ext_\A^{1,2^i} =\Ext_\A^{1,2^i}(\Z/2,\Z/2)$ be the class
corresponding to the generator $Sq^{2^i}\in\A$.
It is also known that $h_i^2\ne0$, $h_i^3\ne0$ for all $i\ge0$ (\cite{adams_non-existence_1960}).
The classes $h_0$, $h_1$, $h_2$, and $h_3$ are known to
detect homotopy elements
$2\iota\in\pi_0^S$, $\eta\in\pi_1^S$, $\nu\in\pi_3^S$, and $\sigma\in\pi_7^S$, where
$\eta$, $\nu$, and $\sigma$ are the three Hopf classes.
So $h_0^2$, $h_1^2$, $h_2^2$, and $h_3^2$ detect respectively the homotopy elements
$4\iota=(2\iota)^2$, $\eta^2$, $\nu^2$, and $\sigma^2$, and
$h_0^3$, $h_1^3$, $h_2^3$, and $h_3^3$ detect respectively the homotopy elements
$8\iota=(2\iota)^3$, $\eta^3$, $\nu^3$, and $\sigma^3$.
(All of these homotopy elements are non-zero in $\pi^S_\ast$
(\cite{mimura_n_20-th_1963,mimura_generalized_1965}))
J.\ F.\ Adams also proves in \cite{adams_non-existence_1960} that
for each $i>3$,
$h_i$ does not detect homotopy elements in $\pi_\ast^S$
by showing that 
\begin{align}
  \label{eq:intro:d2-hi}
  d_2(h_i) &=h_{i-1}^2h_0\ne0
  \quad\mbox{for $i\ge4$ in the ASS for $\pi_\ast^S$.}
\end{align}
So it is not immediate whether $h_i^2$ detect homotopy elements for $i\ge4$.
Any homotopy element in $\pi_{2^{i+1}-2}^S$ for $i\ge1$ which is
detected by $h_i^2$ is denoted by $\theta_i$.
$\theta_i$, if exists, is known as a Kervaire invariant element
\cite{kervaire_manifold_1960,browder_kervaire_1969}.
Thus, $\eta^2=\theta_1$,
$\nu^2=\theta_2$, and $\sigma^2=\theta_3$.
Mahowald and Tangora prove in \cite{mahowald_differentials_1967} that
$h_4^2$ does detect homotopy elements in $\pi_{30}^S$. That is,
there exists a $\theta_4$.
In fact they have shown that $\pi_{30}^S=\Z/2(\theta_4)$.
So $2\theta_4=0$.
Barratt, Jones, and Mahowald prove in \cite{barratt_differentials_1970} that
$h_5^2$ also detects homotopy elements in $\pi_{62}^S$.
Kochman and Mahowald have shown in \cite{kochman_computation_1995} that
there actually exists a $\theta_5$ with $2\theta_5=0$.
Recently, M.\ Hill, M. Hopkins, and D. Ravenel have shown that
$h_i^2$ does not detect homotopy elements for $i\ge7$. The status of $h_6^2$
is still open.

For the family $\{h_i^3\}_{i\ge0}$,
$h_0^3$, $h_1^3$, $h_2^3$, and $h_3^3$ are known to detect homotopy elements
as already mentioned above.
Barratt, Mahowald, and Tangora have shown in \cite{barratt_differentials_1970}
that $h_4^3$ detects homotopy elements in $\pi_{45}^S$.
W.\ H.\ Lin have shown in \cite{lin_differential_1998} that
$h_7^3$ does not detect homotopy elements.
More precisely, Lin proves that $d_4(h_7^3)$ is essentially $h_0^3g_5\ne0$
in the ASS for spheres.
We refer to \cite{lin_differential_1998} for the details of this ``essentiality.''
Here we recall \cite{lin_differential_1998} that
for each $i\ge1$, there is a non-zero class $g_i\in\Ext_\A^{4,2^{i+3}+2^{i+2}}$, and
$\Ext_\A^{4,2^{i+3}+2^{i+2}}=\Z/2(g_i)$.
It is known (\cite{tangora_cohomology_1970}) that
$h_0^2g_1\ne0$, $h_0^2g_2\ne0$, and
$h_0^3g_1=0$, $h_0^3g_2=0$.
Lin has shown in \cite{lin_differential_1998} that
$h_0^3g_5\ne0$ as has just been mentioned above.
If one can show $h_0^3g_{i-2}\ne0$ for $i\ge8$, then
Lin's method in \cite{lin_differential_1998} actually also proves that
$d_4(h_i^3)$ is essentially $h_0^3g_{i-2}$ for all $i\ge8$.
He also remarks in \cite{lin_differential_1998} that his method
does not include the cases $i=5$ and $6$.
In \cite{kan_differential_2001}
Kan has shown that $h_0^3g_3\ne0$ and that $d_4(h_5^3)=h_0^3g_3\ne0$
in the ASS for spheres.

It is the purpose of this paper to show the following.
\begin{thm}\label{thm:main}
  In the mod $2$ Adams spectral sequence for spheres, ${d_4(h_6^3)=h_0^3g_4\ne0}$
  in $\Ext_\A^{7,195}$.
\end{thm}
This together with
\cite{mimura_n_20-th_1963,mimura_generalized_1965,barratt_differentials_1970,lin_differential_1998,kan_differential_2001}
mentioned above completely settle the differentials of the $h_i^3$ family.

The plan for this paper goes as follows.
We will prove the following $\Ext$ group result in \autoref{se:Ext}.
\begin{note}
  \label{eq:intro:h03g4-non-zero}
  $h_0^3g_4\ne0$ in $\Ext_\A^{7,195}$.
\end{note}
\noindent
Granting this, in \autoref{se:boundary} we will outline the proof of the following result.
\begin{thm}\label{thm:intro:boundary}
  $h_0^3g_4$ is a boundary in the ASS for spheres.
\end{thm}
In \autoref{se:Ext}
we will also
prove the following $\Ext$ group results.
\begin{prop}\label{thm:intro:Ext}
  \begin{enumerate}
    \item $\Ext_\A^{s,s+189}=0$ for $0\le s\le2$.
    \item $\Ext_\A^{s,s+189}=\Z/2(h_0^{s-3}h_6^3)$ for $3\le s\le4$.
    \item $\Ext_\A^{5,194}=\Z/2(h_0^2h_6^3)\dsum\Z/2(h_1g_4)$.
  \end{enumerate}
\end{prop}
In \autoref{se:pi} we will compute the
stable homotopy groups $\pi_{62}^S(P^{62}_{46})$ and $\pi_{62}^S(P^{62}_{47})$
($P^l_k$ are the stunted projective spaces)
which we need to complete the proof of \autoref{thm:intro:boundary}.
In \autoref{se:not-hit} we will complete the proof of
\autoref{thm:intro:boundary}, and
also prove the following results in the ASS for spheres.
\begin{prop}
  \label{thm:dr-h6}
  In the ASS for spheres, $d_r(h_6^3)=0\in\Ext_\A^{3+r,191+r}$ for $2\le r\le 3$.
\end{prop}
\begin{prop}
  \label{thm:d2-h1g4}
  In the ASS for spheres, $d_2(h_1g_4)\ne h_0^3g_4$.
\end{prop}
\begin{proof}[Proof of \autoref{thm:main}]
  By \autoref{thm:intro:boundary} $d_r(\alpha)=h_0^3g_4$ for some
  $\alpha\in\Ext_\A^{7-r,196-r}$, $2\le r\le7$. From
  \ref{thm:dr-h6}, \ref{thm:d2-h1g4}, and
  from the fact that $h_0$ is an infinite cycle in the ASS for spheres,
  it is easy to see that
  $d_2(h_0^2h_6^3) =0$,
  $d_2(h_1g_4) \ne h_0^3g_4$,
  $d_3(h_0h_6^3) =0$.
  So by \autoref{thm:intro:Ext},
  the only possibility is $d_4(h_6^3)=h_0^3g_4$.
  This completes the proof of \autoref{thm:main}.
\end{proof}

\newpage
\section{Outline of the proof of \autoref{thm:intro:boundary}}\label{se:boundary}
We recall again that
in the stable homotopy groups $\pi_\ast^S=\pi_\ast^S(S^0)$,
there exist homotopy elements
$\theta_4\in\pi_{30}^S=\Z/2$ and $\theta_5\in\pi_{62}^S$
detected respectively by $h_4^2$ and $h_5^2$ in the ASS for $\pi_\ast^S$,
with $2\theta_4=0$, $2\theta_5=0$
(\cite{mahowald_differentials_1967,barratt_differentials_1970,kochman_computation_1995}).
Throughout this paper we fix a $\theta_5$ with $2\theta_5=0$.
We also assume that there exists a $\theta_6 \in\pi_{126}^S$.
(Recall that it is still not known whether $\theta_6$ exists or not.)
Our proof of the main theorem is more conveniently presented
when $\theta_6$ is assumed to be existent.
The method of our proof can be modified to work without the knowledge of $\theta_6$, and
this will be given in \autoref{se:t6-not-exist} of the paper.
Again, in what follows, we fix a $\theta_6$.

A key step in Kan's proof of the result $d_4(h_5^3)=h_0^3g_3$ in his thesis \cite{kan_differential_2001} is
to show that $\theta_4\theta_5\in\pi_{92}^S$ is non-zero and
is detected by $h_0^2g_3\in\Ext_\A^{6,98}$.
Our proof of \autoref{thm:intro:boundary} is based on a similar result
which is \autoref{thm:t5t6-detected} below.
To state it, we
recall (\cite{lin_ext_2008}) that for each $i\ge0$,
there is a non-zero class $D_3(i)\in\Ext_\A^{4,2^{i+6}+2^i}$.
Consider the classes $h_0^2g_4$ and $h_5^2D_3(1)$ in $\Ext_\A^{6,194}$.
In \autoref{se:Ext} we will show the following $\Ext$ group result.
\begin{note}
  \label{eq:h_0^2g_4-ne-0}
  $h_0^2g_4\ne0$, $h_5^2D_3(1)\ne0$, and they are linearly independent
  in $\Ext_\A^{6,194}$.
\end{note}
\noindent
The main work of this section is to prove the following theorem.
\begin{thm}\label{thm:t5t6-detected}
  $\theta_5\theta_6$ in $\pi_{188}^S$
  is detected by $h_0^2g_4 +\varepsilon h_5^2D_3(1)\ne0$
  in $\Ext_\A^{6,194}$,
  where $\varepsilon=0\mbox{ or }1$.
\end{thm}
From this we will deduce \autoref{thm:intro:boundary} at the end of this section.
\begin{remark}
  \label{remark:h02g4-not-a-boundary}
  From \autoref{thm:intro:Ext} and \autoref{thm:dr-h6},
  we see $h_0^2g_4+\varepsilon h_5^2D_3(1)$ is not a boundary in the
  ASS for $\pi_\ast^S$.
\end{remark}

\autoref{thm:t5t6-detected} is a consequence of \autoref{thm:t5ht6-detected} later, which is
a result about the stable homotopy group $\pi_{188}^S(P)$ parallel to the result
about $\pi_{188}^S$ in \autoref{thm:t5t6-detected} above. Here
$P$ is the suspension spectrum of the infinite real projective space $P^\infty$.
We need to recall the
Kahn-Priddy theorem (\cite{kahn_applications_1972})
in order to state \autoref{thm:t5ht6-detected}, and also
to see what the ``parallel'' is.
In the next few paragraphs we do this recalling.

For any integer $n\ge0$, the real projective space $P^n$
is the space of all lines in $\R^{n+1}$ through the origin.
There is an inclusion $P^n\subset P^{n+1}$
induced by the natural inclusion $\R^{n+1}\subset\R^{n+2}$.
Regard $P^n$ as a subspace of $P^{n+1}$ with this inclusion for each $n\ge0$, and
let $P^\infty=\U_{n=0}^\infty P^n$.
The suspension spectrum of the space $P^n$ is also denoted by $P^n$, and
the suspension spectrum of the space $P^\infty$ is denoted by $P$
as mentioned above.

For any integer $n\ge0$ there is a space map $\Susp^{n+1}P^n\to{t_n}S^{n+1}$
which is called the \emph{$n$th transfer map}.
These transfer maps are compatible with respect to $n$,
inducing a stable map $P\to{t}S^0$, called the \emph{transfer map}, from the spectrum $P$ to the sphere spectrum $S^0$.
These are described in \eqref{eq:r_n-compatible} through \autoref{defn:t} below.

Given a point $L\in P^n$, which is a line through the origin of $\R^{n+1}$,
let $L^\perp$ be the $n$-dimensional hyperplane through the origin of $\R^{n+1}$ that is perpendicular to $L$,
and let $\R^{n+1}\to{r_L}\R^{n+1}$ be the map which is the reflection with respect to $L^\perp$.
Then $r_L$ is an element of $O(n+1)$, the group of orthogonal transformations
from $\R^{n+1}$ to itself.
Clearly the map $P^n\to{r_n}O(n+1)$ defined by $r_n(L)=r_L$ is continuous for each $n$.
Note that
\begin{align}
  r_n(L)&\in O(n+1)-SO(n+1)\quad  \given L\in P^n,
  \label{eq:rnL-component}
\end{align}
where $SO(n+1)\subset O(n+1)$ is the subgroup of orthogonal transformations
of determinant $1$.
These maps are compatible with respect to $n$ in the sense that the following diagram,
where the vertical maps are obvious inclusions,
is commutative:
\begin{align}
  \begin{split}
    \xymatrix{
    P^n\ar[r]^-{r_n}\ar[d]&   O(n+1)\ar[d]\\
    P^{n+1}\ar[r]^-{r_{n+1}}& O(n+2)
    }
  \end{split}
  \label{eq:r_n-compatible}
\end{align}
$S^{n+1}$ can be regarded as $\R^{n+1}\U\{\infty\}$, the one-point compatification of $\R^{n+1}$.
Then each element in $O(n+1)$ can be thought as a continuous function
from $S^{n+1}=\R^{n+1}\U\{\infty\}$ to itself, sending $\infty$ to $\infty$.
The $(n+1)$-fold loop space $\Loop^{n+1}S^{n+1}$ is the function space
of all the maps from $S^{n+1}$ to itself preserving the base point $\infty$.
Thus we get a natural map $O(n+1)\to{i_n}\Loop^{n+1}S^{n+1}$ for each $n$, and
these maps are compatible with respect to $n$ in the sense that
the following diagram, where the vertical maps are obvious inclusions,
is commutative:
\begin{align}
  \begin{split}
    \xymatrix{
    O(n)\ar[r]^-{i_{n-1}}\ar[d]&  \Loop^nS^n\ar[d]\\
    O(n+1)\ar[r]^-{i_n}&          \Loop^{n+1}S^{n+1}
    }
  \end{split}
  \label{eq:i_n-compatible}
\end{align}
Note that
\begin{note}
  $P^n\to{i_nr_n}\Loop^{n+1}S^{n+1}$
  maps $P^n$ to the ``$-1$'' component of $\Loop^{n+1}S^{n+1}$ by \eqref{eq:rnL-component},
  since each element of $O(n+1)-S(n+1)$ when considered as a self-map of $S^{n+1}$
  is a degree $-1$ map.
  \label{eq:inrn}
\end{note}
\noindent
From \eqref{eq:r_n-compatible} and \eqref{eq:i_n-compatible}
we have the following commutative diagram.
\begin{align*}
  \begin{split}
    \xymatrix{
    P^n\ar[r]^-{r_n}\ar[d]&   O(n+1)\ar[r]^-{i_n}\ar[d]& \Loop^{n+1}S^{n+1}\ar[d]\\
    P^{n+1}\ar[r]^-{r_{n+1}}& O(n+2)\ar[r]^-{i_{n+1}}&    \Loop^{n+2}S^{n+2}
    }
  \end{split}
\end{align*}
\begin{defn}\label{defn:t_n}
  The \emph{$n$th transfer map}
  $\Susp^{n+1}P^n\to{t_n}S^{n+1}$ is
  the $(n+1)$-fold adjoint map to $P^n\to{i_nr_n}\Loop^{n+1}S^{n+1}$.
\end{defn}
It is clear that $t_n$ is also compatible with respect to $n$ in the sense that
the following diagram, where the left vertical map is the obvious inclusion,
is commutative:
\begin{align*}
  \begin{split}
    \xymatrix{
    \Susp^{n+2}P^n\ar[r]^-{\Susp t_n}\ar[d]&  \Susp S^{n+1}\ar@{=}[d]\\
    \Susp^{n+2}P^{n+1}\ar[r]^-{t_{n+1}}&      S^{n+2}
    }
  \end{split}
\end{align*}
\begin{defn}\label{defn:t}
  The family of $n$th transfer maps
  $\{t_n\}$ induces a stable map
  $P\to{t}S^0$ from the suspension spectrum $P$
  to the sphere spectrum $S^0$.
  This map $t$ is called the \emph{transfer map}.
\end{defn}
Consider the induced map $\pi_\ast^S(P)\to{t_\#}\pi_\ast^S$ of $P\to{t}S^0$.
Kahn-Priddy theorem (\cite{kahn_applications_1972})
says that
$t_\#$ is onto in
the $2$-primary components for $\ast>0$.
This (Kahn-Priddy theorem) is seen as follows.
Take the $(n+1)$ loop map
\begin{align}
  \Loop^{n+1}\Susp^{n+1}P^n&\to{\Loop^{n+1}t_n}\Loop^{n+1}S^{n+1}
  \label{defn:loop-t_n}
\end{align}
of the $n$th transfer map $\Susp^{n+1} P^n \to{t_n}S^{n+1}$.
Let $(\Loop^{n+1}S^{n+1})_{-1}$ be the ``$-1$'' component of $\Loop^{n+1}S^{n+1}$
as mentioned in \eqref{eq:inrn}.
Let $\Loop^\infty P=\lim_{n}\Loop^{n+1}\Susp^{n+1}P^n$,
$\Loop^\infty S^0=\lim_{n}\Loop^{n+1}S^{n+1}$ and
$(\Loop^\infty S^0)_{-1}=\lim_{n}(\Loop^{n+1}S^{n+1})_{-1}$. Letting
$n-->\infty$ we get
a space map
\begin{align}
  \Loop^\infty P&\to{\Loop^\infty t}\Loop^\infty S^0.
  \label{defn:loop-t}
\end{align}
The stable homotopy groups $\pi_\ast^S(P)$ (resp.\ $\pi_\ast^S$) is one-to-one correspondent
to the unstable homotopy groups $\pi_\ast(\Loop^\infty P)$ (resp.\ $\pi_\ast(\Loop^\infty S^0)$),
and the group homomorphism $\pi_\ast^S(P)\to{t_\#}\pi_\ast^S$ corresponds to the
induced homomorphism
$\pi_\ast(\Loop^\infty P)\to{(\Loop^\infty t)_\#}\pi_\ast(\Loop^\infty S^0)$. Note that
\begin{note}
  $\pi_\ast(\Loop^\infty S^0)$ is equal to $\pi_\ast( (\Loop^\infty S^0)_{-1})$
  for all $\ast>0$.
\end{note}
\noindent
By Sullivan's theory of localization (\cite{sullivan_genetics_1974}) there is a space
$_2(\Loop^\infty S^0)_{-1}$ whose homotopy groups are $2$-primary, and also a map
$(\Loop^\infty S^0)_{-1}\to{L}{}_2(\Loop^\infty S^0)_{-1}$ such that
$\pi_\ast( (\Loop^\infty S^0)_{-1})\to{L_\#}\pi_\ast(_2(\Loop^\infty S^0)_{-1})$ is a $2$-primary isomorphism for all $\ast$.
Let $\Loop^\infty S^0\to{p}(\Loop^\infty S^0)_{-1}$ be the projection to the ``$-1$'' component.
Then Kahn-Priddy theorem says the following composite
\begin{align*}
  \Loop^\infty P&\to{\Loop^\infty t}\Loop^\infty S^0
  \to{p}(\Loop^\infty S^0)_{-1}
  \to{L}{}_2(\Loop^\infty S^0)_{-1}
\end{align*}
is a retraction map.
That is, there is a map $_2(\Loop^\infty S^0)_{-1}\to{i}\Loop^\infty P$ such that
the composite $_2(\Loop^\infty S^0)_{-1}\to{i}\Loop^\infty P\to{Lp\Loop^\infty t}(\Loop^\infty S^0)_{-1}$
is homotopic to the identity map.
Therefore $\Loop^\infty P \to{\Loop^\infty t}\Loop^\infty S^0$ induces a homomorphism $\pi_\ast^S(P)\to{t_\#}\pi_\ast^S$ that is surjective on the $2$-primary homotopy groups for $\ast>0$.
Summarizing, we have the following result. 
\begin{note}
  Given a $2$-primary homotopy element $\alpha\in\pi_\ast^S$ with order $k\ge1$ for $\ast>0$.
  Then there exists an element $\beta\in\pi_\ast^S(P)$ also with order $k$, such that
  $t_\#(\beta)=\alpha$.
  \label{eq:Kahn-Priddy}
\end{note}

The transfer map $P\to{t}S^0$ also induces an $\Ext$ group map
\begin{align}
  \Ext_\A^{s,t}(P)&\to{t_\ast}\Ext_\A^{s+1,t+1}.
  \label{eq:t_ast}
\end{align}
This $\Ext$ group map \eqref{eq:t_ast} is a realization of the
map $\pi_\ast^S(P)\to{t_\#}\pi_\ast^S$.
The algebraic Kahn-Priddy theorem by Lin \cite{lin_algebraic_1981} says that
$t_\ast$ is also onto for $t-s>0$.
In particular we have the following.
\begin{note}
  \label{eq:hnh}
  For $n\ge1$ there exists a class $\widehat{h}_n\in\Ext_\A^{0,2^n-1}(P)$ such that
  $t_\ast(\widehat{h}_n)=h_n$. So $t_\ast(\widehat{h}_nh_n)=h_n^2$.
  In fact it is easy to see that
  $\Ext_\A^{0,2^n-1}(P)=\Z/2(\widehat{h}_n)$, and
  $\Ext_\A^{1,2^{n+1}-1}(P)=\Z/2(\widehat{h}_nh_n)$,
  for $n\ge1$.
\end{note}

There is a
result in $\pi_\ast^S(P)$
that is parallel to \autoref{thm:t5t6-detected}
by Kahn-Priddy theorem.
To describe the result in $\pi_\ast^S(P)$
which is \autoref{thm:t5ht6-detected} below,
first we note the following observation.
\begin{note}
  \label{eq:theta_n-hat}
  Suppose there is an element $\theta_n\in\pi_{2^{n+1}-2}^S$ detected by $h_n^2$
  in the ASS for $\pi_\ast^S$
  with $2\theta_n=0$.
  Then by \eqref{eq:Kahn-Priddy} and \eqref{eq:hnh} there exists an element
  $\widehat{\theta}_n\in\pi_{2^{n+1}-2}^S(P)$ detected by $\widehat{h}_nh_n$ such that
  $t_\#(\theta_n)=\widehat{\theta}_n$ and $2\widehat{\theta}_n=0$.
  By dimensional reason, $\widehat{\theta}_n\in\pi_{2^{n+1}-2}^S(P)$ can be pulled back to an element
  in $\pi_{2^{n+1}-2}^S(P^{2^{n+1}-2})$, which will also be denoted by $\widehat{\theta}_n$.
\end{note}
\noindent
We already know there exist $\theta_4\in\pi_{30}^S$ and $\theta_5\in\pi_{62}^S$ with
$2\theta_4=0$, $2\theta_5=0$,
so by \eqref{eq:theta_n-hat} there exist $\widehat{\theta}_4\in\pi_{30}^S(P)$ and $\widehat{\theta}_5\in\pi_{62}^S(P)$
with $2\widehat{\theta}_4=0$, $2\widehat{\theta}_5=0$. Also,
these latter $\widehat{\theta}_4$ and $\widehat{\theta}_5$ have pullbacks
$\widehat{\theta}_4\in\pi_{30}^S(P^{30})$ and $\widehat{\theta}_5\in\pi_{62}^S(P^{62})$ respectively.
We will see in a moment (\eqref{eq:2t-nonzero} below) that these pullbacks
$\widehat{\theta}_4$, $\widehat{\theta}_5$ do not have
$2\widehat{\theta}_4=0$, $2\widehat{\theta}_5=0$.
Throughout this paper we fix a $\theta_4\in\pi_{30}^S$, a $\theta_5\in\pi_{62}^S$ with
$2\theta_4=0$, $2\theta_5=0$,
and a $\widehat{\theta}_4\in\pi_{30}^S(P)$, a $\widehat{\theta}_5\in\pi_{62}^S(P)$ with $2\widehat{\theta}_4=0$,
$2\widehat{\theta}_5=0$. We also fix a pullback $\widehat{\theta}_4\in\pi_{30}^S(P^{30})$, and
a pullback $\widehat{\theta}_5\in\pi_{62}^S(P^{62})$.

In \cite{lin_ext_2008} it is shown that
for each $i\ge1$
there is a non-zero class $\widehat{g}_i\in\Ext_\A^{3,2^{i+3}+2^{i+2}-1}(P)$
with $t_\ast(\widehat{g}_i)=g_i\in\Ext_\A^{4,2^{i+3}+2^{i+2}}$ where
$t_\ast$ is as in \eqref{eq:t_ast}.
In \autoref{se:Ext} we will show
that $\widehat{g}_4\in\Ext_\A^{3,191}(P)$ pulls back to a class in $\Ext_\A^{3,191}(P^{62})$,
which is also denoted by $\widehat{g}_4$.
There is a class $\widehat{h}_5\in\Ext_\A^{1,32}(P)$ by \eqref{eq:hnh}.
It is clear that
$\widehat{h}_5\in\Ext_\A^{1,32}(P)$ pulls back uniquely to a class in $\Ext_\A^{1,32}(P^{62})$,
which is also denoted by $\widehat{h}_5$.
We will see in a moment that
$\widehat{g}_4h_0^2,\widehat{h}_5h_5D_3(1) \in\Ext_\A^{5,193}(P^{62})$ are linearly independent.
The following is the result in $\pi_\ast^S(P)$ parallel to \autoref{thm:t5t6-detected}.
\begin{thm}
  \label{thm:t5ht6-detected}
  Let $\widehat{\theta}_5\in\pi_{62}^S(P^{62})$ be as in \eqref{eq:theta_n-hat}.
  In $\pi_{188}^S(P^{62})$,
  $\widehat{\theta}_5\theta_6\ne0$ and is detected in the ASS for $\pi_\ast^S$
  by a class
  \begin{align*}
    \tau
    =\widehat{g}_4h_0^2,
    +\varepsilon\widehat{h}_5h_5D_3(1)
    +\bar{\tau}\ne0
    \quad\mbox{in $\Ext_\A^{5,193}(P^{62})$}
  \end{align*}
  where
  $\varepsilon=0$ or $1$.
  Here
  $\bar{\tau}\in\Ext_\A^{5,193}(P^{62})$ is a class such that
  $t_\ast(\bar{\tau})=0$ in $\Ext_\A^{6,194}$.
  Moreover,
  $t_\ast(\widehat{g}_4h_0^2)=h_0^2g_4$, and
  $t_\ast(\widehat{h}_5h_5D_3(1))= h_5^2D_3(1)$.
\end{thm}
\noindent
We remark that $t_\ast$ in \autoref{thm:t5ht6-detected}
is the restriction to $P^{62}$ of the $t_\ast$ in \eqref{eq:t_ast}.
And also note that from \eqref{eq:h_0^2g_4-ne-0} and \autoref{thm:t5ht6-detected}
we see that
the classes $\widehat{g}_4h_0^2$, $\widehat{h}_5h_5D_3(1)$, and $\bar{\tau}$ are linearly independent in
$\Ext_\A^{5,193}(P^{62})$ (if $\bar{\tau}\ne0$),
because $t_\ast(\widehat{g}_4h_0^2)=h_0^2g_4$, $t_\ast(\widehat{h}_5h_5D_3(1))=h_5^2D_3(1)$,
and $t_\ast(\bar{\tau})=0$.
It is easy to see then that \autoref{thm:t5t6-detected} is a direct consequence of \autoref{thm:t5ht6-detected} and \autoref{remark:h02g4-not-a-boundary}.

To prove \autoref{thm:t5ht6-detected}
we will need to recall a well-known result on the stable liftings
of the stable map $S^{62}\to{2\widehat{\theta}_5}P^{62}$,
and to make an analysis
on the ``Adams' filtration'' of the homotopy element
$\widehat{\theta}_5\theta_6\in\pi_{188}^S(P^{62})$ which is
represented by the composite $S^{188}\to{\theta_6}S^{62}\to{\widehat{\theta}_5}P^{62}$.
These are recalled in the next few paragraphs.

First note that the stable map
$S^{62}\to{2\widehat{\theta}_5}P^{62}$ is non-zero,
because $\widehat{h}_5h_0h_5\in\Ext_\A^{2,64}(P^{62})$ is a non-zero infinite cycle
detecting $2\widehat{\theta}_5$ in the ASS for $\pi_\ast^S(P^{62})$.
In general, we
have the following.
\begin{note}
  If $S^{2^{i+1}-2}\to{\widehat{\theta}_i}P^{2^{i+1}-2}$ is a lifting of a
  Kervaire invariant element
  $S^{2^{i+1}-2}\to{\theta_i}S^0$ with $2\theta_i=0$ ($i\ge4$), then
  the stable map $S^{2^{i+1}-2}\to{2\widehat{\theta}_i}P^{2^{i+1}-2}$
  is non-zero as
  $2\widehat{\theta}_i\in\pi_{2^{i+1}-2}^S(P^{2^{i+1}-2})$ is detected in the ASS for $\pi_\ast^S(P^{2^{i+1}-2})$
  by $\widehat{h}_ih_0h_i$.
  \label{eq:2t-nonzero}
\end{note}
\noindent
For $0\le l\le m$,
we say a stable map $S^m\to{\varphi}P^m$ can be lifted to $P^l$ if
there is a stable map $S^m\to{\psi}P^l$ such that
the composite $S^m\to{\psi}P^l\to{i}P^m$ is $\varphi$
as shown in the following diagram \eqref{eq:lifting},
where $P^l\to{i}P^m$ is the inclusion map.
$\psi$ is said to be a lifting of $\varphi$ to $P^l$.
\begin{align}
  \begin{split}
  \xymatrix{
  & P^l\ar[d]^i\\
  S^m\ar[r]^\varphi\ar[ur]^\psi&  P^m
  }
  \end{split}
  \label{eq:lifting}
\end{align}
The result we are going to use is the following.
\begin{prop}\label{thm:lifting-of-2theta5}
  The stable map
  $S^{62}\to{2\widehat{\theta}_5}P^{62}$ can be lifted to $P^{51}$ and no further.
\end{prop}
The proof of \autoref{thm:lifting-of-2theta5} is based on the theory of vector fields on spheres.
We recall this theory as follows.

For each $n\ge1$, a non-zero vector field on $S^{n-1}$ is a continuous map
$S^{n-1}\to{\sigma}S^{n-1}$ such that $\sigma(v)\perp v$ for all $v\in S^{n-1}$.
A set of non-zero vector fields $\{\sigma_1,\cdots,\sigma_r\}$ on $S^{n-1}$
is said to be linearly independent if the set of vectors
$\{\sigma_1(v),\cdots,\sigma_r(v)\}$ is linearly independent for all $v\in S^{n-1}$.
In the case where there is a linearly independent set $\{\sigma_1,\cdots,\sigma_r\}$
of vector fields on $S^{n-1}$, we say $S^{n-1}$ admits an $r$-field.
\begin{note}\label{defn:rho}
  Given an integer $n\ge2$,
  write $n=2^b(2a+1)$, and let $b=c+4d$, where
  $0\le c\le3$. Let
  $\rho(n)=2^c +8d$.
  It is classical that for each $n\ge2$, $S^{n-1}$ admits a $(\rho(n)-1)$-field.
  (\cite{husemoller_fibre_1994})
\end{note}
\noindent
J.\ F.\ Adams proves the following result in \cite{adams_vector_1962}.
\begin{thm}\label{thm:Adams-vector-1962}
  Let $\rho(n)$ be as defined in \eqref{defn:rho}.
  $S^{n-1}$ does not admit a $\rho(n)$-field.
\end{thm}
For $n\ge1$, the Whitehead square
$[\iota_n,\iota_n]=\omega_n\in\pi_{2n-1}(S^n)$
is the homotopy element of the following composite of maps,
\begin{align*}
  S^{2n-1}\to{\sigma}& S^n\V S^n\to{F} S^n,
\end{align*}
where $\sigma$ is the attaching map of the top cell of $S^n\times S^n=(S^n\V S^n)\U_\sigma e^{2n}$,
and $F$ is the folding map.
We say $\omega_n$ can be desuspened $r$ times
or $\omega_n$ is an $r$-fold suspension (for $r\ge0$),
if there is a map
$S^{2n-1-r}\to{\widehat{\omega}_n}S^{n-r}$ such that
$\left[\Susp^r\widehat{\omega}_n\right]=\omega_n$.
The problem of how far can $\omega_n$ be desuspended
is connected to the vector field problem
on the spheres in the following content
as shown in \cite{james_whitehead_1957}.
\begin{thm}\label{thm:James-Whitehead-1957}
  Given an integer $r\ge0$ with $r\le n$.
  If $S^n$ admits an $r$-field then $\omega_n$ is an $r$-fold suspension.
  Conversely if $\omega_n$ is an $r$-fold suspension,
  and if $n>2r$, then $S^n$ admits an $r$-field.
\end{thm}
For $n\ge1$ with $n\le5$, let $m=2^{n+1}-2$,
consider a space map $\Susp^{m+1}S^m \to{2\widehat{\theta}_n}\Susp^{m+1}P^m$
(which is within the stable range)
that represents $2\widehat{\theta}_n\in\pi_m^S(P^m)$.
We refer to \cite{barratt_kervaire_1983} for the following.
\begin{note}\label{eq:lifting-and-omega_desuspension}
  For $r<\frac{m+2}{2}$,
  $\omega_{m+1}$ can be desuspended $r$ times if and only if
  there is a lifting $\Susp^{m+1}S^m \to{\psi}\Susp^{m+1}P^{m-r}$ of the stable map
  $\Susp^{m+1}S^m \to{2\widehat{\theta}_n}\Susp^{m+1}P^m$.
\end{note}
\noindent
Combining the results in \autoref{thm:Adams-vector-1962}, \autoref{thm:James-Whitehead-1957},
and \eqref{eq:lifting-and-omega_desuspension}, we have the following
\begin{coro}\label{thm:lifting-of-2theta_n-and-rho}
  There is a (stable) lifting
  $S^{2^{n+1}-2}\to{\psi}P^{2^{n+1}-2-k}$
  of $S^{2^{n+1}-2}\to{2\widehat{\theta}_n}P^{2^{n+1}-2}$ if and only if
  $k\le\rho(2^{n+1}) -1$.
\end{coro}
\noindent
When $n=5$, $\rho(64)=12$, so \autoref{thm:lifting-of-2theta5} is a special case
of \autoref{thm:lifting-of-2theta_n-and-rho}.

Next we define
``Adams' filtration'' $AF(\alpha)$ of a homotopy element $\alpha\in\pi_\ast^S(X)$
as follows.

Let $K$ be the Eilenberg-MacLane spectrum $K(\Z/2)$, that is,
$\pi_0^S(K)=\Z/2$ and $\pi_\ast^S(K)=0$ for $\ast\ne0$.
Given a stable complex $X$, recall that the mod $2$ Adams resolution of $X$
consists of the followings:
for each $s\ge0$ there is a spectrum $X_s$,
a map $X_{s+1}\to{p_{s+1}}X_s$,
a set of integers $I_s\subset\N$,
a spectrum $K_s=\V_{q\in I_s}\Susp^q K$, and
a map $X_s\to{f_s}K_s$, such that
$X_0 =X$, and for each $s\ge0$,
the map $X_s \to{f_s}K_s$ induces epimorphism on the mod $2$ cohomology,
and there is a cofibration sequence
\begin{align*}
  X_{s+1}\to{p_{s+1}}  X_s\to{f_s} K_s
\end{align*}
Given a map $Y\to{\varphi}X$, we say $AF(\varphi)$, the Adams' filtration of $\varphi$, is $l$ if
$l\ge0$ and $l$ is the largest integer
such that
there is a lifting $Y\to{\varphi_l}X_l$ of $\varphi$ to $X_l$
in the following tower ($l$ may be $\infty$):
\begin{align*}
  \xymatrix{
  && \vdots\ar[d]\\
  &\Susp^{-1}K_l\ar[r]^{j_{l+1}}&
  X_{l+1}\ar[d]^{p_{l+1}}\ar[r]^{f_{l+1}}&    K_{l+1}\\
  &\Susp^{-1}K_{l-1}\ar[r]^{j_l}&
    X_l\ar[d]^{p_l}\ar[r]^{f_l}&  K_l\\
  &&
    X_{l-1}\ar[d]^{p_{l-1}}\ar[r]^{f_{l-1}}&  K_{l-1}\\
  && \vdots\ar[d]\\
  &\Susp^{-1}K_1\ar[r]^{j_2}&
    X_2\ar[d]^{p_2}\ar[r]^{f_2}&  K_2\\
  &\Susp^{-1}K_0\ar[r]^{j_1}&
    X_1\ar[d]^{p_1}\ar[r]^{f_1}&  K_1\\
  Y\ar[rr]^{\varphi}\ar@/^2pc/[uuuuurr]^{\varphi_l}
  && X_0\ar[r]^{f_0}&  K_0
  }
\end{align*}
Let $Z\to{\psi}Y$ be another map, then it is not difficult to show that
$AF(\varphi\psi)\ge AF(\varphi)$.

In order to prove
\autoref{thm:t5ht6-detected},
we have to prove several propositions first.
For any $0\le k\le n$ let $P^n_k=P^n/P^{k-1}$.
Let $P^{62} \to{q_k}P^{62}_k$
be the collapsing map and
$P^k \to{i_k}P^{62}$ be the inclusion map.
\begin{prop}
  \label{thm:qk2t5:essential}
  For any integer $k\ge1$, the composite of stable maps
  $S^{62} \to{2\iota}S^{62} \to{\widehat{\theta}_5}P^{62} \to{q_k}P^{62}_k$
  is essential if and only if $k\le 51$.
\end{prop}
\begin{proof}
  The composite $S^{62} \to{2\iota}S^{62} \to{\widehat{\theta}_5}P^{62} \to{q_k}P^{62}_k$
  is zero if and only of there is a lifting of $S^{62} \to{2\widehat{\theta}_5}P^{62}$
  to $P^{k-1}$ because $P^{62}_k =P^{62}/P^{k-1}$.
  By \autoref{thm:lifting-of-2theta5} such a lifting exists
  if and only if $k-1\ge 51$.
  So the composite $S^{62} \to{2\iota}S^{62} \to{\widehat{\theta}_5}P^{62} \to{q_k}P^{62}_k$
  is zero if and only if $k\ge52$, or equivalently,
  the composite $S^{62} \to{2\iota}S^{62} \to{\widehat{\theta}_5}P^{62} \to{q_k}P^{62}_k$
  is essential if and only if $k\le51$.
  This proves \autoref{thm:qk2t5:essential}.
\end{proof}
\begin{coro}
  \label{thm:q47-theta5-essential}
  The composite of the stable maps $S^{62} \to{\widehat{\theta}_5}P^{62} \to{q_{47}}P^{62}_{47}$ is essential.
\end{coro}
\begin{proof}
  By \autoref{thm:qk2t5:essential},
  $S^{62} \to{q_{47}(2\widehat{\theta}_5)}P^{62}_{47}$ is essential, so 
  $S^{62} \to{q_{47}\widehat{\theta}_5}P^{62}_{47}$ must also be essential.
\end{proof}

The homotopy element $q_{47}\widehat{\theta}_5\in\pi_{62}^S(P^{62}_{47})$
will be described more clearly in \autoref{thm:P62-46:pi62} through \autoref{thm:q47-theta5-eq-2w}
when we give a concrete result on the group $\pi_{62}^S(P^{62}_{47})$.
For this purpose we are going to obtain a partial result
on the Adams filtration $AF(q_{47}\widehat{\theta}_5)$ in \autoref{thm:q47t5:AF-ge-3}.
To this end
we need to consider a (stable) complex $\widetilde{P}^{62}$
and also a map
$P^{62}\to{\widetilde{q}}\widetilde{P}^{62}$.
These are described in the next paragraph.

Recall the stable homotopy element
$\widehat{\theta}_4\in\pi_{30}^S(P)$ with $2\widehat{\theta}_4=0$ in \eqref{eq:theta_n-hat}.
Consider a pullback $\widehat{\theta}_4\in\pi_{30}^S(P^{31})$ of $\widehat{\theta}_4\in\pi_{30}^S(P)$.
From the cofibration sequence
\begin{align*}
  \cdots->\Susp^{-1}P_{32}->P^{31}->P->P_{32}->\cdots,
\end{align*}
we have $\pi_{30}^S(P^{31})\iso \pi_{30}^S(P)$
since $\pi_{30}^S(P_{32})=0$ and $\pi_{31}^S(P_{32})=0$.
So $2\widehat{\theta}_4=0\in\pi_{30}^S(P^{31})$ also.
Now consider the pullback of $\widehat{\theta}_4\in\pi_{30}^S(P)$ to $\pi_{30}^S(P^{30})$
as defined in \eqref{eq:theta_n-hat} which has $2\widehat{\theta}_4\ne0\in\pi_{30}^S(P^{30})$
by \eqref{eq:2t-nonzero}. The composite
\begin{align*}
  S^{30}&\to{2\iota}S^{30}\to{\widehat{\theta}_4}P^{30}
\end{align*}
is non-zero. Let $P^{30}\to{i_{30}}P^{31}$ be the inclusion map, then the composite
\begin{align*}
  S^{30}&\to{2\iota}  S^{30}\to{\widehat{\theta}_4} P^{30}\to{i_{30}}  P^{31}
\end{align*}
is zero as already shown above.
So $S^{30} \to{i_{30}\widehat{\theta}_4}P^{31}$ has an extension $\bar{h}_5^\prime$:
\begin{align*}
  S^{30}\U_{2\iota}e^{31} &\to{\bar{h}_5^\prime}  P^{31}.
\end{align*}
Let $P^{31}\to{i_{31}}P^{62}$ be the inclusion map, and let $\bar{h}_5$ be the composite
\begin{align}
  \bar{h_5}: S^{30}\U_{2\iota}e^{31}\to{\bar{h}_5^\prime}P^{31}\to{i_{31}}P^{62}.
  \label{eq:h5-bar}
\end{align}
\begin{defn}
  \label{defn:P62t}
  \label{defn:qt}
  Define $\widetilde{P}^{62}$ to be the mapping cone $C_{\bar{h}_5}=P^{62}\U_{\bar{h}_5}C(S^{30}\U_{2\iota}e^{31})$, and
  define $P^{62}\to{\widetilde{q}}\widetilde{P}^{62}=C_{\bar{h}_5}$ to be the inclusion map.
\end{defn}
\noindent
There are two properties of this complex $\widetilde{P}^{62}$,
which are \eqref{eq:h5-bar-iso} and \eqref{eq:qt-h5t-eq-0} below,
that will be important to our proof later.
\begin{note}
  $\bar{h}_5$ induces an isomorphism
  on the mod $2$ homology groups
  $H_{31}(S^{30}\U_{2\iota}e^{31})\to[\iso]{(\bar{h}_5)_\ast} H_{31}(P^{62})$.
  \label{eq:h5-bar-iso}
\end{note}
\noindent
Recall from \eqref{eq:hnh} that $\Ext_\A^{0,31}(P^{62})\iso\Ext_\A^{0,31}(P)=\Z/2(\widehat{h}_5)$.
Then \eqref{eq:h5-bar-iso} implies
\begin{note}
  $\Ext_\A^{0,31}(P^{62})\to{\widetilde{q}_\ast}\Ext_\A^{0,31}(\widetilde{P}^{62})$ has
  $\widetilde{q}_\ast(\widehat{h}_5)=0$.
  \label{eq:qt-h5t-eq-0}
\end{note}

\eqref{eq:h5-bar-iso} is shown as follows.
We have the following diagram,
where $S^{30}\to{\bar{\sigma}}\Susp^{-1}P^{31}$
is the attaching map of the top cell of $\Susp^{-1}P^{32}$, and all other
horizontal maps are cofibrations.
\begin{align*}
  \xymatrix{
  &&S^{30}\ar[r]^{2\iota}\ar[d]^{n\iota}\ar@/_/[lld]_{m\iota}&
    S^{30}\ar[r]\ar[d]^{\widehat{\theta}_4}&
    S^{30}\U_{2\iota}e^{31}\ar[r]\ar[d]^{\bar{h}_5^\prime}&
    S^{31}\ar[d]^{n\iota}\\
  S^{30}\ar[r]^-{\bar{\sigma}}\ar@/_1pc/[rr]_{2\iota}&
    \Susp^{-1}P^{31}\ar[r]^{\Susp^{-1}q}&
    S^{30}\ar[r]^{\sigma}&
    P^{30}\ar[r]^{i_{30}}&
    P^{31}\ar[r]^q\ar[d]^{i_{31}}&
    S^{31}\\
  &&&&P^{62}
  }
\end{align*}
$S^{30} \to{i_{30}\widehat{\theta}_4(2\iota)}P^{31}$ is zero implies that
there is a pullback $S^{30}\to{n\iota}S^{30}$
of $S^{30} \to{2\widehat{\theta}_4}P^{30}$ for some $n\in\Z$, and $n\iota \ne0$ because
$S^{30} \to{2\widehat{\theta}_4}P^{30}$ is non-zero.
This pullback $n\iota$ is subjected to an indeterminancy
$S^{30}->\Susp^{-1}P^{31}\to{\Susp^{-1}q}S^{30}$.
Note that the composite
$S^{30} \to{\bar{\sigma}}\Susp^{-1}P^{31} \to{\Susp^{-1}q}S^{30}$ is $2\iota$.
So $S^{30} \to{\Susp^{-1}q\bar{\sigma}(m\iota)=2m\iota}S^{30}$ is an indeterminancy for any $m\in\Z$.
Thus we can choose a pullback $S^{30}\to{n\iota}S^{30}$ of $S^{30} \to{2\widehat{\theta}_4}P^{30}$
with $n=1$. By this choice,
$H_{31}(S^{30}\U_{2\iota}e^{31};\Z)\to{(\bar{h}_5)_\ast}H_{31}(P^{62};\Z)$
is an isomorphism.
This proves \eqref{eq:h5-bar-iso}.

We can now prove a partial result on the Adams filtration of
$S^{62} \to{q_{47}\theta_5}P^{62}_{47}$.
\begin{prop}
  \label{thm:q47t5:AF-ge-3}
  $AF(q_{47}\widehat{\theta}_5)\ge3$ for the composite
  $S^{62} \to{\theta_5}P^{62} \to{q_{47}}P^{62}_{47}$.
\end{prop}
To prove \autoref{thm:q47t5:AF-ge-3} we will use the following results
\eqref{eq:q47p-phi-eq-q47} and \eqref{eq:boundary:P62t:Ext:2-64}.
Let $\widetilde{P}^{46}$ be the $46$-skeleton of $\widetilde{P}^{62}$.
From \eqref{eq:h5-bar} we see that
$\widetilde{P}^{62}/\widetilde{P}^{46}=P^{62}_{47}$.
Let $\widetilde{P}^{62}\to{q_{47}^\prime}\widetilde{P}^{62}/\widetilde{P}^{46}=P^{62}_{47}$ be the collapsing map.
We have the following diagram consisting of a commutative square together with the map $S^{62}\to{\widehat{\theta}_5}P^{62}$:
\begin{align}
  \begin{split}
  \xymatrix{
  S^{62}\ar[r]^{\widehat{\theta}_5}&
    P^{62}\ar[r]^{q_{47}}\ar[d]^{\widetilde{q}}&  P^{62}_{47}\ar@{=}[d]\\
  & \widetilde{P}^{62}\ar[r]^-{q_{47}^\prime}&    \widetilde{P}^{62}/\widetilde{P}^{46}
  }
  \end{split}
  \label{eq:q47p-phi-eq-q47}
\end{align}
In \autoref{se:Ext}
we will prove the following result
\begin{note}
  \label{eq:boundary:P62t:Ext:2-64}
  $\Ext_\A^{2,64}(\widetilde{P}^{62}) =\Z/2(\bar{\widetilde{e}_{31}h_5}h_0)\dsum\Z/2(\bar{e_{47}}h_0h_4)$, and
  $(q_{47}^\prime)_\ast(\bar{\widetilde{e}_{31}h_5}h_5) =(q_{47}^\prime)_\ast(\bar{e_{47}}h_0h_4) =0$
  in $\Ext_\A^{2,64}(P^{62}_{47})$.
\end{note}
\begin{proof}[Proof of \autoref{thm:q47t5:AF-ge-3}]
  From \eqref{eq:q47p-phi-eq-q47} we see that
  $q_{47}\widehat{\theta}_5=q_{47}^\prime\widetilde{q}\widehat{\theta}_5$.
  $\widehat{\theta}_5$ is detected by $\widehat{h}_5h_5$ in the ASS for $\pi_\ast^S(P^{62})$. By
  \eqref{eq:qt-h5t-eq-0}, $\widetilde{q}_\ast(\widehat{h}_5)=0$.
  This implies $AF(\widetilde{q}\widehat{\theta}_5)\ge2$.
  If $AF(\widetilde{q}\widehat{\theta}_5)=2$, then
  $\widetilde{q}\widehat{\theta}_5$ is detected by some class
  $\tau_1\in\Ext_\A^{2,64}(\widetilde{P}^{62})$.
  By \eqref{eq:boundary:P62t:Ext:2-64},
  $(q_{47}^\prime)_\ast(\tau_1)=0$; so
  $AF(q_{47}^\prime\widetilde{q}\widehat{\theta}_5)\ge3$.
  If $AF(\widetilde{q}\widehat{\theta}_5)\ge3$ then clearly
  $AF(q_{47}^\prime\widetilde{q}\widehat{\theta}_5)\ge3$.
  In both cases we have $AF(q_{47}\widehat{\theta}_5)=AF(q_{47}^\prime\widetilde{q}\widehat{\theta}_5)\ge3$.
  This proves \autoref{thm:q47t5:AF-ge-3}.
\end{proof}

We have to describe the map $S^{62}\to{q_{47}\widehat{\theta}_5}P^{62}_{47}$ more precisely
than the result in \eqref{thm:q47t5:AF-ge-3} in order to prove
\autoref{thm:t5ht6-detected}. First we describe some specific homotopy elements
$y,x,w\in\pi_{62}^S(P^{62}_{46})$, and $v\in\pi_{62}^S(P^{62}_{47})$ in
\eqref{eq:y} through \eqref{eq:v} below. In
\eqref{eq:y}, \eqref{eq:w}, and \eqref{eq:v} the horizontal sequences are the obvious
cofibration sequences.
\begin{note}
  \label{eq:y}
  Define $y\in\pi_{62}^S(P^{62}_{46})$ as follows.
  Recall the Hopf classes $S^1\to{\eta}S^0$ and $S^7\to{\sigma}S^0$.
  Since $2\sigma^2=0$ and $\eta\sigma^2=0$, and since
  $P^{48}_{46}=(S^{46}\V S^{47})\U_{\eta\V2\iota}e^{48}$,
  there is a
  coextension $\widetilde{\sigma^2}$:
  \begin{align*}
    \xymatrix{
    \cdots\ar[r]& S^{46}\V S^{47}\ar[r]^-i&
      P^{48}_{46}\ar[r]^-q& S^{48}\ar[r]^-{\eta\V2\iota}& S^{47}\V S^{48}\ar[r]&  \cdots.\\
    & & & S^{62}\ar[u]_{\sigma^2}\ar[ul]^{\widetilde{\sigma^2}}\\
    }
  \end{align*}
  Let $y$ be the composite $S^{62}\to{\widetilde{\sigma^2}}P^{48}_{46}\to{i_{48}}P^{62}_{46}$,
  where $i_{48}$ is the inclusion map. Then
  $y\in\pi_{62}^S(P^{62}_{46})$.
\end{note}
\begin{note}
  \label{eq:w}
  Define $w\in\pi_{62}^S(P^{62}_{46})$ as follows.
  Consider the subcomplex $S^{47}\U_\sigma e^{55}\subset P^{62}_{46}$.
  Since $\sigma(2\sigma)=0$ there is a non-zero coextension
  $\bar{2\sigma}$:
  \begin{align*}
    \xymatrix{
    \cdots\ar[r]& S^{47}\ar[r]^-{i^\prime}&
      S^{47}\U_\sigma e^{55}\ar[r]^-{q^\prime}&
      S^{55}\ar[r]^\sigma&  S^{48}\ar[r]
    &\cdots.\\
    & & & S^{62}\ar[u]^{2\sigma}\ar[ul]^{\bar{2\sigma}}
    }
  \end{align*}
  Let $w$ be the composite $S^{62}\to{\bar{2\sigma}}S^{47}\U_\sigma e^{55}\to{j}P^{62}_{46}$,
  where $j$ is the inclusion map. Then
  $w\in\pi_{62}^S(P^{62}_{46})$.
\end{note}
\begin{note}
  \label{eq:x}
  Define $x=w+y\in\pi_{62}^S(P^{62}_{46})$. 
\end{note}
\begin{note}
  \label{eq:v}
  Define $v\in\pi_{62}^S(P^{62}_{47})$ as follows.
  The cells in dimensions $47$, $51$, and $53$ in $P^{62}_{47}$ form a
  subcomplex which is $S^{47}\U_\nu e^{51}\U_\eta e^{53}\subset P^{62}_{47}$.
  Recall (\cite{toda_composition_1962})
  there is a stable map $S^8\to{\varepsilon}S^0$ which
  has the properties $\eta\varepsilon\ne0$ and $\eta^2\varepsilon=0$.
  It is not difficult to see that the Toda bracket $\toda{\nu,\eta,\eta\varepsilon}$ is zero.
  So there is a non-zero coextension $\widetilde{\eta\varepsilon}$ of $\eta\varepsilon$:
  \begin{align*}
    \xymatrix{
    \cdots\ar[r]&
    S^{47}\U_\nu e^{51}\ar[r]^-{i^{\prime\prime}}&
    S^{47}\U_\nu e^{51}\U_\eta e^{53}\ar[r]^-{q^{\prime\prime}}& S^{53}\ar[r]&
    S^{48}\U_\nu e^{52}\ar[r]&  \cdots\\
    & & & S^{62}\ar[u]_{\eta\varepsilon}
      \ar[ul]^{\widetilde{\eta\varepsilon}}
    }
  \end{align*}
  Let $v$ be the composite $S^{62}\to{\widetilde{\eta\varepsilon}}S^{47}\U_\nu e^{51}\U_\eta e^{53}\to{j'}P^{62}_{47}$,
  where $j'$ is the inclusion map. Then $v$
  is an element in $\pi_{62}^S(P^{62}_{47})$.
\end{note}
\noindent
We will show in \autoref{se:pi} that the elements
$y,x,w\in\pi_{62}^S(P^{62}_{46})$ and $v\in\pi_{62}^S(P^{62}_{47})$ are non-zero.
In fact, in \autoref{se:pi} we will the following
\autoref{thm:P62-46:pi62} and \autoref{thm:P62-47:pi62}.
\begin{thm}\label{thm:P62-46:pi62}
  $\pi^S_{62}(P^{62}_{46})=\Z/4(y)\dsum\Z/8(x)$.
  Moreover $AF(y)=1$, $AF(w)=1$, $AF(2y)=2$, and $AF(x =y+w)=2$.
\end{thm}
\noindent
Let $P^{62}_{46}\to{q_1}P^{62}_{47}$ be the collapsing map, and
consider the induced homomorphism $\pi_{62}^S(P^{62}_{46})\to{(q_1)_\#}\pi_{62}^S(P^{62}_{47})$.
Let $\widetilde{y}=(q_1)_\#(y)$, $\widetilde{w}=(q_1)_\#(w)$, and
$\widetilde{x}=(q_1)_\#(x=y+w)=\widetilde{y}+\widetilde{w}$ in $\pi_{62}^S(P^{62}_{47})$.
\begin{thm}\label{thm:P62-47:pi62}
  $\pi^S_{62}(P^{62}_{47})=\Z/2(\widetilde{y})\dsum\Z/8(\widetilde{x})\dsum\Z/2(v-2\widetilde{x})$. Moreover
  $AF(\widetilde{y})=1$, $AF(\widetilde{x})=2$, and $AF(v)=4$.
\end{thm}
Now we can describe precisely the map
$q_{47}\widehat{\theta}_5$ considered in \autoref{thm:q47-theta5-essential}.
\begin{prop}
  \label{thm:q47-theta5-eq-2w}
  $q_{47}\widehat{\theta}_5=2\widetilde{x}\in\pi_{62}^S(P^{62}_{47})
  =\Z/2(\widetilde{y})\dsum\Z/8(\widetilde{x})\dsum\Z/2(v-2\widetilde{x})$.
\end{prop}
\begin{proof}
  Let $P^{62}\to{q_{46}}P^{62}_{46}$ be the collapsing map.
  It is clear that $(q_{46})_\ast(\widehat{h}_5)=0$,
  and since $\widehat{\theta}_5$ is detected by $\widehat{h}_5h_5$, it follows that
  $AF(q_{46}\widehat{\theta}_5)\ge2$.
  From \autoref{thm:P62-46:pi62} we thus see that
  $q_{46}\widehat{\theta}_5=ay +bx$ where $a=0$ or $2$ and
  $0\le b\le 7$ (if $a=1$ or $3$ then
  $AF(q_{46}\widehat{\theta}_5)=1$).
  Since $2\widetilde{y} =0$ in $\pi_{62}^S(P^{62}_{47})$,
  we have
  \begin{align*}
    q_{47}\widehat{\theta}_5
    &=q_1q_{46}\widehat{\theta}_5=q_1(ay +bx)
    =a\widetilde{y}+b\widetilde{x}
    =b\widetilde{x}\in\pi_{62}^S(P^{62}_{47})
  \end{align*}
  If $b$ is odd then $AF(q_{47}\widehat{\theta}_5)=2$ by \autoref{thm:P62-47:pi62}.
  But \autoref{thm:q47t5:AF-ge-3}
  says $AF(q_{47}\widehat{\theta}_5)\ge3$ so this is a contradiction.
  Therefore $b=0$, $2$, or $4$.
  By \autoref{thm:q47-theta5-essential}, $q_{47}\widehat{\theta}_5$ is essential so $b\ne0$.
  If $b=4$ then
  $2(q_{47}\widehat{\theta}_5)=2b\widetilde{x}=8\widetilde{x}=0$ in $\pi_{62}^S(P^{62}_{47})$,
  because the order of $\widetilde{x}$ is $8$.
  By \autoref{thm:qk2t5:essential} this is impossible.
  Therefore $b=2$, and $q_{47}\widehat{\theta}_5 =2\widetilde{x}$.
  This proves \autoref{thm:q47-theta5-eq-2w}.
\end{proof}


We are almost ready to prove the main result of this section which is \autoref{thm:t5ht6-detected}.
In \autoref{se:pi} we will prove the following.
\begin{note}
  \label{eq:boundary:2w-detected}
  $2\widetilde{x} \in\pi_{62}^S(P^{62}_{47})$ is detected by a non-zero class
  $\bar{e_{48}h_0h_3^2}\in\Ext_\A^{3,65}(P^{62}_{47})$.
\end{note}
\noindent
Note that $\bar{e_{48}h_0h_3^2}h_6^2\in\Ext_\A^{5,193}(P^{62}_{47})$.
In \autoref{se:Ext} we will obtain the following $\Ext$ groups results
\eqref{eq:boundary:P62:Ext-5-193} through \eqref{eq:boundary:h_0^2g_4}.
\begin{note}
  \label{eq:boundary:P62:Ext-5-193}
  \begin{enumerate}
    \item $\Ext_\A^{5,193}(P^{62})=
      \Z/2(\widehat{h}_5h_5D_3(1))
      \dsum \Z/2(\widehat{h}_4h_0^2c_4)
      \dsum \Z/2(\widehat{h}_4h_1f_3)
      \dsum \Z/2(\gamma_{61}h_7)
      \dsum \Z/2(\widehat{g}_4h_0^2)$.
    \item $\Ext_\A^{s,189+s}(P^{62}_{47}) =0$ for $0\le s\le2$.
    \item $\Ext_\A^{3,192}(P^{62}_{47})=\Z/2(\bar{e_{62}h_1}h_6^2)$.
    \item $\Ext_\A^{3,192}(P^{62}) \to{t_\ast}\Ext_\A^{4,193}$ is given by
      $t_\ast(\bar{e_{62}h_1}h_6^2) =h_0h_6^3 \ne0$.
  \end{enumerate}
\end{note}
\begin{note}
  \label{eq:boundary:h_0^2g_4}
  The map
  $\Ext_\A^{5,193}(P^{62}) \to{t_\ast}\Ext_\A^{6,194}$ as defined in \eqref{eq:t_ast}
  is given by the following.
  $t_\ast(\widehat{g}_4h_0^2)= h_0^2g_4\ne0$, and
  $t_\ast(\widehat{h}_5h_5D_3(1))= h_5^2D_3(1)\ne0$, and
  $t_\ast(\widehat{h}_4h_0^2c_4)= t_\ast(\widehat{h}_4h_1f_3)= t_\ast(\gamma_{61}h_7)= 0$
  in $\Ext_\A^{6,194}$.
\end{note}
\begin{note}
  \label{eq:boundary:e_48h_0h_3^2h_6^2-relation}
  The map
  $\Ext_\A^{5,193}(P^{62}) \to{(q_{47})_\ast}\Ext_\A^{5,193}(P^{62}_{47})$
  is given by the following.
  $(q_{47})_\ast(\widehat{h}_5h_5D_3(1)) =(q_{47})_\ast(\widehat{h}_4h_0^2c_4) =(q_{47})_\ast(\widehat{h}_4h_1f_3) =(q_{47})_\ast(\gamma_{61}h_7) =0$ and
  $(q_{47})_\ast(\widehat{g}_4h_0^2) =\bar{e_{48}h_0h_3^2}h_6^2\ne0$ in $\Ext_\A^{5,193}(P^{62}_{47})$.
\end{note}
\noindent
Finally, in \autoref{se:not-hit} we will prove the following result
from the $\Ext$ group information in \autoref{se:Ext}.
\begin{prop}
  \label{thm:t5ht6:AF-ge-5}
  For the composite $S^{188} \to{\theta_6}S^{62} \to{\widehat{\theta}_5}P^{62}$,
  $AF(\widehat{\theta}_5\theta_6)\ge5$.
\end{prop}
Now we prove \autoref{thm:t5ht6-detected} which is restated below.
\theoremstyle{nonumberplain}
\newtheorem{thm:t5ht6-detected}{\autoref{thm:t5ht6-detected}}
\begin{thm:t5ht6-detected}
  Let $\widehat{\theta}_5\in\pi_{62}^S(P)$ be as in \eqref{eq:theta_n-hat}.
  In $\pi_{188}^S(P^{62})$,
  $\widehat{\theta}_5\theta_6\ne0$ and is detected in the ASS for $\pi_\ast^S$
  by a class
  \begin{align*}
    \tau
    =\widehat{g}_4h_0^2,
    +\varepsilon\widehat{h}_5h_5D_3(1)
    +\bar{\tau}\ne0
    \quad\mbox{in $\Ext_\A^{5,193}(P^{62})$}
  \end{align*}
  where
  $\varepsilon=0$ or $1$.
  Here
  $\bar{\tau}\in\Ext_\A^{5,193}(P^{62})$ is a class such that
  $t_\ast(\bar{\tau})=0$ in $\Ext_\A^{6,194}$.
  Moreover,
  $t_\ast(\widehat{g}_4h_0^2)=h_0^2g_4$, and
  $t_\ast(\widehat{h}_5h_5D_3(1))= h_5^2D_3(1)$.
\end{thm:t5ht6-detected}
\begin{proof}
  By \eqref{thm:q47-theta5-eq-2w} and \eqref{eq:boundary:2w-detected}, $q_{47}\widehat{\theta}_5\in\pi_{62}^S(P^{62}_{47})$
  is detected by $\bar{e_{48}h_0h_3^2}$ in the ASS for $\pi_\ast^S(P^{62}_{47})$.
  By \eqref{eq:boundary:e_48h_0h_3^2h_6^2-relation},
  $\bar{e_{48}h_0h_3^2}h_6^2 =(q_{47})_\ast(\widehat{g}_4h_0^2) \ne0$
  in $\Ext_\A^{5,193}(P^{62}_{47})$.
  $(q_{47})_\ast(\widehat{g}_4h_0^2)$ is not a boundary in the ASS for $\pi_\ast^S(P^{62}_{47})$
  by \eqref{eq:boundary:P62:Ext-5-193}~(2), (3), (4),
  by $t_\ast(\widehat{g}_4h_0^2) =h_0^2g_4\ne0$ in \eqref{eq:boundary:h_0^2g_4}, and by
  the differential $d_2(h_0h_6^3) =h_0^2h_5^2h_6^2 =0$ in the ASS for $\pi_\ast^S$.
  Therefore $q_{47}\widehat{\theta}_5\theta_6\in\pi_{188}^S(P^{62}_{47})$
  is non-zero and is detected by $(q_{47})_\ast(\widehat{g}_4h_0^2)$. It follows that
  $\widehat{\theta}_5\theta_6\ne0$ and
  $AF(\widehat{\theta}_5\theta_6)\le AF(q_{47}\widehat{\theta}_5\theta_6) =5$.
  By \eqref{thm:t5ht6:AF-ge-5}, $AF(\widehat{\theta}_5\theta_6)\ge5$. So
  $AF(\widehat{\theta}_5\theta_6)=5$.
  From this and \eqref{eq:boundary:P62:Ext-5-193}~(1) we infer that $\widehat{\theta}_5\theta_6$ is detected by
  \begin{align*}
    \tau
    =\varepsilon_1\widehat{g}_4h_0^2,
    +\varepsilon_2\widehat{h}_5h_5D_3(1)
    +\varepsilon_3\widehat{h}_4h_0^2c_4
    +\varepsilon_4\widehat{h}_4h_1f_3
    +\varepsilon_5\gamma_{61}h_7
  \end{align*}
  for some $\varepsilon_i$ where $\varepsilon_i=0$ or $1$ for $1\le i\le5$
  and $(\varepsilon_1,\varepsilon_2,\varepsilon_3,\varepsilon_4,\varepsilon_5) \ne(0,0,0,0,0)$.
  From \eqref{eq:boundary:2w-detected} and \eqref{eq:boundary:e_48h_0h_3^2h_6^2-relation} we see
  $\varepsilon_1=1$.
  So let $\bar{\tau}=\varepsilon_3\widehat{h}_4h_0^2c_4 +\varepsilon_4\widehat{h}_4h_1f_3 +\varepsilon_5\gamma_{61}h_7$,
  and
  \begin{align*}
    \tau
    =\widehat{g}_4h_0^2,
    +\varepsilon_2\widehat{h}_5h_5D_3(1)
    +\bar{\tau}\in\Ext_\A^{5,193}(P^{62}).
  \end{align*}
  By \eqref{eq:boundary:h_0^2g_4}, we see
  $t_\ast(\widehat{g}_4h_0^2)=h_0^2g_4$,
  $t_\ast(\widehat{h}_5h_5D_3(1))= h_5^2D_3(1)$, and
  $t_\ast(\bar{\tau})=0$ in $\Ext_\A^{6,194}$.
  This proves \autoref{thm:t5ht6-detected}.
\end{proof}

From \autoref{remark:h02g4-not-a-boundary} we see
\autoref{thm:t5ht6-detected} implies the following theorem.
\newtheorem{thm:t5t6-detected}{\autoref{thm:t5t6-detected}}
\begin{thm:t5t6-detected}
  $\theta_5\theta_6$ in $\pi_{188}^S$
  is detected by $h_0^2g_4 +\varepsilon h_5^2D_3(1)\ne0$
  in $\Ext_\A^{6,194}$,
  where $\varepsilon=0\mbox{ or }1$.
\end{thm:t5t6-detected}

Finally we prove \autoref{thm:intro:boundary}, which is restated below.
\newtheorem{thm:intro:boundary}{\autoref{thm:intro:boundary}}
\begin{thm:intro:boundary}
  $h_0^3g_4$ is a boundary in the ASS for spheres.
\end{thm:intro:boundary}
\begin{proof}
  Recall \cite{lin_ext_2008} that there is a relation $h_0D_3(1)=0$ in $\Ext_\A^{5,131}$. Let
  $h_0^2g_4 +\varepsilon h_5^2D_3(1) \in\Ext_\A^{6,194}$
  be the class that detects $\theta_5\theta_6\in\pi_{188}^S$
  in \autoref{thm:t5t6-detected}, and we have
  \begin{align*}
    h_0^3g_4&=h_0(h_0^2g_4 +\varepsilon h_5^2D_3(1)),
  \end{align*}
  where $h_0^3g_4 \ne0$ in $\Ext_\A^{7,195}$ by \eqref{eq:intro:h03g4-non-zero}.
  Since $h_0$ detects $2\iota$, and
  since $h_0^2g_4 +\varepsilon h_5^2D_3(1)$ detects $\theta_5\theta_6$
  by \autoref{thm:t5t6-detected},
  $h_0^3g_4=h_0(h_0^2g_4 +\varepsilon h_5^2D_3(1))$ must be a boundary because
  $2\theta_5\theta_6=0$.
  This proves \autoref{thm:intro:boundary}.
\end{proof}

\newpage
\section{Calculations of some $\Ext$ groups}\label{se:Ext}
In this section we will describe and compute the $\Ext$ group results
that we already have used
in \autoref{se:intro}, \autoref{se:boundary}, and
will also be needed in \autoref{se:not-hit}.
There are many results, and we organize them as follows.
In \eqref{eq:Ext} and \eqref{eq:P:Ext} we list respectively
all the results about $\Ext_\A^{\ast,\ast}$ and $\Ext_\A^{\ast,\ast}(P)$
that we need,
and in \eqref{eq:P62:Ext} are the results about $\Ext_\A^{\ast,\ast}(P^{62})$.
In \eqref{eq:P62-k:Ext} are the results about $\Ext_\A^{\ast,\ast}(P^l_k)$
for some stunted projective spaces $P^l_k =P^l/P^{k-1}$.
Recall the complex $\widetilde{P}^{62}=P^{62}\U_{\bar{h}_5}C(S^{30}\U_{2\iota}e^{31})$ in \eqref{defn:P62t}.
The subcomplex $P^{16}\subset P^{62}$ is also a subcomplex of $\widetilde{P}^{62}$.
Let $\widetilde{P}^{62}_{17}=\widetilde{P}^{62}/P^{16}$.
In \eqref{eq:P62t:Ext} and \eqref{eq:P62-17t:Ext} are the results about
$\Ext_\A^{\ast,\ast}(\widetilde{P}^{62})$ and $\Ext_\A^{\ast,\ast}(\widetilde{P}^{62}_{17})$
respectively.
These are all the $\Ext$ group results that we will compute in this section.

The cohomology classes that appear in \eqref{eq:Ext} and \eqref{eq:P:Ext},
except $\beta_{187}\in\Ext_\A^{5,192}(P)$ in \eqref{eq:P:Ext}~(8),
are as defined in \cite{lin_ext_2008,chen_determination_2011},
and these will be recalled in a moment:
the classes $h_i,c_i,f_i,g_{i+1},D_3(i),V_i$ for $i\ge0$
in \eqref{eq:Ext} and \eqref{eq:P:Ext}
will be recalled in \eqref{eq:Ext:known},
and the classes $\widehat{h}_i,\widehat{g}_{i+1},\gamma_{61}(i),\widehat{D_3}(i)$
for $i\ge0$
in \eqref{eq:P:Ext}
will be recalled in \eqref{eq:P:Ext:known}.
The class $\beta_{187}$ will be defined in \eqref{eq:P:Ext:beta}.
We will simply write $D_3$ for $D_3(0)$, $\gamma_{61}$ for $\gamma_{61}(0)$, and
$\widehat{D_3}$ for $\widehat{D_3}(0)$.
\begin{note}
  \label{eq:Ext}
  \begin{enumerate}
    \item $\Ext_\A^{4,192}=\Z/2(g_4)$.
    \item $\Ext_\A^{5,193}=\Z/2(h_0g_4)$.
    \item $\Ext_\A^{s,s+189}=0$ for $0\le s\le2$.
    \item $\Ext_\A^{s,s+189}=\Z/2(h_0^{s-3}h_6^3)$, $3\le s\le4$.
    \item $\Ext_\A^{5,194}=\Z/2(h_0^2h_6^3)\dsum\Z/2(h_1g_4)$.
    \item $\Ext_\A^{s,s+190}=0$ for $0\le s\le2$.
    \item $\Ext_\A^{6,193}=\Z/2(h_5V_0)$.
    \item $\Ext_\A^{6,194}=\Z/2(h_0^2g_4)\dsum\Z/2(h_5^2D_3(1))$.
    \item $h_0^3g_4,h_1h_5V_0$ are linearly independent in $\Ext_\A^{7,195}$.
  \end{enumerate}
\end{note}
\noindent
Recall the map $P\to{t}S^0$ and the induced homomorphism
$\Ext_\A^{s,t}(P)\to{t_\ast}\Ext_\A^{s+1,t+1}$
mentioned in \eqref{eq:t_ast}.
\eqref{eq:P:Ext}~($j$)~(a) for $1\le j\le9$ below describe $\Ext_\A^{s,t}(P)$
for some particular $(s,t)$.
\eqref{eq:P:Ext}~($j$)~(b) for $j=1,5,8,9$, \eqref{eq:P:Ext}~($j$)~(c) for $2\le j\le9$ and
\eqref{eq:P:Ext}~(8)~(d) describe $t_\ast(x)\in\Ext_\A^{s+1,t+1}$
for the classes $x$ in \eqref{eq:P:Ext}~($j$)~(a) for $1\le j\le9$.
\begin{note}
  \label{eq:P:Ext}
  \begin{enumerate}
    \item
      \begin{enumerate}
        \item $\Ext_\A^{2,48}(P)=\Z/2(\widehat{h}_4h_0h_5)$,
        \item $t_\ast(\widehat{h}_4h_0h_5)=0$ in $\Ext_\A^{3,49}$.
      \end{enumerate}
    \item
      \begin{enumerate}
        \item $\Ext_\A^{3,64}(P)=\Z/2(\widehat{D_3})$,
        \setcounter{enumii}{2}
        \item $t_\ast(\widehat{D_3})=D_3\ne0$ in $\Ext_\A^{4,65}$.
      \end{enumerate}
    \item
      \begin{enumerate}
        \item $\Ext_\A^{2,64}(P)=\Z/2(\widehat{h}_5h_0h_5)$,
        \setcounter{enumii}{2}
        \item $t_\ast(\widehat{h}_5h_0h_5)=h_0h_5^2\ne0$ in $\Ext_\A^{3,65}$.
      \end{enumerate}
    \item
      \begin{enumerate}
        \item $\Ext_\A^{3,191}(P)=\Z/2(\widehat{g}_4)$,
        \setcounter{enumii}{2}
        \item $t_\ast(\widehat{g}_4)=g_4\ne0$ in $\Ext_\A^{4,192}$.
      \end{enumerate}
    \item
      \begin{enumerate}
        \item $\Ext_\A^{4,192}(P)
          =\Z/2(\widehat{h}_4h_0c_4)
          \dsum \Z/2(\widehat{D_3}h_7)
          \dsum \Z/2(\widehat{g}_4h_0)
          \dsum \Z/2(\widehat{h}_7D_3)$,
        \item $t_\ast(\widehat{h}_4h_0c_4)=0$ in $\Ext_\A^{5,193}$,
        \item $t_\ast(\widehat{D_3}h_7)=t_\ast(\widehat{h}_7D_3)=D_3h_7$,
          $t_\ast(\widehat{g}_4h_0)=h_0g_4$, and
          $D_3h_7,h_0g_4$ are linearly independent in $\Ext_\A^{5,193}$.
      \end{enumerate}
    \item
      \begin{enumerate}
        \item $\Ext_\A^{2,191}(P)=\Z/2(\widehat{h}_5h_5h_7)$,
        \setcounter{enumii}{2}
        \item $t_\ast(\widehat{h}_5h_5h_7)=h_5^2h_7=h_6^3\ne0$ in $\Ext_\A^{3,192}$.
      \end{enumerate}
    \item
      \begin{enumerate}
        \item $\Ext_\A^{3,192}(P)=\Z/2(\widehat{h}_5h_0h_5h_7)\dsum\Z/2(\widehat{h}_6h_0h_6^2)$,
        \setcounter{enumii}{2}
        \item $t_\ast(\widehat{h}_5h_0h_5h_7) =t_\ast(\widehat{h}_6h_0h_6^2)=h_0h_6^3\ne0$ in $\Ext_\A^{4,193}$.
      \end{enumerate}
    \item
      \begin{enumerate}
        \item $\Ext_\A^{5,192}(P) =\Z/2(\beta_{187}) \dsum \Z/2(\widehat{h}_4h_0f_3) \dsum \Z/2(\widehat{h}_5V_0)$,
        \item $t_\ast(\widehat{h}_4h_0f_3)=0$ in $\Ext_\A^{6,193}$,
        \item $t_\ast(\widehat{h}_5V_0)=h_5V_0\ne0$ in $\Ext_\A^{6,193}$,
        \item $t_\ast(\beta_{187})=h_5V_0\ne0$ in $\Ext_\A^{6,193}$.
      \end{enumerate}
    \item
      \begin{enumerate}
        \item $\Ext_\A^{5,193}(P)= \Z/2(\widehat{h}_5h_5D_3(1))
          \dsum \Z/2(\widehat{h}_4h_0^2c_4)
          \dsum \Z/2(\widehat{h}_4h_1f_3)
          \dsum \Z/2(\gamma_{61}h_7)
          \dsum \Z/2(\widehat{g}_4h_0^2)$,
        \item $t_\ast(\widehat{h}_4h_0^2c_4)= t_\ast(\widehat{h}_4h_1f_3)= t_\ast(\gamma_{61}h_7)=0$ in $\Ext_\A^{6,194}$,
        \item $t_\ast(\widehat{h}_5h_5D_3(1))=h_5^2D_3(1)$,
          $t_\ast(\widehat{g}_4h_0^2)=h_0^2g_4$, and
          $h_5^2D_3(1),h_0^2g_4$
          are linearly independent in $\Ext_\A^{6,194}$.
      \end{enumerate}
  \end{enumerate}
\end{note}
\noindent
Let $P^{62}\to{i_{62}}P$ be the inclusion map, and
$\Ext_\A^{s,t}(P^{62})\to{(i_{62})_\ast}\Ext_\A^{s,t}(P)$, $s,t\ge0$,
be the induced homomorphisms.
\eqref{eq:P62:Ext}~(i)~($j$) for $1\le j\le 9$ below
describe the groups $\Ext_\A^{s,t}(P^{62})$ for some $s,t\ge0$.
Each of the classes in $\Ext_\A^{s,t}(P^{62})$ in \eqref{eq:P62:Ext}~(i)~($j$)
is a pullback through $(i_{62})_\ast$ of the class of the same name in $\Ext_\A^{s,t}(P)$
given in \eqref{eq:P:Ext}~($j$)~(a) for $1\le j\le 9$, except the following
two classes:
$\bar{e_{62}h_1}h_6^2 \in\Ext_\A^{3,192}(P^{62})$ in \eqref{eq:P62:Ext}~(i)~(7)
is a pullback of $\widehat{h}_6h_0h_6^2 \in\Ext_\A^{3,192}(P)$ in \eqref{eq:P:Ext}~(7), and
$\bar{e_{62}D_3(1)} \in\Ext_\A^{4,192}(P^{62})$ in \eqref{eq:P62:Ext}~(i)~(5)
is a pullback of $\widehat{h}_7D_3 \in\Ext_\A^{4,192}(P)$ in \eqref{eq:P:Ext}~(5).
How to take these pullbacks will be explained in \eqref{eq:P62:Ext:pullback}.
\renewcommand{\labelenumii}{(\arabic{enumii})}
\begin{note}
  \label{eq:P62:Ext}
  \begin{enumerate}[(i)]
    \item 
      \begin{enumerate}
        \item $\Ext_\A^{2,48}(P^{62})=\Z/2(\widehat{h}_4h_0h_5)$,
        \item $\Ext_\A^{3,64}(P^{62})=\Z/2(\widehat{D_3})$,
        \item $\Ext_\A^{2,64}(P^{62})=\Z/2(\widehat{h}_5h_0h_5)$,
        \item $\Ext_\A^{3,191}(P^{62})=\Z/2(\widehat{g}_4)$,
        \item $\Ext_\A^{4,192}(P^{62})
          =\Z/2(\widehat{h}_4h_0c_4)
          \dsum \Z/2(\widehat{D_3}h_7)
          \dsum \Z/2(\widehat{g}_4h_0)
          \dsum \Z/2(\bar{e_{62}D_3(1)})$,
        \item $\Ext_\A^{2,191}(P^{62})=\Z/2(\widehat{h}_5h_5h_7)$,
        \item $\Ext_\A^{3,192}(P^{62})=\Z/2(\widehat{h}_5h_0h_5h_7)\dsum\Z/2(\bar{e_{62}h_1}h_6^2)$,
        \item $\Ext_\A^{5,192}(P^{62})
          =\Z/2(\beta_{187})
          \dsum \Z/2(\widehat{h}_4h_0f_3)
          \dsum \Z/2(\widehat{h}_5V_0)$,
        \item $\Ext_\A^{5,193}(P^{62})=
          \Z/2(\widehat{h}_5h_5D_3(1))
          \dsum \Z/2(\widehat{h}_4h_0^2c_4)
          \dsum \Z/2(\widehat{h}_4h_1f_3)
          \dsum \Z/2(\gamma_{61}h_7)
          \dsum \Z/2(\widehat{g}_4h_0^2)$.
      \end{enumerate}
    \item $\Ext_\A^{s,t}(P^{62})\to[\iso]{(i_{62})_\ast}\Ext_\A^{s,t}(P)$ for all the $(s,t)$ in \eqref{eq:P62:Ext} (i).
    \item $\Ext_\A^{s,t} \to{t_\ast}\Ext_\A^{s+1,t+1}$
      for all $(s,t)$ in \eqref{eq:P62:Ext}~(i)
      are given by
      \eqref{eq:P:Ext}~($j$)~(b) for $j=1,5,8,9$ and
      \eqref{eq:P:Ext}~($j$)~(c) for $2\le j\le9$, except
      $t_\ast(\bar{e_{62}D_3(1)}) =D_3h_7 \ne0$ and
      $t_\ast(\bar{e_{62}h_1}h_6^2) =h_0h_6^3 \ne0$.
  \end{enumerate}
\end{note}
\renewcommand{\labelenumii}{(\alph{enumii})}
\noindent
The cohomology classes in the following \eqref{eq:P62-k:Ext} will be
denoted by their $\bar{E}_2$ representatives in a spectral sequence $\{\bar{E}_r^{i,s,t}\}_{r\ge1}$, whose group of infinite cycles is $\Ext_\A^{s,t+s}(P^l_k =P^l/P^{k-1})$.
This spectral sequence $\{\bar{E}_r^{i,s,t}\}_{r\ge1}$
will be introduced in \eqref{eq:P62-k:Ext:spectral-sequence}.
The $E_2$ term of $\{\bar{E}_r^{i,s,t}\}_{r\ge1}$
is generated over $\Z/2$ by elements
of the form $e_i\alpha$, $\alpha \in\Ext_\A^{s,t+s}$ ($\alpha$ could be $1 \in\Ext_\A^{0,0}=\Z/2$) and
$e_i$ is a $\Z/2$-basis element in $H_\ast(P^l_k)$.
A permanent cycle represented by $e_i\alpha$ is denoted by $\bar{e_i\alpha}$,
and is a non-zero class in $\Ext_\A^{s,t+s+i}(P^l_k)$.
The cohomology classes in \eqref{eq:P62-k:Ext} will be
defined in \eqref{eq:P62-k:Ext:known}.
\begin{note}
  \label{eq:P62-k:Ext}
  \begin{enumerate}
    \item $\Ext_\A^{1,63}(P^{62}_{17})=\Z/2(\bar{e_{31}}h_5)$.
    \item $\Ext_\A^{2,32}(P^{62}_{17})=\Z/2(\bar{e_{23}}h_0h_3)$.
    \item $\Ext_\A^{3,64}(P^{62}_{17})=\Z/2(\bar{e_{23}}h_0h_3h_5)$,\\
      $\Ext_\A^{3,64}(P^{62}) \to{(q_{17})_\ast}\Ext_\A^{3,64}(P^{62}_{17})$ is given by
      $(q_{17})_\ast(\widehat{D_3}) =\bar{e_{23}}h_0h_3h_5$.
    \item $\Ext_\A^{0,47}(P^{62}_{33})=\Z/2(\bar{e_{47}})$.
    \item $\Ext_\A^{2,48}(P^{62}_{33}) =\Z/2(\bar{e_{39}}h_0h_3)$.
    \item $\Ext_\A^{1,31}(P^{16})
      =\Z/2(\bar{e_{15}}h_4)$.
    \item $\Ext_\A^{2,16}(P^{16})
      =\Z/2(\bar{e_7}h_0h_3)$.
    \item $\Ext_\A^{3,32}(P^{16})
      =\Z/2(\bar{e_{15}}h_0h_3^2)$.
    \item $\Ext_\A^{0,31}(P^{32})
      =\Z/2(\bar{e_{31}})$.
    \item $\Ext_\A^{1,47}(P^{32})
      =\Z/2(\bar{e_{31}}h_4)$.
    \item $\Ext_\A^{2,32}(P^{32})
      =\Z/2(\bar{e_{15}}h_0h_4)$.
    \item $\Ext_\A^{3,48}(P^{32})
      =\Z/2(\bar{e_{31}}h_0h_3^2) 
      \dsum\Z/2(\bar{e_{15}}h_0h_4^2)$.
    \item $\Ext_\A^{s,189+s}(P^{62}_{47}) =0$ for $0\le s\le2$.
    \item $\Ext_\A^{2,64}(P^{62}_{47}) =\Z/2(\bar{e_{54}h_1h_3})$.
    \item $\Ext_\A^{3,65}(P^{62}_{47}) =\Z/2(\bar{e_{48}h_0h_3^2})$.
    \item $\Ext_\A^{3,192}(P^{62}_{47})=\Z/2(\bar{e_{62}h_1}h_6^2)$,\\
      $\Ext_\A^{3,192}(P^{62}) \to{(q_{17})_\ast}\Ext_\A^{3,192}(P^{62}_{47})$ is given by
      $(q_{17})_\ast(\bar{e_{62}h_1}h_6^2) =\bar{e_{62}h_1}h_6^2$.
    \item $\Ext_\A^{5,193}(P^{62}_{47})=
      \Z/2(\bar{e_{48}h_0h_3^2}h_6^2)$,\\
      $\Ext_\A^{5,193}(P^{62}) \to{(q_{47})_\ast}\Ext_\A^{5,193}(P^{62}_{47})$ is given by
      $(q_{47})_\ast(\widehat{h}_5h_5D_3(1))
      =(q_{47})_\ast(\widehat{h}_4h_0^2c_4)
      =(q_{47})_\ast(\widehat{h}_4h_1f_3)
      =(q_{47})_\ast(\gamma_{61}h_7)
      =0$, and
      $(q_{47})_\ast(\widehat{g}_4h_0^2)=\bar{e_{48}h_0h_3^2}h_6^2\ne0$.
  \end{enumerate}
\end{note}
\noindent
Recall the complex $\widetilde{P}^{62}$, the map $P^{62}\to{\widetilde{q}}\widetilde{P}^{62}$, and
the induced homomorphism $\Ext_\A^{\ast,\ast}(P^{62})\to{\widetilde{q}_\ast}\Ext_\A^{\ast,\ast}(\widetilde{P}^{62})$ in \eqref{defn:P62t}.
In \eqref{eq:P62t:Ext} below the cohomology classes
are all in the image of $\Ext_\A^{\ast,\ast}(P^{62})\to{\widetilde{q}_\ast}\Ext_\A^{\ast,\ast}(\widetilde{P}^{62})$
except those generated (using the obvious $\Ext_\A^{\ast,\ast}$-action)
by the following four classes:
$\bar{\widetilde{e}_{31}h_5} \in\Ext_\A^{1,63}(\widetilde{P}^{62})$,
$\bar{e_{47}} \in\Ext_\A^{0,47}(\widetilde{P}^{62})$, and
$\bar{e_{39}h_3} \in\Ext_\A^{1,47}(\widetilde{P}^{62})$, and
$\bar{\widetilde{e}_{31}h_4^2} \in\Ext_\A^{2,63}(\widetilde{P}^{62})$.
These four new classes
will be defined in \eqref{eq:P62t:Ext:new} later.
\begin{note}
  \label{eq:P62t:Ext}
  \begin{enumerate}
    \item $\Ext_\A^{0,47}(\widetilde{P}^{62})
      =\Z/2(\bar{e_{47}})$.
    \item $\Ext_\A^{2,48}(\widetilde{P}^{62})
      = \Z/2(\bar{e_{39}h_3}h_0)
      \dsum \Z/2(\widetilde{q}_\ast(\widehat{h}_4h_0h_5))$.
    \item $\Ext_\A^{3,64}(\widetilde{P}^{62})
      = \Z/2(\bar{\widetilde{e}_{31}h_4^2}h_0)
      \dsum \Z/2(\widetilde{q}_\ast(\widehat{D_3}))$,\\
      $\bar{e_{39}h_3}h_0h_4 =\bar{\widetilde{e}_{31}h_4^2}h_0 +\widetilde{q}_\ast(\widehat{D_3})$ in $\Ext_\A^{3,64}(\widetilde{P}^{62})$.
    \item $\Ext_\A^{1,63}(\widetilde{P}^{62})
      =\Z/2(\bar{\widetilde{e}_{31}h_5})
      \dsum \Z/2(\bar{e_{47}}h_4)$.
    \item $\Ext_\A^{2,64}(\widetilde{P}^{62})=\Z/2(\bar{\widetilde{e}_{31}h_5}h_0)\dsum\Z/2(\bar{e_{47}}h_0h_4)$.
    \item $\Ext_\A^{4,192}(\widetilde{P}^{62})
      = Z/2(\widetilde{q}_\ast(\widehat{h}_4h_0c_4))
      \dsum \Z/2(\widetilde{q}_\ast(\widehat{D_3}h_7))
      \dsum \Z/2(\widetilde{q}_\ast(\widehat{g}_4h_0))
      \dsum \Z/2(\widetilde{q}_\ast(\bar{e_{62}D_3(1)}))
      \dsum \Z/2(\bar{\widetilde{e}_{31}h_4^2}h_0h_7)$.
    \item
      We have the following results on $\Ext_\A^{5,193}(\widetilde{P}^{62})$ (see \eqref{eq:P62:Ext}~(i)~(9)):
      \begin{enumerate}
        \item $\widetilde{q}_\ast(\widehat{h}_5h_5D_3(1)) =0$ in $\Ext_\A^{5,193}(\widetilde{P}^{62})$,
        \item $\widetilde{q}_\ast(\widehat{h}_4h_0^2c_4),
          \widetilde{q}_\ast(\widehat{h}_4h_1f_3),
          \widetilde{q}_\ast(\gamma_{61}h_7),
          \widetilde{q}_\ast(\widehat{g}_4h_0^2),
          \bar{\widetilde{e}_{31}h_4^2}h_0^2h_7$ are linearly independent in $\Ext_\A^{5,193}(\widetilde{P}^{62})$.
      \end{enumerate}
    \item $\Ext_\A^{0,189}(\widetilde{P}^{62})=0$.
    \item $\Ext_\A^{1,190}(\widetilde{P}^{62})=0$.
    \item $\Ext_\A^{2,191}(\widetilde{P}^{62})=\Z/2(\bar{\widetilde{e}_{31}h_5}h_7)\dsum \Z/2(\bar{e_{47}}h_4h_7)$.
    \item $\Ext_\A^{3,192}(\widetilde{P}^{62})
      =\Z/2(\bar{\widetilde{e}_{31}h_5}h_0h_7)
      \dsum \Z/2(\bar{e_{47}}h_4h_7h_0)
      \dsum\Z/2(\widetilde{q}_\ast(\bar{e_{62}h_1}h_6^2))$.
  \end{enumerate}
\end{note}
\noindent
For each $k$ with $1\le k\le62$ let $\widetilde{P}^k$ be the $k$-skeleton of $\widetilde{P}^{62}$, and
let $\widetilde{P}^{62}_k=\widetilde{P}^{62}/\widetilde{P}^{k-1}$.
Then $\widetilde{P}^k =P^k$ and $\widetilde{P}^{62}_k=\widetilde{P}^{62}/P^{k}$ for $1\le k\le29$, and
$\widetilde{P}^{62}_k=P^{62}_k$ for $33\le k\le62$.
Let
$\widetilde{P}^{62} \to{q'_k}\widetilde{P}^{62}_k$ be the pinching map.
Then $P^{62} \to{\widetilde{q}}\widetilde{P}^{62}$ induces a map
$P^{62}_k \to{\widetilde{q}_k}\widetilde{P}^{62}_k$ with
the following commutative diagram.
\begin{align*}
  \xymatrix{
  P^{62}\ar[r]^{\widetilde{q}}\ar[d]^{q_k}& \widetilde{P}^{62}\ar[d]^{q_k'}\\
  P^{62}_k\ar[r]^{\widetilde{q}_k}&         \widetilde{P}^{62}_k
  }
\end{align*}
In \eqref{eq:P62-17t:Ext}~(1), (2), (3) below
we describe $\Ext_\A^{s,t}(\widetilde{P}^{62}_{17})$
for $(s,t) =(1,63),(2,32),(3,64)$ respectively.
For the classes in the image of
$\Ext_\A^{\ast,\ast}(P^{62}_{17}) \to{(\widetilde{q}_{17})_\ast}\Ext_\A^{\ast,\ast}(\widetilde{P}^{62}_{17})$
we denote them as such
in \eqref{eq:P62-17t:Ext}~(1) through \eqref{eq:P62-17t:Ext}~(3).
The classes generated by 
$\bar{\widetilde{e}_{31}} \in\Ext_\A^{0,31}(\widetilde{P}^{62}_{17})$
in \eqref{eq:P62-17t:Ext}~(1) and \eqref{eq:P62-17t:Ext}~(3),
and by $\bar{e_{47}} \in\Ext_\A^{0,47}(\widetilde{P}^{62}_{17})$
in \eqref{eq:P62-17t:Ext}~(1),
are not in the image of $(\widetilde{q}_{17})_\ast$.
These classes will be defined in \eqref{eq:P62-17t:Ext:new} later.
We also describe the map 
$\Ext_\A^{\ast,\ast}(\widetilde{P}^{62}) \to{(q_{17}')_\ast}\Ext_\A^{\ast,\ast}(\widetilde{P}^{62}_{17})$
at the dimensions specified in
\eqref{eq:P62-17t:Ext}~(1) through \eqref{eq:P62-17t:Ext}~(3).
For $k\ge33$, $\widetilde{P}^{62}_k= P^{62}_k$.
In \eqref{eq:P62-17t:Ext}~(4), (5), (6) below
we describe the map
$\Ext_\A^{\ast,\ast}(\widetilde{P}^{62}) \to{(q_k')_\ast}\Ext_\A^{\ast,\ast}(\widetilde{P}^{62}_k =P^{62}_k)$
for $(s,t,k) =(2,48,33),(0,47,33),(2,64,47)$ respectively.
We refer to
\eqref{eq:P62-k:Ext} and \eqref{eq:P62t:Ext} for the classes
in these descriptions.
\begin{note}
  \label{eq:P62-17t:Ext}
  \begin{enumerate}
    \item $\Ext_\A^{1,63}(\widetilde{P}^{62}_{17})
      =\Z/2(\bar{\widetilde{e}_{31}}h_5)
      \dsum \Z/2(\bar{e_{47}}h_4)$,\\
      $\Ext_\A^{1,63}(\widetilde{P}^{62}) \to{(q_{17}')_\ast}\Ext_\A^{1,63}(\widetilde{P}^{62}_{17})$ is given by
      $(q_{17}')_\ast(\bar{\widetilde{e}_{31}h_5}) =\bar{\widetilde{e}_{31}}h_5$ and
      $(q_{17}')_\ast(\bar{e_{47}}h_4) =\bar{e_{47}}h_4$.
    \item $\Ext_\A^{2,32}(\widetilde{P}^{62}_{17})=\Z/2( (\widetilde{q}_{17})_\ast(\bar{e_{23}}h_0h_3))$.
    \item $\Ext_\A^{3,64}(\widetilde{P}^{62}_{17})
      =\Z/2(\bar{\widetilde{e}_{31}}h_0h_4^2)
      \dsum \Z/2( (\widetilde{q}_{17})_\ast(\bar{e_{23}}h_0h_3h_5))$,\\
      $\Ext_\A^{3,64}(\widetilde{P}^{62}) \to{(q_{17}')_\ast}\Ext_\A^{3,64}(\widetilde{P}^{62}_{17})$ is given by
      $(q_{17}')_\ast(\bar{\widetilde{e}_{31}h_4^2}h_0) =\bar{\widetilde{e}_{31}}h_0h_4^2$, and\\
      $(q_{17}')_\ast(\widetilde{q}_\ast(\widehat{D_3})) =(\widetilde{q}_{17})_\ast(\bar{e_{23}}h_0h_3h_5)$.
    \item 
      $\Ext_\A^{2,48}(\widetilde{P}^{62}) \to{(q_{33}')_\ast}\Ext_\A^{2,48}(\widetilde{P}^{62}_{33} =P^{62}_{33})$
      is given by
      $(q_{33}')_\ast(\bar{e_{39}h_3}h_0) =\bar{e_{39}h_3}h_0$, and
      $(q_{33}')_\ast(\widetilde{q}_\ast(\widehat{h}_4h_0h_5)) =0$.
    \item
      $\Ext_\A^{0,47}(\widetilde{P}^{62}) \to{(q_{33}')_\ast}\Ext_\A^{0,47}(\widetilde{P}^{62}_{33} =P^{62}_{33})$
      is given by
      $(q_{33}')_\ast(\bar{e_{47}}) =\bar{e_{47}}$.
    \item $\Ext_\A^{2,64}(\widetilde{P}^{62}) \to{(q_{47}')_\ast}\Ext_\A^{2,64}(\widetilde{P}^{62}_{47} =P^{62}_{47})$ is given by
      $(q_{47}')_\ast(\bar{\widetilde{e}_{31}h_5}h_5)=(q_{47}')_\ast(\bar{e_{47}}h_0h_4)=0$.
  \end{enumerate}
\end{note}
\noindent
This completes the statements of all the $\Ext$ groups results in this section.
\begin{remark}
  \label{remark:P:Ext}
  \upshape
  As will be seen later,
  the results in \eqref{eq:P:Ext}~(8)~(c), \eqref{eq:P:Ext}~(9)~(c)
  on $\Ext_\A^{\ast,\ast}(P)$ imply the results in
  \eqref{eq:Ext} (7), (8) respectively.
\end{remark}
In the first half of this section we will describe and prove
\eqref{eq:Ext} and \eqref{eq:P:Ext}, and
in the second half we will describe and prove
\eqref{eq:P62:Ext} through \eqref{eq:P62-17t:Ext}.

\eqref{eq:Ext}~(1) through \eqref{eq:Ext}~(6) follow from some
well-known results on $\Ext_\A^{\ast,\ast}$ obtained in
\cite{adams_structure_1958,lin_ext_2008,chen_determination_2011}.
These known results are summarized in \autoref{thm:Ext:known} below.
In order to describe the cohomology classes in \eqref{eq:Ext} and
to state \autoref{thm:Ext:known}, we need to recall our main tool
for calculating the $\Ext$ groups in \eqref{eq:Ext},
which is the mod $2$ lambda algebra $\Lambda$ (\cite{bousfield_rm_1966}).
Recall that $\Lambda$ is the bigraded differential algebra generated by $\lambda_i\in\Lambda^{1,i}$ for $i\ge0$
with relation
\begin{align}
  \label{eq:Lambda:Adem-relation}
  \lambda_i\lambda_{2i+1+m}
  &=\Sum_{\nu\ge0}{m-1-\nu \choose \nu}\lambda_{m+i-\nu}\lambda_{2i+1+\nu}
  \quad \mbox{for $m\ge0$.}
\end{align}
The differential $\Lambda\to{\delta}\Lambda$ is given by
\begin{align}
  \label{eq:Lambda:diff}
  \delta(\lambda_i)
  &=\Sum_{\nu\ge0}{i-1-\nu \choose \nu+1}\lambda_{i-1-\nu}\lambda_\nu
  \quad \mbox{on the generators $\lambda_i$.}
\end{align}
Then
\begin{note}
  \label{eq:Lambda:H}
  $H^{s,t}(\Lambda)=H^{s,t}(\Lambda,\delta)=\Ext_\A^{s,t+s}$.
\end{note}
\noindent
A monomial $\lambda_{i_1}\cdots\lambda_{i_s}\in\Lambda$ will simply be denoted by
$\lambda_I$ where $I=(i_1,\cdots,i_s)$ is a sequence of non-negative integers.
\begin{note}
  \label{eq:admissible-basis}
  From \eqref{eq:Lambda:Adem-relation} we see that the set of monomials
  $\lambda_I=\lambda_{i_1}\cdots\lambda_{i_s}$ with $2i_j\ge i_{j+1}$, $1\le j\le s-1$, $s\ge0$ ($\lambda_I=1$ if $s=0$),
  is a basis for $\Lambda$.
  This basis is called the admissible basis and
  the monomials in this basis are said to be admissible monomials.
\end{note}
\noindent
There is an operation $\Lambda\to{Sq^0}\Lambda$ given by
\begin{align}
  \label{eq:Sq0}
  Sq^0(\lambda_{j_1}\cdots\lambda_{j_s})
  &=\lambda_{2j_1+1}\cdots\lambda_{2j_s+1},
\end{align}
where $\lambda_{j_1}\cdots\lambda_{j_s}$ is not necessarily admissible.
This operation is compatible with the relations
\eqref{eq:Lambda:Adem-relation} and commutes with the differential
\eqref{eq:Lambda:diff}, so it induces an $\Ext$ groups map
\begin{align*}
  H^{s,t}(\Lambda)=\Ext_\A^{s,t+s}\to{Sq^0} H^{s,2t+s}(\Lambda)=\Ext_\A^{s,2t+2s}.
\end{align*}
Let
$\Ext_\A^{\ast,\ast}\to{(Sq^0)^i}\Ext_\A^{\ast,\ast}$
denote the composite $\underbrace{Sq^0\cdots Sq^0}_i$ for $i\ge1$, and let
$\Ext_\A^{\ast,\ast}\to{(Sq^0)^0}\Ext_\A^{\ast,\ast}$ be the identity map.

For the following cohomology classes we refer to
\cite{lin_ext_2008,chen_determination_2011}.
\begin{note}
  \label{eq:Ext:known}
  \begin{enumerate}
    \item $h_i=\{\lambda_{2^i-1}=(Sq^0)^i(\lambda_0)\}\in\Ext_\A^{1,2^i}$, $i\ge0$.
    \item $c_i=\{(Sq^0)^i(c_0^\ast =\lambda_2\lambda_3^2)\}\in\Ext_\A^{3,2^{i+3}+2^{i+1}+2^i}$, $i\ge0$.
    \item $d_i=\{(Sq^0)^i(d_0^\ast =\lambda_6\lambda_2\lambda_3^2 +\lambda_4^2\lambda_3^2 +\lambda_2\lambda_4\lambda_5\lambda_3)\}\in\Ext_\A^{4,2^{i+4}+2^{i+1}}$, $i\ge0$.
    \item $e_i=\{(Sq^0)^i(\lambda_8\lambda_3^3 +\lambda_4(\lambda_5^2\lambda_3 +\lambda_7\lambda_3^2) +\lambda_2(\lambda_9\lambda_3^2 +\lambda_3^2\lambda_9))\}\in\Ext_\A^{4,2^{i+4}+2^{i+2}+2^i}$, $i\ge0$.
    \item $f_i=\{(Sq^0)^i(\lambda_4\lambda_0\lambda_7^2 +\lambda_3(\lambda_9\lambda_3^2 +\lambda_3\lambda_5\lambda_7) +\lambda_2\lambda_4\lambda_5\lambda_7)\}\in\Ext_\A^{4,2^{i+4}+2^{i+2}+2^{i+1}}$, $i\ge0$.
    \item $g_{i+1}=\{(Sq^0)^i(\lambda_6\lambda_0\lambda_7^2 +\lambda_5(\lambda_9\lambda_3^2 +\lambda_3^2\lambda_9) +\lambda_3(\lambda_{11}\lambda_3^2 +\lambda_5\lambda_9\lambda_3))\}\in\Ext_\A^{4,2^{i+4}+2^{i+3}}$, $i\ge0$.
    \item $p_i=\{(Sq^0)^i(\lambda_0(\lambda_{19}\lambda_7^2 +\lambda_7^2\lambda_{19}))\}\in\Ext_\A^{4,2^{i+5}+2^{i+2}+2^i}$, $i\ge0$.
    \item $D_3(i)=\{(Sq^0)^i(\lambda_0\lambda_{23}\lambda_7\lambda_{31})\}\in\Ext_\A^{4,2^{i+6}+2^i}$, $i\ge0$.
    \item $p'_i=\{(Sq^0)^i(\lambda_0(\lambda_{39}\lambda_{15}^2 +\lambda_{15}^2\lambda_{39}))\}\in\Ext_\A^{4,2^{i+6}+2^{i+3}+2^i}$, $i\ge0$.
    \item $P^1h_1=\{\lambda_2\lambda_0^3\lambda_7 +\lambda_1(\lambda_2\lambda_4\lambda_1^2 +\lambda_1\lambda_2\lambda_4\lambda_1 +\lambda_2\lambda_1^2\lambda_4)\}\in\Ext_\A^{5,14}$.
    \item $P^1h_2=\{\lambda_4\lambda_0^3\lambda_7 +\lambda_3(\lambda_5\lambda_1^3 +\lambda_2\lambda_4\lambda_1^2 +\lambda_1\lambda_2\lambda_4\lambda_1 +\lambda_2\lambda_1^2\lambda_4) +\lambda_2^2\lambda_0^2\lambda_7 +\lambda_1^2\lambda_2\lambda_0\lambda_7\}\in\Ext_\A^{5,16}$.
    \item $n_i=\{(Sq^0)^i(\lambda_9\lambda_3\lambda_5\lambda_7^2 +\lambda_6\lambda_0\lambda_3\lambda_{15}\lambda_7 +\lambda_3\lambda_5\lambda_1\lambda_{15}\lambda_7)\}\in\Ext_\A^{5,2^{i+5}+2^{i+2}}$, $i\ge0$.
    \item $x_i=\{(Sq^0)^i(\lambda_{14}\lambda_2\lambda_3^2\lambda_{15} +\lambda_{12}\lambda_4\lambda_7^3 +\lambda_8^2\lambda_7^3 +\lambda_6\lambda_2\lambda_3^2\lambda_{23} +\lambda_4^2\lambda_3^2\lambda_{23} +\lambda_2\lambda_4\lambda_1\lambda_{15}^2)\}\in\Ext_\A^{5,2^{i+5}+2^{i+3}+2^{i+1}}$, $i\ge0$.
    \item $D_1(i)=\{(Sq^0)^i(\lambda_4\lambda_7\lambda_{11}\lambda_{15}^2)\}\in\Ext_\A^{5,2^{i+5}+2^{i+4}+2^{i+3}+2^i}$, $i\ge0$.
    \item $H_1(i)=\{(Sq^0)^i(\lambda_{14}\lambda_7\lambda_{11}\lambda_{15}^2 +\lambda_{10}\lambda_{11}^2\lambda_{15}^2 +\lambda_8\lambda_1\lambda_7\lambda_{31}\lambda_{15} +\lambda_6\lambda_7\lambda_3\lambda_{31}\lambda_{15} +\lambda_4\lambda_5\lambda_7\lambda_{31}\lambda_{15})\}\in\Ext_\A^{5,2^{i+6}+2^{i+1}+2^i}$, $i\ge0$.
    \item $Q_3(i)=\{(Sq^0)^i( (\lambda_6\lambda_0\lambda_7^2+\lambda_5(\lambda_9\lambda_3^2 +\lambda_3^2\lambda_9) +\lambda_3\lambda_5\lambda_9\lambda_3 +\lambda_3\lambda_{11}\lambda_3^2)\lambda_{47} +\lambda_5(\lambda_{21}\lambda_{11}\lambda_{15}^2 +\lambda_9\lambda_{15}\lambda_7\lambda_{31}) +\lambda_3\lambda_{23}\lambda_{11}\lambda_{15}^2)\}\in\Ext_\A^{5,2^{i+6}+2^{i+3}}$, $i\ge0$.
    \item $K_i=\{(Sq^0)^i(\lambda_0^2\lambda_{47}\lambda_{15}\lambda_{63} +(\lambda_9\lambda_3^2 +\lambda_3^2\lambda_9)\lambda_{47}\lambda_{63} +\lambda_{21}\lambda_{11}\lambda_{31}^3)\}\in\Ext_\A^{5,2^{i+7}+2^{i+1}}$, $i\ge0$.
    \item $J_i=\{(Sq^0)^i(\lambda_0(\lambda_{19}\lambda_7^2 +\lambda_7^2\lambda_{19})\lambda_{95} +\lambda_0\lambda_{43}\lambda_{23}\lambda_{31}^2 +\lambda_0\lambda_{19}\lambda_{31}\lambda_{15}\lambda_{63})\}\in\Ext_\A^{5,2^{i+7}+2^{i+2}+2^i}$, $i\ge0$.
    \item $T_i=\{(Sq^0)^i(\lambda_0^2(\lambda_{79}\lambda_{31}^2 +\lambda_{31}^2\lambda_{79}) +(\lambda_9\lambda_3^2 +\lambda_3^2\lambda_9)\lambda_{63}^2)\}\in\Ext_\A^{5,2^{i+7}+2^{i+4}+2^{i+1}}$, $i\ge0$.
    \item $V_i=\{(Sq^0)^i(\lambda_0\lambda_{31}\lambda_{47}\lambda_{15}\lambda_{63})\}\in\Ext_\A^{5,2^{i+7}+2^{i+5}+2^i}$, $i\ge0$.
    \item $V'_i=\{(Sq^0)^i(\lambda_0\lambda_{23}\lambda_7\lambda_{31}\lambda_{191} +\lambda_0\lambda_{47}\lambda_{15}\lambda_{127}\lambda_{63})\}\in\Ext_\A^{5,2^{i+8}+2^i}$, $i\ge0$.
    \item $U_i=\{(Sq^0)^i(\lambda_0(\lambda_{39}\lambda_{15}^2 +\lambda_{15}^2\lambda_{39})\lambda_{191} +\lambda_0\lambda_{87}\lambda_{47}\lambda_{63}^2 +\lambda_0\lambda_{39}\lambda_{63}\lambda_{31}\lambda_{127})\} \in\Ext_\A^{5,2^{i+8}+2^{i+3}+2^i}$, $i\ge0$.
  \end{enumerate}
\end{note}
\noindent
\begin{thm}[\cite{adams_structure_1958,lin_ext_2008,chen_determination_2011}]\label{thm:Ext:known}
  \begin{enumerate}[(i)]
    \item The algebra $\Ext_\A^{s,\ast}$ for $s\le 3$ is generated by $h_i\ne0$ and $c_i\ne0$ for $i\ge0$, and
      subject only to the relations $h_ih_{i+1}=0$, $h_ih_{i+2}^2=0$, and $h_{i+1}^3=h_i^2h_{i+2}$ for $i\ge0$.
    \item The subalgebra $E'$ of the algebra $\Ext_\A^{s,\ast}$ for $s\le 4$ generated by
      $h_i$ and $c_i$ for $i\ge0$ is subject only to the relations in (i) together with the relations
      $h_i^2h_{i+3}^2=0$, and $h_jc_i=0$ for $j=i-1,i,i+2$ and $i+3$.
    \item The set of classes $d_i,e_i,f_i,g_{i+1},p_i,D_3(i)$ and $p'_i$ for $i\ge0$ is a
      $\Z/2$-basis for the indecomposable elements in $\Ext_\A^{4,\ast}$.
    \item The subalgebra $E''$ of the algebra $\Ext_\A^{s,\ast}$ for $s\le 5$ generated by 
      $h_1,c_i,d_i,e_i,f_i,g_{i+1},p_i,D_3(i)$ and $p'_i$ for $i\ge0$ is subject only to the relations
      in (i) and (iii) together with the following relations (1) through (28), where $j\ge0$.
      \begin{multicols}{2}
      \begin{enumerate}[(1)]
        \item $h_i^2c_j=0$ for $i=j+1,j+4$,
        \item $h_jh_{j+3}c_{j+2}=0$,
        \item $h_jd_{j+1}=0$,
        \item $h_id_j=0$ for $i=j+3,j+4$,
        \item $h_je_{j+1}=0$,
        \item $h_{j+4}e_j=0$,
        \item $h_if_j=0$ for $i=j+1,j+3,j+4$,
        \item $h_{j+3}g_{j+1}=0$,
        \item $h_jp_{j+1}=0$,
        \item $h_ip_j=0$ for $i=j+1,j+2,j+4,j+5$,
        \item $h_jD_3(j+1)=0$,
        \item $h_iD_3(j)=0$ for $i=j,j+5,j+6$,
        \item $h_jp'_{j+1}=0$,
        \item $h_ip'_j=0$ for $i=j+2,j+3,j+6$,
        \item $h_{j+1}h_{j+4}c_j=h_{j+3}e_j$,
        \item $h_jh_{j+4}c_{j+3}=h_{j+5}p'_j$,
        \item $h_{j+5}^2c_j=h_{j+1}p'_j$,
        \item $h_jd_{j+2}=h_{j+3}D_3(j)$,
        \item $h_{j+2}d_{j+1}=h_{j+4}g_{j+1}$,
        \item $h_{j+1}d_{j+1}=h_jp_j$,
        \item $h_{j+2}d_j=h_je_j$,
        \item $h_{j+1}e_j=h_jf_j$,
        \item $h_{j+2}e_{j+1}=h_{j+1}f_{j+1}=h_j^2c_{j+2}$,
        \item $h_{j+2}e_j=h_jg_{j+1}$,
        \item $h_jf_{j+2}=h_{j+4}p'_j$,
        \item $h_jf_{j+1}=h_{j+3}p_j$,
        \item $h_{j+2}f_j=h_{j+1}g_{j+1}$,
        \item $h_{j+3}g_{j+2}=h_{j+5}g_{j+1}$.\\  
      \end{enumerate}
      \end{multicols}
    \item The set of classes $P^1h_1,P^1h_2,n_i,x_i,D_1(i),H_1(i),Q_3(i),K_i,J_i,T_i,V_i,V'_i$ and $U_i$ for $i\ge0$
      is a $\Z/2$-basis for the indecomposable elements in $\Ext_\A^{5,\ast}$.
  \end{enumerate}
\end{thm}
\noindent
It is not difficult to see that
\autoref{thm:Ext:known} covers \eqref{eq:Ext}~(1) through \eqref{eq:Ext}~(6).

\eqref{eq:P:Ext}~(1) through \eqref{eq:P:Ext}~(7) also follow from some
``already known'' results on $\Ext_\A^{s,\ast}(P)$ for $s\le4$
in \cite{lin_ext_2008,chen_determination_2011}.
These results will be summarized in \autoref{thm:P:Ext:known} in a moment.
To state these results, recall (\cite{cohen_adams_1988}) that
the lambda algebra $\Lambda$ can also be used to calculate
the $\Ext$ groups $\Ext_\A^{\ast,\ast}(X)$ for any stable complex $X$,
in particular for $X=P$. This is described as follows.

Given a stable complex $X$,
the mod $2$ Steenrod algebra $\A$ acts on
the mod $2$ reduced homology group
$ H_\ast(X)$ from the right.
Consider the 
bigraded differential module $ H_\ast(X)\tensor\Lambda$
with bigrading 
$( H_\ast(X)\tensor\Lambda)^{s,t} =\Sum_k  H_k(X)\tensor\Lambda^{s,t-k}$.
For $\alpha\in H_k(X)$, we write
$\alpha\lambda_I$ to denote $\alpha\tensor\lambda_I\in H_\ast(X)\tensor\Lambda$,
and write $\alpha$ to denote $\alpha1=\alpha\tensor1\in H_\ast(X)\tensor\Lambda$.
Then we have the following.
\begin{note}
  \label{eq:XLambda}
  \begin{enumerate}
    \item $ H_\ast(X)\tensor\Lambda$ is a differential $\Lambda$-module with
      differential $ H_\ast(X)\tensor\Lambda\to{\delta} H_\ast(X)\tensor\Lambda$ given by
      \begin{align*}
        \delta(\alpha\lambda_I)
        &=\alpha\delta(\lambda_I) +\Sum_{j\ge0}\alpha Sq^{j+1}\lambda_j\lambda_I.
      \end{align*}
    \item $H^{s,t}( H_\ast(X)\tensor\Lambda) =H^{s,t}( H_\ast(X)\tensor\Lambda,\delta) =\Ext_\A^{s,t+s}(X) =\Ext_\A^{s,t+s}( H_\ast(X),\Z/2)$.
  \end{enumerate}
\end{note}
\noindent
Note that since $ H_\ast(X)\tensor\Lambda$ is a differential $\Lambda$-module,
$\Ext_\A^{\ast,\ast}(X)$ is an $\Ext_\A^{\ast,\ast}$-module
in the obvious way.

To calculate $\Ext_\A^{\ast,\ast}(P)$, when $X=P$,
let $e_k$ denote the generator of $ H_k(P)\iso\Z/2$ for $k\ge1$.
The mod $2$ Steenrod algebra $\A$ acts on $e_k\in H_k(P)$ from the right by
\begin{align}
  \label{eq:P:e_k}
  e_kSq^{\nu+1}&={k-\nu-1 \choose \nu+1}e_{k-\nu-1}.
\end{align}
From \eqref{eq:XLambda} (1) and \eqref{eq:P:e_k} we thus see that the
differential $ H_\ast(P)\tensor\Lambda\to{\delta} H_\ast(P)\tensor\Lambda$ of
the differential $\Lambda$-module $ H_\ast(P)\tensor\Lambda=\dsum_{k\ge1}e_k\Lambda$ is given by
\begin{align}
  \label{eq:P:diff}
  \delta(e_k\lambda_I)
  &=e_k\delta(\lambda_I) +\Sum_{\nu\ge0}
  {k-1-\nu \choose \nu+1}e_{k-\nu-1}\lambda_\nu\lambda_I.
\end{align}
And $H^{s,t}( H_\ast(P)\tensor\Lambda,\delta)=\Ext_\A^{s,t+s}(P)$.

Analogous to the operation $\Lambda\to{Sq^0}\Lambda$ in \eqref{eq:Sq0},
there is also an operation
$ H_\ast(P)\tensor\Lambda\to{Sq^0} H_\ast(P)\tensor\Lambda$
given by
\begin{align}
  \label{eq:P:Sq0}
  Sq^0(e_k\lambda_{j_1}\cdots\lambda_{j_s})
  &=e_{2k+1}\lambda_{2j_1+1}\cdots\lambda_{2j_s+1}.
\end{align}
This operation is compatible with the relations
\eqref{eq:Lambda:Adem-relation} and commutes with the differential
\eqref{eq:P:diff}, so it induces a map
\begin{align*}
  H^{s,t}( H_\ast(P)\tensor\Lambda)=\Ext_\A^{s,t+s}(P)
  \to{Sq^0} H^{s,2t+s+1}( H_\ast(P)\tensor\Lambda)=\Ext_\A^{s,2t+2s+1}(P).
\end{align*}
We also let
$\Ext_\A^{\ast,\ast}(P)\to{(Sq^0)^i}\Ext_\A^{\ast,\ast}(P)$
denote the composite $\underbrace{Sq^0\cdots Sq^0}_i$ for $i\ge1$, and let
$\Ext_\A^{\ast,\ast}(P)\to{(Sq^0)^0}\Ext_\A^{\ast,\ast}(P)$
be the identity map.

Recall the transfer map $P\to{t}S^0$ in \eqref{defn:t}. It is not difficult to
show the following.
\begin{note}
  \label{eq:t-bar}
  $P\to{t}S^0$ induces a differential $\Lambda$-map
  $ H_\ast(P)\tensor\Lambda\to{\bar{t}}\Lambda$ given by
  $\bar{t}(e_k\lambda_I)=\lambda_k\lambda_I$
  for any $I=(i_1,\cdots,i_s)$. It gives rise to an
  $\Ext$ groups map
  \begin{align*}
    \Ext_\A^{s,t}(P)\to{t_\ast}\Ext_\A^{s+1,t+1},\quad
    \mbox{for all $s\ge0,t\ge0$.}
  \end{align*}
\end{note}
\noindent
From \eqref{eq:Sq0}, \eqref{eq:P:Sq0} and \eqref{eq:t-bar}
one easily sees the following.
\begin{note}
  \label{eq:t-bar-Sq0-commute}
  There is a commutative diagram
  \begin{align*}
    \xymatrix{
    \Ext_\A^{s,t}(P)\ar[r]^{t_\ast}\ar[d]^{(Sq^0)^i}& \Ext_\A^{s+1,t+1}\ar[d]^{(Sq^0)^i}\\
    \Ext_\A^{s,2^it+2^i-1}(P)\ar[r]^{t_\ast}&         \Ext_\A^{s+1,2^i(t+1)}
    }
  \end{align*}
  for all $i\ge1$.
\end{note}

A class $\alpha\in\Ext_\A^{s,t}(P)$ is a
decomposable element in the $\Ext_\A^{\ast,\ast}$-module $\Ext_\A^{\ast,\ast}(P)$
if $\alpha=\Sum_i a_i\alpha_i$ for some $a_i\in\Ext_\A^{s_i,t_i}$ with
$s_i>0$ and $\alpha_i\in\Ext_\A^{s_i',t_i'}(P)$.
Otherwise $\alpha$ is an indecomposable element.
In \eqref{eq:P:Ext:all-known} below we recall from
\cite{lin_ext_2008,chen_determination_2011}
all the indecomposable elements in $\Ext_\A^{s,\ast}(P)$ for $s\le4$.
Only the classes $\widehat{h}_i,\widehat{g}_{i+1},\widehat{D_3}(i),\gamma_{61}(i)$
in \eqref{eq:P:Ext:all-known} (i) (1), \eqref{eq:P:Ext:all-known} (ii) (4),
\eqref{eq:P:Ext:all-known} (ii) (6), \eqref{eq:P:Ext:all-known} (iii) (11)
for $i\ge0$ will be explicitly described in terms of cycle representations
in $ H_\ast(P)\tensor\Lambda$
as in \eqref{eq:P:Ext:known} below.
The cycle representations
of other classes will not be recorded here because they will not be
used in this paper.
\begin{note}
  \label{eq:P:Ext:known}
  \begin{enumerate}
    \item $\widehat{h}_i=(Sq^0)^{i-1}(\widehat{h}_1)=\{e_{2^i-1}=(Sq^0)^{i-1}(e_1)\}\in\Ext_\A^{0,2^i-1}(P)$, $i\ge1$.
    \item $\widehat{g}_{i+1}=(Sq^0)^i(\widehat{g}_1)=\{(Sq^0)^i(e_6\lambda_0\lambda_7^2 +e_5(\lambda_9\lambda_3^2 +\lambda_3^2\lambda_9) +e_3\lambda_{11}\lambda_3^2 +e_3\lambda_5\lambda_9\lambda_3)\}\in\Ext_\A^{3,2^{i+4}+2^{i+3}-1}(P)$, $i\ge0$.
    \item $\widehat{D_3}(i)=(Sq^0)^i(\widehat{D_3})=\{(Sq^0)^i(e_{22}\lambda_1\lambda_7\lambda_{31} +e_{16}\lambda_{15}^3 +e_{14}\lambda_9\lambda_7\lambda_{31})\}\in\Ext_\A^{3,2^{i+6}+2^i-1}(P)$, $i\ge0$.
    \item $\gamma_{61}(i) =(Sq^0)^i(\gamma_{61} =\{\gamma_{61}^\ast\}) \in\Ext_\A^{4,2^{i+6}+2^{i+1}-1}(P)$,
      $i\ge0$, where
      \begin{align*}
      \gamma_{61}^\ast
      &=
      \begin{minipage}[t]{11cm} 
      $e_{29}(\lambda_{13}\lambda_5\lambda_7^2 +\lambda_9^2\lambda_7^2 +\lambda_5^2\lambda_{15}\lambda_7)
      +e_{28}\lambda_0(\lambda_{19}\lambda_7^2 +\lambda_7^2\lambda_{19})
      +e_{27}\lambda_9\lambda_3\lambda_{15}\lambda_7
      +e_{26}\lambda_2(\lambda_{19}\lambda_7^2 +\lambda_7^2\lambda_{19})
      +e_{26}\lambda_4\lambda_1\lambda_{15}^2
      +e_{25}\lambda_5\lambda_1\lambda_{15}^2
      +e_{23}\lambda_0^3\lambda_{15}\lambda_{23} +e_{23}\lambda_{17}\lambda_7^3 +e_{23}\lambda_9\lambda_{15}\lambda_7^2
      +e_{23}\lambda_7\lambda_9\lambda_{15}\lambda_7 +e_{23}\lambda_3\lambda_5\lambda_{15}^2
      +e_{22}\lambda_6(\lambda_{19}\lambda_7^2 +\lambda_7^2\lambda_{19})
      +e_{22}\lambda_8\lambda_1\lambda_{15}^2 +e_{22}\lambda_2\lambda_3^2\lambda_{31}
      +e_{21}\lambda_7(\lambda_{19}\lambda_7^2 +\lambda_7^2\lambda_{19})
      +e_{21}(\lambda_{19}\lambda_7^2 +\lambda_7^2\lambda_{19})\lambda_7
      +e_{20}\lambda_0\lambda_{11}\lambda_{15}^2
      +e_{19}\lambda_0^2\lambda_{11}\lambda_{31}
      +e_{19}\lambda_5\lambda_{15}^2\lambda_7 +e_{19}\lambda_7\lambda_5\lambda_{15}^2 +e_{19}\lambda_9\lambda_7\lambda_{11}\lambda_{15}
      +e_{18}\lambda_2\lambda_{11}\lambda_{15}^2
      +e_{15}\lambda_1\lambda_{15}\lambda_{23}\lambda_7 +e_{15}\lambda_8\lambda_0\lambda_{15}\lambda_{23} +e_{15}\lambda_0\lambda_8\lambda_{15}\lambda_{23} +e_{15}\lambda_9\lambda_{15}^2\lambda_7 +e_{15}\lambda_{11}\lambda_{17}\lambda_3\lambda_{15} + e_{15}\lambda_{11}\lambda_{13}\lambda_{15}\lambda_7
      +e_{15}\lambda_{11}\lambda_9\lambda_{11}\lambda_{15}
      +e_{14}\lambda_{10}\lambda_3^2\lambda_{31} +e_{14}\lambda_8\lambda_1\lambda_7\lambda_{31} +e_{14}\lambda_5\lambda_0\lambda_{11}\lambda_{31} +e_{14}\lambda_4\lambda_5\lambda_7\lambda_{31} +e_{14}\lambda_3\lambda_2\lambda_{11}\lambda_{31}
      +e_{14}\lambda_6\lambda_{11}\lambda_{15}^2 +e_{14}\lambda_{14}(\lambda_{19}\lambda_7^2 +\lambda_7^2\lambda_{19}) +e_{14}\lambda_{16}\lambda_1\lambda_{15}^2 +e_{14}\lambda_2\lambda_{15}^3
      +e_{13}\lambda_{15}(\lambda_{19}\lambda_7^2 +\lambda_7^2\lambda_{19})
      +e_{13}\lambda_3\lambda_{27}\lambda_3\lambda_{15} +e_{13}\lambda_3^2\lambda_{35}\lambda_7 +e_{13}\lambda_{15}\lambda_{11}\lambda_{15}\lambda_7 +e_{13}\lambda_0\lambda_6\lambda_{11}\lambda_{31}
      +e_{12}\lambda_8\lambda_{11}\lambda_{15}^2
      +e_{11}\lambda_8\lambda_0\lambda_{11}\lambda_{31} +e_{11}\lambda_0\lambda_8\lambda_{11}\lambda_{31} +e_{11}\lambda_1^3\lambda_{47} +e_{11}\lambda_{13}\lambda_{15}^2\lambda_7 +e_{11}\lambda_{15}\lambda_5\lambda_{15}^2 +e_{11}\lambda_{17}(\lambda_{19}\lambda_7^2 +\lambda_7^2\lambda_{19})
      +e_{11}\lambda_0^2\lambda_{19}\lambda_{31}
      +e_{10}\lambda_{10}\lambda_{11}\lambda_{15}^2
      +e_7\lambda_6^2\lambda_{11}\lambda_{31} +e_7\lambda_4\lambda_0\lambda_{19}\lambda_{31} +e_7\lambda_0\lambda_4\lambda_{19}\lambda_{31} +e_7\lambda_4\lambda_8\lambda_{11}\lambda_{31}
      +e_7\lambda_{11}\lambda_5\lambda_7\lambda_{31} +e_7\lambda_7\lambda_9\lambda_7\lambda_{31} +e_7\lambda_{17}\lambda_3^2\lambda_{31} +e_7\lambda_9\lambda_{11}\lambda_3\lambda_{31} +e_7\lambda_5\lambda_1^2\lambda_{47} +e_7\lambda_{13}\lambda_{15}\lambda_{11}\lambda_{15}
      +e_7\lambda_{11}\lambda_{13}\lambda_{15}^2 +e_7\lambda_{17}\lambda_3\lambda_{19}\lambda_{15} +e_7\lambda_9\lambda_{11}\lambda_{19}\lambda_{15}
      +e_6\lambda_2\lambda_3^2\lambda_{47} +e_5\lambda_1^2\lambda_{23}\lambda_{31} +e_3\lambda_0^2\lambda_{11}\lambda_{47}.
      $
      \end{minipage}
      \end{align*}
  \end{enumerate}
\end{note}
\begin{note}
  \label{eq:P:Ext:all-known}
  \begin{enumerate}[(i)]
    \item
      \begin{enumerate}[(1)]
        \item $\widehat{h}_i =(Sq^0)^{i-1}(\widehat{h}_1)\in\Ext_\A^{0,2^i-1}(P)$, $i\ge1$.
        \item $\widehat{c}_i =(Sq^0)^i(\widehat{c}_0)\in\Ext_\A^{2,2^{i+3}+2^{i+1}+2^i-1}(P)$, $i\ge0$.
      \end{enumerate}
    \item
      \begin{enumerate}[(1)]
        \item $\widehat{d}_i =(Sq^0)^i(\widehat{d}_0)\in\Ext_\A^{3,2^{i+4}+2^{i+1}-1}(P)$, $i\ge0$.
        \item $\widehat{e}_i =(Sq^0)^i(\widehat{e}_0) \in\Ext_\A^{3,2^{i+4}+2^{i+2}+2^i-1}(P)$, $i\ge0$.
        \item $\widehat{f}_i =(Sq^0)^i(\widehat{f}_0) \in\Ext_\A^{3,2^{i+4}+2^{i+2}+2^{i+1}-1}(P)$, $i\ge0$.
        \item $\widehat{g}_{i+1} =(Sq^0)^i(\widehat{g}_1) \in\Ext_\A^{3,2^{i+4}+2^{i+3}-1}(P)$, $i\ge0$.
        \item $\widehat{p}_i =(Sq^0)^i(\widehat{p}_0) \in\Ext_\A^{3,2^{i+5}+2^{i+2}+2^i-1}(P)$, $i\ge0$.
        \item $\widehat{D_3}(i) =(Sq^0)^i(\widehat{D_3}) \in\Ext_\A^{3,2^{i+6}+2^i-1}(P)$, $i\ge0$.
        \item $\widehat{p}_i' =(Sq^0)^i(\widehat{p}_0') \in\Ext_\A^{3,2^{i+6}+2^{i+3}+2^i-1}(P)$, $i\ge0$.
        \item $\alpha_{16}(i) =(Sq^0)^i(\alpha_{16}) \in\Ext_\A^{3,2^{i+4}+2^{i+2}-1}(P)$, $i\ge0$.
        \item $\alpha_{21}(i) =(Sq^0)^i(\alpha_{21}) \in\Ext_\A^{3,2^{i+4}+2^{i+3}+2^i-1}(P)$, $i\ge0$.
        \item $\xi_{31}(i) =(Sq^0)^i(\xi_{31}) \in\Ext_\A^{3,2^{i+5}+2^{i+1}+2^i-1}(P)$, $i\ge0$.
      \end{enumerate}
    \item
      \begin{enumerate}[(1)]
        \item $\widehat{P^1h_1}\in\Ext_\A^{4,13}(P)$.
        \item $\widehat{P^1h_2}\in\Ext_\A^{4,15}(P)$.
        \item $\delta_{30}(i) =(Sq^0)^i(\delta_{30}) \in\Ext_\A^{4,2^{i+5}+2^{i+1}+2^i-1}(P)$, $i\ge0$.
        \item $\widehat{n}_i =(Sq^0)^i(\widehat{n}_0) \in\Ext_\A^{4,2^{i+5}+2^{i+2}-1}(P)$, $i\ge0$.
        \item $\gamma_{37}(i) =(Sq^0)^i(\gamma_{37}) \in\Ext_\A^{4,2^{i+5}+2^{i+3}+2^{i+1}-1}(P)$, $i\ge0$.
        \item $\widehat{x}_i =(Sq^0)^i(\widehat{x}_0) \in\Ext_\A^{4,2^{i+5}+2^{i+3}+2^{i+1}-1}(P)$, $i\ge0$.
        \item $\alpha_{40}(i) =(Sq^0)^i(\alpha_{40}) \in\Ext_\A^{4,2^{i+5}+2^{i+3}+2^{i+2}+2^i-1}(P)$, $i\ge0$.
        \item $\phi_{45}(i) =(Sq^0)^i(\phi_{45}) \in\Ext_\A^{4,2^{i+5}+2^{i+4}+2^{i+1}-1}(P)$, $i\ge0$.
        \item $\gamma_{50}(i) =(Sq^0)^i(\gamma_{50}) \in\Ext_\A^{4,2^{i+5}+2^{i+4}+2^{i+2}+2^{i+1}+2^i-1}(P)$, $i\ge0$.
        \item $\widehat{D_1}(i) =(Sq^0)^i(\widehat{D_1}) \in\Ext_\A^{4,2^{i+5}+2^{i+4}+2^{i+3}+2^i-1}(P)$, $i\ge0$.
        \item $\gamma_{61}(i) =(Sq^0)^i(\gamma_{61}) \in\Ext_\A^{4,2^{i+6}+2^{i+1}-1}(P)$, $i\ge0$.
        \item $\widehat{H_1}(i) =(Sq^0)^i(\widehat{H_1}) \in\Ext_\A^{4,2^{i+6}+2^{i+1}+2^i-1}(P)$, $i\ge0$.
        \item $\gamma_{63}(i) =(Sq^0)^i(\gamma_{63}) \in\Ext_\A^{4,2^{i+6}+2^{i+1}+2^i-1}(P)$, $i\ge0$.
        \item $\phi_{65}(i) =(Sq^0)^i(\phi_{65}) \in\Ext_\A^{4,2^{i+6}+2^{i+2}+2^{i+1}-1}(P)$, $i\ge0$.
        \item $\delta_{67}(i) =(Sq^0)^i(\delta_{67}) \in\Ext_\A^{4,2^{i+6}+2^{i+3}-1}(P)$, $i\ge0$.
        \item $\widehat{Q_3}(i) =(Sq^0)^i(\widehat{Q_3}) \in\Ext_\A^{4,2^{i+6}+2^{i+3}-1}(P)$, $i\ge0$.
        \item $\rho_{80}(i) =(Sq^0)^i(\rho_{80}) \in\Ext_\A^{4,2^{i+6}+2^{i+4}+2^{i+2}+2^i-1}(P)$, $i\ge0$.
        \item $\widehat{K}_i =(Sq^0)^i(\widehat{K}_0) \in\Ext_\A^{4,2^{i+7}+2^{i+1}-1}(P)$, $i\ge0$.
        \item $\gamma_{128}(i) =(Sq^0)^i(\gamma_{128}) \in\Ext_\A^{4,2^{i+7}+2^{i+2}+2^i-1}(P)$, $i\ge0$.
        \item $\delta_{128}(i) =(Sq^0)^i(\delta_{128}) \in\Ext_\A^{4,2^{i+7}+2^{i+2}+2^i-1}(P)$, $i\ge0$.
        \item $\widehat{J}_i =(Sq^0)^i(\widehat{J}_0) \in\Ext_\A^{4,2^{i+7}+2^{i+2}+2^i-1}(P)$, $i\ge0$.
        \item $\gamma_{132}(i) =(Sq^0)^i(\gamma_{132}) \in\Ext_\A^{4,2^{i+7}+2^{i+3}+2^i-1}(P)$, $i\ge0$.
        \item $\gamma_{140}(i) =(Sq^0)^i(\gamma_{140}) \in\Ext_\A^{4,2^{i+7}+2^{i+4}+2^i-1}(P)$, $i\ge0$.
        \item $\widehat{T}_i =(Sq^0)^i(\widehat{T}_0) \in\Ext_\A^{4,2^{i+7}+2^{i+4}+2^{i+1}-1}(P)$, $i\ge0$.
        \item $\delta_{144}(i) =(Sq^0)^i(\delta_{144}) \in\Ext_\A^{4,2^{i+7}+2^{i+4}+2^{i+2}+2^i-1}(P)$, $i\ge0$.
        \item $\gamma_{148}(i) =(Sq^0)^i(\gamma_{148}) \in\Ext_\A^{4,2^{i+7}+2^{i+4}+2^{i+3}+2^i-1}(P)$, $i\ge0$.
        \item $\delta_{156}(i) =(Sq^0)^i(\delta_{156}) \in\Ext_\A^{4,2^{i+7}+2^{i+5}+2^i-1}(P)$, $i\ge0$.
        \item $\widehat{V}_i =(Sq^0)^i(\widehat{V}_0) \in\Ext_\A^{4,2^{i+7}+2^{i+5}+2^i-1}(P)$, $i\ge0$.
        \item $\gamma_{164}(i) =(Sq^0)^i(\gamma_{164}) \in\Ext_\A^{4,2^{i+7}+2^{i+5}+2^{i+3}+2^i-1}(P)$, $i\ge0$.
        \item $\widehat{V}_i' =(Sq^0)^i(\widehat{V}_0') \in\Ext_\A^{4,2^{i+8}+2^i-1}(P)$, $i\ge0$.
        \item $\widehat{U}_i =(Sq^0)^i(\widehat{U}_0) \in\Ext_\A^{4,2^{i+8}+2^{i+3}+2^i-1}(P)$, $i\ge0$.
      \end{enumerate}
  \end{enumerate}
\end{note}
\noindent
We also describe the class $\beta_{187}\in\Ext_\A^{5,192}(P)$ in \eqref{eq:P:Ext}~(8) as follows.
$\beta_{187}$ is an indecomposable element too.
\begin{note}
  \label{eq:P:Ext:beta}
  $\beta_{187}=\{e_8\lambda_3^2\lambda_{47}\lambda_{63}^2\}\in\Ext_\A^{5,192}(P)$.
\end{note}
\noindent
Thus all the classes in \eqref{eq:P:Ext}
can be read from \eqref{eq:Ext:known}, \eqref{eq:P:Ext:known} and \eqref{eq:P:Ext:beta},
noting that $\Ext_\A^{\ast,\ast}(P)$ is an $\Ext_\A^{\ast,\ast}$-module.

Now we state the ``already known'' results on $\Ext_\A^{s,\ast}$ for $s\le4$
from \cite{lin_ext_2008,chen_determination_2011} as follows.
\begin{thm}[\cite{lin_ext_2008,chen_determination_2011}]
  \label{thm:P:Ext:known}
  \begin{enumerate}[(i)]
    \item Modulo indecomposable elements in $\Ext_\A^{3,\ast}(P)$,
      $\Ext_\A^{s,\ast}(P)$ for $s\le3$ is generated by
      $\widehat{h}_i,\widehat{c}_i,h_i,c_i$ subject to
      the relations in \autoref{thm:Ext:known} (i), and
      the following relations.
      \begin{multicols}{2}
      \begin{enumerate}[(1)]
        \item $\widehat{h}_ih_{i-1}=0$,
        \item $\widehat{h}_{i+2}h_i^2=\widehat{h}_{i+1}h_{i+1}^2$,
        \item $\widehat{h}_{i+2}h_{i+2}h_i=0$,
        \item $\widehat{h}_{i+3}h_{i+3}h_i^2=0$,
        \item $\widehat{h}_{i+2}h_{i+3}^2h_i=0$,
        \item $\widehat{h}_ic_j=0$ for $j=i-2,i-3$,
        \item $\widehat{c}_ih_j=0$ for $j=i-1,i,i+2,i+3$,
        \item $\widehat{c}_ih_{i+1}=\widehat{h}_{i+1}c_i\ne0$.
      \end{enumerate}
      \end{multicols}
    \item The elements in \eqref{eq:P:Ext:all-known} (ii)
      is a $\Z/2$-base for
      the submodule of indecomposable elements
      in the $\Z/2$-module $\Ext_\A^{3,\ast}(P)$.
    \item The decomposable elements in $\Ext_\A^{4,\ast}(P)$
      is the $\Z/2$-submodule generated by the following elements.
      \begin{multicols}{2}
      \begin{enumerate}[(1)]
        \item $\widehat{h}_{4+n}h_0^4$, $n\ge0$,
        \item $\widehat{h}_lh_j^3h_k$, $1\le l\le j-1<k-3$,
        \item $\widehat{h}_nh_{j-1}^3h_k$, $3\le j+2\le n\le k$, $j<k-2$,
        \item $\widehat{h}_{k+n}h_{j-1}^3h_{k-1}$, $n\ge1$, $1\le j<k-2$,
        \item $\widehat{e}_{j-2}h_k$, $2\le j<k-2$,
        \item $\widehat{c}_{j-3}h_{j+1}h_k$, $3\le j<k-2$,
        \item $\widehat{h}_lh_j^2h_k^2$, $1\le l\le j<k-3$,
        \item $\widehat{h}_nh_{j-1}^2h_k^2$, $3\le j+2\le n\le k$, $j<k-3$,
        \item $\widehat{c}_{k-2}h_{j-1}^2$, $1\le j<k-3$,
        \item $\widehat{h}_{k+1+n}h_{j-1}^2h_{k-1}^2$, $n\ge1$, $1\le j<k-3$,
        \item $\widehat{c}_{j-2}h_k^2$, $2\le j<k-3$,
        \item $\widehat{h}_lh_j^2h_kh_i$, $1\le l\le j<k-2<i-3$,
        \item $\widehat{h}_nh_{j-1}^2h_kh_i$, $3\le j+2\le n\le k$, $j<k-2<i-3$,
        \item $\widehat{h}_nh_{j-1}^2h_{k-1}h_i$, $k+1\le n\le i$, $1\le j<k-2<i-3$,
        \item $\widehat{h}_{i+n+1}h_{j-1}^2h_{k-1}h_{i-1}$, $n\ge0$, $1\le j<k-2<i-3$,
        \item $\widehat{c}_{j-2}h_kh_i$, $2\le j<k-2<i-3$,
        \item $\widehat{h}_lh_jh_k^2h_i$, $1\le l\le j<k-2<i-4$,
        \item $\widehat{f}_jh_i$, $0\le j<i-5$,
        \item $\alpha_{16}(j)h_i$, $0\le j<i-5$,
        \item $\widehat{g}_{j+1}h_i$, $0\le j<i-5$,
        \item $\widehat{h}_nh_{j-1}h_k^2h_i$, $2\le j+1\le n\le k$, $j<k-2<i-4$,
        \item $\widehat{c}_{k-2}h_{j-1}h_i$, $1\le j<k-3<i-5$,
        \item $\widehat{p}_{j-1}'h_i$, $1\le j<i-5$,
        \item $\widehat{h}_nh_{j-1}h_{k-1}^2h_i$, $k+2\le n\le i$, $1\le j<k-2<i-4$,
        \item $\widehat{h}_{i+1+n}h_{j-1}h_{k-1}^2h_{i-1}$, $n\ge0$, $1\le j<k-2<i-4$,
        \item $\widehat{h}_lh_jh_kh_i^2$, $1\le l\le j<k-1<i-3$,
        \item $\widehat{h}_nh_{j-1}h_kh_i^2$, $2\le j+1\le n\le k$, $j<k-1<i-3$,
        \item $\widehat{f}_kh_{j-1}$, $1\le j<k-1$,
        \item $\alpha_{16}(k)h_{j-1}$, $1\le j<k-1$,
        \item $\widehat{g}_{k+1}h_{j-1}$, $1\le j<k-1$,
        \item $\widehat{h}_nh_{j-1}h_{k-1}h_i^2$, $k+1\le n\le i$, $1\le j<k-1<i-3$,
        \item $\widehat{c}_{i-2}h_{j-1}h_{k-1}$, $1\le j<k-1<i-3$,
        \item $\widehat{p}_{k-1}'h_{j-1}$, $1\le j<k-1$,
        \item $\widehat{h}_{i+1+n}h_{j-1}h_{k-1}h_{i-1}^2$, $n\ge1$, $1\le j<k-1<i-3$,
        \item $\widehat{h}_lh_jh_k^3$, $1\le l\le j<k-3$,
        \item $\widehat{f}_jh_{j+5}$, $j\ge0$,
        \item $\alpha_{16}(j)h_{j+5}$, $j\ge0$,
        \item $\widehat{g}_{j+1}h_{j+5}$, $j\ge0$,
        \item $\widehat{h}_nh_{j-1}h_k^3$, $2\le j+1\le n\le k-1$, $j<k-3$,
        \item $\widehat{c}_{k-3}h_{j-1}h_{k+1}$, $1\le j<k-3$,
        \item $\widehat{e}_{k-2}h_{j-1}$, $1\le j<k-3$,
        \item $\widehat{h}_{k+2+n}h_{j-1}h_{k-1}^3$, $n\ge0$, $1\le j<k-3$,
        \item $\widehat{h}_lh_jh_kh_ih_p$, $1\le l\le j<k-1<i-2<p-3$,
        \item $\widehat{h}_nh_{j-1}h_kh_ih_p$, $2\le j+1\le n\le k$, $j<k-1<i-2<p-3$,
        \item $\widehat{h}_nh_{j-1}h_{k-1}h_ih_p$, $k+1\le n\le i$, $1\le j<k-1<i-2<p-3$,
        \item $\widehat{h}_nh_{j-1}h_{k-1}h_{i-1}h_p$, $i+1\le n\le p$, $1\le j<k-1<i-2<p-3$,
        \item $\widehat{h}_{p+n}h_{j-1}h_{k-1}h_{i-1}h_{p-1}$, $n\ge1$, $1\le j<k-1<i-2<p-3$,
        \item $\widehat{D_3}(k-1)h_{j-1}$, $1\le j<k-1$,
        \item $\widehat{D_3}(j-1)h_p$, $1\le j<p-5$,
        \item $\widehat{d}_jh_j$, $j\ge0$,
        \item $\widehat{h}_lc_jh_{j+1}$, $1\le l\le j$,
        \item $\xi_{31}(j-1)h_j$, $j\ge1$,
        \item $\widehat{h}_{j+3+n}h_jc_{j-1}$, $n\ge0$, $j\ge1$,
        \item $\widehat{h}_lh_jc_k$, $1\le l\le j<k-1$,
        \item $\widehat{h}_nh_{j-1}c_k$, $2\le j+1\le n\le k+1$, $j<k-1$,
        \item $\widehat{h}_{k+3+n}h_{j-1}c_{k-1}$, $n\ge0$, $1\le j<k-1$,
        \item $\xi_{31}(k-1)h_{j-1}$, $1\le j<k-1$,
        \item $\alpha_{21}(k-1)h_{j-1}$, $1\le j<k-1$,
        \item $\widehat{p}_{k-1}h_{j-1}$, $1\le j<k-1$,
        \item $\widehat{d}_kh_{j-1}$, $1\le j<k-1$,
        \item $\widehat{h}_lc_jh_k$, $1\le l\le j+1<k-2$,
        \item $\widehat{d}_jh_k$, $0\le j<k-4$,
        \item $\alpha_{21}(j-1)h_k$, $1\le j<k-3$,
        \item $\widehat{p}_{j-1}h_k$, $1\le j<k-4$,
        \item $\widehat{e}_{j+1}h_{j-1}$, $j\ge1$,
        \item $\widehat{h}_nc_{j-1}h_k$, $4\le j+3\le n\le k$, $j<k-3$,
        \item $\widehat{h}_{k+1+n}c_{j-1}h_{k-1}$, $n\ge0$, $1\le j<k-4$,
        \item $\xi_{31}(j-1)h_k$, $1\le j<k-3$,
        \item $\widehat{h}_ld_j$, $1\le l\le j+2$,
        \item $\widehat{h}_{j+4+n}d_{j-1}$, $n\ge0$, $j\ge1$,
        \item $\xi_{31}(j-1)h_{j+1}$, $j\ge1$,
        \item $\xi_{31}(j-1)h_{j+2}$, $j\ge1$,
        \item $\widehat{h}_le_j$, $1\le l\le j+2$,
        \item $\widehat{h}_{j+4+n}e_{j-1}$, $n\ge0$, $j\ge1$,
        \item $\widehat{h}_lf_j$, $1\le l\le j+2$,
        \item $\widehat{h}_{j+4+n}f_{j-1}$, $n\ge0$, $j\ge1$,
        \item $\widehat{p}_{j-1}'h_{j-1}$, $j\ge1$,
        \item $\widehat{h}_lg_{j+1}$, $1\le l\le j+2$,
        \item $\widehat{h}_{j+4+n}g_j$, $n\ge1$, $j\ge1$,
        \item $\widehat{h}_lp_j$, $1\le l\le j$ and $l=j+3$,
        \item $\widehat{g}_{j+2}h_j$, $j\ge0$,
        \item $\widehat{h}_{j+1}p_{j-1}'$, $j\ge1$,
        \item $\widehat{h}_{j+5+n}p_{j-1}$, $n\ge0$, $j\ge1$,
        \item $\widehat{h}_{j+3}D_3(j-1)$, $j\ge1$,
        \item $\widehat{h}_lD_3(j)$, $1\le l\le j+3$,
        \item $\widehat{h}_{j+6+n}D_3(j-1)$, $n\ge0$, $j\ge1$,
        \item $\widehat{h}_lp_j'$, $1\le l\le j+1$ and $j+3\le l\le j+4$,
        \item $\widehat{g}_{j+3}h_j$, $j\ge0$,
        \item $\widehat{h}_{j+6+n}p_{j-1}'$, $n\ge0$, $j\ge1$.
      \end{enumerate}
      \end{multicols}
    \item The elements in \eqref{eq:P:Ext:all-known} (iii)
      is a $\Z/2$-base for
      the submodule of indecomposable elements
      in the $\Z/2$-module $\Ext_\A^{4,\ast}(P)$.
  \end{enumerate}
\end{thm}
It is not difficult to see that
\autoref{thm:P:Ext:known} covers \eqref{eq:P:Ext}~($j$)~(a) for $1\le j\le 7$.

Let $\Ext_\A^{s,t}(P)\to{t_\ast}\Ext_\A^{s+1,t+1}$ be as in \eqref{eq:t-bar}.
From \eqref{eq:Ext:known}, \eqref{eq:t-bar}, \eqref{eq:t-bar-Sq0-commute}
and \eqref{eq:P:Ext:known}
it is easy to see the following
\eqref{eq:P:Ext:t_ast} (1) through \eqref{eq:P:Ext:t_ast} (3).
\eqref{eq:P:Ext:t_ast} (4)
is proved in \cite{lin_ext_2008,chen_determination_2011}.
\begin{note}
  \label{eq:P:Ext:t_ast}
  $\Ext_\A^{s,t}(P)\to{t_\ast}\Ext_\A^{s+1,t+1}$ has the following properties:
  \begin{multicols}{2}
    \begin{enumerate}
      \item $t_\ast(\widehat{h}_i)=h_i$, $i\ge0$.
      \item $t_\ast(\widehat{g}_{i+1})=g_{i+1}$, $i\ge0$.
      \item $t_\ast(\widehat{D_3}(i))=D_3(i)$, $i\ge0$.
      \item $t_\ast(\gamma_{61}(i))=0$, $i\ge0$.
    \end{enumerate}
  \end{multicols}
\end{note}
\noindent
It is easy to see that the results in
\eqref{eq:P:Ext}~($j$)~(b) for $j=1,5$, \eqref{eq:P:Ext}~($j$)~(c) for $2\le j\le7$,
\eqref{eq:P:Ext}~(8)~(b) and \eqref{eq:P:Ext}~(9)~(b)
follow from \eqref{thm:Ext:known} and \eqref{eq:P:Ext:t_ast}.
The equations in \eqref{eq:P:Ext}~(8)~(c) and \eqref{eq:P:Ext}~(9)~(c) also
follow straightforwardly from \eqref{eq:P:Ext:t_ast}.
The equation $t_\ast(\beta_{187})=h_5V_0$ in \eqref{eq:P:Ext}~(8)~(d)
is proved as follows.
By \eqref{eq:Ext:known}, $V_0=\{\lambda_0\lambda_{31}\lambda_{47}\lambda_{15}\lambda_{63}\}$.
So $h_5V_0=V_0h_5=\{\lambda_0\lambda_{31}\lambda_{47}\lambda_{15}\lambda_{63}\lambda_{31}\}$.
By \eqref{eq:Lambda:Adem-relation}, \eqref{eq:t-bar} and \eqref{eq:P:Ext:beta} we have
$t_\ast(\beta_{187})=\{\lambda_8\lambda_3^2\lambda_{47}\lambda_{63}^2
=\lambda_0\lambda_{15}\lambda_7^2\lambda_{63}\lambda_{95}\}$.
Straightforward calculations show that
$\delta(\lambda_0\lambda_{23}\lambda_7\lambda_{63}\lambda_{95})
=\lambda_0\lambda_{15}\lambda_7^2\lambda_{63}\lambda_{95} +\lambda_0\lambda_{31}\lambda_{47}\lambda_{15}\lambda_{63}\lambda_{31}$.
So $t_\ast(\beta_{187})=h_5V_0$.
Nontriviality and linear independence as claimed in \eqref{eq:P:Ext}~(8)~(c),
\eqref{eq:P:Ext}~(8)~(d), and \eqref{eq:P:Ext}~(9)~(c) will be proved later.

We recall
the algebraic Kahn-Priddy theorem of Lin mentioned in \autoref{se:boundary} as follows.
\begin{thm}[algebraic Kahn-Priddy theorem]\label{thm:algebraic-Kahn-Priddy}
  $\Ext_\A^{s,t}(P)\to{t_\ast}\Ext_\A^{s+1,t+1}$ is onto
  for all $t,s$ with $t-s>0$.
\end{thm}
From \autoref{thm:algebraic-Kahn-Priddy},
it is easy to see that
the $\Ext$ group results in \eqref{eq:Ext}~(7), (8)
follow from \eqref{eq:P:Ext}~(8) and \eqref{eq:P:Ext}~(9) (cf. \autoref{remark:P:Ext}).
Thus to complete the proof of \eqref{eq:Ext} and \eqref{eq:P:Ext}
it suffices to prove \eqref{eq:Ext}~(9),
\eqref{eq:P:Ext}~(8)~(a) and \eqref{eq:P:Ext}~(9)~(a) (which we will simply say ``the remaining \eqref{eq:P:Ext}~(a)''), and
\eqref{eq:P:Ext}~(8)~(c) and \eqref{eq:P:Ext}~(9)~(c) (which we will simply say ``the remaining \eqref{eq:P:Ext}~(c)'').

For this purpose we need to recall how \autoref{thm:algebraic-Kahn-Priddy} is proved.
First of all define $\Lambda(n)\subset\Lambda$ as follows.
\begin{note}
  \label{eq:Ln}
  For each $n\ge0$, let $\Lambda(n)$ be the $\Z/2$-submodule of $\Lambda$
  generated by the admissible monomials $\lambda_{i_1}\cdots\lambda_{i_s}$ with
  $i_1\le n$.
\end{note}
\noindent
It is not difficult to see that each $\Lambda(n)$ is a subcomplex of $\Lambda$.
So we have an increasing filtration
$1=\Lambda(0)\subset\Lambda(1)\subset\Lambda(2)\subset\cdots\subset\Lambda(n)\subset\Lambda(n+1)\subset\cdots$
of subcomplexes of $\Lambda$.
The differential $\Lambda$-module $ H_\ast(P)\tensor\Lambda$
also has an increasing filtration of subcomplexes
$F(0)=0\subset F(1)\subset F(2)\subset\cdots\subset F(n)\subset F(n+1)\subset\cdots$
defined as follows.
\begin{note}
  \label{eq:Fn}
  $F(n)=\Sum_{1\le k\le n} H_k(P)\tensor\Lambda$, for $n\ge1$.
\end{note}
\noindent
It is easy to see the following.
\begin{note}
  \label{eq:Fi-mod-eq-L}
  $F(i)/F(i-1)\iso\Susp^i\Lambda$.
\end{note}
\noindent
In \cite{lin_algebraic_1981} Lin has shown the following.
\begin{note}
  \label{eq:phi}
  There is a chain map $\Lambda\to{\phi} H_\ast(P)\tensor\Lambda$ that goes
  from $\Lambda^{s+1,t-s}$ to $( H_\ast(P)\tensor\Lambda)^{s,t-s}$, and so
  induces an $\Ext$ group map $\Ext_\A^{s+1,t+1}\to{\phi_\ast}\Ext_\A^{s,t}(P)$ for $t-s\ge0$.
  The chain map $\phi$ has the following properties.
  \begin{enumerate}
    \item $\phi(\Lambda(i))\subset F(i)$ for all $i\ge0$.
    \item $\phi(\lambda_{i_1}\lambda_{i_2}\cdots\lambda_{i_s})\equiv e_{i_1}\lambda_{i_2}\cdots\lambda_{i_s}\mod{F(i_1-1)}$
      for any admissible monomial $\lambda_{i_1}\cdots\lambda_{i_s}$ with $i_1\ge1$.
    \item Let $\psi$ be the composite $\Lambda\to{\phi} H_\ast(P)\tensor\Lambda\to{\bar{t}}\Lambda$
      where $\bar{t}$ is as in \eqref{eq:t-bar}. Then for each admissible monomial
      $\lambda_{i_1}\lambda_{i_2}\cdots\lambda_{i_s}\in\Lambda$ with $i_1\ge1$ there is the relation
      \begin{align*}
        \psi(\lambda_{i_1}\lambda_{i_2}\cdots\lambda_{i_s})
        \equiv\lambda_{i_1}\lambda_{i_2}\cdots\lambda_{i_s}\mod{\Lambda(i_1-1)}.
      \end{align*}
  \end{enumerate}
\end{note}
\noindent
It is easy to see that
\autoref{thm:algebraic-Kahn-Priddy}
follows from \eqref{eq:t-bar-Sq0-commute} and \eqref{eq:phi} (3).
To complete our proofs for \eqref{eq:Ext} and \eqref{eq:P:Ext},
we are going to use a spectral sequence $\{E_r^{i,s,t}\}_{r\ge1}$ for
computing $\Ext_\A^{\ast,\ast}(P)$ described in \eqref{eq:P:Ext:spectral-sequence} below
to prove the remaining \eqref{eq:P:Ext}~(a).
Then we will use these results, and the map $\phi_\ast$ in \eqref{eq:phi} to deduce
\eqref{eq:Ext}~(9) and the remaining \eqref{eq:P:Ext}~(c).

The spectral sequence for computing
$\Ext_\A^{s,t+s}(P)=H^{s,t}( H_\ast(P)\tensor\Lambda)$ is
the one defined by the filtration $\{F(i)\}_{i\ge0}$
for $ H_\ast(P)\tensor\Lambda$ in \eqref{eq:Fn}.
From \eqref{eq:Fi-mod-eq-L} we see that
\begin{align}
  H^{s,t}(F(i)/F(i-1))&\iso H^{s,t}(\Susp^i\Lambda)\iso\Susp^i\Ext_\A^{s,t+s-i}.
\end{align}
So we have the following.
\begin{note}
  \label{eq:P:Ext:spectral-sequence}
  The filtration $\{F(i)\}_{i\ge0}$ gives rise to a spectral sequence $\{E_r^{i,s,t}\}_{r\ge1}$ with
  $E_1^{i,s,t}= H^{s,t}(F(i)/F(i-1))\iso\Susp^i\Ext_\A^{s,t+s-i}$,
  $i\ge1$,
  and $\dsum_{i\ge1}E_\infty^{i,s,t}\iso\Ext_\A^{s,t+s}(P)$ as $\Z/2$-modules.
  The differential $d_r$ of this spectral sequence goes from $E_r^{i,s,t}$ to $E_r^{i-r,s+1,t-1}$.
\end{note}
\noindent
We will simply write
$E_r^{\ast,s,\ast}\to{d_r}E_r^{\ast,s+1,\ast}$
to indicate that we are considering the differentials for a fixed $s$
and for all $i,t,r$.
From \eqref{eq:P:Ext:spectral-sequence} we see that
for a given $s>1$,
if one knows completely the $\Ext$ groups $\Ext_\A^{s,\ast}$ and
has computed all the differentials
$E_r^{\ast,s-1,\ast}\to{d_r}E_r^{\ast,s,\ast}$ and
$E_r^{\ast,s,\ast}\to{d_r}E_r^{\ast,s+1,\ast}$, then
one can determine the $\Ext$ groups $\Ext_\A^{s,\ast}(P)$.
This is done in \cite{cohen_adams_1988}, \cite{lin_ext_2008} and \cite{chen_determination_2011},
where
the differentials $E_r^{\ast,s,\ast}\to{d_r}E_r^{\ast,s+1,\ast}$
with $0\le s\le2$, with $s=3$ and with $s=4$ are computed respectively.
To deduce the remaining \eqref{eq:P:Ext}~(a)
we have to compute the differentials $E_r^{i,s,t}\to{d_r}E_r^{i-r,s+1,t-1}$
for $(s,t) =(4,188), (4,189)$ and for $(s,t) =(5,187), (5,188)$, all $i,r$.
These will be described in
\eqref{eq:P:Ext:5-192} and \eqref{eq:P:Ext:5-193}
below.
As will be seen later in \eqref{eq:Ext:7-195},
to deduce \eqref{eq:Ext} (9) we also need to compute
the above differential for $(s,t) =(4,190)$ and for $(s,t) =(5,189)$.
This will be described in \eqref{eq:P:Ext:5-194} below.

To state
\eqref{eq:P:Ext:5-192}, \eqref{eq:P:Ext:5-193} and \eqref{eq:P:Ext:5-194}
we use the following conventions.
By \eqref{eq:P:Ext:spectral-sequence},
if $\alpha$ is a basis element in $\Ext_\A^{s,\ast}$
then $e_i\alpha=e_i\tensor\alpha$ is a basis element in $E_1^{i,s,\ast}$.
We will write $e_i\alpha->e_{i-r}\beta$,
where $e_{i-r}\beta$ is some basis element in $E_1^{i-r,s+1,\ast}$,
to indicate that
both $e_i\alpha$ and $e_{i-r}\beta$ survive to $E_r^{\ast,\ast,\ast}$
and $d_r(e_i\alpha)=e_{i-r}\beta$ in the spectral sequence.
In \eqref{eq:P:Ext:5-192}, \eqref{eq:P:Ext:5-193} and \eqref{eq:P:Ext:5-194} below
there are two families of non-trivial differentials:
the differentials in
($j$) (a) for $j=\ref{eq:P:Ext:5-192},\ref{eq:P:Ext:5-193},\ref{eq:P:Ext:5-194}$
(which will be refered to as the ``(a)-differentials'') that go from $E_r^{\ast,4,\ast}$ to $E_r^{\ast,5,\ast}$,
and the differentials in
($j$) (b) for $j=\ref{eq:P:Ext:5-192},\ref{eq:P:Ext:5-193},\ref{eq:P:Ext:5-194}$
(which will be refered to as the ``(b)-differentials'') that go from $E_r^{\ast,5,\ast}$ to $E_r^{\ast,6,\ast}$.
Any basis element $e_i\alpha$ in $E_1^{\ast,5,\ast}$ that
is not in the targets of the (a)-differentials nor in the sources of the (b)-differentials
is a cycle that survives to $E_\infty^{\ast,5,\ast}$.
These permanent cycles are listed in
($j$) (c) for $j=\ref{eq:P:Ext:5-192},\ref{eq:P:Ext:5-193},\ref{eq:P:Ext:5-194}$
(which will be refered to as the ``(c)-cycles'').
For each (c)-cycle $e_i\alpha$, we will write
$e_i\alpha<->\beta$ to mean that $\beta\in\Ext_\A^{\ast,\ast}(P)$
is the cohomology class represented by this cycle.

All of the (a)-differentials
are covered by \cite{chen_determination_2011} ((4.1.1) through (4.1.289) in \cite{chen_determination_2011}).
In \eqref{eq:P:Ext:5-192}, \eqref{eq:P:Ext:5-193} and \eqref{eq:P:Ext:5-194} below
we simply list these differentials and omit the calculations.
Computing the (b)-differentials is the main work of this section from which
we are going to deduce the remaining \eqref{eq:P:Ext}~(a) and \eqref{eq:P:Ext}~(c).
These calculations are lengthy and will not be given in this paper.
Also, the targets of the (b)-differentials are of the form $e_i\beta$
where $\beta\in\Ext_\A^{6,\ast}$.
So to state the (b)-differentials
we need to describe some indecomposable elements in $\Ext_\A^{6,\ast}$
(elements which are not sums of products of elements in \eqref{eq:Ext:known}).
These elements are listed in \eqref{eq:Ext:6:new} below.
\begin{note}
  \label{eq:Ext:6:new}
  \begin{enumerate}
    \item $\bar{\lambda_{19}D_1}=\{\bar{\lambda_{19}D_1}^\ast =\lambda_{19}D_1^\ast +\lambda_{15}\lambda_1\lambda_2\lambda_3^2\lambda_{47} +(\lambda_{14}\lambda_9\lambda_7 +\lambda_{11}\lambda_{12}\lambda_7 +\lambda_7\lambda_{16}\lambda_7 +\lambda_7\lambda_7\lambda_{16})\lambda_{11}\lambda_{15}^2
      \}$.
    \item $\bar{\lambda_{39}D_1(1)} =Sq^0(\bar{\lambda_{19}D_1})$.
    \item $\bar{\lambda_{40}D_1(1)} =\{\lambda_2\lambda_3^2(\lambda_{79}\lambda_{31}^2 +\lambda_{31}^2\lambda_{79})\}$.
    \item $A =\{\lambda_9(\lambda_4\lambda_7\lambda_{11}\lambda_{15}^2) +\lambda_5\lambda_1\lambda_2\lambda_3^2\lambda_{47} +\lambda_3\lambda_0^3\lambda_{11}\lambda_{47}\}$.
    \item $A' =\{\lambda_{20}\lambda_0^2c_2^\ast +\lambda_{18}\lambda_2\lambda_0c_2^\ast +\lambda_{17}\lambda_1\lambda_2c_2^\ast +\lambda_{15}\lambda_4\lambda_7\lambda_5\lambda_{15}^2 + \lambda_{14}\lambda_6\lambda_0c_2^\ast +\lambda_{13}(\lambda_5\lambda_2 +\lambda_3\lambda_4)c_2^\ast +\lambda_{12}(\lambda_8\lambda_0 +\lambda_0\lambda_8)c_2^\ast + \lambda_{11}(\lambda_8\lambda_9(\lambda_{19}\lambda_7^2 +\lambda_7^2\lambda_{19}) +\lambda_3\lambda_6c_2^\ast +\lambda_0\lambda_7\lambda_{13}\lambda_{15}^2) + \lambda_{12}f_0^\ast\lambda_{31} +\lambda_{10}(\lambda_2\lambda_4\lambda_{15}^3 +\lambda_{10}\lambda_0c_2^\ast +\lambda_2\lambda_8c_2^\ast) + \lambda_9(\lambda_9\lambda_2c_2^\ast +\lambda_1\lambda_{10}c_2^\ast) + \lambda_7(\lambda_4\lambda_7\lambda_{13}\lambda_{15}^2 +\lambda_7\lambda_6c_2^\ast +\lambda_5\lambda_8c_2^\ast) + \lambda_6(\lambda_6\lambda_8c_2^\ast +\lambda_{10}\lambda_0\lambda_{15}^3 +\lambda_9^2\lambda_3^2 +\lambda_9\lambda_3^2\lambda_9 +\lambda_7\lambda_{11}\lambda_3^2\lambda_{31} +\lambda_6\lambda_4\lambda_{15}^3 +\lambda_7\lambda_5\lambda_9\lambda_3\lambda_{31})\}$. 
    \item $D_2=\{\lambda_0^4\lambda_{11}\lambda_{47}\}$. 
    \item $t_0= \{\lambda_5(\lambda_9\lambda_3\lambda_5\lambda_7^2 +\lambda_6\lambda_0\lambda_3\lambda_{15}\lambda_7 +\lambda_3\lambda_5\lambda_1\lambda_{15}\lambda_7) +\lambda_3(\lambda_{11}\lambda_3\lambda_5\lambda_7^2+\lambda_5\lambda_1\lambda_{13}\lambda_7^2+\lambda_3\lambda_5\lambda_1^2\lambda_{23})\}$.
    \item $q_0= \{\lambda_2\lambda_0^3\lambda_{15}^2 +\lambda_1( (\lambda_5\lambda_1^3+\lambda_2\lambda_4\lambda_1^2+\lambda_1\lambda_2\lambda_4\lambda_1+\lambda_2\lambda_1^2\lambda_4)\lambda_{23} +\lambda_6\lambda_0\lambda_3\lambda_{15}\lambda_7 +\lambda_2\lambda_{10}\lambda_5\lambda_7^2 +\lambda_2\lambda_4\lambda_7\lambda_3\lambda_{15} +\lambda_9\lambda_3\lambda_5\lambda_7^2)\}$.
  \end{enumerate}
\end{note}
\noindent
We also need to show the following.
\begin{note}
  \label{eq:P:Ext:nontriviality}
  Let $S$ be the set of all the classes $\beta$ such that
  $e_i\beta$ appears in the targets of a (b)-differential
  in \eqref{eq:P:Ext:5-192} through \eqref{eq:P:Ext:5-194}.
  Then $S$ is linearly independent in $\Ext_\A^{6,\ast}$.
\end{note}
\noindent
\eqref{eq:P:Ext:nontriviality} will be proved in \eqref{eq:Ext:6:nonzero} later.

Finally we make some remarks on the (c)-cycles.
From \eqref{eq:P:Ext:spectral-sequence},
a permanent cycle in the spectral sequence
represents a class in $\Ext_\A^{\ast,\ast}(P)$.
One needs some further calculations
to describe the class represented by a permanent cycle in terms of
the cohomology classes in \eqref{eq:Ext:known} and \eqref{eq:P:Ext:known}.
These calculations will not be given in this paper.
Instead, we illustrate our calculations for two examples of these.
The first example is $e_{14}h_1f_3$
in \eqref{eq:P:Ext:5-192} (c) (2)
which represents the cohomology class $\widehat{h}_4h_0f_3\in\Ext_\A^{5,192}(P)$.
From \eqref{eq:Ext:known},
$f_3=\{f_3^\ast\}\in\Ext_\A^{4,176}$ where
\begin{align*}
  f_3^\ast
  &=(Sq^0)^3(f_0^\ast)
  = \lambda_{39}\lambda_7\lambda_{63}^2 +\lambda_{31}(\lambda_{79}\lambda_{31} +\lambda_{31}\lambda_{47}\lambda_{63}) +\lambda_{23}\lambda_{39}\lambda_{47}\lambda_{63}.
\end{align*}
From \eqref{eq:P:Ext:known} we have
$\widehat{h}_4=\{e_{15}\}$,
so $\widehat{h}_4h_0=\{e_{15}\lambda_0\}=\{e_{14}\lambda_1 +e_{12}\lambda_3 +e_8\lambda_7\}$,
and $\widehat{h}_4h_0f_3=\{(e_{14}\lambda_1 +e_{12}\lambda_3 +e_8\lambda_7)f_3^\ast\}$.
It follows that the permanent cycle $e_{14}h_1f_3$ represents the class
$\widehat{h}_4h_0f_3$.

The other example is the following.
\begin{note}
  \label{eq:e47T0-eq-g4h02}
  The cycle $e_{47}T_0$ represents the cohomology class
  $\widehat{g}_4h_0^2\in\Ext_\A^{5,193}(P)$
  in the the spectral sequence of $\Ext_\A^{\ast,\ast}(P)$.
\end{note}
\noindent
The result \eqref{eq:e47T0-eq-g4h02} will be used to prove \eqref{eq:P62-k:Ext}~(17).
It is explained as follows.
From \eqref{eq:P:Ext:known} we have $\widehat{g}_4=\{\widehat{g}_4^\ast\}\in\Ext_\A^{3,191}(P)$ where
\begin{align*}
  \widehat{g}_4^\ast
  &=(Sq^0)^3(\widehat{g}_1^\ast)
  = e_{55}\lambda_7\lambda_{63}^2 +e_{47}(\lambda_{79}\lambda_{31}^2 +\lambda_{31}^2\lambda_{79})
    +e_{31}\lambda_{47}\lambda_{79}\lambda_{31},
\end{align*}
so $\widehat{g}_4h_0^2 =\{\widehat{g}_4^\ast\lambda_0^2\}$.
From \eqref{eq:Ext:known} we have $T_0=\{T_0^\ast\}\in\Ext_\A^{5,146}$ where
\begin{align*}
  T_0^\ast
  &=\lambda_0^2(\lambda_{79}\lambda_{31}^2 +\lambda_{31}^2\lambda_{79}) +(\lambda_9\lambda_3^2 +\lambda_3^2\lambda_9)\lambda_{63}^2.
\end{align*}
The infinite cycle $e_{47}T_0$ represents $\widehat{g}_4h_0\in\Ext_\A^{5,193}(P)$ because
$\delta(B_0)\equiv \widehat{g}_4^\ast\lambda_0^2 +e_{47}T_0^\ast \mod{F(46)}$
where
\begin{align*}
  B_0
  &=
  \begin{minipage}[t]{13.5cm}
  $e_{55}\lambda_7\lambda_{63}\lambda_{64}\lambda_0 +e_{55}\lambda_7\lambda_{63}\lambda_0\lambda_{64} +e_{55}\lambda_7\lambda_{64}\lambda_0\lambda_{63} +e_{55}\lambda_7\lambda_0\lambda_{64}\lambda_{63}
  +e_{55}\lambda_8\lambda_0\lambda_{63}^2 +e_{55}\lambda_0\lambda_8\lambda_{63}^2 +e_{56}\lambda_0\lambda_7\lambda_{63}^2
  +e_{54}\lambda_2\lambda_7\lambda_{63}^2
  +e_{57}\lambda_3^2\lambda_{63}^2 +e_{52}\lambda_4\lambda_7\lambda_{63}^2 +e_{51}\lambda_5\lambda_7\lambda_{63}^2 +
  e_{47}(\lambda_{79}\lambda_{31}\lambda_{32}\lambda_0 +\lambda_{79}\lambda_{31}\lambda_0\lambda_{32} +\lambda_{79}\lambda_{32}\lambda_0\lambda_{31}
  +\lambda_{79}\lambda_0\lambda_{32}\lambda_{31}
  +\lambda_{80}\lambda_0\lambda_{31}^2 +\lambda_{63}\lambda_{16}\lambda_{32}\lambda_{31} +\lambda_{64}\lambda_{16}\lambda_{31}^2 +\lambda_0\lambda_{80}\lambda_{31}^2
  +\lambda_{31}\lambda_{47}\lambda_{64}\lambda_0 +\lambda_{31}\lambda_{47}\lambda_0\lambda_{64}
  +\lambda_{31}\lambda_{48}\lambda_0\lambda_{63} +\lambda_{31}\lambda_0\lambda_{48}\lambda_{63}
  +\lambda_0\lambda_{32}\lambda_{47}\lambda_{63} +\lambda_{32}\lambda_0\lambda_{47}\lambda_{63} +\lambda_9\lambda_7\lambda_{63}^2
  +\lambda_7^2\lambda_{64}^2 +\lambda_7\lambda_8\lambda_{64}\lambda_{63}
  +\lambda_7\lambda_8\lambda_{63}\lambda_{64} +\lambda_8\lambda_7\lambda_{64}\lambda_{63}
  +\lambda_8^2\lambda_{63}^2 +\lambda_8\lambda_7\lambda_{63}\lambda_{64}
  +\lambda_{47}\lambda_{31}\lambda_{32}^2 +\lambda_{48}\lambda_{31}\lambda_{32}\lambda_{31} +\lambda_{48}\lambda_{31}^2\lambda_{32}^2 +\lambda_{47}\lambda_{32}^2\lambda_{31}
  +\lambda_{47}\lambda_{32}\lambda_{31}\lambda_{32}
  +\lambda_{31}^2\lambda_{16}\lambda_{64} +\lambda_{31}\lambda_{32}\lambda_{15}\lambda_{64} +\lambda_{32}\lambda_{31}\lambda_{15}\lambda_{64} +\lambda_{32}\lambda_{31}\lambda_{16}\lambda_{63}
  +\lambda_{32}^2\lambda_{15}\lambda_{63})$.
  \end{minipage}
\end{align*}
It follows that the infinite cycle $e_{47}T_0$ represents $\widehat{g}_4h_0^2$.

Now we describe the differentials
$E_r^{i,4,188}\to{d_r}E_r^{i-r,5,187}\to{d_r}E_r^{i-2r,6,186}$,
$E_r^{i,4,189}\to{d_r}E_r^{i-r,5,188}\to{d_r}E_r^{i-2r,6,187}$, and
$E_r^{i,4,190}\to{d_r}E_r^{i-r,5,189}\to{d_r}E_r^{i-2r,6,188}$ in
\eqref{eq:P:Ext:5-192}, \eqref{eq:P:Ext:5-193}, \eqref{eq:P:Ext:5-194} below respectively
(in the spectral sequence \eqref{eq:P:Ext:spectral-sequence} for $\Ext_\A^{\ast,\ast}(P)$).
\renewcommand{\labelenumi}{(\alph{enumi})}
\renewcommand{\labelenumii}{(\arabic{enumii})}
\renewcommand{\labelenumiii}{(\roman{enumiii})}
\begin{note}
  \label{eq:P:Ext:5-192}
  \begin{enumerate}
    \item
      \begin{multicols}{2}
\begin{enumerate}
  \item $e_{12}h_2c_4->e_{11}h_0h_2c_4$
  \item $e_{16}f_3->e_{15}h_0f_3$
  \item $e_{16}h_4^3h_7->e_{15}h_0h_4^3h_7$
  \item $e_{20}c_2h_7->e_{19}h_0c_2h_7$
  \item $e_{24}h_2^2h_5h_7->e_{20}h_2^3h_5h_7$
  \item $e_{24}e_3->e_{23}h_0e_3$
  \item $e_{28}h_1^2h_5h_7->e_{26}h_1^3h_5h_7$
  \item $e_{30}h_0^2h_5h_7->e_{29}h_0^3h_5h_7$
  \item $e_{46}p'_1->e_{31}V_0$
  \item $e_{42}c_1h_7->e_{38}h_2c_1h_7$
  \item $e_{54}h_1h_3h_6^2->e_{40}h_3^3h_7$
  \item $e_{44}h_1^2h_4h_7->e_{42}h_1^3h_4h_7$
  \item $e_{44}p_2->e_{43}h_0p_2$
  \item $e_{46}h_0^2h_4h_7->e_{45}h_0^3h_4h_7$
  \item $e_{48}d_3->e_{47}h_0d_3$
  \item $e_{53}c_0h_7->e_{51}h_1c_0h_7$
  \item $e_{56}h_2^2h_6^2->e_{52}h_2^3h_6^2$
  \item $e_{54}h_0^2h_3h_7->e_{53}h_0^3h_3h_7$
  \item $e_{60}h_1^2h_6^2->e_{58}h_1^3h_6^2$
  \item $e_{62}h_0^2h_6^2->e_{61}h_0^3h_6^2$
  \item $e_{102}h_1c_3->e_{78}D_1(1)$
  \item $e_{84}c_2h_6->e_{83}h_0c_2h_6$
  \item $e_{88}h_4c_3->e_{87}h_0h_4c_3$
  \item $e_{119}p'_0->e_{94}h_1g_3$
  \item $e_{96}g_3->e_{95}h_0g_3$
  \item $e_{106}c_1h_6->e_{102}h_2c_1h_6$
  \item $e_{104}f_2->e_{103}h_0f_2$
  \item $e_{108}h_1^2h_4h_6->e_{106}h_1^3h_2h_4h_6$
  \item $e_{108}e_2->e_{107}h_0e_2$
  \item $e_{110}h_0^2h_4h_6->e_{109}h_0^3h_4h_6$
  \item $e_{118}p_1->e_{111}h_4D_3$
  \item $e_{117}c_0h_6->e_{115}h_1c_0h_6$
  \item $e_{118}h_0^2h_3h_6->e_{117}h_0^3h_3h_6$
  \item $e_{120}d_2->e_{119}h_0d_2$
  \item $e_{124}h_1^2h_5->e_{122}h_1^3h_5^2$
  \item $e_{126}h_0^2h_5^2->e_{125}h_0^3h_5^2$
  \item $e_{138}c_1h_5->e_{134}h_2c_1h_5$
  \item $e_{140}h_3c_2->e_{135}D_1$
  \item $e_{155}p_0->e_{142}h_1g_2$
  \item $e_{144}h_0g_2->e_{143}h_0g_2$
  \item $e_{150}e_1->e_{146}h_2e_1$
  \item $e_{148}f_1->e_{147}h_0f_1$
  \item $e_{149}c_0h_5->e_{147}h_1c_0h_5$
  \item $e_{150}h_0^2h_3h_5->e_{149}h_0^3h_3h_5$
  \item $e_{156}d_1->e_{154}h_1d_1$
  \item $e_{158}h_0^2h_4^2->e_{157}h_0^3h_4^2$
  \item $e_{165}c_0h_4->e_{163}h_1c_0h_4$
  \item $e_{171}e_0->e_{166}h_1g_1$
  \item $e_{168}g_1->e_{167}h_0g_1$
  \item $e_{170}f_0->e_{169}h_0f_0$
  \item $e_{174}d_0->e_{173}h_0d_0$
  \item $e_{188}h_0^4->e_{187}h_0^5$
\end{enumerate}
      \end{multicols}
    \item
      \begin{multicols}{2}
\begin{enumerate}
  \item $e_{10}c_1h_5h_7->e_6h_2c_1h_5h_7$
  \item $e_{12}h_1^2c_4->e_{10}h_1^3c_4$
  \item $e_{12}h_5p_2->e_{11}h_0h_5p_2$
  \item $e_{12}h_3c_2h_7->e_7D_1h_7$
  \item $e_{14}h_0^2c_4->e_{13}h_0^3c_4$
  \item $e_{16}h_5d_3->e_{15}h_0h_5d_3$
  \item $e_{16}g_2h_7->e_{15}h_0g_2h_7$
  \item $e_{20}f_1h_7->e_{19}h_0f_1h_7$
  \item $e_{21}c_0h_5h_7->e_{19}h_1c_0h_5h_7$
  \item $e_{22}h_1e_3->e_8h_3^2e_3$
  \item $e_{22}e_1h_7->e_{18}h_2e_1h_7$
  \item $e_{22}h_0^2h_3h_5h_7->e_{21}h_0^3h_3h_5h_7$
  \item $e_{24}x_2->e_{23}h_0x_2$
  \item $e_{27}p_0h_7->e_{14}h_1g_2h_7$
  \item $e_{28}d_1h_7->e_{26}h_1d_1h_7$
  \item $e_{30}h_0^2h_4^2h_7->e_{29}h_0^3h_4^2h_7$
  \item $e_{32}h_4d_3->e_{31}h_0h_4d_3$
  \item $e_{37}c_0h_4h_7->e_{35}h_1c_0h_4h_7$
  \item $e_{40}h_3d_3->e_{20}h_2x_2$
  \item $e_{40}g_1h_7->e_{39}h_0g_1h_7$
  \item $e_{42}f_0h_7->e_{41}h_0f_0h_7$
  \item $e_{42}c_1h_6^2->e_{38}h_2c_1h_6^2$
  \item $e_{43}e_0h_7->e_{38}h_2f_0h_7$
  \item $e_{44}p'_1h_1->e_{42}h_0^2p_2$
  \item $e_{46}h_1d_3->e_{30}h_1h_4d_3$
  \item $e_{46}d_0h_7->e_{45}h_0d_0h_7$
  \item $e_{46}T_0->e_{45}h_0T_0$
  \item $e_{48}Q_3(1)->e_{47}h_0Q_3(1)$
  \item $e_{48}n_2->e_{47}h_0n_2$
  \item $e_{53}c_0h_6^2->e_{51}h_1c_0h_6^2$
  \item $e_{54}h_0^2h_3h_6^2->e_{53}h_0^3h_3h_6^2$
  \item $e_{54}h_3D_3(1)->e_{46}h_3^2D_3(1)$
  \item $e_{58}H_1(1)->e_{54}h_2H_1(1)$
  \item $e_{58}h_2D_3(1)->e_{51}h_3J_0$
  \item $e_{59}J_0->e_{46}h_1Q_3(1)$
  \item $e_{60}h_0^4h_7->e_{59}h_0^5h_7$
  \item $e_{62}K_0->e_{61}h_0K_0$
  \item $e_{80}g_2h_6->e_{79}h_0g_2h_6$
  \item $e_{84}f_1h_6->e_{83}h_0f_1h_6$
  \item $e_{86}e_1h_6->e_{82}h_2e_1h_6$
  \item $e_{88}h_4f_2->e_{87}h_0h_4f_2$
  \item $e_{91}p_0h_6->e_{78}h_1g_2h_6$
  \item $e_{92}d_1h_6->e_{90}h_1d_1h_6$
  \item $e_{92}h_2g_3->e_{91}h_0h_2g_3$
  \item $e_{94}h_0^2h_5^3->e_{93}h_0^3h_5^3$
  \item $e_{100}h_2f_2->e_{95}h_4^2D_3$
  \item $e_{101}c_0h_4h_6->e_{99}h_1c_0h_4h_6$
  \item $e_{102}h_0^2c_3->e_{101}h_0^3c_3$
  \item $e_{102}h_1f_2->e_{86}h_0^2h_4c_3$
  \item $e_{104}h_2e_2->e_{100}h_2^2e_2$
  \item $e_{104}g_1h_6->e_{103}h_0g_1h_6$
  \item $e_{106}f_0h_6->e_{105}h_0f_0h_6$
  \item $e_{107}e_0h_6->e_{102}h_1g_1h_6$
  \item $e_{108}x_1->e_{107}h_0x_1$
  \item $e_{110}d_0h_6->e_{109}h_0d_0h_6$
  \item $e_{112}h_3d_2->e_{111}h_0h_3d_2$
  \item $e_{116}h_1p_1->e_{106}h_1x_1$
  \item $e_{117}c_0h_5^2->e_{115}h_1c_0h_5^2$
  \item $e_{118}h_0p_0'->e_{117}h_0^2p_0'$
  \item $e_{120}n_1->e_{119}h_0n_1$
  \item $e_{120}Q_3->e_{119}h_0Q_3$
  \item $e_{123}h_2D_3->e_{115}\bar{\lambda_{19}D_1}$ 
  \item $e_{124}h_0^4h_6->e_{123}h_0^5h_6$
  \item $e_{125}h_1D_3->e_{118}h_1Q_3$
  \item $e_{125}H_1->e_{123}h_1H_1$
  \item $e_{136}g_1h_5->e_{135}h_0g_1h_5$
  \item $e_{138}f_0h_5->e_{137}h_0f_0h_5$
  \item $e_{139}e_0h_5->e_{134}h_1g_1h_5$
  \item $e_{140}h_2g_2->e_{139}h_0h_2g_2$
  \item $e_{142}d_0h_5->e_{141}h_0d_0h_5$
  \item $e_{148}h_1e_1->e_{146}h_1^2e_1$
  \item $e_{150}x_0->e_{149}h_0x_0$
  \item $e_{152}h_2d_1->e_{148}h_2^2d_1$
  \item $e_{156}n_0->e_{150}t_0$
  \item $e_{156}h_0^4h_5->e_{155}h_0^5h_5$
  \item $e_{164}h_2g_1->e_{163}h_0h_2g_1$
  \item $e_{170}h_0e_0->e_{169}h_0^2e_0$
  \item $e_{172}h_0^4h_4->e_{171}h_0^5h_4$
  \item $e_{172}h_1d_0->e_{170}h_1^2d_0$
  \item $e_{176}P^1h_2->e_{175}h_0P^1h_2$
  \item $e_{178}P^1h_1->e_{172}h_0^2d_0$
\end{enumerate}
      \end{multicols}
    \item
      \begin{multicols}{2}
\begin{enumerate}
  \item $e_8h_2^2c_4<->\beta_{187}$
  \item $e_{14}h_1f_3<->\widehat{h}_4h_0f_3$
  \item $e_{30}D_3(1)h_5<->\widehat{h}_5V_0$
\end{enumerate}
      \end{multicols}
  \end{enumerate}
\end{note}
\begin{note}
  \label{eq:P:Ext:5-193}
  \begin{enumerate}
    \item
      \begin{multicols}{2}
\begin{enumerate}
  \item $e_{21}c_2h_7->e_{13}h_3c_2h_7$
  \item $e_{16}h_0c_4->e_{15}h_0^2c_4$
  \item $e_{17}f_3->e_{15}h_1f_3$
  \item $e_{25}h_2^2h_5h_7->e_{21}h_2^3h_5h_7$
  \item $e_{24}h_0h_3h_5h_7->e_{23}h_0^2h_3h_5h_7$
  \item $e_{25}e_3->e_{23}h_1e_3$
  \item $e_{28}h_0h_2h_5h_7->e_{27}h_1^3h_5h_7$
  \item $e_{32}h_0h_4^2h_7->e_{31}h_0^2h_4^2h_7$
  \item $e_{45}p_2->e_{31}D_3(1)h_5$
  \item $e_{48}h_0h_3^2h_7->e_{38}c_0h_4h_7$
  \item $e_{44}h_0h_2h_4h_7->e_{43}h_1^3h_4h_7$
  \item $e_{49}d_3->e_{47}h_1d_3$
  \item $e_{54}c_0h_7->e_{47}d_0h_7$
  \item $e_{57}h_2^2h_6^2->e_{53}h_2^3h_6^2$
  \item $e_{56}h_0h_3h_6^2->e_{55}h_0^2h_3h_6^2$
  \item $e_{60}h_0h_2h_6^2->e_{59}h_1^3h_6^2$
  \item $e_{62}h_0^3h_7->e_{61}h_0^4h_7$
  \item $e_{85}c_2h_6->e_{77}h_3c_2h_6$
  \item $e_{89}h_4c_3->e_{79}D_1(1)$
  \item $e_{109}h_1^2h_4h_6->e_{63}K_0$
  \item $e_{96}h_0h_5^3->e_{95}h_0^2h_5^3$
  \item $e_{97}g_3->e_{95}h_1g_3$
  \item $e_{104}h_0c_3->e_{103}h_0^2c_3$
  \item $e_{105}f_2->e_{103}h_1f_2$
  \item $e_{107}c_1h_6->e_{103}h_2c_1h_6$
  \item $e_{108}h_0h_2h_4h_6->e_{107}h_0^2h_2h_4h_6$
  \item $e_{109}e_2->e_{101}h_3e_2$
  \item $e_{109}h_1^2h_4h_6->e_{63}K_0$
  \item $e_{118}c_0h_6->e_{111}d_0h_6$
  \item $e_{119}p_1->e_{93}h_2g_3$
  \item $e_{120}p_0'->e_{119}h_0p_0'$
  \item $e_{121}d_2->e_{117}h_2d_2$
  \item $e_{124}h_0h_2h_5^2->e_{123}h_0^2h_2h_5^2$
  \item $e_{126}h_0^3h_6->e_{125}h_0^4h_6$
  \item $e_{128}D_3->e_{126}h_1D_3$
  \item $e_{139}c_1h_5->e_{135}h_2c_1h_5$
  \item $e_{145}g_2->e_{143}h_1g_2$
  \item $e_{148}h_0c_2->e_{147}h_0^2c_2$
  \item $e_{149}f_1->e_{141}h_3f_1$
  \item $e_{150}c_0h_5->e_{143}d_0h_5$
  \item $e_{151}e_1->e_{139}f_0h_5$
  \item $e_{156}p_0->e_{155}h_0p_0$
  \item $e_{157}d_1->e_{126}H_1$
  \item $e_{158}h_0^3h_5->e_{157}h_0^4h_5$
  \item $e_{169}g_1->e_{167}h_1g_1$
  \item $e_{171}f_0->e_{165}h_2g_1$
  \item $e_{172}e_0->e_{171}h_0e_0$
  \item $e_{174}h_0^3h_4->e_{173}h_0^4h_4$
  \item $e_{180}h_1c_0->e_{177}P^1h_2$
  \item $e_{182}h_0^3h_3->e_{179}P^1h_1$
\end{enumerate}
      \end{multicols}
    \item
      \begin{multicols}{2}
\begin{enumerate}
  \item $e_{11}c_1h_5h_7->e_7h_2c_1h_5h_7$
  \item $e_{12}h_0h_2c_4->e_{11}h_1^3c_4$
  \item $e_{16}h_0f_3->e_{15}h_0^2f_3$
  \item $e_{17}h_3e_3->e_9h_3^2e_3$
  \item $e_{17}g_2h_7->e_{15}h_1g_2h_7$
  \item $e_{20}h_0c_2h_7->e_{19}h_0^2c_2h_7$
  \item $e_{21}f_1h_7->e_{13}h_3f_1h_7$
  \item $e_{22}c_0h_5h_7->e_{15}d_0h_5h_7$
  \item $e_{23}e_1h_7->e_{11}f_0h_5h_7$
  \item $e_{24}h_0e_3->e_{23}h_0^2e_3$
  \item $e_{25}x_2->e_{23}h_1x_2$
  \item $e_{28}p_0h_7->e_{27}h_0p_0h_7$
  \item $e_{30}h_0^3h_5h_7->e_{29}h_0^4h_5h_7$
  \item $e_{32}V_0->e_{31}h_0V_0$
  \item $e_{33}g_3h_6->e_{31}h_1g_3h_6$
  \item $e_{41}h_3d_3->e_{21}h_2x_2$
  \item $e_{41}g_1h_7->e_{39}h_1g_1h_7$
  \item $e_{43}f_0h_7->e_{37}h_2g_1h_7$
  \item $e_{44}h_0p_2->e_{43}h_0^2p_2$
  \item $e_{44}e_0h_7->e_{43}h_0e_0h_7$
  \item $e_{45}h_1p_1'->e_{30}h_1V_0$
  \item $e_{46}h_0^3h_4h_7->e_{45}h_0^4h_4h_7$
  \item $e_{48}h_0d_3->e_{47}h_0^2d_3$
  \item $e_{49}Q_3(1)->e_{47}h_1Q_3(1)$
  \item $e_{49}n_2->e_{47}h_1n_2$
  \item $e_{52}h_1c_0h_7->e_{49}P^1h_2h_7$
  \item $e_{54}h_0^3h_3h_7->e_{51}P^1h_1h_7$
  \item $e_{54}c_0h_6^2->e_{47}d_0h_6^2$
  \item $e_{55}h_3D_3(1)->e_{39}\bar{\lambda_{39}D_1(1)}$ 
  \item $e_{59}H_1(1)->e_{55}h_2H_1(1)$
  \item $e_{59}h_2D_3(1)->e_{45}h_2Q_3(1)$
  \item $e_{60}J_0->e_{59}h_0J_0$
  \item $e_{62}h_0^3h_6^2->e_{61}h_0^4h_6^2$
  \item $e_{81}g_2h_6->e_{79}h_1g_2h_6$
  \item $e_{84}h_0c_2h_6->e_{83}h_0^2c_2h_6$
  \item $e_{85}f_1h_6->e_{77}h_3f_1h_6$
  \item $e_{88}h_0h_4c_3->e_{87}h_0^2h_4c_3$
  \item $e_{92}p_0h_6->e_{91}h_0p_0h_6$
  \item $e_{93}d_1h_6->e_{62}H_1h_6$
  \item $e_{96}h_0g_3->e_{95}h_0^2g_3$
  \item $e_{104}h_0f_2->e_{103}h_0^2f_2$
  \item $e_{105}h_2e_2->e_{101}h_2^2e_2$
  \item $e_{105}g_1h_6->e_{103}h_1g_1h_6$
  \item $e_{107}f_0h_6->e_{101}h_2g_1h_6$
  \item $e_{108}e_0h_6->e_{107}h_0e_0h_6$
  \item $e_{108}h_0e_2->e_{107}h_0^2e_2$
  \item $e_{109}x_1->e_{93}h_1^2g_3$
  \item $e_{110}h_0^3h_4h_6->e_{109}h_0^4h_4h_6$
  \item $e_{112}h_4D_3->e_{110}h_1h_4D_3$
  \item $e_{113}h_3d_2->e_{105}d_1h_4^2$
  \item $e_{116}c_0h_1h_6->e_{113}P^1h_2h_6$
  \item $e_{118}h_0^3h_3h_6->e_{115}P^1h_1h_6$
  \item $e_{121}Q_3->e_{119}h_1Q_3$
  \item $e_{121}n_1->e_{109}t_1$
  \item $e_{124}h_2D_3->e_{117}h_2Q_3$
  \item $e_{126}h_0^3h_5^2->e_{125}h_0^4h_5^2$
  \item $e_{136}D_1->e_{129}D_2$
  \item $e_{137}g_1h_5->e_{135}h_1g_1h_5$
  \item $e_{140}e_0h_5->e_{139}h_0e_0h_5$
  \item $e_{144}h_0g_2->e_{143}h_0^2g_2$
  \item $e_{148}h_0f_1->e_{147}h_0^2f_1$
  \item $e_{148}h_1c_0h_5->e_{145}P^1h_2h_5$
  \item $e_{149}h_3d_1->e_{140}h_0h_2g_2$
  \item $e_{150}h_0^3h_3h_5->e_{147}P^1h_1h_5$
  \item $e_{151}x_0->e_{138}h_1e_0h_5$
  \item $e_{153}h_2d_1>e_{149}h_2^2d_1$
  \item $e_{158}h_0^3h_4^2->e_{155}q_0$
  \item $e_{164}h_1c_0h_4->e_{161}h_2^2g_1$
  \item $e_{168}h_0g_1->e_{167}h_0^2g_1$
  \item $e_{170}h_1e_0->e_{164}h_2^2e_0$
  \item $e_{173}h_1d_0->e_{171}h_1^2d_0$
  \item $e_{174}h_0d_0->e_{173}h_0^2d_0$
  \item $e_{188}h_0^5->e_{187}h_0^6$
\end{enumerate}
      \end{multicols}
    \item
      \begin{multicols}{2}
\begin{enumerate}
  \item $e_9h_2^2c_4<->\widehat{h}_5h_5D_3(1)$
  \item $e_{13}h_1^2c_4<->\widehat{h}_4h_0^2c_4$
  \item $e_{13}h_2f_3<->\widehat{h}_4h_1f_3$
  \item $e_{29}d_1h_7<->\gamma_{61}h_7$
  \item $e_{47}T_0<->\widehat{g}_4h_0^2$
\end{enumerate}
      \end{multicols}
  \end{enumerate}
\end{note}
\begin{note}
  \label{eq:P:Ext:5-194}
  \begin{enumerate}
    \item
      \begin{multicols}{2}
\begin{enumerate}
  \item $e_2g_4->e_1h_0g_4$
  \item $e_{14}h_2c_4->e_{13}h_0h_2c_4$
  \item $e_{16}h_1c_4->e_{14}h_1^2c_4$
  \item $e_{18}f_3->e_{17}h_0f_3$
  \item $e_{18}h_4^3h_7->e_{17}h_0h_4^3h_7$
  \item $e_{22}c_2h_7->e_{21}h_0c_2h_7$
  \item $e_{24}h_1h_3h_5h_7->e_{22}h_2^3h_5h_7$
  \item $e_{26}h_2^2h_5h_7->e_{23}c_0h_5h_7$
  \item $e_{26}e_3->e_{25}h_0e_3$
  \item $e_{32}h_1h_4^2h_7->e_{29}p_0h_7$
  \item $e_{32}h_0^2h_5h_7->e_{31}h_0^3h_5h_7$
  \item $e_{44}c_1h_7->e_{30}d_1h_7$
  \item $e_{46}p_2->e_{45}h_0p_2$
  \item $e_{48}h_0^2h_4h_7->e_{47}h_0^3h_4h_7$
  \item $e_{48}p_1'->e_{46}h_1p_1'$
  \item $e_{50}d_3->e_{49}h_0d_3$
  \item $e_{56}h_0^2h_3h_7->e_{55}h_0^3h_3h_7$
  \item $e_{56}h_1h_3h_6^2->e_{54}h_1^2h_3h_6^2$
  \item $e_{58}h_2^2h_6^2->e_{55}c_0h_6^2$
  \item $e_{64}h_0^2h_6^2->e_{63}h_0^3h_6^2$
  \item $e_{64}D_3(1)->e_{61}J_0$
  \item $e_{86}c_2h_6->e_{85}h_0c_2h_6$
  \item $e_{90}h_4c_3->e_{89}h_0h_4c_3$
  \item $e_{96}h_1h_5^3->e_{93}p_0h_6$
  \item $e_{98}g_3->e_{97}h_0g_3$
  \item $e_{104}h_1c_3->e_{102}h_1^2c_3$
  \item $e_{106}h_3^3h_6->e_{103}c_0h_4h_6$
  \item $e_{106}f_2->e_{105}h_0f_2$
  \item $e_{108}c_1h_6->e_{94}d_1h_6$
  \item $e_{110}e_2->e_{109}h_0e_2$
  \item $e_{112}h_0^2h_4h_6->e_{111}h_0^3h_4h_6$
  \item $e_{120}h_0^2h_3h_6->e_{119}h_0^3h_3h_6$
  \item $e_{120}p_1->e_{118}h_1p_1$
  \item $e_{121}p_0'->e_{119}p_0'h_1$
  \item $e_{122}d_2->e_{121}h_0d_2$
  \item $e_{128}h_0^2h_5^2->e_{127}h_0^3h_5^2$
  \item $e_{129}D_3->e_{127}h_1D_3$
  \item $e_{140}c_1h_5->e_{125}h_2D_3$
  \item $e_{142}h_3c_2->e_{127}H_1$
  \item $e_{146}g_2->e_{145}h_0g_2$
  \item $e_{150}f_1->e_{149}h_0f_1$
  \item $e_{152}e_1->e_{150}h_1e_1$
  \item $e_{152}h_0^2h_3h_5->e_{151}h_0^3h_3h_5$
  \item $e_{156}h_0^2h_2h_5->e_{124}h_1^3h_5^2$
  \item $e_{160}h_0^2h_4^2->e_{159}h_0^3h_4^2$
  \item $e_{168}h_2c_1->e_{158}n_0$
  \item $e_{170}g_1->e_{169}h_0g_1$
  \item $e_{172}h_1^2h_4->e_{141}e_0h_5$
  \item $e_{172}f_0->e_{171}h_0f_0$
  \item $e_{173}e_0->e_{165}h_3e_0$
  \item $e_{176}d_0->e_{175}h_0d_0$
  \item $e_{181}h_1c_0->e_{174}h_1d_0$
  \item $e_{190}h_0^4->e_{189}h_0^5$
\end{enumerate}
      \end{multicols}
    \item
      \begin{multicols}{2}
\begin{enumerate}
  \item $e_{10}h_2c_4->e_7c_0c_4$
  \item $e_{14}h_5p_2->e_{13}h_0h_5p_2$
  \item $e_{16}h_0^2c_4->e_{15}h_0^3c_4$
  \item $e_{16}h_1f_3->e_{14}h_1^2f_3$
  \item $e_{18}g_2h_7->e_{17}h_0g_2h_7$
  \item $e_{22}f_1h_7->e_{21}h_0f_1h_7$
  \item $e_{24}e_1h_7->e_{22}h_1e_1h_7$
  \item $e_{24}h_0^2h_3h_5h_7->e_{23}h_0^3h_3h_5h_7$
  \item $e_{24}h_1e_3->e_{22}h_1^2e_3$
  \item $e_{32}h_5D_3(1)->e_{29}h_5J_0$
  \item $e_{32}h_0^2h_4^2h_7->e_{31}h_0^3h_4^2h_7$
  \item $e_{33}V_0->e_{31}h_1V_0$
  \item $e_{34}g_3h_6->e_{33}h_0g_3h_6$
  \item $e_{40}h_2c_1h_7->e_{30}n_0h_7$
  \item $e_{42}h_2p_2->e_{39}\bar{\lambda_{40}D_1(1)}$
  \item $e_{42}g_1h_7->e_{41}h_0g_1h_7$
  \item $e_{44}f_0h_7->e_{43}h_0f_0h_7$
  \item $e_{44}h_1^3h_4h_7->e_{31}e_0h_5h_7$
  \item $e_{45}e_0h_7->e_{38}h_3e_0h_7$
  \item $e_{48}h_1d_3->e_{45}h_4J_0$
  \item $e_{48}d_0h_7->e_{47}h_0d_0h_7$
  \item $e_{48}T_0->e_{47}h_0T_0$
  \item $e_{50}Q_3(1)->e_{49}h_0Q_3(1)$
  \item $e_{50}n_2->e_{49}h_0n_2$
  \item $e_{53}h_1c_0h_7->e_{46}h_1d_0h_7$
  \item $e_{56}h_0^2h_3h_6^2->e_{55}h_0^3h_3h_6^2$
  \item $e_{56}h_3D_3(1)->e_{53}h_3J_0$
  \item $e_{62}h_0^4h_7->e_{61}h_0^5h_7$
  \item $e_{64}K_0->e_{63}h_0K_0$
  \item $e_{78}c_2h_3h_6->e_{63}H_1h_6$
  \item $e_{82}g_2h_6->e_{81}h_0g_2h_6$
  \item $e_{86}f_1h_6->e_{85}h_0f_1h_6$
  \item $e_{88}e_1h_6->e_{86}h_1e_1h_6$
  \item $e_{90}h_4f_2->e_{89}h_0h_4f_2$
  \item $e_{94}h_2g_3->e_{93}h_0h_2g_3$
  \item $e_{96}h_1g_3->e_{94}h_1^2g_3$
  \item $e_{96}h_0^2h_5^3->e_{95}h_0^3h_5^3$
  \item $e_{104}h_0^2c_3->e_{103}h_0^3c_3$
  \item $e_{104}h_2c_1h_6->e_{94}h_6n_0$
  \item $e_{104}h_1f_2->e_{102}h_1^2f_2$
  \item $e_{106}h_2e_2->e_{102}h_2^2e_2$
  \item $e_{106}g_1h_6->e_{105}h_0g_1h_6$
  \item $e_{108}f_0h_6->e_{107}h_0f_0h_6$
  \item $e_{109}e_0h_6->e_{101}h_3e_0h_6$
  \item $e_{110}x_1->e_{109}h_0x_1$
  \item $e_{112}d_0h_6->e_{111}h_0d_0h_6$
  \item $e_{113}h_4D_3->e_{111}h_1h_4D_3$
  \item $e_{114}h_3d_2->e_{113}h_0h_3d_2$
  \item $e_{117}h_1c_0h_6->e_{110}h_1d_0h_6$
  \item $e_{120}h_0p_0'->e_{119}h_0^2p_0'$
  \item $e_{122}n_1->e_{121}h_0n_1$
  \item $e_{122}Q_3->e_{121}h_0Q_3$
  \item $e_{126}h_0^4h_6->e_{125}h_0^5h_6$
  \item $e_{136}h_2c_1h_5->e_{126}h_5n_0$
  \item $e_{137}D_1->e_{127}A$
  \item $e_{138}g_1h_5->e_{137}h_0g_1h_5$
  \item $e_{140}f_0h_5->e_{139}h_0f_0h_5$
  \item $e_{142}h_2g_2->e_{141}h_0h_2g_2$
  \item $e_{144}h_1g_2->e_{136}h_1g_1h_5$
  \item $e_{144}d_0h_5->e_{143}h_0d_0h_5$
  \item $e_{148}h_0^2c_2->e_{127}A'$
  \item $e_{149}h_1c_0h_5->e_{142}h_1d_0h_5$
  \item $e_{152}x_0->e_{151}h_0x_0$
  \item $e_{154}h_2d_1->e_{150}h_2^2d_1$
  \item $e_{156}h_1d_1->e_{125}h_1H_1$
  \item $e_{158}h_0^4h_5->e_{157}h_0^5h_5$
  \item $e_{166}h_2g_1->e_{165}h_0h_2g_1$
  \item $e_{168}h_1g_1->e_{162}h_2^2g_1$
  \item $e_{172}h_0e_0->e_{171}h_0^2e_0$
  \item $e_{174}h_0^4h_4->e_{173}h_0^5h_4$
  \item $e_{178}P^1h_2->e_{177}h_0P^1h_2$
  \item $e_{180}P^1h_1->e_{178}h_1P^1h_1$
\end{enumerate}
      \end{multicols}
    \item
      \begin{multicols}{2}
\begin{enumerate}
  \item $e_1D_3h_7<->\widehat{h}_1D_3h_7$
  \item $e_{12}c_1h_5h_7<->\xi_{31}h_5h_7$
  \item $e_{14}h_3c_2h_7<->\widehat{H}_1h_7$
  \item $e_{28}h_1^3h_5h_7<->\widehat{h}_5h_0^3h_5h_7$
  \item $e_{60}H_1(1)<->\widehat{h}_7H_1$
  \item $e_{60}h_2D_3(1)<->\widehat{h}_7h_1D_3$
  \item $e_{60}h_1^3h_6^2<->\widehat{h}_6h_0^3h_6^2$
\end{enumerate}
      \end{multicols}
  \end{enumerate}
\end{note}
\renewcommand{\labelenumi}{(\arabic{enumi})}
\renewcommand{\labelenumii}{(\alph{enumii})}
\noindent
The remaining \eqref{eq:P:Ext}~(a) (\eqref{eq:P:Ext}~(8)~(a), \eqref{eq:P:Ext}~(9)~(a))
follow from
the results in \eqref{eq:P:Ext:5-192} (c) and \eqref{eq:P:Ext:5-193} (c) above.
This proves the remaining \eqref{eq:P:Ext}~(a).

To prove the remaining \eqref{eq:P:Ext}~(c) (\eqref{eq:P:Ext}~(8)~(c), \eqref{eq:P:Ext}~(9)~(c)),
\eqref{eq:P:Ext:nontriviality}, and \eqref{eq:Ext} (9)
we recall
in \eqref{eq:leading-term} and \eqref{eq:Ext:type} below
from \cite{lin_ext_2008,chen_determination_2011}
an important notion about the spectral sequence $\{E_r^{i,s,t}\}$ in \eqref{eq:P:Ext:spectral-sequence}.
\begin{note}
  \label{eq:leading-term}
  Given a (possibly zero) class $\alpha\in\Ext_\A^{s,t}$ with $t-s>0$ and $s\ge1$.
  Let $x$ be a cycle in $\Lambda^{s,t-s}$ representing the class $\alpha$
  with $x\ne0$.
  Since $t-s>0$ we have
  $x\equiv\lambda_j\lambda_I\mod{\Lambda(j-1)}$ for some $j\ge1$,
  where $\lambda_I\in\Lambda^{s-1,t-s-j}$ is a nonzero chain
  such that for any admissible monomial $\lambda_{j_1}\cdots\lambda_{j_{s-1}}$
  appearing in the admissible expansion of $\lambda_I$,
  $\lambda_j\lambda_{j_1}\cdots\lambda_{j_{s-1}}$ is admissible.
\end{note}
\noindent
Since $x$ is a cycle and $\Lambda(j-1)$ is a subcomplex of $\Lambda$,
$\lambda_I$ in \eqref{eq:leading-term} is actually a cycle.
Recall the map $\Ext_\A^{s+1,t+1}\to{\phi_\ast}\Ext_\A^{s,t}(P)$
in \eqref{eq:phi}.
It is easy to see
from \eqref{eq:phi} (2)
that
\begin{align*}
  \phi(x)\equiv e_j\lambda_I\mod{F(j-1)}.
\end{align*}
Since $\lambda_I$ is a cycle,
$e_j\{\lambda_I\}$ is an element in $E_1^{\ast,s-1,\ast}$
in the spectral sequence considered in \eqref{eq:P:Ext:spectral-sequence}.
And since $\phi(x)$ is a cycle,
$e_j\{\lambda_I\}$ is an infinite cycle in the spectral sequence
which may not be a nontrivial one.
To show that $\alpha\ne0$
it suffices to show that $\phi_\ast(\alpha)\ne0$.
It is not difficult to see that there are two cases for showing this.
\begin{note}
  \label{eq:Ext:type}
  \begin{enumerate}
    \item If $e_j\{\lambda_I\}$ is non-zero in the spectral sequence, then
      $\alpha$ is a non-zero class.
      In this case
      we will use $\bar{\phi}(\alpha)$ to denote $e_j\{\lambda_I\}$,
      and use the correspondence
      \begin{align*}
        \alpha=\{x\}=\{\bar{\lambda_j\lambda_I}\}->e_j\{\lambda_I\}=\bar{\phi}(\alpha),
      \end{align*}
      to indicate that the class $\phi_\ast(\alpha)$ is represented by the infinite cycle
      $e_j\{\lambda_I\}=\bar{\phi}(\alpha)$
      (in the spectral sequence in \eqref{eq:P:Ext:spectral-sequence} for $\Ext_\A^{\ast,\ast}(P)$).
    \item If $e_j\{\lambda_I\}$ is a boundary in the spectral sequence, but
      $\phi(x)$ is homologous to a cycle $z$
      such that $z\equiv e_k\lambda_J\mod{F(k-1)}$
      for some $k$ with $1\le k<j$ and some $\lambda_J\in\Lambda^{s-1,\ast}$,
      and such that $e_k\{\lambda_J\}$ is non-zero in the spectral sequence,
      then $\alpha$ is also a non-zero class.
      In this case
      we will use $\widetilde{\phi}(\alpha)$ to denote $e_k\{\lambda_J\}$,
      and use the correspondence
      \begin{align*}
        \alpha=\{x\}=\{\bar{\lambda_j\lambda_I}\}->e_k\{\lambda_J\}=\widetilde{\phi}(\alpha).
      \end{align*}
      to indicate that the class $\phi_\ast(\alpha)$ is represented by the infinite cycle
      $e_k\{\lambda_J\}=\widetilde{\phi}(\alpha)$.
  \end{enumerate}
\end{note}

\eqref{eq:Ext:6-193} and \eqref{eq:Ext:6-194} below
give the correspondences (defined in \eqref{eq:Ext:type})for the elements
in the remaining \eqref{eq:P:Ext}~(c).
\begin{note}
  \label{eq:Ext:6-193}
\begin{enumerate}
  \item $h_5V_0->e_{30}D_3(1)h_5=\bar{\phi}(h_5V_0)$
\end{enumerate}
\end{note}
\begin{note}
  \label{eq:Ext:6-194}
\begin{enumerate}
  \item $h_0^2g_4->e_{47}T_0=\bar{\phi}(h_0^2g_4)$
  \item $h_5^2D_3(1)->e_9h_2^2c_4=\widetilde{\phi}(h_5^2D_3(1))$ 
\end{enumerate}
\end{note}
\noindent
\eqref{eq:Ext:6-193} shows that $\phi_\ast(h_5V_0)\ne0$,
and \eqref{eq:Ext:6-194} shows that
$\phi_\ast(h_0^2g_4),\phi_\ast(h_5^2D_3(1))$ are linearly independent.
This proves the remaining \eqref{eq:P:Ext}~(c).

To prove \eqref{eq:P:Ext:nontriviality}, we list
in \eqref{eq:Ext:6:nonzero} below
each $\beta\in\Ext_\A^{6,\ast}$ such that
$e_i\beta$ appears in the targets of the (b)-differentials
in \eqref{eq:P:Ext:5-192} and \eqref{eq:P:Ext:5-193}.
The calculations for these correspondences
are lengthy and will not be given in this paper.
\begin{note}
  \label{eq:Ext:6:nonzero}
\begin{enumerate}
  \item $A->e_9D_1 =\bar{\phi}(A)$
  \item $A'->e_{20}h_0^2c_2 =\bar{\phi}(A')$
  \item $\bar{\lambda_{19}D_1}->e_{19}D_1=\bar{\phi}(\bar{\lambda_{19}D_1})$ 
  \item $\bar{\lambda_{39}D_1(1)}->e_{39}D_1(1)=\bar{\phi}(\bar{\lambda_{39}D_1(1)})$ 
  \item $\bar{\lambda_{40}D_1(1)}->e_{40}D_1(1)=\bar{\phi}(\bar{\lambda_{40}D_1(1)})$
  \item $D_2->e_{28}h_0^3h_4^2=\bar{\phi}(D_2)$
  \item $q_0->e_2h_0^2h_4^2=\widetilde{\phi}(q)$ 
  \item $t_0->e_5n_0=\widetilde{\phi}(t_0)$

  \item $h_0^3h_3h_5h_7 ->e_{84}h_1^3h_4h_6=\bar{\phi}(h_0^3h_3h_5h_7)$
  \item $h_0^3h_3h_6^2  ->e_{53}h_1^3h_4h_6=\bar{\phi}(h_0^3h_3h_6^2)$
  \item $h_0^3h_5^3     ->e_{28}h_1^3h_5^2=\bar{\phi}(h_0^3h_4^2h_6)$
  \item $h_0^3h_4^2h_7  ->e_{92}h_1^3h_5^2=\bar{\phi}(h_0^3h_4^2h_7)$
  \item $h_0^4h_4h_6    ->e_{40}h_0^3h_3h_5=\widetilde{\phi}(h_0^4h_4h_6)$
  \item $h_0^4h_4h_7    ->e_{104}h_0^3h_3h_5=\widetilde{\phi}(h_0^4h_4h_7)$
  \item $h_0^4h_5^2     ->e_{24}h_0^3h_3h_5=\widetilde{\phi}(h_0^4h_5^2)$
  \item $h_0^4h_5h_7    ->e_{88}h_0^3h_3h_6=\widetilde{\phi}(h_0^4h_5h_7)$
  \item $h_0^4h_6^2     ->e_{56}h_0^3h_3h_6=\widetilde{\phi}(h_0^4h_6^2)$
  \item $h_0^5h_4       ->e_6P^1h_1=\widetilde{\phi}(h_0^5h_4)$
  \item $h_0^5h_5       ->e_{22}P^1h_1=\widetilde{\phi}(h_0^5h_5)$
  \item $h_0^5h_6       ->e_{54}P^1h_1=\widetilde{\phi}(h_0^5h_6)$
  \item $h_0^5h_7       ->e_{118}P^1h_1=\widetilde{\phi}(h_0^5h_7)$

  \item $h_0^2c_2h_6    ->e_{17}h_1^2c_3=\bar{\phi}(h_0^2c_2h_6)$
  \item $h_0^2c_2h_7    ->e_{81}h_1^2c_3=\bar{\phi}(h_0^2c_2h_7)$
  \item $h_0^2h_4c_3    ->e_{13}h_1^2c_3=\bar{\phi}(h_0^2h_4c_3)$
  \item $h_0^3c_3       ->e_{20}h_1^3h_5^2=\bar{\phi}(h_0^3c_3)$
  \item $h_0^3c_4       ->e_{44}h_1^3h_6^2=\bar{\phi}(h_0^3c_4)$
  \item $h_1^3c_4       ->e_{41}h_2^3h_6^2=\bar{\phi}(h_1^3c_4)$
  \item $h_1c_0h_5^2    ->e_{18}h_2c_1h_5=\bar{\phi}(h_1c_0h_5^2)$
  \item $h_1c_0h_6^2    ->e_{50}h_2c_1h_6=\bar{\phi}(h_1c_0h_6^2)$
  \item $h_2c_1h_6^2    ->e_{37}h_3c_2h_6=\bar{\phi}(h_2c_1h_6^2)$
  \item $h_1c_0h_4h_6   ->e_{34}h_2c_1h_5=\bar{\phi}(h_1c_0h_4h_6)$
  \item $h_1c_0h_4h_7   ->e_{98}h_2c_1h_5=\bar{\phi}(h_1c_0h_4h_7)$
  \item $h_1c_0h_5h_7   ->e_{82}h_2c_1h_6=\bar{\phi}(h_1c_0h_5h_7)$   
  \item $h_2c_1h_5h_7   ->e_{69}h_3c_2h_6=\bar{\phi}(h_2c_1h_5h_7)$
  \item $c_0c_4         ->e_{36}c_1h_6^2=\bar{\phi}(c_0c_4)$          
  \item $h_0^2d_3=h_4K_0 +D_1(1)h_5   ->e_{15}K_0=\widetilde{\phi}(e_{15}K_0)$
  \item $h_0^2d_0=h_2P^1h_2           ->e_3P^1h_2=\bar{\phi}(h_2P^1h_2)$
  \item $h_1^2d_0       ->e_1h_1d_0=\widetilde{\phi}(h_1^2d_0)$
  \item $h_2^2d_1       ->e_3h_2d_1=\widetilde{\phi}(h_2^2d_1)$
  \item $d_0h_6^2       ->e_{45}d_1h_6=\bar{\phi}(d_0h_6^2)$
  \item $h_0d_0h_5      ->e_{12}h_1d_1=\bar{\phi}(h_0d_0h_5)$
  \item $h_0d_0h_6      ->e_{44}h_1d_1=\bar{\phi}(h_0d_0h_6)$
  \item $h_0d_0h_7      ->e_{108}h_1d_1=\bar{\phi}(h_0d_0h_7)$
  \item $h_1d_1h_6=h_0p_0h_6          ->e_{25}h_2d_2=\bar{\phi}(h_1d_1h_6)$
  \item $h_1d_1h_7=h_0p_0h_7          ->e_{89}h_2d_2=\bar{\phi}(h_1d_1h_7)$
  \item $h_1d_0h_5      ->e_{11}h_2d_1=\bar{\phi}(h_1d_0h_5)$
  \item $h_1d_0h_6      ->e_{43}h_2d_1=\bar{\phi}(h_1d_0h_6)$
  \item $h_1d_0h_7      ->e_{107}h_2d_1=\bar{\phi}(h_1d_0h_7)$
  \item $d_0h_5h_7      ->e_{77}d_1h_6=\bar{\phi}(d_0h_5h_7)$
  \item $h_0h_3d_2      ->e_4h_2d_2=\widetilde{\phi}(h_0h_3d_2)$
  \item $h_1h_4d_3      ->e_9h_3d_3=\widetilde{\phi}(h_1h_4d_3)$
  \item $h_0h_4d_3      ->e_{14}h_1d_3=\widetilde{\phi}(h_0h_4d_3)$
  \item $h_0h_5d_3      ->e_{30}h_1d_3=\bar{\phi}(h_0h_5d_3)$

  \item $h_0^2e_0 =h_0h_2d_0          ->e_2h_1d_0=\bar{\phi}(h_0h_2d_0)$
  \item $h_1^2e_1=h_1h_3d_1 =h_0^2f_1=h_0h_3p_0  ->e_3x_0=\widetilde{\phi}(h_1h_3d_1)$
  \item $h_2^2e_2=h_2h_4d_2 =h_1^2f_2=h_1h_4p_1  ->e_7x_1=\widetilde{\phi}(h_2h_4d_2)$
  \item $h_3^2e_3=h_3h_5d_3 =h_2^2f_3=h_2h_5p_2  ->e_{15}x_2=\widetilde{\phi}(h_3h_5d_3)$

  \item $h_0h_3g_1=h_2^2e_0           ->e_2h_1g_1=\bar{\phi}(h_2^2e_0)$

  \item $h_0^2e_2 =h_1x_1             ->e_{27}h_2c_1h_5=\bar{\phi}(h_0^2e_2)$       
  \item $h_0^2e_3       ->e_{61}f_1h_6=\widetilde{\phi}(h_0^2e_3)$    
  \item $h_1^2e_3 =h_2x_2             ->e_{55}h_3c_2h_6=\bar{\phi}(h_1^2e_3)$       
  \item $h_0e_0h_5=h_2d_0h_5          ->e_9h_3d_1=\bar{\phi}(h_0e_0h_5)$
  \item $h_0e_0h_6=h_2d_0h_6          ->e_{41}h_3d_1=\bar{\phi}(h_0e_0h_6)$
  \item $h_0e_0h_7=h_2d_0h_7          ->e_{105}h_3d_1=\bar{\phi}(h_0e_0h_7)$
  \item $h_1e_1h_6=h_3d_1h_6          ->e_{19}h_4d_2=\bar{\phi}(h_1e_1h_6)$
  \item $h_1e_1h_7=h_3d_1h_7          ->e_{85}h_4d_2=\bar{\phi}(h_1e_1h_7)$
  \item $h_2e_1h_6=h_1f_1h_6          ->e_{17}h_3e_2=\bar{\phi}(h_2e_1h_6)$
  \item $h_2e_1h_7=h_1f_1h_7          ->e_{81}h_3e_2=\bar{\phi}(h_2e_1h_7)$
  \item $h_3e_0h_6      ->e_{34}c_1h_2h_5=\bar{\phi}(h_3e_0h_6)$
  \item $h_3e_0h_7      ->e_{98}c_1h_2h_5=\bar{\phi}(h_3e_0h_7)$
  \item $e_0h_5h_7      ->e_{74}e_1h_6=\bar{\phi}(e_0h_5h_7)$

  \item $h_0^2f_2=h_0h_4p_0'          ->e_{14}h_1p_0'=\bar{\phi}(h_0^2f_2)$
  \item $h_0^2f_3=h_5T_0              ->e_{31}T_0=\bar{\phi}(h_5T_0)$
  \item $h_1^2f_3=h_1h_5p_1'          ->e_{29}h_2p_1'=\bar{\phi}(h_1^2f_3)$
  \item $h_0f_0h_5=h_1e_0h_5          ->e_8h_2e_1=\bar{\phi}(h_0f_0h_5)$
  \item $h_0f_0h_6=h_1e_0h_6          ->e_{40}h_2e_1=\bar{\phi}(h_0f_0h_6)$
  \item $h_0f_0h_7=h_1e_0h_7          ->e_{104}h_2e_1=\bar{\phi}(h_0f_0h_7)$
  \item $h_0h_4f_2      ->e_{14}h_1f_2=\bar{\phi}(h_0h_4f_2)$
  \item $h_1h_4f_2      ->e_{13}h_2f_2=\bar{\phi}(h_1h_4f_2)$
  \item $h_0f_1h_6=h_3p_0h_6          ->e_{18}p_1h_4=\bar{\phi}(h_0f_1h_6)$
  \item $h_0f_1h_7=h_3p_0h_7          ->e_{82}p_1h_4=\bar{\phi}(h_0f_1h_7)$
  \item $f_0h_5h_7      ->e_{73}f_1h_6=\bar{\phi}(f_0h_5h_7)$
  \item $h_2f_0h_7 =h_1g_1h_7         ->e_{101}h_3f_1=\bar{\phi}(h_2f_0h_7)$
  \item $h_3f_1h_6 =h_2g_2h_6 =h_2g_3h_4 =e_1h_3^2  ->e_{11}h_4f_2=\bar{\phi}(h_3f_1h_6)$
  \item $h_3f_1h_7=Sq^0(h_2f_0h_6)    ->e_{75}h_4f_2=\bar{\phi}(h_3f_1h_7)$

  \item $h_0^2g_1=h_2^2d_0=h_0h_2e_0  ->e_2h_1e_0=\widetilde{\phi}(h_2^2d_0)$
  \item $h_0^2g_2=h_3x_0              ->e_7x_0=\bar{\phi}(h_3x_0)$
  \item $h_0^2g_3                     ->e_{22}h_1p_0'=\widetilde{\phi}(h_0^2g_3)$
  \item $h_1^2g_3=Sq^0(h_0^2g_2)      ->e_{15}x_1=\bar{\phi}(h_4x_1)$
  \item $h_4P^1h_2=h_2^2g_1           ->e_3h_2g_1=\bar{\phi}(h_2^2g_1)$
  \item $h_0g_1h_5                    ->e_6h_1g_2=\bar{\phi}(h_0g_1h_5)$
  \item $h_1h_4g_3 =h_1g_2h_6=e_1h_3h_6             ->e_{13}h_2g_3=\bar{\phi}(h_1h_4g_3)$
  \item $h_1g_2h_7=h_3e_1h_7          ->e_{77}h_4e_2=\bar{\phi}(h_1g_2h_7)$
  \item $h_0g_1h_6=h_2e_0h_6          ->e_{38}e_1h_3=\bar{\phi}(h_0g_1h_6)$
  \item $h_0g_1h_7=h_2e_0h_7          ->e_{102}e_1h_3=\bar{\phi}(h_0g_1h_7)$
  \item $h_0g_2h_6                    ->e_{14}h_1g_3=\bar{\phi}(h_0g_2h_6)$
  \item $h_0g_2h_7                    ->e_{78}h_1g_3=\bar{\phi}(h_0g_2h_7)$
  \item $h_1g_3h_6=h_1h_4d_3          ->e_9h_3d_3=\widetilde{\phi}(h_1h_4d_3)$
  \item $h_0g_3h_6=h_0h_4d_3          ->e_{14}h_1d_3=\widetilde{\phi}(h_0h_4d_3)$
  \item $h_0h_2g_2=h_0h_3f_1          ->e_6h_1f_1=\widetilde{\phi}(h_0h_2g_2)$
  \item $h_1g_1h_5                    ->e_5h_2g_2=\bar{\phi}(h_1g_1h_5)$
  \item $h_1g_1h_6=h_2f_0h_6          ->e_{37}h_3f_1=\bar{\phi}(h_1g_1h_6)$
  \item $h_1g_1h_7=h_2f_0h_7          ->e_{101}h_3f_1=\bar{\phi}(h_1g_1h_7)$
  \item $h_2g_1h_6                    ->e_{35}h_3g_2=\bar{\phi}(h_2g_1h_6)$
  \item $h_2g_1h_7                    ->e_{99}h_3g_2=\bar{\phi}(h_2g_1h_7)$
  \item $h_0h_2g_1                    ->e_2h_1g_2=\bar{\phi}(h_0h_2g_1)$
  \item $h_0h_2g_3=h_0h_4e_2          ->e_{12}h_2e_2=\bar{\phi}(h_0h_2g_3)$

  \item $h_0^2p_2                     ->e_{49}h_4e_2=\widetilde{\phi}(h_0^2p_2)$  
  \item $h_0h_5p_2                    ->e_{28}h_2p_2=\widetilde{\phi}(h_0h_5p_2)$

  \item $h_4^2D_3                     ->e_4h_1^2c_3=\widetilde{\phi}(h_4^2D_3)$
  \item $h_3^2D_3(1)                  ->e_7h_3D_3(1)=\widetilde{\phi}(h_3^2D_3(1))$
  \item $h_1h_4D_3                    ->e_1h_4D_3=\widetilde{\phi}(h_1h_4D_3)$

  \item $h_0^2p_0'                    ->e_{24}h_3e_1=\widetilde{\phi}(h_0^2p_0')$ 




  \item $h_0J_0->e_{57}h_1p_1=\bar{\phi}(h_0J_0)$
  \item $h_3J_0->e_{50}h_4p_0=\bar{\phi}(h_3J_0)$
  \item $h_4J_0->e_{15}J_0=\widetilde{\phi}(h_4J_0)$
  \item $h_5J_0->e_{30}H_1(1)=\widetilde{\phi}(h_5J_0)$

  \item $h_0T_0->e_{76}h_1^3h_5^2=\bar{\phi}(h_0T_0)$

  \item $h_0K_0->e_{44}h_0^2h_2h_4h_6=\bar{\phi}(h_0K_0)$

  \item $h_0n_1 =h_4D_1 =h_2^2D_3 ->e_3h_3D_3 =\widetilde{\phi}(h_4D_1)$
  \item $h_0n_2 ->e_{30}D_1(1) =\widetilde{\phi}(h_0n_2)$
  \item $h_1n_2 =h_5D_1(1) =h_3^2D_3(1) ->e_7h_4D_3(1) =\widetilde{\phi}(h_5D_1(1))$

  \item $h_0Q_3->e_{22}h_1g_2=\bar{\phi}(h_0Q_3)$
  \item $h_0Q_3(1)->e_{46}h_1g_3=\bar{\phi}(h_0Q_3(1))$
  \item $h_1Q_3->e_{21}h_2g_2=\bar{\phi}(h_1Q_3)$
  \item $h_1Q_3(1)->e_{45}h_2g_2=\bar{\phi}(h_1Q_3(1))$
  \item $h_2Q_3(1)->e_{43}h_3g_3=\bar{\phi}(h_2Q_3(1))$

  \item $h_0x_0->e_{13}h_1c_0h_4=\bar{\phi}(h_0x_0)$
  \item $h_1x_1->e_{27}h_2c_1h_5=\bar{\phi}(h_1x_1)$
  \item $h_2x_2->e_{55}h_3c_2h_6=\bar{\phi}(h_2x_2)$
  \item $h_0x_1 =h_2h_4D_3->e_3h_4D_3 =\widetilde{\phi}(h_2h_4D_3)$
  \item $h_1x_2 =h_3h_5D_3(1)->e_7h_5D_3(1) =\widetilde{\phi}(h_3h_5D_3(1))$

  \item $h_1V_0->e_{28}h_2D_3(1)=\bar{\phi}(h_1V_0)$
  \item $h_0V_0=h_5K_0->e_{15}T_0=\widetilde{\phi}(h_5K_0)$

  \item $h_1H_1->e_{11}D_1=\widetilde{\phi}(h_1H_1)$
  \item $h_2H_1(1)->e_{23}D_1(1)=\widetilde{\phi}(h_2H_1(1))$
\end{enumerate}
\end{note}
\noindent
It is not difficult to see that \eqref{eq:Ext:6:nonzero} implies \eqref{eq:P:Ext:nontriviality}.

To prove \eqref{eq:Ext} (9),
we give in \eqref{eq:Ext:7-195} below
the correspondences (defined in \eqref{eq:Ext:type}) for $h_0^3g_4,h_1h_5V_0\in\Ext_\A^{7,195}$.
\begin{note}
  \label{eq:Ext:7-195}
  \begin{enumerate}
    \item $h_0^3g_4->e_{45}e_0h_6^2=\widetilde{\phi}(h_0^3g_4)$
    \item $h_1h_5V_0->e_{29}h_2V_0=\bar{\phi}(h_1h_5V_0)$
  \end{enumerate}
\end{note}
\noindent
A check againt the differentials in \eqref{eq:P:Ext:5-194} (b)
shows that $e_{45}e_0h_6^2,e_{29}h_2V_0$ are not boundaries in
the spectral sequence.
It follows that $e_{45}e_0h_6^2,e_{29}h_2V_0$ are infinite cycles
which represent non-trivial classes in $\Ext_\A^{6,194}(P)$,
and these classes are linearly independent.
This proves \eqref{eq:Ext} (9).

This completes the proofs of \eqref{eq:Ext} and \eqref{eq:P:Ext}.

The remainder of this section is devoted to proving
\eqref{eq:P62:Ext}, \eqref{eq:P62-k:Ext}, \eqref{eq:P62t:Ext}, and \eqref{eq:P62-17t:Ext}
(simply ``\eqref{eq:P62:Ext} through \eqref{eq:P62-17t:Ext}'').

First we prove \eqref{eq:P62:Ext}.
For $0\le l\le m$ we can
calculate $\Ext_\A^{s,t}(P^m_l)$ using the
spectral sequence in \eqref{eq:P:Ext:spectral-sequence}
with the following modifications.
The mod $2$ reduced homology group
$ H_\ast(P^m_l)$ has the same right $\A$-module
structure as that of $ H_\ast(P)$ given in \eqref{eq:P:e_k}
except that $e_k=0$ for $k<l$ and for $k>m$.
Therefore we have the following.
\begin{note}
  \label{eq:P62-k:Ext:spectral-sequence}
  Let
  $\bar{F}(0)=0\subset\bar{F}(l)\subset\bar{F}(l+1)\cdots\subset\bar{F}(m)= H_\ast(P^m_l)$
  be the increasing filtration of subcomplexes given by
  $\bar{F}(n)=\Sum_{l\le k\le n} H_k(P^m_l)\tensor\Lambda$, $l\le n\le m$,
  then $\{\bar{F}(i)\}_{l\le i\le m}$ defines a spectral sequence
  $\{\bar{E}_r^{i,s,t}\}_{r\ge1}$
  with
  $\bar{E}_1^{i,s,t}\iso\Susp^i\Ext_\A^{s,t+s-i}$ and
  $\dsum_{i\ge1}\bar{E}_\infty^{i,s,t}\iso\Ext_\A^{s,t+s}(P^m_l)$,
  which is identical to the one in \eqref{eq:P:Ext:spectral-sequence} except that
  $\bar{E}_r^{i,s,t}=0$ for $i<l$ and for $i>m$.
\end{note}
\noindent
Computing the spectral sequence $\{\bar{E}_r^{i,s,t}\}_{r\ge1}$
is equivalent to computing $\{E_r^{i,s,t}\}_{r\ge1}$
(put $e_k=0$ for $k<l$ and for $k>m$ for $\bar{E}_r^{i,s,t}$.)
Thus the results
\eqref{eq:P62:Ext}~(i)~(8), (9) follow from the spectral sequence calculations in
\eqref{eq:P:Ext:5-192} and \eqref{eq:P:Ext:5-193}.
It is also easy to see the following.
\begin{note}
  \label{eq:P62:Ext:pullback}
  Let $P^m\to{i_m}P$ be the inclusion map and
  $\Ext_\A^{s,t}(P^m)\to{(i_m)_\ast}\Ext_\A^{s,t}(P)$
  be the induced homomorphism.
  If $\alpha\in\Ext_\A^{s,t}(P)$ is a non-zero class
  represented in the spectral sequence $\{E_r^{i,s,t}\}_{r\ge0}$ by
  an infinite cycle $e_k\beta$ with $k\le m$, then
  $\alpha$ can be pulled back to a non-zero class $\bar{\alpha}\in\Ext_\A^{s,t,}(P^m)$
  represented in the spectral sequence $\{\bar{E}_r^{i,s,t}\}_{r\ge0}$ also by
  the infinite cycle $e_k\beta$.
\end{note}
\noindent
From \eqref{eq:P:Ext:known} and \eqref{eq:P62:Ext:pullback}, it is not difficult to see that
each class in \eqref{eq:P:Ext},
except $\widehat{h}_6h_0h_6^2\in\Ext_\A^{3,192}(P)$ and $\widehat{h}_7D_3\in\Ext_\A^{4,192}(P)$,
can be pulled back to $\Ext_\A^{\ast,\ast}(P^{62})$.
The classes $\widehat{h}_6h_0h_6^2,\widehat{h}_7D_3$ can also be pulled back
to $\Ext_\A^{\ast,\ast}(P^{62})$
because of the following cycle representations for these classes.
\begin{note}
  \begin{enumerate}
    \item $e_{62}h_1h_6^2<->\widehat{h}_6h_0h_6^2$
    \item $e_{62}D_3(1)<->\widehat{h}_7D_3$
  \end{enumerate}
\end{note}
\noindent
Thus all the classes in \eqref{eq:P:Ext} for $\Ext_\A^{\ast,\ast}(P)$ can be pulled back
to $\Ext_\A^{\ast,\ast}(P^{62})$ via $\Ext_\A^{s,t}(P^{62})\to{(i_{62})_\ast}\Ext_\A^{s,t}(P)$.
And the pullbacks of $\widehat{h}_6h_0h_6^2 \in\Ext_\A^{3,192}(P)$ and $\widehat{h}_7D_3 \in\Ext_\A^{4,192}(P)$ are denoted by
$\bar{e_{62}h_1}h_6^2 \in\Ext_\A^{3,192}(P^{62})$ and $\bar{e_{62}D_3(1)} \in\Ext_\A^{3,192}(P^{62})$ respectively.

To prove \eqref{eq:P62:Ext}~(i)~($j$) for $1\le j\le9$ (simply \eqref{eq:P62:Ext}~(i)) and \eqref{eq:P62:Ext}~(ii),
consider the cofibration sequence
$\Susp(P_{m+1}=P/P^m)->P^m\to{i_m}P$ and the induced long exact sequence
of $\Ext$ groups.
\begin{align*}
  \cdots->\Ext_\A^{s-1,t}(P_{m+1})\to{\delta} \Ext_\A^{s,t}(P^m)
  \to{(i_m)_\ast} \Ext_\A^{s,t}(P)
  \to{} \Ext_\A^{s,t}(P_{m+1})
  ->\cdots.
\end{align*}
Since each of the classes in $\Ext_\A^{s,t}(P^{62})$ in \eqref{eq:P62:Ext}
is a pullback from $\Ext_\A^{s,t}(P)$ via $(i_m)_\ast$,
we need to show that for $m=62$ and for
each of the $(s,t)$ in $\eqref{eq:P62:Ext}$~(i)~(1) through (9),
the map $\Ext_\A^{s-1,t}(P_{63})\to{\delta}\Ext_\A^{s,t}(P^{62})$ is zero.
This is equivalent to saying that
in the spectral sequence $\{E_r^{i,s,t-s}\}_{r\ge1}$ for $\Ext_\A^{s,t}(P)$,
each of the differentials
$E_r^{i,s-1,t-s}\to{d_r}E_r^{i-r,s,t-s}$
that starts with the basis element $e_i\alpha$ with $i>62$
must either be a zero differential or hit a target $e_{i-r}\beta$ with $i-r>62$.
A check against the list of differentials in \cite{lin_ext_2008}
shows that this is indeed the case,
This proves \eqref{eq:P62:Ext}~(i) and \eqref{eq:P62:Ext}~(ii).
\eqref{eq:P62:Ext}~(iii) follows easily from 
\eqref{eq:P62:Ext}~(i) and \eqref{eq:P62:Ext}~(ii).
This proves \eqref{eq:P62:Ext}.

The result \eqref{eq:P62-k:Ext}~(1) through (16) also follows from the list of the differentials in \cite{lin_ext_2008}
once we describe cycle representations for the classes in \eqref{eq:P62-k:Ext} as follows.
\begin{note}
  \label{eq:P62-k:Ext:known}
  \begin{enumerate}
    \item $\bar{e_{23}} =\{e_{23}\} \in\Ext_\A^{0,23}(P^{62}_{17})$.
    \item $\bar{e_{31}} =\{e_{31}\} \in\Ext_\A^{0,31}(P^{62}_{17}) =\Ext_\A^{0,31}(P^{32})$.
    \item $\bar{e_{39}} =\{e_{39}\} \in\Ext_\A^{0,39}(P^{62}_{33})$.
    \item $\bar{e_{47}} =\{e_{47}\} \in\Ext_\A^{0,47}(P^{62}_{33})$.
    \item $\bar{e_{15}} =\{e_{15}\} \in\Ext_\A^{0,15}(P^{16}) =\Ext_\A^{0,15}(P^{32})$.
    \item $\bar{e_7} =\{e_7\} \in\Ext_\A^{0,7}(P^{16})$.
    \item $\bar{e_{62}h_1} =\{e_{62}\lambda_1 +e_{60}\lambda_3 +e_{56}\lambda_7 +e_{48}\lambda_{15} +e_{32}\lambda_{31}\} \in\Ext_\A^{1,64}(P^{62}_{47})$.
    \item $\bar{e_{48}h_0h_3^2} =\{e_{48}\lambda_0\lambda_7^2 +e_{47}(\lambda_9\lambda_3^2 +\lambda_3^2\lambda_9)\} \in\Ext_\A^{3,65}(P^{62}_{47})$.
    \item $\bar{e_{54}h_1h_3} =\{e_{54}\lambda_1\lambda_7 +e_{48}\lambda_7^2\} \in\Ext_\A^{2,64}(P^{62}_{47})$.
  \end{enumerate}
\end{note}
\noindent
To prove
\begin{enumerate}
  \item[\eqref{eq:P62-k:Ext}~(17)\quad]
    $(q_{47})_\ast(\widehat{g}_4h_0^2) =\bar{e_{48}h_0h_3^2}h_6^2\ne0$
    in $\Ext_\A^{5,193}(P^{62}_{47})$,
\end{enumerate}
we recall from \eqref{eq:e47T0-eq-g4h02} that
$\widehat{g}_4h_0^2 \in\Ext_\A^{5,193}(P)$ is represented by
$e_{47}T_0$ in the spectral sequence for $\Ext_\A^{\ast,\ast}(P)$.
Also recall from \eqref{eq:Ext:known} (19) that
\begin{align*}
  T_0
  &=\{T_0^\ast =\lambda_0^2(\lambda_{79}\lambda_{31}^2 +\lambda_{31}^2\lambda_{79}) +(\lambda_9\lambda_3^2 +\lambda_3^2\lambda_9)\lambda_{63}^2\} \in\Ext_\A^{5,146}.
\end{align*}
It follows that $(q_{47})_\ast(\widehat{g}_4h_0^2) =\{e_{47}T_0^\ast\} \in\Ext_\A^{5,193}$.
By \eqref{eq:P62-k:Ext:known} (11),
$\bar{e_{48}h_0h_3^2} =\{e_{48}\lambda_0\lambda_7^2 +e_{47}(\lambda_9\lambda_3^2 +\lambda_3^2\lambda_9)\}$. So $\bar{e_{48}h_0h_3^2}h_6^2 =(q_{47})_\ast(\widehat{g}_4h_0^2)$ because
$\delta(e_{48}\lambda_0(\lambda_{79}\lambda_{31}^2 +\lambda_{31}^2\lambda_{79})) =e_{48}\lambda_0\lambda_7^2\lambda_{63}^2 +e_{47}(\lambda_9\lambda_3^2 +\lambda_3^2\lambda_9)\lambda_{63}^2 +e_{47}T_0^\ast$ in $H_\ast(P^{62}_{47})\tensor\Lambda$.

This completes the proof of \eqref{eq:P62-k:Ext}.

Next we prove \eqref{eq:P62t:Ext}.
To describe cohomology classes in $\Ext_\A^{s,t}(X)$ for $X=\widetilde{P}^{62}$
we need to study $ H_\ast(X=\widetilde{P}^{62})$.
Recall from
\eqref{defn:P62t} that
$\widetilde{P}^{62}$ is the mapping cone $C_{\bar{h}_5}$ of a map
$S^{30}\U_{2\iota}e^{31}\to{\bar{h}_5}P^{62}$
which is an extension of the map $S^{30}\to{\widehat{\theta}_4}P^{62}$, and that
$P^{62}\to{\widetilde{q}}\widetilde{P}^{62}=P^{62}\U_{\bar{h}_5}C(S^{30}\U_{2\iota}e^{31})$
is the inclusion map.
Consider the following long exact sequence.
\begin{align*}
  \cdots-> H_{31}(S^{30}\U_{2\iota}e^{31})&\to{(\bar{h}_5)_\ast}
   H_{31}(P^{62})\to{\widetilde{q}_\ast}  H_{31}(\widetilde{P}^{62})\\
  &\quad
  \to{\Delta} H_{30}(S^{30}\U_{2\iota}e^{31})\to{(\bar{h}_5)_\ast}
   H_{30}(P^{62})\to{\widetilde{q}_\ast}  H_{30}(\widetilde{P}^{62})
  ->\cdots
\end{align*}
By \eqref{eq:h5-bar-iso},
$ H_{31}(S^{30}\U_{2\iota}e^{31})\to{(\bar{h}_5)_\ast} H_{31}(P^{62})$
is an isomorphism.
$ H_{30}(S^{30}\U_{2\iota}e^{31})\to{(\bar{h}_5)_\ast} H_{30}(P^{62})$
is zero because $S^{30}\to{\widehat{\theta}_4}P$ is detected by $\widehat{h}_4h_4$ in the ASS of $\pi_\ast^S(P)$.
So from the above long exact sequence we see the following.
\begin{note}
  \begin{enumerate}
    \item $\Z/2\iso H_k(P^{62})\to[\iso]{\widetilde{q}_\ast} H_k(\widetilde{P}^{62})\iso\Z/2$ for $k\ne31$.
    \item $\Z/2\iso H_{31}(P^{62})\to{\widetilde{q}_\ast=0} H_{31}(\widetilde{P}^{62})\iso\Z/2$.
  \end{enumerate}
  \label{eq:qt_ast}
\end{note}
\noindent
To describe the right $\A$-module $ H_\ast(\widetilde{P}^{62})$,
denote by $\widetilde{e}_{31}$
the generator of $ H_{31}(\widetilde{P}^{62})\iso\Z/2$, and
denote by $e_k$
the generator of $ H_k(\widetilde{P}^{62})\iso\Z/2$
for $1\le k\le62$ with $k\ne31$. Thus
\begin{align}
  \label{eq:P62t:Hk}
   H_k(\widetilde{P}^{62})&=
  \begin{cases}
    \Z/2(e_k) &\mbox{for $1\le k\le 62$ and $k\ne31$.}\\
    \Z/2(\widetilde{e}_{31})  &\mbox{for $k=31$}\\
    0         &\mbox{otherwise}
  \end{cases}
\end{align}
The mod $2$ Steenrod algebra $\A$ acts on
$e_k$ with $1\le k\le 62$ and $k\ne31$
the same way as given in \eqref{eq:P:e_k} but with $e_{31}=0$, and
the definition
of $S^{30}\U_{2\iota}e^{31}\to{\bar{h}_5}P^{62}$ in \eqref{eq:h5-bar}
shows that
$\A$ acts on $\widetilde{e}_{31}\in H_{31}(\widetilde{P}^{62})$
from the right by
\begin{align*}
  \widetilde{e}_{31}Sq^{16}=e_{15}.
\end{align*}
So for an element of the form
$e_k\lambda_I\in H_\ast(\widetilde{P}^{62})\tensor\Lambda$
with $1\le k\le62$ and $k\ne31$,
the differential $\delta(e_k\lambda_I)$ is as given
by formula \eqref{eq:P:diff} but with $e_{31}=0$.
And for an element of the form $\widetilde{e}_{31}\lambda_I\in H_{31}(\widetilde{P}^{62})\tensor\Lambda$,
$\delta(\widetilde{e}_{31}\lambda_I)$ is given by
\begin{align}
  \label{eq:P62t:diff}
  \delta(\widetilde{e}_{31}\lambda_I)
  &=\widetilde{e}_{31}\delta(\lambda_I) +e_{15}\lambda_{15}\lambda_I.
\end{align}
And $H^{s,t}( H_\ast(\widetilde{P}^{62})\tensor\Lambda,\delta)=\Ext_\A^{s,t+s}(\widetilde{P}^{62})$.

In \eqref{eq:P62t:Ext:new} below we define the cohomology classes that appears in
\eqref{eq:P62t:Ext}.
By \eqref{eq:Lambda:diff}, \eqref{eq:XLambda}, and \eqref{eq:P62t:diff},
each chain $\Sum_j e_{k_j}\lambda_{I_j}$ of $\{\Sum_j e_{k_j}\lambda_{I_j}\}$
in the following list is easily seen to be a cycle.
\begin{note}
  \label{eq:P62t:Ext:new}
  \begin{enumerate}
    \item $\bar{\widetilde{e}_{31}h_5}=\{\widetilde{e}_{31}\lambda_{31}\}\in\Ext_\A^{1,63}(\widetilde{P}^{62})$.
    \item $\bar{e_{47}}=\{e_{47}\}\in\Ext_\A^{0,47}(\widetilde{P}^{62})$.
    \item $\bar{e_{39}h_3}=\{e_{39}\lambda_7 +e_{23}\lambda_{23} +\widetilde{e}_{31}\lambda_{15}\}\in\Ext_\A^{1,47}(\widetilde{P}^{62})$.
    \item $\bar{\widetilde{e}_{31}h_4^2}=\{\widetilde{e}_{31}\lambda_{15}^2 +e_{23}\lambda_7\lambda_{31}\}\in\Ext_\A^{2,63}(\widetilde{P}^{62})$.
  \end{enumerate}
\end{note}

To prove \eqref{eq:P62t:Ext}
we use the induced homomorphism
$ H_\ast(P^{62})\to{\widetilde{q}_\ast} H_\ast(\widetilde{P}^{62})$.
From \eqref{eq:qt_ast}
we see $\ker(\widetilde{q}_\ast)=\Z/2(e_{31})$ and
$\coker(\widetilde{q}_\ast)=\Z/2(\widetilde{e}_{31})$,
so there is the following exact sequence
\begin{align}
  \label{eq:P62t:Ext:qt_ast-ker-coker}
  \begin{split}
  0&->\Z/2(e_{31})=\ker(\widetilde{q}_\ast)
  \to{i_1}  H_\ast(P^{62})\\
  &\qquad
  \to{\widetilde{q}_\ast}  H_\ast(\widetilde{P}^{62})
  \to{q_2} \Z/2(\widetilde{e}_{31})=\coker(\widetilde{q}_\ast)
  ->0,
  \end{split}
\end{align}
where $\Z/2(e_{31})\to{i_1} H_\ast(P^{62})$ is
given by $i_1(e_{31})=e_{31}$,
and $ H_\ast(\widetilde{P}^{62})\to{q_2}\Z/2(\widetilde{e}_{31})$ is
by $q_2(\widetilde{e}_{31})=\widetilde{e}_{31}$.
Let $M=\frac{ H_\ast(P^{62})}{\Z/2(e_{31})}$.
Let $ H_\ast(P^{62})\to{q_1}M=\frac{ H_\ast(P^{62})}{\Z/2(e_{31})}$ be the quotient map,
and $M\to{i_2} H_\ast(\widetilde{P}^{62})$ be the inclusion map,
then we have the following horizontal and oblique short exact sequences.
\begin{align}
  \label{eq:P62t:Ext:short-exact-sequence}
  \begin{split}
  \xymatrix{
  & & & & 0\ar[dl]\\
  0\ar[r]&  \Z/2(e_{31})\ar[r]^{i_1}&
   H_\ast(P^{62})\ar[r]^{q_1}\ar[d]^{\widetilde{q}_\ast}&
  M\ar[r]\ar[ld]^{i_2}&  0\\
  & &  H_\ast(\widetilde{P}^{62})\ar[dl]^{q_2}\\
  & \Z/2(\widetilde{e}_{31})\ar[dl]\\
  0
  }
  \end{split}
\end{align}
All the maps in \eqref{eq:P62t:Ext:short-exact-sequence},
except $ H_\ast(P^{62})\to{\widetilde{q}_\ast} H_\ast(\widetilde{P}^{62})$,
are algebraic homomorphisms, meaning that they are not induced by topological maps.
It follows that \eqref{eq:P62t:Ext:short-exact-sequence} induces the following
long exact sequences of $\Ext$ groups
(we denote by $\Ext_\A^{s,t}(M)$ the $\Ext$ group $\Ext_\A^{s,t}(M,\Z/2)$):
\begin{align}
  \label{eq:P62t:Ext:long-exact-sequence}
  \begin{split}
  \xymatrix{
  & & & & \Ext_\A^{s-1,t-31}\ar[dl]_-{\Delta_2}\\
  \cdots\ar[r]&  \Ext_\A^{s,t-31}\ar[r]^-{(i_1)_\ast}&
  \Ext_\A^{s,t}(P^{62})\ar[r]^{(q_1)_\ast}\ar[d]^{\widetilde{q}_\ast}&
  \Ext_\A^{s,t}(M)\ar[r]^-{\Delta_1}\ar[ld]^-{(i_2)_\ast}& \Ext_\A^{s+1,t-31}\ar[r]&  \cdots\\
  & & \Ext_\A^{s,t}(\widetilde{P}^{62})\ar[dl]^-{(q_2)_\ast}\\
  & \Ext_\A^{s,t-31}
  }
  \end{split}
\end{align}
And from this we see:
\begin{note}
  \label{eq:P62t:Ext:knowledge}
  To calculate $\Ext_\A^{s,t}(\widetilde{P}^{62})$ from
  the exact sequences in \eqref{eq:P62t:Ext:long-exact-sequence}
  we need the following information:
  \begin{inparaenum}[(i)]
    \item $\Ext_\A^{s,t-31}$,
    \item $\Ext_\A^{s,t}(P^{62})$,
    \item $\Ext_\A^{s+1,t-31}$, and
    \item $\Ext_\A^{s-1,t-31}$.
  \end{inparaenum}
  We also need to know the maps
  $\Ext_\A^{s,t}(M) \to{\Delta_1}\Ext_\A^{s+1,t-31}$ and
  $\Ext_\A^{s-1,t-31} \to{\Delta_2}\Ext_\A^{s,t}(M)$.
\end{note}

Now we prove the $\Ext$ group results in \eqref{eq:P62t:Ext}.
The equation $\bar{e_{39}h_3}h_0h_4 =\bar{\widetilde{e}_{31}h_4}h_0 +\widetilde{q}_\ast(\widetilde{D_3})$ in \eqref{eq:P62t:Ext}~(3)
will be proved later.
To see
\begin{enumerate}
  \item[\eqref{eq:P62t:Ext}~(8)\quad] $\Ext_\A^{0,189}(\widetilde{P}^{62})=0$,
  \item[\eqref{eq:P62t:Ext}~(9)\quad] $\Ext_\A^{1,190}(\widetilde{P}^{62})=0$,
\end{enumerate}
we note from \autoref{thm:Ext:known} and \autoref{thm:P:Ext:known}
that the $\Ext$ groups in \eqref{eq:P62t:Ext:knowledge}~(i) through (iv) are all zero
for $(s,t) =(0,189),(1,190)$.
Therefore $\Ext_\A^{0,189}(\widetilde{P}^{62})=0$ and $\Ext_\A^{1,190}(\widetilde{P}^{62})=0$.
This proves \eqref{eq:P62t:Ext}~(8) and \eqref{eq:P62t:Ext}~(9).
It is not difficult to see
\begin{enumerate}
  \item[\eqref{eq:P62t:Ext}~(1)\quad]
    $\Ext_\A^{0,47}(\widetilde{P}^{62})
    =\Z/2(\bar{e_{47}})$,
  \item[\eqref{eq:P62t:Ext}~(4)\quad]
    $\Ext_\A^{1,63}(\widetilde{P}^{62})
    =\Z/2(\bar{\widetilde{e}_{31}h_5})
    \dsum \Z/2(\bar{e_{47}}h_4)$,
\end{enumerate}
from the differential
$H_\ast(\widetilde{P}^{62})\tensor\Lambda \to{\delta}H_\ast(\widetilde{P}^{62})\tensor\Lambda$
as given by \eqref{eq:P62t:diff}.
To see \eqref{eq:P62t:Ext}~(2),(3),(5),(6),(10),(11),
we will use the knowledge of
\eqref{eq:P62t:Ext:knowledge}~(i), (iii), and (iv) as given by \eqref{eq:Ext:known},
and that of \eqref{eq:P62t:Ext:knowledge}~(ii) as given by \eqref{eq:P62:Ext}.
In \eqref{eq:P62t:Ext:info} below we collect this information
in \eqref{eq:P62t:Ext:info}~($j$)~(i), (ii), (iii), (iv) for each $j=2,3,5,6,10,11$.
These will imply \eqref{eq:P62t:Ext:info}~($j$)~(v) and \eqref{eq:P62t:Ext:info}~($j$)~(vi), and
note that \eqref{eq:P62t:Ext:info}~($j$)~(vi) is precisely \eqref{eq:P62t:Ext}~($j$), $j=2,3,5,6,10,11$.
\begin{note}
  \label{eq:P62t:Ext:info}
  \renewcommand{\labelenumii}{(\roman{enumii})}
  \begin{enumerate}
    \setcounter{enumi}{1}
    \item
      \begin{enumerate}
        \item $\Ext_\A^{2,17}=\Z/2(h_0h_4)$,
        \item $\Ext_\A^{2,48}(P^{62})=\Z/2(\widehat{h}_4h_0h_5)$,
        \item $\Ext_\A^{3,17}=\Z/2(h_0h_3^2)$,
        \item $\Ext_\A^{1,17}=0$,
        \item $\Ext_\A^{2,48}(M)=\Z/2( (q_1)_\ast(\widehat{h}_4h_0h_5))$,
        \item $\Ext_\A^{2,48}(\widetilde{P}^{62})=\Z/2(\bar{e_{39}h_3}h_0)
          \dsum\Z/2( q_\ast(\widehat{h}_4h_0h_5))$.
      \end{enumerate}
    \item
      \begin{enumerate}
        \item $\Ext_\A^{3,33}=\Z/2(h_0h_4^2)$,
        \item $\Ext_\A^{3,64}(P^{62})=\Z/2(\widehat{D_3})$,
        \item $\Ext_\A^{4,33}=0$,
        \item $\Ext_\A^{2,33}=\Z/2(h_0h_5)$,
        \item $\Ext_\A^{3,64}(M)=\Z/2( (q_1)_\ast(\widehat{D_3}))$,
        \item $\Ext_\A^{3,64}(\widetilde{P}^{62})=\Z/2(\bar{\widetilde{e}_{31}h_4^2}h_0)
          \dsum \Z/2(\widetilde{q}_\ast(\widehat{D_3}))$.
      \end{enumerate}
    \setcounter{enumi}{4}
    \item \begin{enumerate}
        \item $\Ext_\A^{2,33}=\Z/2(h_0h_5)$,
        \item $\Ext_\A^{2,64}(P^{62})=\Z/2(\widehat{h}_5h_0h_5)$,
        \item $\Ext_\A^{3,33}=\Z/2(h_0h_4^2)$,
        \item $\Ext_\A^{1,33}=0$,
        \item $\Ext_\A^{2,64}(M)=\Z/2(\bar{e_{47}}h_0h_4)$,
        \item $\Ext_\A^{2,64}(\widetilde{P}^{62})
          =\Z/2(\bar{\widetilde{e}_{31}h_5}h_0)
          \dsum\Z/2(\bar{e_{47}}h_0h_4)$.
      \end{enumerate}
    \item \begin{enumerate}
        \item $\Ext_\A^{4,161}=\Z/2(h_0h_4^2h_7)$,
        \item $\Ext_\A^{4,192}(P^{62})
          =\Z/2(\widehat{h}_4h_0c_4)
          \dsum \Z/2(\widehat{D_3}h_7)
          \dsum \Z/2(\widehat{g}_4h_0)
          \dsum \Z/2(\bar{e_{62}D_3(1)})$,
        \item $\Ext_\A^{5,161}=\Z/2(V_0)$,
        \item $\Ext_\A^{3,161}=h_0h_5h_7$,
        \item $\Ext_\A^{4,192}(M)
          =\Z/2( (q_1)_\ast(\widehat{h}_4h_0c_4))
          \dsum \Z/2( (q_1)_\ast(\widehat{D_3}h_7))
          \dsum \Z/2( (q_1)_\ast(\widehat{g}_4h_0))
          \dsum \Z/2( (q_1)_\ast(\bar{e_{62}D_3(1)}))$,
        \item $\Ext_\A^{4,192}(\widetilde{P}^{62})=
          Z/2(\widetilde{q}_\ast(\widehat{h}_4h_0c_4))
          \dsum \Z/2(\widetilde{q}_\ast(\widehat{D_3}h_7))
          \dsum \Z/2(\widetilde{q}_\ast(\widehat{g}_4h_0))
          \dsum \Z/2(\widetilde{q}_\ast(\bar{e_{62}D_3(1)}))
          \dsum \Z/2(\bar{\widetilde{e}_{31}h_4^2}h_0h_7)$.
      \end{enumerate}
    \setcounter{enumi}{9}
    \item
      \begin{enumerate}
        \item $\Ext_\A^{2,160}=\Z/2(h_5h_7)$,
        \item $\Ext_\A^{2,191}(P^{62})=\Z/2(\widehat{h}_5h_5h_7)$,
        \item $\Ext_\A^{3,160}=\Z/2(h_4^2h_7)$,
        \item $\Ext_\A^{1,160}=0$,
        \item $\Ext_\A^{2,191}(M)=\Z/2(\bar{e_{47}}h_4h_7)$,
        \item $\Ext_\A^{2,191}(\widetilde{P}^{62})
          =\Z/2(\bar{\widetilde{e}_{31}h_5}h_7)
          \dsum \Z/2(\bar{e_{47}}h_4h_7)$.
      \end{enumerate}
    \item
      \begin{enumerate}
        \item $\Ext_\A^{3,161}=\Z/2(h_0h_5h_7)$,
        \item $\Ext_\A^{3,192}(P^{62})=\Z/2(\widehat{h}_5h_0h_5h_7)\dsum\Z/2(\bar{e_{62}h_1}h_6^2)$,
        \item $\Ext_\A^{4,161}=\Z/2(h_0h_4^2h_7)$,
        \item $\Ext_\A^{2,161}=0$,
        \item $\Ext_\A^{3,192}(M)=\Z/2( (q_1)_\ast(\bar{e_{62}h_1}h_6^2))
          \dsum\Z/2(\bar{e_{47}h_4h_0h_7})$,
        \item $\Ext_\A^{3,192}(\widetilde{P}^{62})
          =\Z/2(\bar{\widetilde{e}_{31}h_5}h_0h_7)
          \dsum \Z/2(\bar{e_{47}}h_4h_7h_0)
          \dsum\Z/2(\widetilde{q}_\ast(\bar{e_{62}h_1}h_6^2))$.
      \end{enumerate}
  \end{enumerate}
  \renewcommand{\labelenumii}{(\alph{enumii})}
\end{note}
\noindent
We explain how we obtain \eqref{eq:P62t:Ext:info}~(2)~(v) and \eqref{eq:P62t:Ext:info}~(2)~(vi)
from \eqref{eq:P62t:Ext:info}~(2)~(i), (ii), (iii), (iv)
as an illustration of such deductions.
From the definition of $\Z/2(e_{31})\to{i_1} H_\ast(P^{62})$
we have $(i_1)_\ast(h_0h_4)=\widehat{h}_5h_0h_4=0$ in $\Ext_\A^{2,48}(P^{62})$ and
$(i_1)_\ast(h_0h_3^2)=\widehat{h}_5h_0h_3^2=\widehat{h}_4h_4^2h_0\ne0$ in $\Ext_\A^{3,48}(P^{62})$.
From this we deduce \eqref{eq:P62t:Ext:info}~(2)~(v).
Since $\Delta_2(h_0h_4)=\widehat{h}_4h_4^2h_0=\widehat{h}_5h_3^2h_0=0$ in $\Ext_\A^{3,48}(M)$, $h_0h_4\in\Ext_\A^{2,17}$ is pulled back to a class in $\Ext_\A^{2,48}(\widetilde{P}^{62})$ via $(q_2)_\ast$.
From \eqref{eq:P62t:Ext:new} we see this class
is $\bar{e_{39}h_3}h_0$.
Thus $\Ext_\A^{2,48}(\widetilde{P}^{62})=\Z/2(\bar{e_{39}h_3}h_0) \dsum\Z/2( q_\ast(\widehat{h}_4h_0h_5))$
as claimed in \eqref{eq:P62t:Ext:info}~(2)~(vi).
\eqref{eq:P62t:Ext:info}~($j$)~(vi) for $j=3,5,6,10,11$
are similarly deduced.
This proves \eqref{eq:P62t:Ext}~(2), (3), (5), (6), (10), (11).
To see
\begin{enumerate}
  \item[\eqref{eq:P62t:Ext}~(7)~(a)\quad]
    $\widetilde{q}_\ast(\widehat{h}_5h_5D_3(1)) =0$ in $\Ext_\A^{5,193}(\widetilde{P}^{62})$,
  \item[\eqref{eq:P62t:Ext}~(7)~(b)\quad]
    $\widetilde{q}_\ast(\widehat{h}_4h_0^2c_4),
    \widetilde{q}_\ast(\widehat{h}_4h_1f_3),
    \widetilde{q}_\ast(\gamma_{61}h_7),
    \widetilde{q}_\ast(\widehat{g}_4h_0^2),
    \bar{\widetilde{e}_{31}h_4^2}h_0^2h_7$ are linearly independent in $\Ext_\A^{5,193}(\widetilde{P}^{62})$,
\end{enumerate}
we also use the diagram \eqref{eq:P62t:Ext:long-exact-sequence}
with following the knowledge of \eqref{eq:P62t:Ext:knowledge}~(i) and (iv) for $(s,t)=(5,193)$
as given by \eqref{eq:Ext:known},
and that of \eqref{eq:P62t:Ext:knowledge}~(iv) for $(s,t)=(5,193)$
as given by \eqref{eq:P62:Ext}~(i)~(9).
\begin{note}
  \label{eq:P62t:Ext:5-193:info}
  \begin{enumerate}[(i)]
    \item $\Ext_\A^{5,162} =\Z/2(h_0^2h_4^2h_7) \dsum\Z/2(h_5D_3(1))$.
    \item $\Ext_\A^{5,193}(P^{62})=
      \Z/2(\widehat{h}_5h_5D_3(1))
      \dsum \Z/2(\widehat{h}_4h_0^2c_4)
      \dsum \Z/2(\widehat{h}_4h_1f_3)
      \dsum \Z/2(\gamma_{61}h_7)
      \dsum \Z/2(\widehat{g}_4h_0^2)$.
    \setcounter{enumi}{3}
    \item $\Ext_\A^{4,162} =\Z/2(h_1h_4^2h_7) \dsum\Z/2(h_0^2h_5h_7)$.
  \end{enumerate}
\end{note}
\noindent
The map $\Ext_\A^{5,162} \to{(i_1)_\ast}\Ext_\A^{5,193}(P^{62})$ is given by
$(i_1)_\ast(h_5D_3(1)) =\widehat{h}_5h_5D_3(1) \in\Ext_\A^{5,193}(P^{62})$, and
$(i_1)_\ast(h_0^2h_4^2h_7) =\widehat{h}_5h_0^2h_4^2h_7 =0$ in $\Ext_\A^{5,193}(P^{62})$.
So we see $(q_1)_\ast(\widehat{h}_5h_5D_3(1)) =0$ and therefore
$\widetilde{q}_\ast(\widehat{h}_5h_5D_3(1)) =(i_2)_\ast(q_1)_\ast(\widehat{h}_5h_5D_3(1)) =0$.
This proves \eqref{eq:P62t:Ext}~(7)~(a).
We also see that
$(q_1)_\ast(\widehat{h}_4h_0^2c_4), (q_1)_\ast(\widehat{h}_4h_1f_3), (q_1)_\ast(\gamma_{61}h_7), (q_1)_\ast(\widehat{g}_4h_0^2)$ are linearly independent in $\Ext_\A^{5,193}(M)$.
The map $\Ext_\A^{4,162} \to{\Delta_2}\Ext_\A^{5,193}(M)$ is given by
$\Delta_2(h_1h_4^2h_7) =\widehat{h}_4h_4h_1h_4^2h_7 =\widehat{h}_4h_1h_3^2h_5h_7 =0$ and
$\Delta_2(h_0^2h_5h_7) =\widehat{h}_4h_4h_0^2h_5h_7 =0$ in $\Ext_\A^{5,193}(M)$,
so the map $\Ext_\A^{5,193}(M) \to{(i_2)_\ast}\Ext_\A^{5,193}(\widetilde{P}^{62})$
is a monomorphism.
From \eqref{eq:P62t:Ext:5-193:info}~(i)
we also have $h_0^2h_4^2h_7 \in\Ext_\A^{5,162}$ with
$\Delta_2(h_0^2h_4^2h_7) =\widehat{h}_4h_4h_0^2h_4^2h_7 =\widehat{h}_4h_0^2h_3^2h_5h_7 =0$
in $\Ext_\A^{6,193}(M)$, so
$h_0^2h_4^2h_7 \in\Ext_\A^{5,162}$ has a pullback
$\bar{\widetilde{e}_{31}h_4^2}h_0^2h_7 \in\Ext_\A^{5,193}(\widetilde{P}^{62})$.
From these it is not difficult to see that
$\widetilde{q}_\ast(\widehat{h}_4h_0^2c_4)$, $\widetilde{q}_\ast(\widehat{h}_4h_1f_3)$, $\widetilde{q}_\ast(\gamma_{61}h_7)$, $\widetilde{q}_\ast(\widehat{g}_4h_0^2)$, $\bar{\widetilde{e}_{31}h_4^2}h_0^2h_7$ are linearly independent in $\Ext_\A^{5,193}(\widetilde{P}^{62})$.
This proves \eqref{eq:P62t:Ext}~(7)~(b).

The equation
\begin{enumerate}
  \item[\eqref{eq:P62t:Ext}~(3)\quad]
    $\bar{e_{39}h_3}h_0h_4 =\bar{\widetilde{e}_{31}h_4^2}h_0 +\widetilde{q}_\ast(\widehat{D_3})$ in $\Ext_\A^{3,64}(\widetilde{P}^{62})$.
\end{enumerate}
is proved by the following \eqref{eq:P62t:Ext:equation}.
\begin{note}
  \label{eq:P62t:Ext:equation}
  \begin{enumerate}
    \item From \eqref{eq:P:Ext:known} (3),
      $\widehat{D_3} =\{e_{22}\lambda_1\lambda_7\lambda_{31} +e_{16}\lambda_{15}^3 +e_{14}\lambda_9\lambda_7\lambda_{31}\} \in\Ext_\A^{3,64}(P)$, and
      $\{e_{22}\lambda_1\lambda_7\lambda_{31} +e_{16}\lambda_{15}^3 +e_{14}\lambda_9\lambda_7\lambda_{31}\} =\{e_{23}\lambda_7\lambda_{31}\lambda_0 +e_{15}(\lambda_{15}^2\lambda_{16} +\lambda_{15}\lambda_{16}\lambda_{15} +\lambda_{16}\lambda_{15}^2) +e_{16}\lambda_{15}^3\}$. So\\
      $\widetilde{q}_\ast(\widehat{D_3}) =\{e_{23}\lambda_7\lambda_{31}\lambda_0 +e_{15}(\lambda_{15}^2\lambda_{16} +\lambda_{15}\lambda_{16}\lambda_{15} +\lambda_{16}\lambda_{15}^2) +e_{16}\lambda_{15}^3\}$.
    \item It is not difficult to see that
      $\bar{\widetilde{e}_{31}h_4^2h_0} =\{\widetilde{e}_{31}\lambda_{15}^2\lambda_0 +e_{15}(\lambda_{15}^2\lambda_{16} +\lambda_{15}\lambda_{16}\lambda_{15} +\lambda_{16}\lambda_{15}^2) +e_{16}\lambda_{15}^3\}$.
    \item It is not difficult to see that
      $\bar{e_{39}h_3} =\{e_{39}\lambda_7 +e_{23}\lambda_{23} +\widetilde{e}_{31}\lambda_{15}\}$. So\\
      $\bar{e_{39}h_3}h_0h_4 =\bar{e_{39}h_3}h_4h_0 =\{(e_{39}\lambda_7 +e_{23}\lambda_{23} +\widetilde{e}_{31}\lambda_{15}) \lambda_{15}\lambda_0\} =\{e_{23}\lambda_7\lambda_{31}\lambda_0 +\widetilde{e}_{31}\lambda_{15}^2\lambda_0\}$.
    \item Then from (1), (2), and (3) we have
      $\widetilde{q}_\ast(\widehat{D_3}) +\bar{\widetilde{e}_{31}h_4^2h_0} =\{e_{23}\lambda_7\lambda_{31}\lambda_0 +\widetilde{e}_{31}\lambda_{15}^2\lambda_0\} =\bar{e_{39}h_3}h_0h_4$.
  \end{enumerate}
\end{note}

This completes the proof of \eqref{eq:P62t:Ext}.

We proceed to prove
\eqref{eq:P62-17t:Ext}~(1), (2), (3).
In \eqref{eq:P62-17t:Ext:new} below we define some cohomology classes in $\Ext_\A^{\ast,\ast}(\widetilde{P}^{62}_{17})$.
\begin{note}
  \label{eq:P62-17t:Ext:new}
  \begin{enumerate}
    \item $\bar{\widetilde{e}_{31}} =\{\widetilde{e}_{31}\} \in\Ext_\A^{0,31}(\widetilde{P}^{62}_{17})$.
    \item $\bar{e_{47}} =\{e_{47}\} \in\Ext_\A^{0,47}(\widetilde{P}^{62}_{17})$.
  \end{enumerate}
\end{note}
\noindent
$P^{62}_{17}=P^{62}/P^{16}\to{\widetilde{q}_{17}}\widetilde{P}^{62}_{17}=\widetilde{P}^{62}/P^{16}$ is
induced by $P^{62}\to{\widetilde{q}}\widetilde{P}^{62}$,
so we can modify diagram \eqref{eq:P62t:Ext:long-exact-sequence}
for $P^{62}$ to $P^{62}_{17}$ as follows,
\begin{align}
  \label{eq:P62-17t:Ext:long-exact-sequence}
  \begin{split}
  \xymatrix{
  & & & & \Ext_\A^{s-1,t-31}\ar[dl]_-{\Delta'_2}\\
  \cdots\ar[r]&  \Ext_\A^{s,t-31}\ar[r]^-{(i_1)_\ast}&
  \Ext_\A^{s,t}(P^{62}_{17})\ar[r]^{(q_1)_\ast}\ar[d]^{(\widetilde{q}_{17})_\ast}&
  \Ext_\A^{s,t}(M_{17})\ar[r]^-{\Delta'_1}\ar[ld]^-{(i_2)_\ast}& \Ext_\A^{s+1,t-31}\ar[r]&  \cdots\\
  & & \Ext_\A^{s,t}(\widetilde{P}^{62}_{17})\ar[dl]^-{(q_2)_\ast}\\
  & \Ext_\A^{s,t-31}
  }
  \end{split}
\end{align}
where $M_{17}=\frac{ H_\ast(P^{62}_{17})}{\Z/2(e_{31})}$.
To calcluate \eqref{eq:P62-17t:Ext} (1), (2), (3)
we use the same method as in the proof of \eqref{eq:P62t:Ext}.
The information in
\eqref{eq:P62-17t:Ext:2-32} ($j$) (i), (ii), (iii), (iv) for $1\le j\le 3$ below are
collected from \eqref{thm:Ext:known} and
\eqref{eq:P62-k:Ext}~(1), (2), (3).
From these we derive \eqref{eq:P62-17t:Ext:2-32} ($j$) (v) and \eqref{eq:P62-17t:Ext:2-32} ($j$) (vi), and note
that \eqref{eq:P62-17t:Ext:2-32} ($j$) (vi) are precisely \eqref{eq:P62-17t:Ext} ($j$)
for $1\le j\le 3$.
\begin{note}
  \label{eq:P62-17t:Ext:2-32}
  \renewcommand{\labelenumii}{(\roman{enumii})}
  \begin{enumerate}
    \item
      \begin{enumerate}
        \item $\Ext_\A^{1,32} =\Z/2(h_5)$,
        \item $\Ext_\A^{1,63}(P^{62}_{17}) =\Z/2(\bar{e_{31}}h_5)$,
        \item $\Ext_\A^{2,32} =\Z/2(h_4^2)$,
        \item $\Ext_\A^{0,32} =0$,
        \item $\Ext_\A^{1,63}(M_{17}) =\Z/2(\bar{e_{47}}h_4)$,
        \item $\Ext_\A^{1,63}(\widetilde{P}^{62}_{17}) =\Z/2(\bar{\widetilde{e}_{31}}h_5)
          \dsum \Z/2(\bar{e_{47}}h_4)$.
      \end{enumerate}
    \item 
      \begin{enumerate}
        \item $\Ext_\A^{2,1}=0$,
        \item $\Ext_\A^{2,32}(P^{62}_{17})=\Z/2(\bar{e_{23}}h_0h_3)$,
        \item $\Ext_\A^{3,1}=0$,
        \item $\Ext_\A^{1,1}=\Z/2(h_0)$,
        \item $\Ext_\A^{2,32}(M_{17})=\Z/2( (q_1)_\ast(\bar{e_{23}}h_0h_3))$,
        \item $\Ext_\A^{2,32}(\widetilde{P}^{62}_{17})
          =\Z/2( (\widetilde{q}_{17})_\ast(\bar{e_{23}}h_0h_3))$.
      \end{enumerate}
    \item
      \begin{enumerate}
        \item $\Ext_\A^{3,33} =\Z/2(h_0h_4^2)$,
        \item $\Ext_\A^{3,64}(P^{62}_{17})=\Z/2(\bar{e_{23}}h_0h_3h_5)$,
        \item $\Ext_\A^{4,33} =0$,
        \item $\Ext_\A^{2,33} =\Z/2(h_0h_5)$,
        \item $\Ext_\A^{3,64}(M_{17}) =\Z/2( (q_1)_\ast(\bar{e_{23}}h_0h_3h_5))$,
        \item $\Ext_\A^{3,64}(\widetilde{P}^{62}_{17})
          =\Z/2(\bar{\widetilde{e}_{31}}h_0h_4^2)
          \dsum \Z/2( (\widetilde{q}_{17})_\ast(\bar{e_{23}}h_0h_3h_5))$.
      \end{enumerate}
  \end{enumerate}
  \renewcommand{\labelenumii}{(\alph{enumii})}
\end{note}
\noindent
This proves the $\Ext$ group results in \eqref{eq:P62-17t:Ext}~(1), (2), (3).
It is not difficult to see that
the map $\Ext_\A^{1,63}(\widetilde{P}^{62}) \to{(q_{17}')_\ast}\Ext_\A^{1,63}(\widetilde{P}^{62}_{17})$
is as given in \eqref{eq:P62-17t:Ext}~(1), and
the map $\Ext_\A^{3,64}(\widetilde{P}^{62}) \to{(q_{17}')_\ast}\Ext_\A^{3,64}(\widetilde{P}^{62}_{17})$
is as given in \eqref{eq:P62-17t:Ext}~(3).
Note that
$\bar{\widetilde{e}_{31}}h_4 =\{\widetilde{e}_{31}\lambda_{15}\} \in\Ext_\A^{1,63}(\widetilde{P}^{62}_{17})$ is a decomposable element because
there is a class $\bar{\widetilde{e}_{31}} =\{\widetilde{e}_{31}\} \in\Ext_\A^{0,31}(\widetilde{P}^{62}_{17})$
(but $\widetilde{e}_{31}$ is not a cycle in $H_\ast(\widetilde{P}^{62})\tensor\Lambda$).
This proves \eqref{eq:P62-17t:Ext}~(1), (2), (3).
Finally, it is not difficult to see
\begin{enumerate}[]
  \item[\eqref{eq:P62-17t:Ext}~(4)\quad]
    $\Ext_\A^{2,48}(\widetilde{P}^{62}) \to{(q_{33}')_\ast}\Ext_\A^{2,48}(\widetilde{P}^{62}_{33} =P^{62}_{33})$
    is given by
    $(q_{33}')_\ast(\bar{e_{39}h_3}h_0) =\bar{e_{39}h_3}h_0$,
    $(q_{33}')_\ast(\widetilde{q}_\ast(\widehat{h}_4h_0h_5)) =0$,
  \item[\eqref{eq:P62-17t:Ext}~(5)\quad]
    $\Ext_\A^{0,47}(\widetilde{P}^{62}) \to{(q_{33}')_\ast}\Ext_\A^{0,47}(\widetilde{P}^{62}_{33} =P^{62}_{33})$
    is given by
    $(q_{33}')_\ast(\bar{e_{47}}) =\bar{e_{47}}$,
  \item[\eqref{eq:P62-17t:Ext}~(6)\quad]
    $\Ext_\A^{2,64}(\widetilde{P}^{62}) \to{(q_{47}')_\ast}\Ext_\A^{2,64}(\widetilde{P}^{62}_{47} =P^{62}_{47})$ is given by
    $(q_{47}')_\ast(\bar{\widetilde{e}_{31}h_5}h_5)=(q_{47}')_\ast(\bar{e_{47}}h_0h_4)=0$,
\end{enumerate}
from the cycle representations of these classes.

This completes the proof of \eqref{eq:P62-17t:Ext}.

\newpage
\section{Proofs of \autoref{thm:P62-46:pi62} and \autoref{thm:P62-47:pi62}}\label{se:pi}
In this section we prove 
\autoref{thm:P62-46:pi62} and \autoref{thm:P62-47:pi62} which are
about the stable homotopy groups $\pi_{62}^S(P^{62}_{46})$ and $\pi_{62}^S(P^{62}_{47})$.
The precise statements of these theorems will be recalled later separately
when we come to the prove them.
Our tool for proving
\autoref{thm:P62-46:pi62} and \autoref{thm:P62-47:pi62}
is the Adams spectral sequence for $\pi_\ast^S(P^m_l)$ with $1\le l< m$.
In the first half of this section
we will compute
\begin{note}
  \label{eq:se-pi:Ext}
  The $\Ext$ groups $\Ext_\A^{s,t}(P^{62}_{46})$ and
  $\Ext_\A^{s,t}(P^{62}_{47})$ for $s\ge0$, $61\le t-s\le 63$.
\end{note}
\noindent
In the second half
we will then prove \autoref{thm:P62-46:pi62} and \autoref{thm:P62-47:pi62}.

To compute \eqref{eq:se-pi:Ext}
we will use a long exact sequence in $\Ext$ groups as in \eqref{eq:X:Ext:long-exact-sequence}
below.
Let $Y\to{f}X\to{g}Z$ be a cofibration sequence, and suppose
$0->H_\ast(Y)\to{f_\ast}H_\ast(X)\to{g_\ast}H_\ast(Z)->0$ is an exact sequence
in mod $2$ homology. Then
we have a long exact sequence of $\Ext$ groups for all $s,t$.
\begin{align}
  \label{eq:X:Ext:long-exact-sequence}
  \begin{split}
  \cdots
  -> \Ext_\A^{s-1,t}(Z)
  \to{\delta} \Ext_\A^{s,t}(Y)
  &\to{f_\ast} \Ext_\A^{s,t}(X)\\
  &\quad\to{g_\ast} \Ext_\A^{s,t}(Z)
  \to{\delta} \Ext_\A^{s+1,t}(Y)
  ->\cdots
  \end{split}
\end{align}
Then to compute $\Ext_\A^{s,t}(X)$ for certain $s,t$ it suffices to
know ${\ker(\Ext_\A^{s,t}(Z)\to{\delta}\Ext_\A^{s+1,t}(Y))}$ and
${\coker(\Ext_\A^{s-1,t}(Z)\to{\delta}\Ext_\A^{s,t}(Y))}$.

The computations of \eqref{eq:se-pi:Ext} will be accomplished in the following steps.
\begin{note}
  \label{eq:se-pi:Ext:plan}
  \begin{enumerate}
    \item Compute $\Ext_\A^{s,t}(M_2 =S^0\U_{2\iota}e^1)$
      for $s\ge0$, $0\le t-s\le2$ and $13\le t-s\le 17$
      with the cofibration sequence $S^0->M_2 =S^0\U_{2\iota}e^1->S^1$.
    \item Compute $\Ext_\A^{s,t}(P_{15})$ for $s\ge0$, $15\le t-s\le17$ with
      $\Ext_\A^{s,t}(M_2 =S^0\U_{2\iota}e^1)$ for $s\ge0$, $0\le t-s\le2$ in (1).
    \item Compute $\Ext_\A^{s,t}(P^{14})$ for $s\ge0$, $13\le t-s\le16$
      with the cofibration sequence $P^{14}->P->P_{15}$.
    \item Compute $\Ext_\A^{s,t}(P^{62}_{47})$ for $s\ge0$, $61\le t-s\le63$
      with the cofibration sequence
      $P^{48}_{47}\homo\Susp^{47}M_2->P^{62}_{47}->P^{62}_{49}$.
    \item Compute $\Ext_\A^{s,t}(P^{48}_{46})$ for $s\ge0$, $60\le t-s\le63$
      with the cofibration sequence $S^{46}->P^{48}_{46}->P^{48}_{47}\homo\Susp^{47}M_2$.
    \item Compute $\Ext_\A^{s,t}(P^{62}_{46})$ for $s\ge0$, $61\le t-s\le63$
      with the cofibration sequence $P^{48}_{46}->P^{62}_{46}->P^{62}_{49}$.
  \end{enumerate}
\end{note}

We begin with \eqref{eq:se-pi:Ext:plan} (1).
Since $\Ext_\A^{s,t}(S^l)\iso\Ext_\A^{s,t-l}$ for all $s,t,l$,
from \eqref{eq:X:Ext:long-exact-sequence}
we have the following long exact sequence.
\begin{align}
  \label{eq:M2:Ext:long-exact-sequence}
  \begin{split}
  \cdots
  -> \Ext_\A^{s-1,t-1}
  \to{(2\iota)_\ast} \Ext_\A^{s,t}
  &-> \Ext_\A^{s,t}(M_2 =S^0\U_{2\iota}e^1)\\
  &\quad-> \Ext_\A^{s,t-1}
  \to{(2\iota)_\ast} \Ext_\A^{s+1,t}
  ->\cdots.
  \end{split}
\end{align}
Here the group homomorphisms $\Ext_\A^{s-1,t-1} \to{(2\iota)_\ast}\Ext_\A^{s,t}$ and 
$\Ext_\A^{s,t-1} \to{(2\iota)_\ast}\Ext_\A^{s+1,t}$
are described as follows.
\begin{note}
  \label{eq:2iota-ast}
  Since $h_0\in\Ext_\A^{1,1}$
  corresponds to $Sq^1\in\A$ and $Sq^1$ detects the stable map $S^0\to{2\iota}S^0$,
  it follows that the induced homomorphism
  $\Ext_\A^{s,t}\to{(2\iota)_\ast}\Ext_\A^{s+1,t+1}$ is given by
  \begin{align*}
    (2\iota)_\ast(x)&=h_0 x,\quad\mbox{for all $x\in\Ext_\A^{s,t}$ and for all $s,t$.}
  \end{align*}
\end{note}
\noindent
We recall from
\cite{tangora_cohomology_1970} the $\Ext$ groups
$\Ext_\A^{s,t}$ for $s\ge0$, $0\le t-s\le17$ as in \autoref{fig:Ext}.
Here, and in all figures of this section,
the horizontal axis will be the $(t-s)$-axis and
the vertical axis will be the $s$-axis.
From \eqref{eq:M2:Ext:long-exact-sequence}, \eqref{eq:2iota-ast} and \autoref{fig:Ext}
it is not difficult to deduce the structure of
$\Ext_\A^{s,t}(M_2)$ for $s\ge0$, $0\le t-s\le17$.
The results are given in \autoref{fig:M2:Ext}.
We need to explain the notation for the classes in \autoref{fig:M2:Ext}, and also
the $h_i$-extensions, $i=0,1$ among some of these classes.
We recall from \eqref{eq:P62-k:Ext:spectral-sequence} that there is a spectral sequence
$\{\bar{E}_r^{i,s,t}\}_{r\ge1}$ with
$\bar{E}_1^{i,s,t}\iso\Susp^i\Ext_\A^{s,t+s-i}$ for all $i,s,t$ and
$\dsum_{i\ge1}\bar{E}_\infty^{i,s,t}\iso\Ext_\A^{s,t+s}(P^m_l)$ for all $s,t$.
A cohomology class represented by
$e_i\alpha\in\bar{E}_1^{i,s,t}$ for $l\le i\le m$ and $\alpha\in\Ext_\A^{s,t+s-i}$
in this spectral sequence
will be denoted by $\bar{e_i\alpha}\in\Ext_\A^{s,t+s}(P^m_l)$.
In particular, since
$M_2 =S^0\U_{2\iota}e^1 \homo\Susp^{-l}(P^{l+1}_l =S^l\U_{2\iota}e^{l+1})$
(or equivalently, $P^{l+1}_l =\Susp^l M_2$)
for $l=2k+1>0$,
in \autoref{fig:M2:Ext} we display
$\Ext_\A^{s,t}(P^{l+1}_l)$ for $s\ge0$, $l\le t-s\le l+5$ and $l+13\le t-s\le l+17$
for all $l=2k+1$, from which
one can read $\Ext_\A^{s,t}(M_2)$ via $\Ext_\A^{s,t}(M_2) \iso\Ext_\A^{s,t+l}(P^{l+1}_l)$.
And the cohomology classes in \autoref{fig:M2:Ext} are denoted by
the same notations of the classes in $\Ext_\A^{s,t}(P^m_l)$ just mentioned.
The cases that we will use are
$\Ext_\A^{s,t}(P^{16}_{15})$ (for $l=15$) and $\Ext_\A^{s,t}(P^{48}_{47})$ (for $l=47$).
\begin{figure}
  \begin{center}
\begin{sseq}[grid=crossword,entrysize=7.5mm,xlabelstep=2,ylabelstep=4]{0...17}{10}
  \ssmoveto 0 0     \ssdropcirc \ssname{1}
  \ssmove 0 1       \ssdropcirc \ssstroke
                                \ssname{h0}
  \ssmove 0 1       \ssdropcirc \ssstroke
  \ssmove 0 1       \ssdropcirc \ssstroke
  \ssmove 0 1       \ssdropcirc \ssstroke
  \ssmove 0 1       \ssdropcirc \ssstroke
  \ssmove 0 1       \ssdropcirc \ssstroke
  \ssmove 0 1       \ssdropcirc \ssstroke
  \ssmove 0 1       \ssdropcirc \ssstroke
  \ssmove 0 1       \ssdropcirc \ssstroke
  \ssmove 0 1       \ssdropcirc \ssstroke
  \ssgoto{1}
  \ssmove 1 1       \ssdropcirc \ssstroke
                                \ssname{h1}\ssdroplabel[D]{h_1}
  \ssmove 1 1       \ssdropcirc \ssstroke
  \ssmove 1 1       \ssdropcirc \ssstroke
  \ssmove 0 {-1}    \ssdropcirc \ssstroke
  \ssmove 0 {-1}    \ssdropcirc \ssstroke
                                \ssname{h2}\ssdroplabel[D]{h_2}

  \ssmoveto 6 2     \ssdropcirc \ssname{h02}\ssdroplabel[D]{h_2^2}

  \ssmoveto 7 1     \ssdropcirc \ssname{h3}\ssdroplabel[D]{h_3}
  \ssmove 0 1       \ssdropcirc \ssstroke
  \ssmove 0 1       \ssdropcirc \ssstroke
  \ssmove 0 1       \ssdropcirc \ssstroke
  \ssgoto{h3}
  \ssmove 1 1       \ssdropcirc \ssstroke
  \ssmove 1 1       \ssdropcirc \ssstroke

  \ssmoveto 8 3     \ssdropcirc \ssname{c0}\ssdroplabel[D]{c_0}
  \ssmove 1 1       \ssdropcirc \ssstroke

  \ssmoveto 9 5     \ssdropcirc \ssname{Ph1}\ssdroplabel[D]{P^1h_1}
  \ssmove 1 1       \ssdropcirc \ssstroke
  \ssmove 1 1       \ssdropcirc \ssstroke
  \ssmove 0 {-1}    \ssdropcirc \ssstroke
  \ssmove 0 {-1}    \ssdropcirc \ssstroke
                                \ssname{Ph2}\ssdroplabel[D]{P^1h_2}

  \ssmoveto {14} 2  \ssdropcirc \ssname{h32}\ssdroplabel[D]{h_3^2}
  \ssmove 0 1       \ssdropcirc \ssstroke

  \ssmoveto {14} 4  \ssdropcirc \ssname{d0}\ssdroplabel[D]{d_0}
  \ssmove 0 1       \ssdropcirc \ssstroke
  \ssmove 0 1       \ssdropcirc \ssstroke
  \ssgoto{d0}
  \ssmove 1 1       \ssdropcirc \ssstroke
  \ssmove 1 1       \ssdropcirc \ssstroke
  \ssmove 1 1       \ssdropcirc \ssstroke
  \ssmove 0 {-1}    \ssdropcirc \ssstroke
  \ssmove 0 {-1}    \ssdropcirc \ssstroke
  \ssmove 0 {-1}    \ssdropcirc \ssstroke
                                \ssname{e0}\ssdroplabel[D]{e_0}
  \ssmove 1 1       \ssdropcirc \ssstroke
  \ssmove 0 {-1}    \ssdropcirc \ssstroke
                                \ssname{f0}\ssdroplabel[D]{f_0}

  \ssmoveto {15} 1  \ssdropcirc \ssname{h4}\ssdroplabel[D]{h_4}
  \ssmove 0 1       \ssdropcirc \ssstroke
  \ssmove 0 1       \ssdropcirc \ssstroke
  \ssmove 0 1       \ssdropcirc \ssstroke
  \ssmove 0 1       \ssdropcirc \ssstroke
  \ssmove 0 1       \ssdropcirc \ssstroke
  \ssmove 0 1       \ssdropcirc \ssstroke
  \ssmove 0 1       \ssdropcirc \ssstroke
  \ssgoto{h4}
  \ssmove 1 1       \ssdropcirc \ssstroke
  \ssmove 1 1       \ssdropcirc \ssstroke
  \ssmove 1 1       \ssdropcirc \ssstroke
  \ssmove 0 {-1}    \ssdropcirc \ssstroke
  \ssmove 0 {-1}    \ssdropcirc \ssstroke
                                \ssname{h2h4}\ssdroplabel[D]{h_2h_4}

  \ssmoveto {16} 7  \ssdropcirc \ssname{Pc0}\ssdroplabel[D]{P^1c_0}
  \ssmove 1 1       \ssdropcirc \ssstroke

  \ssmoveto {17} 9  \ssdropcirc \ssname{P2h1}\ssdroplabel[D]{P^2h_1}
  \ssmove 1 1       \ssdropcirc \ssstroke
  \ssmove 1 1       \ssdropcirc \ssstroke
  \ssmove 0 {-1}    \ssdropcirc \ssstroke
  \ssmove 0 {-1}    \ssdropcirc \ssstroke
                                \ssname{P2h2}\ssdroplabel[D]{P^2h_2}
\end{sseq}
  \end{center}
  \caption{$\Ext_\A^{s,t}$, $s\ge0$, $13\le t-s\le17$}
  \label{fig:Ext}
\end{figure}
\begin{figure}
  \begin{center}
\begin{sseq}[grid=crossword,entrysize=8mm,xlabelstep=2,ylabelstep=4,
  xlabels={l;l+1;l+2;l+3;l+4;l+5,l+13;l+14;l+15;l+16;l+17}
  ]{0...2,13...17}{10}
  \ssmoveto 0 0     \ssdropcirc \ssname{e1}
  \ssmove 1 1       \ssdropcirc \ssstroke
                                \ssname{e1h1}\ssdroplabel[D]{\bar{e_l}h_1}
  \ssmove 1 1       \ssdropcirc \ssstroke
  \ssmove 0 {-1}    \ssdropcirc \ssstroke[dashed]
                                \ssname{e2h1}\ssdroplabel[D]{\bar{e_{l+1}h_1}}
  \ssmove 1 1       \ssdropcirc \ssstroke
  \ssmove 1 1       \ssdropcirc \ssstroke

  \ssmoveto 3 1     \ssdropcirc \ssname{e1h2}\ssdroplabel[D]{\bar{e_l}h_2}

  \ssmoveto 6 2     \ssdropcirc \ssname{e1h22}\ssdroplabel[D]{\bar{e_l}h_2^2}

  \ssmoveto 8 4     \ssdropcirc \ssname{e2h03h3}
  \ssmove 1 1       \ssdropcirc \ssstroke[dashed]
                                \ssname{e1Ph1}
  \ssmove 1 1       \ssdropcirc \ssstroke
  \ssmove 0 {-1}    \ssdropcirc \ssstroke[dashed]
                                \ssname{e2Ph1}
  \ssmove 1 1       \ssdropcirc \ssstroke
  \ssmove 1 1       \ssdropcirc \ssstroke

  \ssmoveto 7 2     \ssdropcirc \ssname{e2h22}\ssdroplabel[D]{\bar{e_{l+1}h_2^2}}
  \ssmove 1 1       \ssdropcirc \ssstroke[dashed]
                                \ssname{e1c0}
  \ssmove 1 1       \ssdropcirc \ssstroke
  \ssmove 0 {-1}    \ssdropcirc \ssstroke[dashed]
                                \ssname{e2h1h3}
  \ssmove 1 1       \ssdropcirc \ssstroke
  \ssmove 1 1       \ssdropcirc \ssstroke[dashed]

  \ssmoveto 7 1     \ssdropcirc \ssname{e1h3}
  \ssmove 1 1       \ssdropcirc \ssstroke
  \ssmove 1 1       \ssdropcirc \ssstroke
  \ssmove 0 {-1}    \ssdropcirc \ssstroke[dashed]
                                \ssname{e2h1h3}
  \ssmove 1 1       \ssdropcirc \ssstroke

  \ssmoveto {15} 1  \ssdropcirc \ssname{e1h4}\ssdroplabel[D]{\bar{e_l}h_4}
  \ssmove 1 1       \ssdropcirc \ssstroke
  \ssmove 1 1       \ssdropcirc \ssstroke
  \ssmove 0 {-1}    \ssdropcirc \ssstroke[dashed]
                                \ssname{e2h1h4}\ssdroplabel[D]{\bar{e_{l+1}h_1}h_4}
  \ssmove 1 1       \ssdropcirc \ssstroke

  \ssmoveto {14} 2  \ssdropcirc \ssname{e1h32}\ssdroplabel[D]{\bar{e_l}h_3^2}

  \ssmoveto {15} 3  \ssdropcirc \ssname{e2h0h32}\ssdroplabel[D]{\bar{e_{l+1}h_0h_3^2}}

  \ssmoveto {14} 4  \ssdropcirc \ssname{e1d0}\ssdroplabel[D]{\bar{e_l}d_0}
  \ssmove 1 1       \ssdropcirc \ssstroke
  \ssmove 1 1       \ssdropcirc \ssstroke
  \ssmove 0 {-1}    \ssdropcirc \ssstroke[dashed]
                                \ssname{e2h1d0}\ssdroplabel[D]{\bar{e_{l+1}h_1}d_0}
  \ssmove 1 1       \ssdropcirc \ssstroke
  \ssmove 1 1       \ssdropcirc \ssstroke

  \ssmoveto {15} 6  \ssdropcirc \ssname{e2h02d0}\ssdroplabel[L]{\bar{e_{l+1}h_0^2d_0}}
  \ssmove 1 1       \ssdropcirc \ssstroke[dashed]
                                \ssname{e1Pc0}\ssdroplabel[L]{\bar{e_l}P^1c_0}
  \ssmove 1 1       \ssdropcirc \ssstroke
  \ssmove 0 {-1}    \ssdropcirc \ssstroke[dashed]
                                \ssname{e2Pc0}
  \ssmove 1 1       \ssdropcirc \ssstroke

  \ssmoveto {16} 8  \ssdropcirc \ssname{e2h07h4}\ssdroplabel[L]{\bar{e_{l+1}h_0^7h_4}}
  \ssmove 1 1       \ssdropcirc \ssstroke[dashed]
                                \ssname{e1P2h1}\ssdroplabel[L]{\bar{e_l}P^2h_1}
  \ssmove 1 1       \ssdropcirc \ssstroke
\end{sseq}
  \end{center}
  \caption{$\Ext_\A^{s,t}(P^{l+1}_l)$, $l=2k+1$, $s\ge0$, $l\le t-s\le l+5$ and $l+13\le t-s\le l+17$}
  \label{fig:M2:Ext}
\end{figure}

We also note some of the group extensions in $\Ext_\A^{\ast,\ast}(P^{l+1}_l)$.
In \autoref{fig:M2:Ext}
there are the following relations,
which are denoted by dashed lines.
\begin{note}
  \label{eq:M2:Ext:extension}
  \begin{enumerate}
    \item $\bar{e_{l+1}h_1}h_0 =\bar{e_l}h_1^2$ in $\Ext_\A^{2,l+4}(P^{l+1}_l)$.
    \item $\bar{e_{l+1}h_0^2d_0}h_1 =\bar{e_l}P^1c_0$ in $\Ext_\A^{7,l+23}(P^{l+1}_l)$.
    \item $\bar{e_{l+1}h_0^7h_4}h_1 =\bar{e_l}P^2h_1$ in $\Ext_\A^{9,l+26}(P^{l+1}_l)$.
    \item $\bar{e_{l+1}P^1c_0}h_0 =\bar{e_l}h_1P^1c_0$ in $\Ext_\A^{8,l+25}(P^{l+1}_l)$.
  \end{enumerate}
\end{note}
\noindent
Note that the relations
$\bar{e_{l+1}h_1}h_0h_4 =\bar{e_l}h_1^2h_4$ in $\Ext_\A^{3,l+20}(P^{l+1}_l)$, and
$\bar{e_{l+1}h_1}h_0d_0 =\bar{e_l}h_1^2d_0$ in $\Ext_\A^{6,l+22}(P^{l+1}_l)$
in \autoref{fig:M2:Ext}
also follow from \eqref{eq:M2:Ext:extension} (1).
We will only prove the relation \eqref{eq:M2:Ext:extension} (1) below
because the other two relations in \eqref{eq:M2:Ext:extension} will not be used in this paper.
For $l=2k+1>0$ it is not difficult to see
the chain representations of the following classes in $\Ext_\A^{\ast,\ast}(P^{l+1}_l)$.
\begin{note}
  \label{eq:M2:Ext}
  \begin{enumerate}
    \item $\bar{e_l} =\{e_l\} \in\Ext_\A^{0,l}(P^{l+1}_l)$.
    \item $\bar{e_{l+1}h_1} =\{e_{l+1}\lambda_1\} \in\Ext_\A^{1,l+3}(P^{l+1}_l)$.
  \end{enumerate}
\end{note}
\noindent
The relation \eqref{eq:M2:Ext:extension} (1)
follows from \eqref{eq:M2:Ext} (1), (2), and from
$\delta(e_{l+1}\lambda_2) =e_{l+1}\lambda_1\lambda_0 +e_l\lambda_1^2$
in $H_\ast(P^{l+1}_l)\tensor\Lambda$.

This completes the computations for \eqref{eq:se-pi:Ext:plan} (1).

To compute
$\Ext_\A^{s,t}(P_{15})$ for $s\ge0$, $15\le t-s\le17$ in
\eqref{eq:se-pi:Ext:plan} (2),
note that by dimensional reasons,
$\Ext_\A^{s,t}(P_{15}) \iso\Ext_\A^{s,t}(P^{18}_{15})$ for $s\ge0$, $15\le t-s\le17$.
Since $P^{18}_{15} =(S^{15}\U_{2\iota}e^{16})\U_{\bar{\eta}}C(S^{16}\U_{2\iota}e^{17})$,
where $S^{16}\U_{2\iota}e^{17}\to{\bar{\eta}}S^{15}$ is an extension of the map
$S^{16}\to{\eta}S^{15}$, there is a
the cofibration sequence
$\Susp^{15}M_2 \iso S^{15}\U_{2\iota}e^{16}->P^{18}_{15}->\Susp^{17}M_2 \iso S^{17}\U_{2\iota}e^{18}$
such that $0-> H_\ast(\Susp^{15}M_2)-> H_\ast(P^{18}_{15})-> H_\ast(\Susp^{17}M_2)->0$ is exact.
Thus by \eqref{eq:X:Ext:long-exact-sequence} we have the following long exact sequence.
\begin{align}
  \label{eq:P-15:Ext:long-exact-sequence}
  \begin{split}
  \cdots
  -> \Ext_\A^{s-1,t-17}(M_2)
  \to{\delta} &\Ext_\A^{s,t-15}(M_2)
  -> \Ext_\A^{s,t}(P^{18}_{15})\\
  &\quad-> \Ext_\A^{s,t-17}(M_2)
  \to{\delta} \Ext_\A^{s+1,t-15}(M_2)
  ->\cdots.
  \end{split}
\end{align}
From \eqref{eq:P-15:Ext:long-exact-sequence},
\autoref{fig:M2:Ext} (put $l=15$ and $l=17$)
and the differentials
$\delta(\bar{e_{17}}) =\bar{e_{15}}h_1$, $\delta(\bar{e_{18}h_1}) =0$
in the spectral sequence for $\Ext_\A^{\ast,\ast}(P^{18}_{15})$,
we deduce the structure of $\Ext_\A^{s,t}(P^{18}_{15})$ for $s\ge0$, $15\le t-s\le17$.
The results are given in \autoref{fig:P18-15:Ext}.
This completes the computations for \eqref{eq:se-pi:Ext:plan} (2).

We proceed to compute
$\Ext_\A^{s,t}(P^{14})$ for $s\ge0$, $13\le t-s\le16$ in \eqref{eq:se-pi:Ext:plan} (3).
From the cofibration sequence $P^{14}\to{i_{14}}P\to{q_{15}}P_{15}$, and
from \eqref{eq:X:Ext:long-exact-sequence} we have the following long exact sequence.
\begin{align}
  \label{eq:P14:Ext:long-exact-sequence}
  \begin{split}
  \cdots
  -> \Ext_\A^{s-1,t}(P)
  \to{(q_{15})_\ast} &\Ext_\A^{s-1,t}(P_{15})
  \to{\delta} \Ext_\A^{s,t}(P^{14})\\
  &\quad-> \Ext_\A^{s,t}(P)
  \to{(q_{15})_\ast} \Ext_\A^{s,t}(P_{15})
  ->\cdots.
  \end{split}
\end{align}
We recall from \cite{whitehead_recent_1970} the $\Ext$ groups
$\Ext_\A^{s,t}(P)$ for $s\ge0$, $13\le t-s\le16$
as in \autoref{fig:P:Ext}.
From \eqref{eq:P14:Ext:long-exact-sequence},
\autoref{fig:P18-15:Ext} and \autoref{fig:P:Ext}
it is not difficult to deduce $\Ext_\A^{s,t}(P^{14})$ for $s\ge0$, $13\le t-s\le16$,
noting that
the map $\Ext_\A^{s,t}(P)\to{(q_{15})_\ast}\Ext_\A^{s,t}(P_{15})$ is given by
$(q_{15})_\ast(\widehat{h}_4) =\bar{e_{15}} \in\Ext_\A^{0,15}(P_{15})$,
and $\bar{e_{16}h_1} \in\Ext_\A^{1,18}(P_{15})$ is not in the image of $(q_{15})_\ast$.
These results are given in \autoref{fig:P14:Ext}.
This completes the computations of \eqref{eq:se-pi:Ext:plan} (3).
\begin{figure}
  \begin{minipage}{.4\textwidth}
    \begin{center}
\begin{sseq}[grid=crossword,entrysize=10mm,xlabelstep=1,ylabelstep=4]{15...17}{9}
  \ssmoveto {15} 0  \ssdropcirc \ssname{e15}\ssdroplabel[D]{\bar{e_{15}}}

  \ssmoveto {17} 1  \ssdropcirc \ssname{e16h1}\ssdroplabel[D]{\bar{e_{16}h_1}}
  \ssmove 1 1       \ssdropcirc \ssstroke
\end{sseq}
    \end{center}
    \caption{$\Ext_\A^{s,t}(P^{18}_{15})$, $s\ge0$, $15\le t-s\le 17$}
    \label{fig:P18-15:Ext}
  \end{minipage}
  \qquad
  \begin{minipage}{.4\textwidth}
    \begin{center}
\begin{sseq}[grid=crossword,entrysize=10mm,xlabelstep=1,ylabelstep=4]{13...16}{9}
  \ssmoveto {14} 1  \ssdropbull \ssname{h3hh3}\ssdroplabel[D]{\widehat{h}_3h_3}
  \ssmove 0 1       \ssdropbull \ssstroke

  \ssmove 0 1       \ssdropbull \ssname{d0h}\ssdroplabel[D]{\widehat{d}_0}
  \ssmove 0 1       \ssdropbull \ssstroke
  \ssmove 0 1       \ssdropbull \ssstroke
  \ssgoto{d0h}
  \ssmove 1 1       \ssdropbull \ssstroke
  \ssmove 1 1       \ssdropbull \ssstroke
  \ssmove 1 1       \ssdropbull \ssstroke

  \ssmoveto {15} 0  \ssdropbull \ssname{h4h}\ssdroplabel[RD]{\widehat{h}_4}
  \ssmove 0 1       \ssdropbull \ssstroke
  \ssmove 0 1       \ssdropbull \ssstroke
  \ssmove 0 1       \ssdropbull \ssstroke
  \ssmove 0 1       \ssdropbull \ssstroke
  \ssmove 0 1       \ssdropbull \ssstroke
  \ssmove 0 1       \ssdropbull \ssstroke
  \ssmove 0 1       \ssdropbull \ssstroke
  \ssgoto{h4h}
  \ssmove 1 1       \ssdropbull \ssstroke
  \ssmove 1 1       \ssdropbull \ssstroke

  \ssmoveto {15} 2  \ssdropbull \ssname{h1hh32}\ssdroplabel[D]{\widehat{h}_1h_3^2}

  \ssmoveto {16} 1  \ssdropbull \ssname{h1hh4}\ssdroplabel[D]{\widehat{h}_1h_4}
  \ssmove 1 1       \ssdropbull \ssstroke

  \ssmoveto {16} 4  \ssdropbull \ssname{a16}\ssdroplabel[D]{\alpha_{16}}

  \ssmoveto {16} 6  \ssdropbull \ssname{pc0h}\ssdroplabel[D]{P^1\widehat{c_0}}
  \ssmove 1 1       \ssdropbull \ssstroke
\end{sseq}
    \end{center}
    \caption{$\Ext_\A^{s,t}(P)$, $s\ge0$, $13\le t-s\le16$}
    \label{fig:P:Ext}
  \end{minipage}
\end{figure}
\begin{figure}
  \begin{minipage}{.4\textwidth}
    \begin{center}
\begin{sseq}[grid=crossword,entrysize=10mm,xlabelstep=1,ylabelstep=4]{13...16}{9}
  \ssmoveto {14} 1  \ssdropbull \ssname{h3hh3}\ssdroplabel[LD]{\widehat{h}_3h_3}
  \ssmove 0 1       \ssdropbull \ssstroke

  \ssmove 0 1       \ssdropbull \ssname{d0h}\ssdroplabel[D]{\widehat{d}_0}
  \ssmove 0 1       \ssdropbull \ssstroke
  \ssmove 0 1       \ssdropbull \ssstroke
  \ssgoto{d0h}
  \ssmove 1 1       \ssdropbull \ssstroke
  \ssmove 1 1       \ssdropbull \ssstroke
  \ssmove 1 1       \ssdropbull \ssstroke

  \ssmoveto {15} 1  \ssdropbull \ssname{h4hh0}\ssdroplabel[D]{\widehat{h}_4h_0}
  \ssmove 0 1       \ssdropbull \ssstroke
  \ssmove 0 1       \ssdropbull \ssstroke
  \ssmove 0 1       \ssdropbull \ssstroke
  \ssmove 0 1       \ssdropbull \ssstroke
  \ssmove 0 1       \ssdropbull \ssstroke
  \ssmove 0 1       \ssdropbull \ssstroke

  \ssmoveto {16} 1  \ssdropbull \ssname{h4hh1}\ssdroplabel[D]{\widehat{h}_4h_1}
  \ssmove 1 1       \ssdropbull \ssstroke

  \ssmoveto {15} 2  \ssdropbull \ssname{h1hh32}\ssdroplabel[D]{\widehat{h}_1h_3^2}

  \ssmoveto {16} 1  \ssdropbull \ssname{h1hh4}\ssdroplabel[D]{\widehat{h}_1h_4}
  \ssmove 1 1       \ssdropbull \ssstroke

  \ssmoveto {16} 4  \ssdropbull \ssname{a16}\ssdroplabel[D]{\alpha_{16}}

  \ssmoveto {16} 6  \ssdropbull \ssname{pc0t}\ssdroplabel[D]{P^1\widehat{c_0}}
  \ssmove 1 1       \ssdropbull \ssstroke

  \ssmoveto {16} 2  \ssdropbull \ssname{e14h12}\ssdroplabel[U]{\bar{e_{14}h_1^2}}
\end{sseq}
    \end{center}
    \caption{$\Ext_\A^{s,t}(P^{14})$, $s\ge0$, $13\le t-s\le16$}
    \label{fig:P14:Ext}
  \end{minipage}
  \qquad
  \begin{minipage}{.4\textwidth}
    \begin{center}
\begin{sseq}[grid=crossword,entrysize=10mm,xlabelstep=1,ylabelstep=4]{61...63}{9}

  \ssmoveto {61} 4  
  \ssmove 1 1 
  \ssmove 1 1 
  \ssmove 0 {-1}  \ssdropcirc         \ssname{e48h1d0}

  \ssmoveto {62} 1  \ssdropcirc       \ssname{e47h4} \ssdroplabel[D]{\bar{e_{47}}h_4}

  \ssmoveto {62} 3  \ssdropcirc       \ssname{e48h0h32} \ssdroplabel[L]{\bar{e_{48}h_0h_3^2}}

  \ssmoveto {62} 6  \ssdropcirc       \ssname{e48h02d0} \ssdroplabel[L]{\bar{e_{48}h_0^2d_0}}

  \ssmoveto {62} 1  
  \ssmove 0 1 \ssdropbull             \ssname{e55h0h3} \ssdroplabel[L]{\bar{e_{55}h_0h_3}}
  \ssline[dashed]{0}{1}

  \ssmoveto {62} 3 
  \ssmove 0 1 \ssdropbull             \ssname{e53h1c0} \ssdroplabel[L]{\bar{e_{53}h_1c_0}}
  \ssmove 0 1 \ssdropbull \ssstroke   \ssname{e51P1h2} \ssdroplabel[L]{\bar{e_{51}P^1h_2}}
  \ssline[dashed,curve=1mm]{0}{-2}\ssgoto{e51P1h2}
  \ssline[dashed]{0}{1}

  \ssmoveto {63} 1  \ssdropbull       \ssname{e63h0}
  \ssmove 0 1 \ssdropbull \ssstroke
  \ssmove 0 1 \ssdropbull \ssstroke
  \ssmove 0 1 \ssdropbull \ssstroke
  \ssmove 0 1 \ssdropbull \ssstroke
  \ssmove 0 1 \ssdropbull \ssstroke
  \ssmove 0 1 \ssdropbull \ssstroke

  \ssmoveto {63} 2  \ssdropbull \ssname{e49h32}\ssdroplabel[R]{\bar{e_{49}h_3^2}}

  \ssmoveto {63} 8  \ssdropcirc \ssname{e48h07h4}
\end{sseq}
    \end{center}
    \caption{$\Ext_\A^{s,t}(P^{62}_{47})$, $s\ge0$, $61\le t-s\le63$}
    \label{fig:P62-47:Ext}
  \end{minipage}
\end{figure}

To compute $\Ext_\A^{s,t}(P^{62}_{47})$ for $s\ge0$, $61\le t-s\le63$ for \eqref{eq:se-pi:Ext:plan} (4),
we use the the cofibration sequence $P^{48}_{47}->P^{62}_{47}->P^{62}_{49}$.
Since $H^\ast(P^{62}_{49})\iso\Susp^{48} H^\ast(P^{14})$
as left $\A$-modules,
$\Ext_\A^{s,t}(P^{62}_{49}) \iso\Ext_\A^{s,t-48}(P^{14})$ for all $s,t$.
So by \eqref{eq:X:Ext:long-exact-sequence} we have
the following long exact sequence.
\begin{align}
  \label{eq:P62-47:Ext:long-exact-sequence}
  \begin{split}
  \cdots-> \Ext_\A^{s-1,t-48}(P^{14})
  \to{\delta} &\Ext_\A^{s,t}(P^{48}_{47})
  -> \Ext_\A^{s,t}(P^{62}_{47})\\
  &\quad-> \Ext_\A^{s,t-48}(P^{14})
  \to{\delta} \Ext_\A^{s+1,t}(P^{48}_{47})
  ->\cdots
  \end{split}
\end{align}
From \eqref{eq:P62-47:Ext:long-exact-sequence},
\autoref{fig:M2:Ext} (put $l=47$) and \autoref{fig:P14:Ext},
it is not difficult to deduce $\Ext_\A^{s,t}(P^{62}_{47})$ for $s\ge0$, $61\le t-s\le 63$.
The results are given in \autoref{fig:P62-47:Ext}.
The group extensions in \autoref{fig:P62-47:Ext} are given in the following \eqref{eq:P62-47:Ext:extension}.
\begin{note}
  \label{eq:P62-47:Ext:extension}
  \begin{enumerate}
    \item $\bar{e_{55}h_0h_3}h_0=\bar{e_{48}h_0h_3^2}$ by the following (a), (b), and (c).
      \begin{enumerate}
        \item $\bar{e_{55}h_0h_3} =\{\bar{e_{55}h_0h_3}^\ast\}$ where
          $\bar{e_{55}h_0h_3}^\ast \equiv e_{55}\lambda_0\lambda_7 +e_{48}\lambda_7^2 \mod{F(47)}$,
        \item $\bar{e_{48}h_0h_3^2} =\{\bar{e_{48}h_0h_3^2}^\ast\}$ where
          $\bar{e_{48}h_0h_3^2}^\ast \equiv e_{48}\lambda_7^2\lambda_0 \mod{F(47)}$,
        \item $\bar{e_{55}h_0h_3}^\ast\lambda_0 +\bar{e_{48}h_0h_3^2}^\ast
          \equiv\delta(e_{55}\lambda_0\lambda_8 +e_{56}\lambda_0\lambda_7 +e_{54}\lambda_2\lambda_7 +e_{57}\lambda_3^2 +e_{52}\lambda_4\lambda_7 +e_{51}\lambda_5\lambda_7) \mod{F(47)}$.
      \end{enumerate}
    \item The extension from $\bar{e_{48}h_0h_3^2}$ to $\bar{e_{51}P^1h_2}$
      in \autoref{fig:P62-47:Ext} is nontrivial and
      will be explained in the proof of \autoref{thm:P62-46:pi62} later.
    \item $\bar{e_{53}h_1c_0}h_0 =\bar{e_{51}P^1h_2}$ by the following.
      \begin{enumerate}[(a)]
        \item $P^1h_2 =\toda{h_1,h_1c_0,h_0} \in\Ext_\A^{5,16}$.
      \end{enumerate}
    \item $\bar{e_{51}P^1h_2}h_0 =\bar{e_{48}h_0^2d_0}$ because the following (a), (b), (c), and (d).
      \begin{enumerate}[(a)]
        \item Recall the following
          \begin{enumerate}
            \item[\eqref{eq:Ext:known}~(3)\quad]
              $d_0 =\{d_0^\ast =\lambda_6\lambda_2\lambda_3^2 +\lambda_4^2\lambda_3^2 +\lambda_2\lambda_4\lambda_5\lambda_3)\} \in\Ext_\A^{4,18}$
          \end{enumerate}
        \item $\bar{e_{51}P^1h_2} =\{\bar{e_{51}P^1h_2}^\ast\}$ where
          $\bar{e_{51}P^1h_2}^\ast\equiv e_{51}\lambda_0^4\lambda_{11} +e_{48}\lambda_0d_0^\ast \mod{F(47)}$,
        \item $\bar{e_{48}h_0^2d_0} =\{\bar{e_{48}h_0^2d_0}^\ast\}$ where
          $\bar{e_{48}h_0^2d_0}^\ast\equiv e_{48}\lambda_0d_0^\ast\lambda_0 \mod{F(47)}$,
        \item $\bar{e_{51}P^1h_2}^\ast\lambda_0 +\bar{e_{48}h_0^2d_0}^\ast
          \equiv\delta(
          e_{51}\lambda_0^4\lambda_{12}
          +e_{52}\lambda_0^4\lambda_{11}
          +e_{51}\lambda_0^2\lambda_6\lambda_3^2
          +e_{50}(\lambda_2\lambda_0^3\lambda_{11} +\lambda_4\lambda_1\lambda_2\lambda_3^2)
          +e_{49}(\lambda_1\lambda_2\lambda_0^2\lambda_{11} +\lambda_2\lambda_4\lambda_2\lambda_3^2 +\lambda_5\lambda_1\lambda_2\lambda_3^2)
          )\mod{F(47)}$.
      \end{enumerate}
  \end{enumerate}
\end{note}

To compute
$\Ext_\A^{s,t}(P^{48}_{46})$ for $s\ge0$, $60\le t-s\le63$
for \eqref{eq:se-pi:Ext:plan} (5),
consider the cofibration sequence
$S^{46} ->P^{48}_{46}=S^{46}\U_{\widetilde{\eta}}C(S^{46})\U_{2\iota}e^{47}) ->S^{47}\U_{2\iota}e^{48}$.
By \eqref{eq:X:Ext:long-exact-sequence} this induces the following long exact sequence of $\Ext$ groups.
\begin{align}
  \label{eq:P48-46:Ext:long-exact-sequence}
  \begin{split}
  \cdots
  -> &\Ext_\A^{s-1,t}(S^{47}\U_{2\iota}e^{48})
  \to{\delta} \Ext_\A^{s,t}(S^{46})
  -> \Ext_A^{s,t}(P^{48}_{46})\\
  &\qquad\qquad-> \Ext_\A^{s,t}(S^{47}\U_{2\iota}e^{48})
  \to{\delta} \Ext_\A^{s+1,t}(S^{46})
  -> \cdots.
  \end{split}
\end{align}
From \eqref{eq:P48-46:Ext:long-exact-sequence},
\autoref{fig:Ext} and \autoref{fig:M2:Ext} (put $l=47$)
it is not difficult
to deduce the structure of $\Ext_\A^{s,t}(P^{48}_{46})$ for $s\ge0$, $61\le t-s\le63$.
These results are given in \autoref{fig:P48-46:Ext}.
The only group extension in \autoref{fig:P48-46:Ext} that will be used in this paper
is explained in the following.
\begin{note}
  \label{eq:P48-46:Ext:extension}
  $\bar{e_{47}}h_4h_0 =\bar{e_{47}}h_0h_4 =\bar{e_{46}}h_1h_4 =\bar{e_{46}}h_4h_1$ because
    $e_{47}\lambda_0 +e_{46}\lambda_1 \equiv\delta(e_{48}) \mod{F(45)}$.
\end{note}
\noindent
This completes the computation for \eqref{eq:se-pi:Ext:plan}~(5).

From \autoref{fig:P14:Ext} (put $l=47$) and \autoref{fig:P48-46:Ext}
it is not difficult to deduce $\Ext_\A^{s,t-s}(P^{62}_{46})$ for $s\ge0$, $61\le t-s\le63$.
These results are given in \autoref{fig:P62-46:Ext}.
The group extensions given in \autoref{fig:P62-46:Ext}
are obtained by the same computations in
\eqref{eq:P62-47:Ext:extension}~(1) through (4), except the followings.
\begin{note}
  \label{eq:P62-46:Ext:extension}
  \begin{enumerate}
    \item $\bar{e_{47}}h_4h_0 =\bar{e_{46}h_4}h_1$ follows from \eqref{eq:P48-46:Ext:extension}.
    \item $\bar{e_{48}h_0^2d_0}h_0 =\bar{e_{46}}P^1c_0$ because the following.
      \begin{enumerate}[(a)]
        \item $P^1c_0 =\toda{h_1,h_0^2d_0,h_0}$.
      \end{enumerate}
  \end{enumerate}
\end{note}
\noindent
This completes the computations for \eqref{eq:se-pi:Ext:plan} (6).
\begin{figure}
  \begin{minipage}{.4\textwidth}
    \begin{center}
\begin{sseq}[grid=crossword,entrysize=10mm,xlabelstep=1,ylabelstep=4]{61...63}{9}
  \ssmoveto {61} 1  \ssdropast \ssname{e46h4}\ssdroplabel[D]{\bar{e_{46}}h_4}
  \ssmove 0 1       \ssdropast \ssstroke
  \ssmove 0 1       \ssdropast \ssstroke
  \ssmove 0 1       \ssdropast \ssstroke
  \ssmove 0 1       \ssdropast \ssstroke
  \ssmove 0 1       \ssdropast \ssstroke
  \ssmove 0 1       \ssdropast \ssstroke
  \ssmove 0 1       \ssdropast \ssstroke
  \ssgoto{e46h4}
  \ssmove 1 1       \ssdropast \ssstroke

  \ssmoveto {61} 5  \ssdropast  \ssname{e46h1d0}
  \ssmove 0 {-1}    \ssdropcirc \ssstroke[dashed]\ssname{e47d0}
  \ssmove 1 1       \ssdropcirc \ssstroke
  \ssmove 1 1       \ssdropcirc \ssstroke

  \ssmoveto {63} 6  \ssdropast  
  \ssmove 0 {-1}  \ssdropast  \ssstroke
  \ssmove 0 {-1}  \ssdropast  \ssstroke

  \ssmoveto {62} 7  \ssdropast        \ssname{e46pc0}

  \ssmoveto {61} 2  \ssdropcirc       \ssname{e47h32}


  \ssmoveto {62} 1  \ssdropcirc       \ssname{e47h4}  \ssdroplabel[D]{\bar{e_{47}}h_4}
  \ssline[dashed]{0}{1}               \ssgoto{e47h4}
  \ssmove 1 1 \ssdropcirc \ssstroke

  \ssmoveto {62} 3  \ssdropcirc       \ssname{e48h0h32}

  \ssmoveto {62} 6  \ssdropcirc       \ssname{e48h02d0}
  \ssmove 1 1 \ssdropcirc \ssstroke
\end{sseq}
    \end{center}
    \caption{$\Ext_\A^{s,t}(P^{48}_{46})$, $s\ge0$, $61\le t-s\le 63$}
    \label{fig:P48-46:Ext}
  \end{minipage}
  \qquad
  \begin{minipage}{.4\textwidth}
    \begin{center}
\begin{sseq}[grid=crossword,entrysize=10mm,xlabelstep=1,ylabelstep=4]{61...63}{9}
  \ssmoveto {61} 1  \ssdropast  \ssname{e46h4}\ssdroplabel[D]{\bar{e_{46}h_4}}
  \ssmove 0 1       \ssdropast  \ssstroke
  \ssmove 0 1       \ssdropast  \ssstroke
  \ssmove 0 1       \ssdropast  \ssstroke
  \ssmove 0 1       \ssdropast  \ssstroke
  \ssmove 0 1       \ssdropast  \ssstroke
  \ssmove 0 1       \ssdropast  \ssstroke
  \ssmove 0 1       \ssdropast  \ssstroke
  \ssgoto{e46h4}
  \ssmove 1 1       \ssdropast  \ssstroke

  \ssmoveto {63} 6  \ssdropast
  \ssmove 0 {-1}    \ssdropast  \ssstroke
  \ssmove 0 {-1}    \ssdropast  \ssstroke

  \ssmoveto {62} 7  \ssdropast  \ssname{e46pc0}\ssdroplabel[D]{\bar{e_{46}}P^1c_0}

  \ssmoveto {62} 1  \ssdropcirc \ssname{e47h4}\ssdroplabel[D]{\bar{e_{47}}h_4}
  \ssline[dashed]{0}{1}

  \ssmoveto {62} 3  \ssdropcirc \ssname{e48h0h32}\ssdroplabel[D]{\bar{e_{48}h_0h_3^2}}

  \ssmoveto {62} 6  \ssdropcirc \ssname{e48h02d0}\ssdroplabel[D]{\bar{e_{48}h_0^2d_0}}
  \ssline[dashed]{0}{1}

  \ssmoveto {62} 1
  \ssmove 0 1       \ssdropbull \ssname{e55h0h3}\ssdroplabel[D]{\bar{e_{55}h_0h_3}}
  \ssline[dashed]{0}{1}

  \ssmoveto {62} 3
  \ssmove 0 1
  \ssmove 0 1       \ssdropbull \ssname{e51P1h2}\ssdroplabel[D]{\bar{e_{51}P^1h_2}}
  \ssline[dashed]{0}{-2}
  \ssgoto{e51P1h2}
  \ssline[dashed]{0}{1}

  \ssmoveto {63} 1  \ssdropbull \ssname{e63h0}
  \ssmove 0 1       \ssdropbull \ssstroke
  \ssmove 0 1       \ssdropbull \ssstroke
  \ssmove 0 1       \ssdropbull \ssstroke
  \ssmove 0 1       \ssdropbull \ssstroke
  \ssmove 0 1       \ssdropbull \ssstroke
  \ssmove 0 1       \ssdropbull \ssstroke

  \ssmoveto {63} 2  \ssdropbull \ssname{e49h32}\ssdroplabel[R]{\bar{e_{49}h_3^2}}
\end{sseq}
    \end{center}
    \caption{$\Ext_\A^{s,t}(P^{62}_{46})$, $s\ge0$, $61\le t-s\le63$}
    \label{fig:P62-46:Ext}
  \end{minipage}
\end{figure}
\begin{figure}
  \begin{minipage}{.4\textwidth}
    \begin{center}
\begin{sseq}[grid=crossword,entrysize=10mm,xlabelstep=1,ylabelstep=4]{61...63}{9}
  \ssmoveto {61} 1  \ssdropast  \ssname{e46h4}
  \ssmove 0 1       \ssdropast  \ssstroke
  \ssmove 0 1       \ssdropast  \ssstroke
  \ssmove 0 1       \ssdropast  \ssstroke
  \ssmove 0 1       \ssdropast  \ssstroke
  \ssmove 0 1       \ssdropast  \ssstroke
  \ssmove 0 1       \ssdropast  \ssstroke
  \ssmove 0 1       \ssdropast  \ssstroke
  \ssgoto{e46h4}
  \ssmove 1 1       \ssdropast  \ssstroke

  \ssmoveto {63} 6  \ssdropast
  \ssmove 0 {-1}    \ssdropast  \ssstroke
  \ssmove 0 {-1}    \ssdropast  \ssstroke

  \ssmoveto {62} 7  \ssdropast  \ssname{e46pc0}\ssdroplabel[D]{\bar{e_{46}}P^1c_0}

  \ssmoveto {62} 1  \ssdropcirc \ssname{e47h4}\ssdroplabel[D]{\bar{e_{47}}h_4}
  \ssline[dashed]{0}{1}

  \ssmoveto {62} 3  \ssdropcirc \ssname{e48h0h32}\ssdroplabel[D]{\bar{e_{48}h_0h_3^2}}

  \ssmoveto {62} 6  \ssdropcirc \ssname{e48h02d0}\ssdroplabel[D]{\bar{e_{48}h_0^2d_0}}
  \ssline[dashed]{0}{1}

  \ssmoveto {62} 1
  \ssmove 0 1       \ssdropbull \ssname{e55h0h3}\ssdroplabel[D]{\bar{e_{55}h_0h_3}}
  \ssline[dashed]{0}{1}

  \ssmoveto {62} 3
  \ssmove 0 1
  \ssmove 0 1       \ssdropbull \ssname{e51P1h2}\ssdroplabel[D]{\bar{e_{51}P^1h_2}}
  \ssline[dashed]{0}{-2}
  \ssgoto{e51P1h2}
  \ssline[dashed]{0}{1}

  \ssmoveto {63} 0  \ssdropbull \ssname{e63}
  \ssmove 0 1       \ssdropbull \ssstroke
  \ssmove 0 1       \ssdropbull \ssstroke
  \ssmove 0 1       \ssdropbull \ssstroke
  \ssmove 0 1       \ssdropbull \ssstroke
  \ssmove 0 1       \ssdropbull \ssstroke
  \ssmove 0 1       \ssdropbull \ssstroke
  \ssmove 0 1       \ssdropbull \ssstroke

  \ssmoveto {63} 2  \ssdropbull \ssname{e49h32}\ssdroplabel[R]{\bar{e_{49}h_3^2}}
\end{sseq}
    \end{center}
    \caption{$\Ext_\A^{s,t}(P_{46})$, $s\ge0$, $61\le t-s\le63$}
    \label{fig:P-46:Ext}
  \end{minipage}
  \qquad
  \begin{minipage}{.4\textwidth}
    \begin{center}
\begin{sseq}[grid=crossword,entrysize=10mm,xlabelstep=1,ylabelstep=4]{62...64}{5}
  \ssmoveto {61} 1  \ssdropast \ssname{e46h4}\ssdroplabel[D]{\bar{e_{46}}h_4}
  \ssmove 0 1       \ssdropast \ssstroke
  \ssmove 0 1       \ssdropast \ssstroke
  \ssmove 0 1       \ssdropast \ssstroke
  \ssmove 0 1       \ssdropast \ssstroke
  \ssmove 0 1       \ssdropast \ssstroke
  \ssmove 0 1       \ssdropast \ssstroke
  \ssmove 0 1       \ssdropast \ssstroke
  \ssgoto{e46h4}
  \ssmove 1 1       \ssdropast \ssstroke

  \ssmoveto {61} 5  \ssdropast  \ssname{e46h1d0}
  \ssmove 0 {-1}    \ssdropcirc \ssstroke[dashed]\ssname{e47d0}

  \ssmoveto {63} 6  \ssdropast  
  \ssmove 0 {-1}  \ssdropast  \ssstroke
  \ssmove 0 {-1}  \ssdropast  \ssstroke

  \ssmoveto {62} 7  \ssdropast        \ssname{e46pc0}
  \ssmoveto {61} 2  \ssdropcirc       \ssname{e47h32}

  \ssmoveto {62} 1  \ssdropcirc       \ssname{e47h4}  \ssdroplabel[D]{\bar{e_{47}}h_4}
  \ssline[dashed]{0}{1}               \ssgoto{e47h4}

  \ssmoveto {62} 3  \ssdropcirc       \ssname{e48h0h32}

  \ssmoveto {62} 6  \ssdropcirc       \ssname{e48h02d0}
  \ssmove 1 1 \ssdropcirc \ssstroke

  \ssmoveto {64} 2  \ssdropcirc \ssname{e49h0h4}\ssdroplabel[D]{\bar{e_{49}h_0}h_4}
  \ssmove 0 1       \ssdropcirc \ssstroke
  \ssmove 0 1       \ssdropcirc \ssstroke
  \ssmove 0 1       \ssdropcirc \ssstroke

  \ssmoveto {63} 2  \ssdropcirc \ssname{e49h32}\ssdroplabel[D]{\bar{e_{49}h_3^2}}
  \ssmove 0 1       \ssdropcirc \ssstroke
\end{sseq}
    \end{center}
    \caption{$\Ext_\A^{s,t}(P^{49}_{46})$, $0\le s\le4$, $62\le t-s\le64$}
    \label{fig:P49-46:Ext}
  \end{minipage}
\end{figure}

In order to prove \autoref{thm:P62-46:pi62} and \ref{thm:P62-47:pi62},
we first prove two differentials and also show an infinite cycle
in the ASS for $\pi_\ast^S(P^{62}_{46})$ and also in the ASS for $\pi_\ast^S(P^{62}_{47})$.
The ranges of these two spectral sequences of our discussion are shown
in \autoref{fig:P62-47:Ext} and \autoref{fig:P62-46:Ext}.
The two differentials are stated in \autoref{thm:P62-46:ASS:d5}.
The infinite cycle is given in \autoref{thm:P62-46:ASS:e49h32}.
We also need the following $\Ext$ group information:
$\Ext_\A^{s,t}(P_{46})$ for $s\ge0$ and $61\le t-s\le63$, and
$\Ext_\A^{s,t}(P^{49}_{46})$ for $0\le s\le4$ and $62\le t-s\le64$.
These results are given in \autoref{fig:P-46:Ext} and \autoref{fig:P49-46:Ext}.
The computations of these groups
follow the same procedure as the computation of, for example,
$\Ext_\A^{s,t}(M_2)$ in \autoref{fig:M2:Ext}, and
will not be repeated here.

To prove \autoref{thm:P62-46:ASS:d5} we need \autoref{thm:2q46t5:detected}.
We recall 
the stable map $S^{62} \to{\widehat{\theta}_5}P^{62}$ in \eqref{eq:theta_n-hat}.
For $1\le k\le62$ let
$P^{62} \to{q_k}P^{62}_k$ be the collapsing map and let 
$P^k \to{i_k}P^{62}$ be the inclusion map.
Let $P^{62}_{46} \to{q'}P^{62}_{47}$ and $P^{62}_{46} \to{q''}P^{62}_{49}$
also be the collapsing maps.
Recall the following \eqref{thm:qk2t5:essential} and \eqref{thm:q47t5:AF-ge-3} in \autoref{se:boundary}.
\theoremstyle{nonumberplain}
\newtheorem{thm:qk2t5:essential}{\autoref{thm:qk2t5:essential}} 
\begin{thm:qk2t5:essential*}
  The composite
  $S^{62} \to{2\iota}S^{62} \to{\widehat{\theta}_5}P^{62} \to{q_k}P^{62}_k$
  is essential
  if and only if $1\le k\le51$.
\end{thm:qk2t5:essential*}
\theoremstyle{nonumberplain}
\newtheorem{thm:q47t5:AF-ge-3}{\autoref{thm:q47t5:AF-ge-3}} 
\begin{thm:q47t5:AF-ge-3*}
  $AF(q_{47}\widehat{\theta}_5)\ge3$ for the composite
  $S^{62} \to{\theta_5}P^{62} \to{q_{47}}P^{62}_{47}$.
\end{thm:q47t5:AF-ge-3*}
\begin{lemm}
  \label{thm:2q46t5:detected}
  The stable homotopy element
  $2q_{46}\widehat{\theta}_5 \in\pi_{62}^S(P^{62}_{46})$ is detected by
  $\bar{e_{51}P^1h_2} \in\Ext_\A^{5,67}(P^{62}_{46})$
  in the ASS for $\pi_\ast^S(P^{62}_{46})$.
\end{lemm}
\begin{proof}
  Refer to \autoref{fig:P62-46:Ext} for the ASS for $\pi_\ast^S(P^{62}_{46})$
  in the range of our discussion.
  $\widehat{\theta}_5 \in\pi_{62}^S(P^{62})$ is detected by
  $\widehat{h}_5h_5 \in\Ext_\A^{1,63}(P^{62})$, so
  $AF(\widehat{\theta}_5) =1$
  in the ASS for $\pi_\ast^S(P^{62})$.
  It is clear that $(q_{46})_\ast(\widehat{h}_5) =0$.
  So $AF(q_{46}\widehat{\theta}_5)\ge2$,
  and this implies that $AF(2q_{46}\widehat{\theta}_5)\ge3$ for the ASS for $\pi_\ast^S(P^{62}_{46})$.
  From \autoref{fig:P62-46:Ext} we see that
  $2q_{46}\widehat{\theta}_5$ must be detected by
  one of the following four $\Ext$ group classes:
  $\bar{e_{48}h_0h_3^2} \in\Ext_\A^{3,65}(P^{62}_{46})$,
  $\bar{e_{51}P^1h_2} \in\Ext_\A^{5,67}(P^{62}_{46})$,
  $\bar{e_{48}h_0^2d_0} \in\Ext_\A^{6,68}(P^{62}_{46})$, or
  $\bar{e_{46}P^1c_0} \in\Ext_\A^{7,69}(P^{62}_{46})$.
  If $AF(2q_{46}\widehat{\theta}_5)\ge6$, that is,
  if $2q_{46}\widehat{\theta}_5$ is detected by
  $\bar{e_{48}h_0^2d_0}$ or $\bar{e_{46}P^1c_0}$, then
  \begin{align}
    \label{eq:2q49t5}
    AF(2q_{49}\widehat{\theta}_5 &=q''(2q_{46}\widehat{\theta}_5))\ge6
    \quad\mbox{in the ASS for $\pi_\ast^S(P^{62}_{49})$.}
    \tag{*}
  \end{align}
  Since $\Ext_\A^{s,62+s}(P^{62}_{49}) \iso\Ext_\A^{s,14+s}(P^{14})$,
  from \autoref{fig:P14:Ext} we see that $\Ext_\A^{s,62+s}(P^{62}_{49}) =0$ for $s\ge6$.
  This together with \eqref{eq:2q49t5} imply
  $2q_{49}\widehat{\theta}_5 =0$ in $\pi_{62}^S(P^{62}_{49})$.
  But this contradicts \autoref{thm:qk2t5:essential} (for $k=49$).
  Therefore $AF(2q_{46}\widehat{\theta}_5)\le5$ in the ASS for $\pi_\ast^S(P^{62}_{46})$.
  It follows that
  $S^{62} \to{2q_{46}\widehat{\theta}_5}P^{62}_{46}$, which is non-zero by \autoref{thm:lifting-of-2theta5},
  must be detected by
  $\alpha =\bar{e_{48}h_0h_3^2}$ or $\bar{e_{51}P^1h_2}$ in \autoref{fig:P62-46:Ext}.
  From \autoref{fig:P62-47:Ext} we see that
  $(q')_\ast(\alpha) \ne0$
  in $\Ext_\A^{\ast,\ast}(P^{62}_{47})$,
  which is also denoted by $\bar{e_{48}h_0h_3^2}$ and $\bar{e_{51}P^1h_2}$ respectively,
  where $P^{62}_{46} \to{q'}P^{62}_{47}$ is the collapsing map.
  Comparing \autoref{fig:P62-46:Ext} and \autoref{fig:P62-47:Ext}
  we see that
  $\Ext_\A^{s,63+s}(P^{62}_{46}) \iso\Ext_\A^{s,63+s}(P^{62}_{47})$
  for $0\le s\le3$, so
  if $\alpha$ is not a boundary in the ASS for $\pi_\ast^S(P^{62}_{46})$,
  then $(q')_\ast(\alpha)$ is not a boundary in the ASS for $\pi_\ast^S(P^{62}_{47})$ either.
  Therefore $2q_{47}\widehat{\theta}_5 =q'(2q_{46}\widehat{\theta}_5)$
  is detected by $(q')_\ast(\alpha)$ ($=\bar{e_{48}h_0h_3^2}$ or $\bar{e_{51}P^1h_2}$)
  in the ASS for $\pi_\ast^S(P^{62}_{47})$.
  By \autoref{thm:q47t5:AF-ge-3},
  $AF(q_{47}\widehat{\theta}_5)\ge3$.
  So $AF(2q_{47}\widehat{\theta}_5)\ge4$.
  Therefore $2q_{47}\widehat{\theta}_5$
  must be detected by 
  $\bar{e_{51}P^1h_2} \in\Ext_\A^{5,67}(P^{62}_{47})$.
  So $2q_{46}\widehat{\theta}_5$
  has to be detected by 
  $\bar{e_{51}P^1h_2} \in\Ext_\A^{5,67}(P^{62}_{46})$.
  This completes the proof of \autoref{thm:2q46t5:detected}.
\end{proof}
\begin{prop}
  \label{thm:P62-46:ASS:d5}
  $d_5(\bar{e_{63}h_0}) =\bar{e_{48}h_0^2d_0}$ and
  $d_5(\bar{e_{63}h_0^2}) =\bar{e_{46}}P^1c_0$
  in the ASS for $\pi_\ast^S(P^{62}_{46})$
  and the ASS for $\pi_\ast^S(P^{62}_{47})$.
\end{prop}
\begin{proof}
  By \autoref{thm:2q46t5:detected},
  $2q_{46}\widehat{\theta}_5 =q_{46}(2\widehat{\theta}_5)$
  is detected by $\bar{e_{51}P^1h_2}$ in the ASS for $\pi_\ast^S(P^{62}_{46})$.
  Consider the inclusion map $P^{62}_{46} \to{i_{62}}P_{46}$.
  Recall that $2\widehat{\theta}_5 =0$ in $\pi_\ast^S(P)$, so that
  $i_{62}q_{46}(2\widehat{\theta}_5) =0$ in $\pi_\ast^S(P_{46})$.
  From \autoref{fig:P-46:Ext} we see
  $(i_{62})_\ast(\bar{e_{51}P^1h_2})\ne0$ in $\Ext_\A^{5,67}(P_{46})$,
  which is also denoted by $\bar{e_{51}P^1h_2}$.
  So $\bar{e_{51}P^1h_2}$ must be a boundary in the ASS for $\pi_\ast^S(P_{46})$.
  Comparing \autoref{fig:P62-46:Ext} and \autoref{fig:P-46:Ext},
  we see $\Ext_\A^{s,63+s}(P_{46}) \iso\Ext_\A^{s,63+s}(P^{62}_{46})$
  for $1\le s\le3$ and
  $\Ext_\A^{0,63}(P_{46}) =\Z/2(\bar{e_{63}})$.
  Thus $d_5(\bar{e_{63}}) =\bar{e_{51}P^1h_2}$
  in the ASS for $\pi_\ast^S(P_{46})$.
  From \autoref{fig:P62-46:Ext},
  $\bar{e_{51}P^1h_2}h_0 =\bar{e_{48}h_0^2d_0}$ and
  $\bar{e_{51}P^1h_2}h_0^2 =\bar{e_{48}h_0^2d_0}h_0 =\bar{e_{46}}P^1c_0$, so it follows that
  $d_5(\bar{e_{63}h_0}) =\bar{e_{48}h_0^2d_0}$ and
  $d_5(\bar{e_{63}h_0^2}) =\bar{e_{46}}P^1c_0$
  in the ASS for $\pi_\ast^S(P_{46})$.
  Comparing \autoref{fig:P62-47:Ext}, \autoref{fig:P62-46:Ext}, and \autoref{fig:P-46:Ext}
  we see that these two differentials
  also hold in the ASS for $\pi_\ast^S(P^{62}_{46})$
  and the ASS for $\pi_\ast^S(P^{62}_{47})$.
  This proves \autoref{thm:P62-46:ASS:d5}.
\end{proof}

We proceed to show an infinite cycle
in the ASS for $\pi_\ast^S(P^{62}_{46})$ and the ASS for $\pi_\ast^S(P^{62}_{47})$.
\begin{prop}
  \label{thm:P62-46:ASS:e49h32}
  $\bar{e_{49}h_3^2}$ is an infinite cycle
  in the ASS for $\pi_\ast^S(P^{62}_{46})$
  and the ASS for $\pi_\ast^S(P^{62}_{47})$.
\end{prop}
\begin{proof}
  Consider the complex $P^{49}_{46}=(S^{46} \V S^{47})\U_{(\eta,2\iota)}e^{48} \U_\eta e^{49}$.
  We have the following cofibration sequence.
  \begin{align*}
    \xymatrix{
    \cdots\ar[r]&
    P^{49}_{46} =(S^{46} \V S^{47})\U_{(\eta,2\iota)}e^{48} \U_\eta e^{49} \ar[r]&
    S^{49}\ar[r]&
    (S^{47} \V S^{48})\U_{(\eta,2\iota)}e^{49}\ar[r]&
    \cdots\\
    && S^{63}\ar[u]^{\sigma^2}\ar[ul]^{\widetilde{\sigma^2}}&
    }
  \end{align*}
  Since $\eta\sigma^2 =0$ there is a coextension
  $S^{63} \to{\widetilde{\sigma^2}}P^{49}_{46}$ to
  $S^{63} \to{\sigma^2}S^{49}$.
  Since $AF(\sigma^2)=2$ in the ASS for $\pi_\ast^S$, we have
  $AF(\widetilde{\sigma^2})\le2$ in the ASS for $\pi_\ast^S(P^{49}_{46})$.
  From \autoref{fig:P49-46:Ext} we see
  $S^{63} \to{\widetilde{\sigma^2}}P^{49}_{46}$ is detected by
  $\bar{e_{49}h_3^2}$ in the ASS for $\pi_\ast^S(P^{49}_{46})$.
  So $\bar{e_{49}h_3^2} \in\Ext_\A^{2,65}(P^{49}_{46})$ is an infinite cycle
  and in the ASS for $\pi_\ast^S(P^{49}_{46})$.
  It follows that
  $\bar{e_{49}h_3^2}$ is also an infinite cycle
  in the ASS for $\pi_\ast^S(P^{62}_{46})$
  and also in the ASS for $\pi_\ast^S(P^{62}_{47})$.
  This completes the proof of \autoref{thm:P62-46:ASS:e49h32}.
\end{proof}

We proceed to prove \autoref{thm:P62-46:pi62} and \autoref{thm:P62-47:pi62}.

To describe the stable homotopy groups
$\pi_{62}^S(P^{62}_{46})$ and $\pi_{62}^S(P^{62}_{47})$ in \autoref{thm:P62-46:pi62} and \autoref{thm:P62-47:pi62},
we need to recall the stable homotopy elements
$y,w,x\in\pi_{62}^S(P^{62}_{46})$ and $v\in\pi_{62}^S(P^{62}_{47})$ defined in
\eqref{eq:y},
\eqref{eq:w},
\eqref{eq:x}, and \eqref{eq:v} respectively.
In Propositions \ref{thm:y}, \ref{thm:w}, \ref{thm:x}, and \ref{thm:v} below 
we recall their definitions and show that they are essential.
Recall the Hopf classes $S^1\to{\eta}S^0$, $S^3\to{\nu}S^0$ and $S^7\to{\sigma}S^0$.
To show that the stable homotopy elements $y,w,x$ and $v$ are essential
we need the information of the following $\Ext$ group structures:
$\Ext_\A^{s,t}(S^{47}\U_\sigma e^{55})$ for $s\ge0$ and $61\le t-s\le 63$,
$\Ext_\A^{s,t}(S^{47}\U_\nu e^{55}\U_\eta e^{53})$ for $s\ge0$ and $61\le t-s\le 63$, and
$\Ext_\A^{s,t}(S^{47}\U_{2\iota} e^{48}\U_\nu e^{51}\U_\eta e^{53})$ for $0\le s\le4$ and $61\le t-s\le63$.
These are given
\autoref{fig:Csigma:Ext} and \autoref{fig:Cnueta:Ext}.
The computations of these groups
are analogous to that of
$\Ext_\A^{s,t}(M_2)$ in \autoref{fig:M2:Ext}, and
will not be repeated here.
\begin{figure}
  \begin{minipage}{.4\textwidth}
    \begin{center}
\begin{sseq}[grid=crossword,entrysize=10mm,xlabelstep=1,ylabelstep=4]{61...63}{5}
  \ssmoveto {47} 0  \ssdropcirc \ssname{e47}\ssdroplabel[D]{\bar{e_{47}}}
  \ssmove 0 1       \ssdropcirc \ssstroke
  \ssmove 0 1       \ssdropcirc \ssstroke
  \ssmove 0 1       \ssdropcirc \ssstroke
  \ssmove 0 1       \ssdropcirc \ssstroke
  \ssmove 0 1       \ssdropcirc \ssstroke
  \ssmove 0 1       \ssdropcirc \ssstroke
  \ssmove 0 1       \ssdropcirc \ssstroke
  \ssmove 0 1       \ssdropcirc \ssstroke
  \ssmove 0 1       \ssdropcirc \ssstroke
  \ssmove 0 1       \ssdropcirc \ssstroke
  \ssgoto{e47}
  \ssmove 1 1       \ssdropcirc \ssstroke
  \ssmove 1 1       \ssdropcirc \ssstroke
  \ssmove 1 1       \ssdropcirc \ssstroke
  \ssmove 0 {-1}    \ssdropcirc \ssstroke
  \ssmove 0 {-1}    \ssdropcirc \ssstroke
                                \ssname{e47h2}\ssdroplabel[D]{\bar{e_{47}}h_2}

  \ssmoveto {53} 2  \ssdropcirc \ssname{e47h22}\ssdroplabel[D]{\bar{e_{47}h_2^2}}

  \ssmoveto {55} 3  \ssdropcirc \ssname{e47c0}\ssdroplabel[D]{\bar{e_{47}}c_0}
  \ssmove 1 1       \ssdropcirc \ssstroke

  \ssmoveto {55} 4  \ssdropcirc \ssname{e55h03}\ssdroplabel[D]{\bar{e_{55}h_0^4}}
  \ssmove 0 1       \ssdropcirc \ssstroke
  \ssmove 0 1       \ssdropcirc \ssstroke
  \ssmove 0 1       \ssdropcirc \ssstroke
  \ssmove 0 1       \ssdropcirc \ssstroke
  \ssmove 0 1       \ssdropcirc \ssstroke
  \ssmove 0 1       \ssdropcirc \ssstroke

  \ssmoveto {56} 5  \ssdropcirc \ssname{e47P1h1}\ssdroplabel[D]{\bar{e_{47}}P^1h_1}
  \ssmove 1 1       \ssdropcirc \ssstroke
  \ssmove 1 1       \ssdropcirc \ssstroke
  \ssmove 0 {-1}    \ssdropcirc \ssstroke
  \ssmove 0 {-1}    \ssdropcirc \ssstroke
                                \ssname{e47P1h2}\ssdroplabel[D]{\bar{e_{47}}P^1h_2}

  \ssmoveto {58} 1  \ssdropcirc \ssname{e55h2}\ssdroplabel[D]{\bar{e_{55}}h_2}
  \ssmove 0 1       \ssdropcirc \ssstroke
  \ssmove 0 1       \ssdropcirc \ssstroke

  \ssmoveto {61} 2  \ssdropcirc \ssname{e55h22}\ssdroplabel[D]{\bar{e_{55}h_2^2}}

  \ssmoveto {61} 4  \ssdropcirc \ssname{e47d0}\ssdroplabel[D]{\bar{e_{47}}d_0}
  \ssmove 0 1       \ssdropcirc \ssstroke
  \ssmove 0 1       \ssdropcirc \ssstroke
  \ssgoto{e47d0}
  \ssmove 1 1       \ssdropcirc \ssstroke

  \ssmoveto {62} 3  \ssdropcirc \ssname{e55h02h3}\ssdroplabel[LD]{\bar{e_{55}h_0^2h_3}}
  \ssmove 0 1       \ssdropcirc \ssstroke

  \ssmoveto {62} 1  \ssdropcirc \ssname{e47h4}\ssdroplabel[D]{\bar{e_{47}}h_4}
  \ssmove 0 1       \ssdropcirc \ssstroke
  \ssmove 0 1       \ssdropcirc \ssstroke
  \ssmove 0 1       \ssdropcirc \ssstroke
  \ssmove 0 1       \ssdropcirc \ssstroke
  \ssmove 0 1       \ssdropcirc \ssstroke
  \ssmove 0 1       \ssdropcirc \ssstroke
  \ssmove 0 1       \ssdropcirc \ssstroke
  \ssgoto{e47h4}
  \ssmove 1 1       \ssdropcirc \ssstroke
  \ssmove 1 1       \ssdropcirc \ssstroke
  \ssmove 1 1       \ssdropcirc \ssstroke

  \ssmoveto {63} 3  \ssdropcirc \ssname{e55c0}
  \ssmove 1 1       \ssdropcirc \ssstroke

  \ssmoveto {63} 7  \ssdropcirc \ssname{e47pc0}\ssdroplabel[D]{\bar{e_{47}}P^1c_0}
  \ssmove 1 1       \ssdropcirc \ssstroke

  \ssmoveto {63} 2  \ssdropcirc \ssname{e55h1h3}\ssdroplabel[D]{\bar{e_{55}h_1h_3}}
  \ssmove 1 1       \ssdropcirc \ssstroke
  \ssmove 1 1       \ssdropcirc \ssstroke
\end{sseq}
    \end{center}
    \caption{$\Ext_\A^{s,t}(S^{47}\U_\sigma e^{55})$, $0\le s\le4$, $61\le t-s\le63$}
    \label{fig:Csigma:Ext}
  \end{minipage}
  \qquad
  \begin{minipage}{.4\textwidth}
    \begin{center}
\begin{sseq}[grid=crossword,entrysize=10mm,xlabelstep=1,ylabelstep=4]{61...63}{5}
  \ssmoveto {47} 0  \ssdropcirc \ssname{e47}\ssdroplabel[D]{e_{47}}
  \ssmove 0 1       \ssdropcirc \ssstroke
  \ssmove 0 1       \ssdropcirc \ssstroke
  \ssmove 0 1       \ssdropcirc \ssstroke
  \ssmove 0 1       \ssdropcirc \ssstroke
  \ssmove 0 1       \ssdropcirc \ssstroke
  \ssmove 0 1       \ssdropcirc \ssstroke
  \ssmove 0 1       \ssdropcirc \ssstroke
  \ssmove 0 1       \ssdropcirc \ssstroke
  \ssmove 0 1       \ssdropcirc \ssstroke
  \ssmove 0 1       \ssdropcirc \ssstroke
  \ssmove 0 1       \ssdropcirc \ssstroke
  \ssmove 0 1       \ssdropcirc \ssstroke
  \ssgoto{e47}
  \ssmove 1 1       \ssdropcirc \ssstroke
  \ssmove 1 1       \ssdropcirc \ssstroke

  \ssmoveto {51} 3  \ssdropcirc \ssname{e51h03}\ssdroplabel[D]{\bar{e_{51}h_0^3}}
  \ssmove 0 1       \ssdropcirc \ssstroke
  \ssmove 0 1       \ssdropcirc \ssstroke
  \ssmove 0 1       \ssdropcirc \ssstroke
  \ssmove 0 1       \ssdropcirc \ssstroke
  \ssmove 0 1       \ssdropcirc \ssstroke
  \ssmove 0 1       \ssdropcirc \ssstroke
  \ssmove 0 1       \ssdropcirc \ssstroke
  \ssmove 0 1       \ssdropcirc \ssstroke
  \ssmove 0 1       \ssdropcirc \ssstroke
  \ssmove 0 1       \ssdropcirc \ssstroke

  \ssmoveto {53} 1  \ssdropcirc \ssname{e53h0}\ssdroplabel[D]{\bar{e_{53}h_0}}
  \ssmove 0 1       \ssdropcirc \ssstroke
  \ssmove 0 1       \ssdropcirc \ssstroke
  \ssmove 0 1       \ssdropcirc \ssstroke
  \ssmove 0 1       \ssdropcirc \ssstroke
  \ssmove 0 1       \ssdropcirc \ssstroke
  \ssmove 0 1       \ssdropcirc \ssstroke
  \ssmove 0 1       \ssdropcirc \ssstroke
  \ssmove 0 1       \ssdropcirc \ssstroke
  \ssmove 0 1       \ssdropcirc \ssstroke
  \ssmove 0 1       \ssdropcirc \ssstroke
  \ssmove 0 1       \ssdropcirc \ssstroke
  \ssmove 0 1       \ssdropcirc \ssstroke
  \ssmove 0 1       \ssdropcirc \ssstroke

  \ssmoveto {54} 2  \ssdropcirc \ssname{e51h0h2}\ssdroplabel[D]{\bar{e_{51}h_0h_2}}

  \ssmoveto {54} 1  \ssdropcirc \ssname{e47h3}\ssdroplabel[D]{\bar{e_{47}}h_3}
  \ssmove 0 1       \ssdropcirc \ssstroke
  \ssmove 0 1       \ssdropcirc \ssstroke
  \ssmove 0 1       \ssdropcirc \ssstroke
  \ssgoto{e47h3}
  \ssmove 1 1       \ssdropcirc \ssstroke

  \ssmoveto {55} 3  \ssdropcirc \ssname{e47c0}\ssdroplabel[D]{\bar{e_{47}}c_0}
  \ssmove 1 1       \ssdropcirc \ssstroke

  \ssmoveto {56} 5  \ssdropcirc \ssname{e47P1h1}\ssdroplabel[D]{\bar{e_{47}}P^1h_1}
  \ssmove 1 1       \ssdropcirc \ssstroke
  \ssmove 1 1       \ssdropcirc \ssstroke
  \ssmove 0 {-1}    \ssdropcirc \ssstroke
  \ssmove 0 {-1}    \ssdropcirc \ssstroke
                                \ssname{e47P1h2}\ssdroplabel[D]{\bar{e_{47}}P^1h_2}

  \ssmoveto {56} 1  \ssdropcirc \ssname{e53h2}\ssdroplabel[D]{\bar{e_{53}h_2}}
  \ssmove 0 1       \ssdropcirc \ssstroke
  \ssmove 0 1       \ssdropcirc \ssstroke

  \ssmoveto {58} 1  \ssdropcirc \ssname{e51h3}\ssdroplabel[D]{\bar{e_{51}h_3}}
  \ssmove 0 1       \ssdropcirc \ssstroke
  \ssmove 0 1       \ssdropcirc \ssstroke
  \ssmove 0 1       \ssdropcirc \ssstroke

  \ssmoveto {59} 3  \ssdropcirc \ssname{e51c0}\ssdroplabel[D]{\bar{e_{51}c_0}}

  \ssmoveto {59} 2  \ssdropcirc \ssname{e53h22}\ssdroplabel[D]{\bar{e_{53}h_2^2}}

  \ssmoveto {60} 2  \ssdropcirc \ssname{e53h0h3}\ssdroplabel[D]{\bar{e_{53}h_0h_3}}
  \ssmove 0 1       \ssdropcirc \ssstroke
  \ssmove 0 1       \ssdropcirc \ssstroke

  \ssmoveto {60} 5  \ssdropcirc \ssname{e51P1h1}\ssdroplabel[D]{\bar{e_{51}P^1h_1}}

  \ssmoveto {62} 1  \ssdropcirc \ssname{e47h4}\ssdroplabel[D]{\bar{e_{47}}h_4}
  \ssmove 0 1       \ssdropcirc \ssstroke
  \ssmove 0 1       \ssdropcirc \ssstroke
  \ssmove 0 1       \ssdropcirc \ssstroke
  \ssmove 0 1       \ssdropcirc \ssstroke
  \ssmove 0 1       \ssdropcirc \ssstroke
  \ssmove 0 1       \ssdropcirc \ssstroke
  \ssmove 0 1       \ssdropcirc \ssstroke
  \ssgoto{e47h4}
  \ssmove 1 1       \ssdropcirc \ssstroke
  \ssmove 1 1       \ssdropcirc \ssstroke

  \ssmoveto {61} 2  \ssdropcirc \ssname{e47h32}\ssdroplabel[D]{\bar{e_{47}}h_3^2}
  \ssmove 0 1       \ssdropcirc \ssstroke

  \ssmoveto {62} 4  \ssdropcirc \ssname{e53h1c0}\ssdroplabel[R]{\bar{e_{53}h_1c_0}}

  \ssmoveto {62} 3  \ssdropcirc \ssname{e53h12h3}\ssdroplabel[R]{\bar{e_{53}h_1^2h_3}}

  \ssmoveto {61} 4  \ssdropcirc \ssname{e47d0}\ssdroplabel[D]{\bar{e_{47}}d_0}
  \ssmove 0 1       \ssdropcirc \ssstroke
  \ssgoto{e47d0}
  \ssmove 1 1       \ssdropcirc \ssstroke
  \ssmove 1 1       \ssdropcirc \ssstroke

  \ssmoveto {62} 6  \ssdropcirc \ssname{e51h0P1h2}

  \ssmoveto {63} 7  \ssdropcirc \ssname{e47P1c0}\ssdroplabel[D]{\bar{e_{47}}P^1c_0}
  \ssmove 1 1       \ssdropcirc \ssstroke

  \ssmoveto {64} 5  \ssdropcirc \ssname{e53P1h2}\ssdroplabel[D]{\bar{e_{53}P^1h_2}}
  \ssmove 0 1       \ssdropcirc \ssstroke
  \ssmove 0 1       \ssdropcirc \ssstroke

  \ssmoveto {64} 4  \ssdropcirc \ssname{e47e0}\ssdroplabel[D]{\bar{e_{47}}e_0}

  \ssmoveto {64} 9  \ssdropcirc \ssname{e47P2h1}\ssdroplabel[D]{\bar{e_{47}P^2h_1}}
\end{sseq}
    \end{center}
    \caption{$\Ext_\A^{s,t}(W_1 =S^{47}\U_\nu e^{51}\U_\eta e^{53})$, $0\le s\le4$, $61\le t-s\le63$}
    \label{fig:Cnueta:Ext}
  \end{minipage}
\end{figure}
\begin{figure}
  \begin{center}
\begin{sseq}[grid=crossword,entrysize=10mm,xlabelstep=1,ylabelstep=4]{61...63}{5}
  \ssmoveto {62} 1  \ssdropcirc \ssname{e47h4}\ssdroplabel[D]{\bar{e_{47}}h_4}
  \ssmove 1 1       \ssdropcirc \ssstroke
  \ssmove 1 1       \ssdropcirc \ssstroke

  \ssmoveto {61} 2  \ssdropcirc \ssname{e47h32}\ssdroplabel[D]{\bar{e_{47}}h_3^2}

  \ssmoveto {62} 4  \ssdropcirc \ssname{e53h1c0}\ssdroplabel[D]{\bar{e_{53}h_1c_0}}

  \ssmoveto {62} 3  \ssdropcirc \ssname{e53h12h3}\ssdroplabel[D]{\bar{e_{53}h_1^2h_3}}

  \ssmoveto {61} 4  \ssdropcirc \ssname{e47d0}\ssdroplabel[D]{\bar{e_{47}}d_0}
  \ssmove 1 1       \ssdropcirc \ssstroke
  \ssmove 1 1       \ssdropcirc \ssstroke

  \ssmoveto {62} 6  \ssdropcirc \ssname{e51h0P1h2}\ssdroplabel[L]{\bar{e_{51}h_0P^1h_2}}

  \ssmoveto {63} 7  \ssdropcirc \ssname{e47P1c0}\ssdroplabel[D]{\bar{e_{47}}P^1c_0}
  \ssmove 1 1       \ssdropcirc \ssstroke
\end{sseq}
  \end{center}
  \caption{$\Ext_\A^{s,t}(W_2 =S^{47}\U_{2\iota} e^{48}\U_\nu e^{51}\U_\eta e^{53})$, $0\le s\le4$, $61\le t-s\le63$}
  \label{fig:Cnueta2iota:Ext}
\end{figure}
\begin{prop}
  \label{thm:y}
  Let $P^{48}_{46}\to{i_{48}}P^{62}_{46}$ be the inclusion map.
  Let $S^{62}\to{\widetilde{\sigma^2}}P^{48}_{46}$ be some coextension
  of $S^{62}\to{\sigma^2}S^{48}$ as shown in the following diagram.
  \begin{align*}
    \xymatrix{
    \cdots\ar[r]&
    P^{48}_{46}\ar[r]^-{q_{48}}& S^{48}\ar[r]^-{\eta\V2\iota}& S^{47}\V S^{48}\ar[r]&  \cdots.\\
    & & S^{62}\ar[u]_{\sigma^2}\ar[ul]^{\widetilde{\sigma^2}}\\
    }
  \end{align*}
  Let $y$ be the composite $S^{62}\to{\widetilde{\sigma^2}}P^{48}_{46}\to{i_{48}}P^{62}_{46}$.
  Then $y$ is detected by $\bar{e_{47}}h_4\in\Ext_\A^{1,63}(P^{62}_{46})$ in the ASS for $\pi_\ast^S(P^{62}_{46})$.
\end{prop}
\begin{proof}
  $S^{62}\to{\sigma^2}S^{48}$ is detected by $h_3^2 \in\Ext_\A^{2,16}$ in the ASS for $\pi_\ast^S$, so $AF(\sigma^2)=2$.
  Thus $AF(\widetilde{\sigma^2})\le2$ in the ASS for $\pi_\ast^S(P^{48}_{46})$.
  From \autoref{fig:P48-46:Ext}, the nonzero stable homotopy class $\widetilde{\sigma^2}$ must be detected by
  $\bar{e_{47}}h_4 \in\Ext_\A^{1,63}(P^{48}_{46})$ or $\bar{e_{47}}h_0h_4 \in\Ext_\A^{2,64}(P^{48}_{46})$.
  Suppose $\widetilde{\sigma^2}$ is detected by $\bar{e_{47}}h_0h_4$ in the ASS for $\pi_\ast^S(P^{48}_{46})$,
  then $AF(\widetilde{\sigma}^2) =2$.
  It is not difficult to see that $(q_{48})_\ast(\bar{e_{47}}h_0h_4)=0$ in $\Ext_\A^{2,64}$, so it follows that
  $AF(\sigma^2 =q_{48}(\widetilde{\sigma^2}))\ge3$ in the ASS for $\pi_\ast^S$,
  but we have seen that $AF(\sigma^2) =2$ so this is a contradiction.
  Therefore $\widetilde{\sigma^2}$ is detected by $\bar{e_{47}}h_4$ in the ASS for $\pi_\ast^S(P^{48}_{46})$.
  From \autoref{fig:P62-46:Ext}, $\bar{e_{47}}h_4 \ne0$ in $\Ext_\A^{1,63}(P^{62}_{46})$
  is not a boundary in the ASS for $\pi_\ast^S(P^{62}_{46})$,
  so $y=i\widetilde{\sigma^2}\in\pi_{62}^S(P^{62}_{46})$ is essential and
  is detected by $\bar{e_{47}}h_4\in\Ext_\A^{1,63}(P^{62}_{46})$ in the ASS for $\pi_\ast^S(P^{62}_{46})$.
  This proves \autoref{thm:y}.
\end{proof}
\begin{prop}
  \label{thm:w}
  Let $S^{47}\U_\sigma e^{55}\to{j}P^{62}_{46}$ be the inclusion map.
  Let $S^{62}\to{\bar{2\sigma}}S^{47}\U_\sigma e^{55}$ be some coextension of
  $S^{62}\to{2\sigma}S^{55}$ as shown in the following diagram.
  \begin{align*}
    \xymatrix{
    \cdots\ar[r]&
      S^{47}\U_\sigma e^{55}\ar[r]^-{q_{55}}&
      S^{55}\ar[r]^\sigma&  S^{48}\ar[r]& \cdots\\
    & & S^{62}\ar[u]^{2\sigma}\ar[ul]^{\bar{2\sigma}}
    }
  \end{align*}
  Let $w$ be the composite $S^{62}\to{\bar{2\sigma}}S^{47}\U_\sigma e^{55}\to{j}P^{62}_{46}$.
  Then $w$ is detected by $\bar{e_{47}}h_4 \in\Ext_\A^{1,63}(P^{62}_{46})$
  in the ASS for $\pi_\ast^S(P^{62}_{46})$.
\end{prop}
\begin{proof}
  $S^{62} \to{2\sigma}S^{55}$ is detected by $h_0h_3$ in the ASS for $\pi_\ast^S$, so $AF(2\sigma) =2$.
  Thus $AF(\bar{2\sigma})\le2$ in the ASS for $\pi_\ast^S(S^{47}\U_\sigma e^{55})$.
  From \autoref{fig:Csigma:Ext}, the nonzero stable homotopy class $\bar{2\sigma}$ must be detected by
  $\bar{e_{47}}h_4 \in\Ext_\A^{1,63}(S^{47}\U_\sigma e^{55})$ or
  $\bar{e_{47}}h_0h_4 \in\Ext_\A^{2,64}(S^{47}\U_\sigma e^{55})$.
  Suppose $\bar{2\sigma}$ is detected by $\bar{e_{47}}h_0h_4$ in the ASS for $\pi_\ast^S(S^{47}\U_\sigma e^{55})$,
  then $AF(\bar{2\sigma}) =2$.
  It is not difficult to see that $(q_{55})_\ast(e_{47}h_0h_4) =0$ in $\Ext_\A^{2,64}$, so it follows that
  $AF(2\sigma =q_{55}\bar{2\sigma})\ge3$,
  but we already have shown that $AF(2\sigma) =2$ so this is a contradiction.
  Therefore $\bar{2\sigma}$ is detected by $\bar{e_{47}}h_4$ in the ASS for $\pi_\ast^S(S^{47}\U_\sigma e^{55})$.
  From \autoref{fig:P62-46:Ext}, $\bar{e_{47}}h_4 \ne0$ in $\Ext_\A^{1,63}(P^{62}_{46})$
  is not a boundary in the ASS for $\pi_\ast^S(P^{62}_{46})$,
  so $w=j\bar{2\sigma}\in\pi_{62}^S(P^{62}_{46})$ is essential and
  is detected by $\bar{e_{47}}h_4\in\Ext_\A^{1,63}(P^{62}_{46})$ in the ASS for $\pi_\ast^S(P^{62}_{46})$.
  This proves \autoref{thm:w}.
\end{proof}
\begin{prop}
  \label{thm:x}
  $x=y+w\in\pi_{62}^S(P^{62}_{46})$ is detected by $\bar{e_{55}h_0h_3} +\varepsilon\bar{e_{47}}h_0h_4 \in\Ext_\A^{2,64}(P^{62}_{46})$, $\varepsilon=0$ or $1$.
\end{prop}
\begin{proof}
  Let $P^{62}_{46} \to{q_{49}}P^{62}_{49}$ be the quotient map and
  consider the composite $S^{62}\to{x=y+w}P^{62}_{46}\to{q_{49}}P^{62}_{49}$.
  From the definition of $y$ it is not difficult to see that
  $q_{49}y=0$ in $\pi_{62}^S(P^{62}_{49})$, as $q_{49}$ pinches $P^{48}_{46}$.
  It is also not difficult to see that $S^{55}$ is a subcomplex of $P^{62}_{49}$.
  Let $S^{55} \to{i_{55}}P^{62}_{49}$ be the inclusion map. Then from the definition of $x$ and $w$,
  $q_{49}x =q_{49}(y+w) =q_{49}w$ is the composite $S^{62} \to{2\sigma}S^{55} \to{i_{55}}P^{62}_{49}$.
  $S^{62} \to{2\sigma}S^{55}$ is detected by $h_0h_3$ in the ASS for $\pi_\ast^S$.
  Recall that, since
  $ H_\ast(P^{62}_{49}) \iso H_{\ast-48}(P^{14})$
  as left-$\A$ modules, $\Ext_\A^{s,t}(P^{62}_{49}) \iso\Ext_\A^{s,t-48}(P^{14})$.
  $\bar{e_{55}h_0h_3} \in\Ext_\A^{2,64}(P^{62}_{49})$ corresponds to
  $\widehat{h}_3h_0h_3 \in\Ext_\A^{2,16}(P^{14})$ under this isomorphism.
  So from \autoref{fig:P14:Ext}, we deduce that
  $\bar{e_{55}h_0h_3} \in\Ext_\A^{2,64}(P^{62}_{49})$ is not a boundary in the ASS for $\pi_\ast^S(P^{62}_{49})$
  and $\Ext_\A^{2,64}(P^{62}_{49}) =\Z/2(\bar{e_{55}h_0h_3})$. So
  $q_{49}x =i_{55}(2\sigma) \in\pi_{62}^S(P^{62}_{49})$ is detected by $\bar{e_{55}h_0h_3} \in\Ext_\A^{2,64}(P^{62}_{49})$,
  and $AF(q_{49}x =i_{55}(2\sigma)) =2$ in the ASS for $\pi_\ast^S(P^{62}_{49})$.
  It follows that $AF(x)\le2$ in the ASS for $\pi_\ast^S(P^{62}_{46})$.
  Since both $y$ and $w$ are detected by $\bar{e_{47}}h_4$ in the ASS for $\pi_\ast^S(P^{62}_{46})$,
  $AF(x =y+w)\ge2$.
  Therefore $AF(x) =2$.
  From \autoref{fig:P62-46:Ext}
  $x$ must be detected by $\varepsilon_1\bar{e_{55}h_0h_3} +\varepsilon_2\bar{e_{47}}h_0h_4\in\Ext_\A^{2,64}(P^{62}_{46})$
  for $\varepsilon_i=0$ or $1$, $i=1,2$.
  Because $q_{49}x$ is detected by $\bar{e_{55}h_0h_3}$,
  $x$ must be detected by $\bar{e_{55}h_0h_3} +\varepsilon_2\bar{e_{47}}h_0h_4$ where $\varepsilon_2=0$ or $1$
  (that is, $\varepsilon_1 =1$).
  This proves \autoref{thm:x}.
\end{proof}

We are now ready to prove \autoref{thm:P62-46:pi62}, which
is restated as follows, where
$y,w$ and $x$ are as in Propositions \ref{thm:y}, \ref{thm:w} and \ref{thm:x} respectively.
\theoremstyle{nonumberplain}
\newtheorem{thm:P62-46:pi62}{Theorem \ref{thm:P62-46:pi62}} 
\begin{thm:P62-46:pi62*}
  $\pi^S_{62}(P^{62}_{46})=\Z/4(y)\dsum\Z/8(x)$.
  Moreover $AF(y)=1$, $AF(w)=1$, $AF(2y)=2$, and $AF(x =y+w)=2$.
\end{thm:P62-46:pi62*}
\begin{proof}[Proof of \autoref{thm:P62-46:pi62}]
  We refer to \autoref{fig:P62-46:Ext} for our discussions on the Adams spectral sequence for $\pi_\ast^S(P^{62}_{46})$
  from which we want to show $\pi_{62}^S(P^{62}_{46}) =\Z/2(y) \dsum\Z/2(x)$.
  \begin{enumerate}
    \item 
      We have shown in Propositions \ref{thm:y}, \ref{thm:w} and \ref{thm:x} that
      $y$ and $w$ are detected by $\bar{e_{47}}h_4$ and
      $x$ is detected by $\bar{e_{55}h_0h_3} +\varepsilon\bar{e_{47}}h_0h_4$, so
      $AF(y) =AF(w) =1$ and $AF(x)=2$.
      Note that this implies that $\bar{e_{47}}h_0h_4$ is
      an infinite cycle in the ASS for $\pi_\ast^S(P^{62}_{46})$.
      From \autoref{fig:P62-46:Ext}, we see
      $\bar{e_{47}}h_0h_4$ is not a boundary. Since $y$ is
      detected by $\bar{e_{47}}h_4$ it follows that $2y$ is
      detected by $\bar{e_{47}}h_0h_4$.
      Therefore $AF(2y)=2$.
    \item
      To show $4y=0$ in $\pi_{62}^S(P^{62}_{46})$, consider the following diagram, where
      the composite
      $S^{62} \to{y}P^{48}_{46} \to{q_{48}}S^{48}$ is $\sigma^2$, and
      $S^{46}\to{\sigma_{46}}S^{46}\to{i_{46}}P^{47}_{46}\to{q_{47}}S^{47}$,
      $S^{47}\to{\sigma_{47}}S^{47}\to{i_{47}}P^{48}_{47}\to{q_{48}}S^{48}$ are
      cofibration sequences.
      \begin{align*}
        \xymatrix{
        & & S^{46}\ar[r]^{\sigma_{46}}& S^{46}\ar[d]^{i_{46}}\\
        & & S^{47}\ar[r]^{\sigma_{47}}& P^{47}_{46}\ar[d]^{i_{47}}\ar[r]^{q_{47}}& S^{47}\\
        S^{62}\ar[r]^{2\iota}
        &S^{62}\ar[r]^{2\iota}\ar@{.>}[urr]^{\bar{2y}}\ar@/^1pc/[uurr]^{\widetilde{2y}}
        & S^{62}\ar[r]^y& P^{48}_{46}\ar[r]^{q_{48}}& S^{48}
        }
      \end{align*}
      Since $S^{62} \to{2\iota}S^{62} \to{y}P^{48}_{46} \to{q_{48}}S^{48}$
      is $2\sigma^2 =0$, the composite $S^{62} \to{2\iota}S^{62} \to{y}P^{48}_{46}$
      can be lifted to $S^{62} \to{\bar{2y}}P^{47}_{46}$.
      Since $S^{62} \to{\bar{2y}}P^{47}_{46} \to{q_{47}}S^{47}$
      is the Toda bracket $\toda{2\iota,\sigma^2,2\iota} =\eta\sigma^2 =0$,
      the map $S^{62} \to{\bar{2y}}P^{47}_{46}$
      can be lifted to $S^{62} \to{\widetilde{2y}}S^{46}$, and
      the composite $S^{62} \to{2\iota}S^{62} \to{\widetilde{2y}}S^{46}$
      is contained in the Toda bracket
      \begin{align*}
        2\toda{\eta,\sigma^2,2\iota}=\toda{2\iota,\eta,\sigma^2}2\iota.
      \end{align*}
      We have $\toda{2\iota,\eta,\sigma^2} =0$ and this is proved as follows.
      $\toda{2\iota,\eta,\sigma^2}\subset\pi_{16}^S =\Z/2(\{h_1h_4\})\dsum\Z/2(\{P^1c_0\})$, where
      $\{h_1h_4\},\{P^1c_0\}$ denotes the stable homotopy elements detected by these classes
      in the ASS for spheres. Since
      $AF(\toda{2\iota,\eta,\sigma^2})\ge3$, and since
      $\eta\toda{2\iota,\eta,\sigma^2}= \toda{\eta,2\iota,\eta}\sigma^2=\{2\nu,0\}\sigma^2=\{0\}$
      but $\eta\{P^1c_0\}\ne0$, we see that $\toda{2\iota,\eta,\sigma^2} =0$.
      This proves that $4y =0$, and therefore the order of $y$ is $4$.
    \item In the ASS for $\pi_\ast^S(P^{62}_{46})$,
      $d_5(\bar{e_{63}h_0}) =\bar{e_{48}h_0^2d_0}$, and
      $d_5(\bar{e_{63}h_0}h_0) =\bar{e_{46}}P^1c_0$ by \autoref{thm:P62-46:ASS:d5}, so the classes
      $\bar{e_{48}h_0^2d_0},\bar{e_{46}}P^1c_0$
      do not represent homotopy classes.
      It follows that the order of $x$ is at most $8$.
    \item
      To show that $2x \ne0$ in $\pi_\ast^S(P^{62}_{46})$, recall from
      \eqref{eq:P62-47:Ext:extension} (1) that
      there is an extension $\bar{e_{55}h_0h_3}h_0 =\bar{e_{48}h_0h_3^2}$ and that
      $\bar{e_{47}}h_0^2h_4 =0$ in $\Ext_\A^{3,65}(P^{62}_{46})$.
      Since $x$ is detected by $\bar{e_{55}h_0h_3} +\varepsilon\bar{e_{47}}h_0h_4$,
      $(\bar{e_{55}h_0h_3} +\varepsilon\bar{e_{47}}h_0h_4)h_0
      =\bar{e_{55}h_0h_3}h_0
      =\bar{e_{48}h_0h_3^2}$ is an infinite cycle in the ASS for $\pi_\ast^S(P^{62}_{46})$,
      and from \autoref{fig:P62-46:Ext} and (3) we see
      $\bar{e_{48}h_0h_3^2}$ is not a boundary.
      So $2x$ is detected by $\bar{e_{48}h_0h_3^2} =\bar{e_{55}h_0h_3}h_0$
      in the ASS for $\pi_\ast^S(P^{62}_{46})$.
    \item
      To show that $4x \ne0$, note that $4y =0$ by (2), so
      $4x =4(w+y) =4w$.
      Consider the following diagram, where
      $S^{62} \to{\bar{2\sigma}}S^{47}\U_\sigma e^{55}$ is
      the coextension as in \autoref{thm:w},
      the map $S^{47}\U_\sigma e^{55} \to{i'}P^{55}_{46}$ is the inclusion
      of the subcomplex $S^{47}\U_\sigma e^{55}$, and
      for each $51\le l\le55$,
      $S^l \to{f_l}P^l_{46} \to{i_l}P^{l+1}_{46} \to{q_{l+1}}S^{l+1}$ is a cofibration sequence.
      \begin{align*}
        \xymatrix{
        & & & S^{51}\ar[r]^{\sigma_{51}}& P^{51}_{46}\ar[d]^{i_{51}}\ar[r]^{q_{51}}& S^{51}\\
        & & & & \vdots\ar[d]\\
        & & & S^{54}\ar[r]^{\sigma_{54}}& P^{54}_{46}\ar[d]^{i_{54}}\ar[r]^{q_{54}}& S^{54}\\
        & & & S^{55}\ar[r]^{\sigma_{55}}& P^{55}_{46}\ar[d]^{i_{55}}\ar[r]^{q_{55}}& S^{55}\\
        & & & S^{56}\ar[r]^{\sigma_{56}}& P^{56}_{46}\ar[d]^{i_{56}}\ar[r]^{q_{56}}& S^{56}\\
        & & & & \vdots\ar[d]\\
        S^{62}\ar[r]^{2\iota}\ar@/^2pc/[uuuuuurrr]^{\widetilde{4w}}&
        S^{62}\ar[r]^{2\iota}&
        S^{62}\ar[rr]^w\ar@/^2pc/[uuur]^\sigma
        \ar@/^2pc/[uuurr]_-(0.3){\bar{w} =(i'(\bar{2\sigma}))}
        \ar@/^5pc/[uuuurr]^{\bar{w}''}& &
        P^{62}_{46}
        }
      \end{align*}
      $S^{62} \to{\bar{w} =i'(\bar{2\sigma})}P^{55}_{46}$ is a lifting of
      $S^{62} \to{w}P^{62}_{46}$.
      Note that from the cell structure of $P^{56}_{46}$,
      $q_{55}\sigma_{55} =2\iota$.
      So by choosing another lifting
      $S^{62} \to{\bar{w}' =\bar{w} +\sigma_{55}\sigma}P^{55}_{46}$
      of $S^{62} \to{w}P^{62}_{46}$, we have
      $q_{55}\bar{w}' =0$, and $S^{62} \to{\bar{w}'}P^{55}_{46}$
      can be lifted to $S^{62} \to{\bar{w}''}P^{54}_{46}$
      with $q_{54}\bar{w}'' =\eta\sigma$.
      From this and following similar arguments as in (1), it is not difficult
      to see that
      $S^{62} \to{4\iota}S^{62} \to{w}P^{62}_{46}$ has a lifting
      $S^{62} \to{\widetilde{4w}}P^{51}_{46}$ with
      $q_{51}(\widetilde{4w}) =\toda{\eta,\eta^2\sigma,2\iota}$.
      The Toda bracket $\toda{\eta,\eta^2\sigma,2\iota}$
      is detected by $P^1h_2 \in\Ext_\A^{5,16}$ in the ASS for $\pi_\ast^S$.
      By \autoref{thm:P62-46:ASS:e49h32} and (3),
      $\bar{e_{51}P^1h_2} \in\Ext_\A^{5,67}(P^{62}_{46})$
      is not a boundary
      in the ASS for $\pi_\ast^S(P^{62}_{46})$.
      It follows that $4x$ is
      detected by $\bar{e_{51}P^1h_2}$
      in the ASS for $\pi_\ast^S(P^{62}_{46})$.
      This together with (3) proves that
      the order of $x$ is $8$.
  \end{enumerate}
  This completes the proof of \autoref{thm:P62-46:pi62}.
\end{proof}

In order to prove \autoref{thm:P62-47:pi62} we need to recall the
stable homotopy element $v\in\pi_{62}^S(P^{62}_{47})$ in \eqref{eq:v}. In
\autoref{thm:v} below we recall its definition and show that it is essential.
\begin{prop}
  \label{thm:v}
  $W_1 =S^{47}\U_\nu e^{51}\U_\eta e^{53}$ is a subcomplex of $P^{62}_{47}$.
  Let $W_1 \to{j'}P^{62}_{47}$ be the inclusion map.
  Let $S^{62} \to{\widetilde{\eta\varepsilon}}W_1$
  be some coextension of $S^{62} \to{\eta\varepsilon}S^{53}$ as in the following diagram.
  \begin{align*}
    \xymatrix{
    \cdots\ar[r]&
      W_1 =S^{47}\U_\nu e^{51}\U_\eta e^{53}\ar[r]^-{q_4}&
      S^{53}\ar[r]&
      S^{48}\U_\nu e^{52}\ar[r]&
      \cdots\\
    & & S^{62}\ar[u]_{\eta\varepsilon}
      \ar[ul]^{\widetilde{\eta\varepsilon}}
    }
  \end{align*}
  Let $v$ be the composite
  $S^{62} \to{\widetilde{\eta\varepsilon}}W_1 \to{j'}P^{62}_{47}$.
  Then $v$ is detected by $\bar{e_{53}h_1c_0} \in\Ext_\A^{4,66}(P^{62}_{47})$
  in the ASS for $\pi_\ast^S(P^{62}_{47})$.
\end{prop}
\begin{proof}
  $S^{62} \to{\eta\varepsilon}S^{53}$ is detected by $h_1c_0 \in\Ext_\A^{4,13}$
  in the ASS for $\pi_\ast^S$, so
  $AF(\widetilde{\eta\varepsilon})\le4$ in
  the ASS for $\pi_\ast^S(W_1)$.
  From \autoref{fig:Cnueta:Ext} we see that
  $\Ext_\A^{1,63}(W_1) =0$,
  $\Ext_\A^{2,64}(W_1) =\Z/2(\bar{e_{47}}h_4)$, and
  $\Ext_\A^{3,65}(W_1) =\Z/2(\bar{e_{47}}h_0h_4)$.
  It is well-known that
  in the ASS for $\pi_\ast^S$,
  $d_2(h_4) =h_0h_3^2$ and $d_3(h_0h_4) =h_0d_0$,
  so in the ASS for $\pi_\ast^S(S^{47}\U_\nu e^{51}\U_\eta e^{53})$,
  $d_2(\bar{e_{47}}h_4) =\bar{e_{47}}h_0h_3^2 \ne0$ and
  $d_3(\bar{e_{47}}h_0h_4) =\bar{e_{47}}h_0d_0 \ne0$.
  It follows that $AF(\widetilde{\eta\varepsilon})\ge3$.

  Consider the subcomplex
  $W_2 =S^{47}\U_{2\iota} e^{48}\U_\nu e^{51}\U_\eta e^{53}\subset P^{62}_{47}$. Let
  $W_1 \to{i_5}W_2 \to{j_1}P^{62}_{47}$ be inclusion maps and let
  $W_2 \to{q_5}S^{53}$ be the collapsing map.
  Then since $AF(\widetilde{\eta\varepsilon})\ge3$ in the ASS for $\pi_\ast^S(W_1)$ and
  since $AF(q_5(i_5\widetilde{\eta\varepsilon})) =\eta\varepsilon)=4$, we see
  $3\le AF(i_5\widetilde{\eta\varepsilon})\le4$ in the ASS for $\pi_\ast^S(W_2)$.
  If $AF(i_5\widetilde{\eta\varepsilon}) =3$ in the ASS for $\pi_\ast^S(W_2)$ then from
  \autoref{fig:Cnueta2iota:Ext} we see that
  $i_5\widetilde{\eta\varepsilon}$ must be detected by
  $\bar{e_{53}h_1^2h_3} \in\Ext_\A^{3,65}(W_2)$.
  But it is not difficult to see that
  $(q_5)_\ast(\bar{e_{53}h_1^2h_3}) =h_1^2h_3\ne0$ in $\Ext_\A^{3,12}$,
  so it follows that
  $AF(q_5(i_5\widetilde{\eta\varepsilon}) =q_4\widetilde{\eta\varepsilon} =\eta\varepsilon) =3$,
  and this is a contradiction because
  $AF(\eta\varepsilon) =4$.
  Therefore $AF(i_5\widetilde{\eta\varepsilon}) =4$.
  From \autoref{fig:Cnueta2iota:Ext} we see that
  $i_5\widetilde{\eta\varepsilon}$ must be detected by
  $\bar{e_{53}h_1c_0} \in\Ext_\A^{3,65}(W_2)$.
  From \autoref{fig:P62-47:Ext},
  \autoref{thm:P62-46:ASS:d5}, and \autoref{thm:P62-46:ASS:e49h32},
  $\bar{e_{53}h_1c_0} \ne0$ is not a boundary in the ASS for $\pi_\ast^S(P^{62}_{47})$,
  so $v =j_1(i_5\widetilde{\eta\varepsilon})$ is detected by
  $\bar{e_{53}h_1c_0} \in\Ext_\A^{4,66}(P^{62}_{47})$.
  This proves \autoref{thm:v}.
\end{proof}
Now we prove \autoref{thm:P62-47:pi62}, which
is restated as follows.
\theoremstyle{nonumberplain}
\newtheorem{thm:P62-47:pi62}{Theorem \ref{thm:P62-47:pi62}} 
\begin{thm:P62-47:pi62*}
  Let $S^{62} \to{\widetilde{y}}P^{62}_{47}$ be the composite
  $S^{62} \to{y}P^{62}_{46} \to{q_1}P^{62}_{47}$,
  $S^{62} \to{\widetilde{w}}P^{62}_{47}$ be the composite
  $S^{62} \to{w}P^{62}_{46} \to{q_1}P^{62}_{47}$, and
  $S^{62} \to{\widetilde{x} =\widetilde{y} +\widetilde{w}}P^{62}_{47}$. Then
  \begin{align*}
    \pi^S_{62}(P^{62}_{47})
    &=\Z/2(\widetilde{y})\dsum\Z/8(\widetilde{x})\dsum\Z/2(v-2\widetilde{x}).
  \end{align*}
  And $AF(\widetilde{y})=1$, $AF(\widetilde{x})=2$, and $AF(v)=4$.
\end{thm:P62-47:pi62*}
\begin{proof}[Proof of \autoref{thm:P62-47:pi62}]
  \begin{enumerate}
    \item
      By \autoref{thm:v}, $AF(v) =4$ in the ASS for $\pi_\ast^S(P^{62}_{47})$.
      From \autoref{thm:y} we see that $\bar{e_{47}}h_4 \in\Ext_\A^{1,63}(P^{62}_{46})$
      detects the stable homotopy class $y\in\pi_{62}^S(P^{62}_{46})$.
      From \autoref{fig:P62-47:Ext}, $\bar{e_{47}}h_4$ is not a boundary
      in the ASS for $\pi_\ast^S(P^{62}_{47})$,
      so $\bar{e_{47}}h_4$ detects $\widetilde{y} =q_1y$, and $AF(\widetilde{y}) =1$.
    \item
      $2\widetilde{y} =0$ because
      we have shown in the proof of \autoref{thm:P62-46:pi62} that
      $S^{62}\to{2\iota} S^{62}\to{y}P^{62}_{48}$ can be lifted to
      $S^{46}\subset P^{48}_{46}$.
    \item
      We already have shown in the proof of \autoref{thm:P62-46:pi62} that
      $\bar{e_{55}h_0h_3} \in\Ext_\A^{2,64}(P^{62}_{46})$,
      $\bar{e_{48}h_0h_3^2} \in\Ext_\A^{3,65}(P^{62}_{46})$,
      $\bar{e_{51}P^1h_2} \in\Ext_\A^{5,67}(P^{62}_{46})$
      detect the stable homotopy classes $x,2x,4x\in\pi_{62}^S(P^{62}_{46})$
      respectively.
      From \autoref{fig:P62-47:Ext}, \autoref{thm:P62-46:ASS:d5}, and
      \autoref{thm:P62-46:ASS:e49h32},
      the classes
      $\bar{e_{55}h_0h_3},\bar{e_{48}h_0h_3^2},\bar{e_{51}P^1h_2}$ are not boundaries
      in the ASS for $\pi_\ast^S(P^{62}_{47})$,
      so $\bar{e_{55}h_0h_3}$ detects $\widetilde{x} =q_1x$ and $AF(\widetilde{x}) =2$,
      $\bar{e_{48}h_0h_3^2}$ detects $2\widetilde{x} =q_1(2x)$, and
      $\bar{e_{51}P^1h_2}$ detects $4\widetilde{x} =q_1(4x)$.

    \item We have shown in \autoref{thm:v} that
      $\bar{e_{53}h_1c_0} \in\Ext_\A^{4,66}(P^{62}_{47})$ detects $v \in\pi_{62}^S(P^{62}_{47})$.
      In \eqref{eq:P62-47:Ext:extension}~(3) we have seen $\bar{e_{53}h_1c_0}h_0 =\bar{e_{51}P^1h_2}$.
      We have seen that $\bar{e_{51}P^1h_2}$ is not a boundary in
      the ASS for $\pi_\ast^S(P^{62}_{47})$.
      Therefore $2v$ is detected by $\bar{e_{55}P^1h_2} \in\Ext_\A^{5,67}(P^{62}_{47})$.
    \item We have $d_5(e_{63}h_0)=\bar{e_{48}h_0^2d_0}$ by \autoref{thm:P62-46:ASS:d5}. So from
      \autoref{fig:P62-47:Ext}, $2v =4\widetilde{x}$ and $4v=8\widetilde{x}=0$.
  \end{enumerate}
  From (3), (4), and (5) it follows that
  $AF(v-2\widetilde{x}) =3$,
  and that $2(v-2\widetilde{x}) =2v-4\widetilde{x} =0$.
  So the order of $v-2\widetilde{x}$ is $2$.
  From $AF(\widetilde{y}) =1$,
  $AF(\widetilde{x}) =2$, $AF(2\widetilde{x}) =3$, $AF(4\widetilde{x}) =5$, and
  $AF(v-2\widetilde{x}) =3$,
  it is not difficult to see that $\widetilde{y},\widetilde{x},v-2\widetilde{x}$
  are linearly independent.
  This completes the proof of \autoref{thm:P62-47:pi62}.
\end{proof}

\newpage
\section{Proof of \autoref{thm:dr-h6}, \autoref{thm:d2-h1g4}, and \autoref{thm:t5ht6:AF-ge-5}}\label{se:not-hit}
In this section we prove the technical results
\autoref{thm:dr-h6}, \autoref{thm:d2-h1g4}, and \autoref{thm:t5ht6:AF-ge-5}.
These are the remaining theorems in \autoref{se:intro} and \autoref{se:boundary}
that need to be proved to complete the proof of \autoref{thm:main}.
\theoremstyle{nonumberplain}
\newtheorem{thm:dr-h6}{\autoref{thm:dr-h6}} 
\begin{thm:dr-h6*}
  In the ASS for spheres, $d_r(h_6^3)=0\in\Ext_\A^{3+r,191+r}$ for $2\le r\le 3$.
\end{thm:dr-h6*}
\newtheorem{thm:d2-h1g4}{\autoref{thm:d2-h1g4}}
\begin{thm:d2-h1g4*}
  In the ASS for spheres, $d_2(h_1g_4)\ne h_0^3g_4$.
\end{thm:d2-h1g4*}
\newtheorem{thm:t5ht6:AF-ge-5}{\autoref{thm:t5ht6:AF-ge-5}} 
\begin{thm:t5ht6:AF-ge-5*}
  The composite
  $S^{188} \to{\theta_6}S^{62} \to{\widehat{\theta}_5} P$ has $AF=5$ in the Adams spectral sequence for $P =P^\infty$ ($AF=$ Adams filtration).
\end{thm:t5ht6:AF-ge-5*}
The order of our proving these propositions is
\ref{thm:t5ht6:AF-ge-5}, \ref{thm:dr-h6}, and \ref{thm:d2-h1g4}.
Before we start proving \autoref{thm:t5ht6:AF-ge-5},
we first recall the following differential in the ASS for $\pi_\ast^S$
mentioned in \eqref{eq:intro:d2-hi}.
\begin{align}
  \label{eq:d2-hi}
  d_2(h_i) &=h_{i-1}^2h_0\ne0
  \quad\mbox{in the ASS for $\pi_\ast^S$.}
\end{align}
We remark that it is not difficult to see the following differential \eqref{eq:d2-hih}
in the ASS for $\pi_\ast^S(P)$
from \eqref{eq:d2-hi}, \eqref{eq:t-bar}, \eqref{eq:P:Ext:t_ast},
and the structure
of $\Ext_\A^{s,t}(P)$ for $0\le s\le3$ as given in \autoref{thm:P:Ext:known}.
\begin{align}
  \label{eq:d2-hih}
  d_2(\widehat{h}_i) &=\widehat{h}_{i-1}h_0h_{i-1}\ne0
  \quad\mbox{in the ASS for $\pi_\ast^S(P)$.}
\end{align}

To prove \autoref{thm:t5ht6:AF-ge-5} we begin by proving
the following result in \autoref{thm:qt5t6:AF-ge-5}.
Recall the complex $\widetilde{P}^{62}$ and
the map $P^{62} \to{\widetilde{q}}\widetilde{P}^{62}$
in \eqref{defn:P62t} which has the property that
$\widetilde{q}_\ast(\widehat{h}_5 =\bar{e_{31}}) =0$ for
$\widehat{h}_5 \in\Ext_\A^{0,31}(P^{62}) =\Z/2$.
\begin{lemm}\label{thm:qt5t6:AF-ge-5}
  The composite
  $S^{188} \to{\theta_6}S^{62} \to{\widehat{\theta}_5} P^{62}\to{\widetilde{q}} \widetilde{P}^{62}$
  has $AF\ge5$
  in the ASS for $\pi_\ast^S(\widetilde{P}^{62})$.
\end{lemm}
\begin{proof}
  Recall \eqref{eq:theta_n-hat} that $\widehat{\theta}_5$ is detected by $\widehat{h}_5h_5$,
  so $AF(\widehat{\theta}_5) =1$.
  Since $\widetilde{q}_\ast(\widehat{h}_5)=0$,
  $AF(\widetilde{q}\widehat{\theta}_5)\ge2$.
  Consider the following two cases.
  \begin{enumerate}

    \item 
      If $AF(\widetilde{q}\widehat{\theta}_5)=2$, then
      $\widetilde{q}\widehat{\theta}_5 \in\pi_{62}^S(\widetilde{P}^{62})$
      is detected by a non-zero class
      $\alpha \in\Ext_\A^{2,64}(\widetilde{P}^{62}) =\Z/2(\bar{\widetilde{e}_{31}h_3}h_0)
      \dsum\Z/2(\bar{e_{47}}h_0h_4)$ (\eqref{eq:P62t:Ext} (2)). So
      \begin{align*}
        \alpha
        &=\varepsilon_1 \bar{\widetilde{e}_{31}h_5}h_0
        +\varepsilon_2 \bar{e_{47}}h_0h_4,
        \quad\mbox{$\varepsilon_1,\varepsilon_2=0$ or $1$ and $(\varepsilon_1,\varepsilon_2) \ne(0,0)$.}
      \end{align*}
      It is not difficult to see that
      $\bar{\widetilde{e}_{31}h_5} =\{\widetilde{e}_{31}\lambda_{31}\} \in\Ext_\A^{1,63}(\widetilde{P}^{62})$ and
      $\bar{e_{47}} =\{e_{47}\} \in\Ext_\A^{0,47}(\widetilde{P}^{62})$. Thus
      \begin{align*}
        (\bar{\widetilde{e}_{31}h_5}h_0)h_6^2
          &=\bar{\widetilde{e}_{31}h_5}h_6^2h_0
          =\{\widetilde{e}_{31}\lambda_{31}\lambda_{63}^2\lambda_0\}=0,\\
        (\bar{e_{47}}h_0h_4)h_6^2
          &=\bar{e_{47}}h_0h_4h_6^2 =0
          \quad\mbox{(since $h_4h_6^2 =0$).}
      \end{align*}
      So $\alpha h_6^2=0$. Since $\theta_6$ is detected by $h_6^2$
      in the ASS for $\pi_\ast^S$, we see
      $AF(\widetilde{q}\widehat{\theta}_5\theta_6)\ge5$.

    \item
      If $AF(\widetilde{q}\widehat{\theta}_5)\ge3$
      then clearly $AF(\widetilde{q}\widehat{\theta}_5\theta_6)\ge3+2 =5$.

  \end{enumerate}
  This proves the \autoref{thm:qt5t6:AF-ge-5}.
\end{proof}

The scheme to prove \autoref{thm:t5ht6:AF-ge-5} is as follows.
Since $\widehat{\theta}_5$, $\theta_6$
are detected by $\widehat{h}_5h_5\in\Ext_\A^{1,63}(P^{62})$ and
$h_6^2\in\Ext_\A^{2,128}$
in the ASS for $\pi_\ast^S(P^{62})$ and for $\pi_\ast^S$
respectively, and since $h_5h_6=0$ in $\Ext_\A^{\ast,\ast}$,
it follows that $AF(\widehat{\theta}_5\theta_6)\ge4$.
If $AF(\widehat{\theta}_5\theta_6)=4$ then
$\widehat{\theta}_5\theta_6$ must be detected by a non-zero class
$\tau\in\Ext_\A^{4,192}(P^{62})$. 
Consider the induced homomorphism
$\Ext_\A^{\ast,\ast}(P^{62}) \to{\widetilde{q}_\ast}\Ext_\A^{\ast,\ast}(\widetilde{P}^{62})$.
By \autoref{thm:qt5t6:AF-ge-5},
$AF(\widetilde{q}\widehat{\theta}_5\theta_6)=5$.
Thus
either $q_\ast(\tau)=0$ in $\Ext_\A^{4,192}(\widetilde{P}^{62})$
or $q_\ast(\tau)$ is a boundary in the ASS for $\pi_\ast^S(\widetilde{P}^{62})$.
To show that neither of these two cases is possible
we need the following \autoref{thm:P62t:Ext:2-191:d2}.
To state it, we recall the following.
\begin{enumerate}
  \item[\eqref{eq:P62t:Ext}~(4)\quad]
    $\Ext_\A^{2,191}(\widetilde{P}^{62})
      = \Z/2(\bar{\widetilde{e}_{31}h_5}h_7) \dsum\Z/2(\bar{e_{47}}h_4h_7)$.
  \item[\eqref{eq:P62t:Ext}~(6)\quad]
    $\Ext_\A^{4,192}(\widetilde{P}^{62})
      = Z/2(\widetilde{q}_\ast(\widehat{h}_4h_0c_4))
      \dsum \Z/2(\widetilde{q}_\ast(\widehat{D_3}h_7))
      \dsum \Z/2(\widetilde{q}_\ast(\widehat{g}_4h_0))
      \dsum \Z/2(\widetilde{q}_\ast(\bar{e_{62}D_3(1)}))
      \dsum \Z/2(\bar{\widetilde{e}_{31}h_4^2}h_0h_7)$.
\end{enumerate}
\begin{lemm}
  \label{thm:P62t:Ext:2-191:d2}
  In the ASS for $\pi_\ast^S(\widetilde{P}^{62})$,
  \begin{enumerate}
    \item
      $d_2(\bar{\widetilde{e}_{31}h_5}h_7)
      =\widetilde{q}_\ast(\widehat{D_3}h_7)
      +\bar{\widetilde{e}_{31}h_0h_4^2}h_7$,
    \item
      $d_2(\bar{e_{47}}h_4h_7)
      =\widetilde{q}_\ast(\widehat{D_3}h_7)
      +\bar{\widetilde{e}_{31}h_0h_4^2}h_7$.
  \end{enumerate}
\end{lemm}
The proof of this \autoref{thm:P62t:Ext:2-191:d2} is rather lengthy.
We will postpond this proof to the end of this section.
Now we prove \autoref{thm:t5ht6:AF-ge-5}.
\begin{proof}[Proof of \autoref{thm:t5ht6:AF-ge-5}]
  We want to show the composite
  $S^{188} \to{\theta_6}S^{62} \to{\widehat{\theta}_5}P^{62}$ has $AF\ge5$.
  We already know that $AF(\widehat{\theta}_5\theta_6)\ge4$ (since $\widehat{h}_5h_5h_6^2 =0$).
  Recall
  \begin{enumerate}
    \item[\eqref{eq:P62:Ext}~(i)~(5)\quad]
      $\Ext_\A^{4,192}(P^{62})
        =\Z/2(\widehat{h}_4h_0c_4)
        \dsum \Z/2(\widehat{D_3}h_7)
        \dsum \Z/2(\widehat{g}_4h_0)
        \dsum \Z/2(\bar{e_{62}D_3(1)})$.
  \end{enumerate}
  Suppose $AF(\widehat{\theta}_5\theta_6) =4$. Then
  \begin{note}
    \label{eq:tau}
    $\widehat{\theta}_5\theta_6$ is detected by
    $\tau
      =\varepsilon_1\widehat{h}_4h_0c_4
      +\varepsilon_2\widehat{D_3}h_7
      +\varepsilon_3\widehat{g}_4h_0
      +\varepsilon_4\widehat{h}_7D_3$, where
      $\varepsilon_i=0$ or $1$ for $1\le i\le 4$, and
    $(\varepsilon_1,\varepsilon_2,\varepsilon_3,\varepsilon_4) \ne(0,0,0,0)$.
  \end{note}
  \noindent
  We are going to show that this would lead to a contradiction.
  Recall
  \begin{enumerate}
    \item[\eqref{eq:P62t:Ext}~(6)\quad]
      $\Ext_\A^{4,192}(\widetilde{P}^{62})
        = Z/2(\widetilde{q}_\ast(\widehat{h}_4h_0c_4))
        \dsum \Z/2(\widetilde{q}_\ast(\widehat{D_3}h_7))
        \dsum \Z/2(\widetilde{q}_\ast(\widehat{g}_4h_0))
        \dsum \Z/2(\widetilde{q}_\ast(\widehat{h}_7D_3))  
        \dsum \Z/2(\bar{\widetilde{e}_{31}h_4^2}h_0h_7)$
  \end{enumerate}
  Let $P^{62} \to{\widetilde{q}}\widetilde{P}^{62}$ be as in \autoref{thm:qt5t6:AF-ge-5}.
  Then from \eqref{eq:tau} we have:
  \begin{note}
    \label{eq:qtau}
    $\widetilde{q}_\ast(\tau)
    =\varepsilon_1\widetilde{q}_\ast(\widehat{h}_4h_0c_4)
    +\varepsilon_2\widetilde{q}_\ast(\widehat{D_3}h_7)
    +\varepsilon_3\widetilde{q}_\ast(\widehat{g}_4h_0)
    +\varepsilon_4\widetilde{q}_\ast(\widehat{h}_7D_3)$,
    where
    $\varepsilon_i=0$ or $1$ for $1\le i\le 4$, and
    $(\varepsilon_1,\varepsilon_2,\varepsilon_3,\varepsilon_4) \ne(0,0,0,0)$.
    $\widetilde{q}_\ast(\tau)$ is non-zero in $\Ext_\A^{4,192}(\widetilde{P}^{62})$ by
    \eqref{eq:P62t:Ext}~(6) and the fact that
    $(\varepsilon_1,\varepsilon_2,\varepsilon_3,\varepsilon_4) \ne(0,0,0,0)$.
  \end{note}
  \noindent
  By \autoref{thm:qt5t6:AF-ge-5},
  $\widetilde{q}_\ast(\tau)$ has to be a boundary
  (if $\widetilde{q}_\ast(\tau)$ is not a boundary then $AF(\widetilde{q}\widehat{\theta}_5\theta_6) =4$).
  So $d_r(\beta) =\widetilde{q}_\ast(\tau)$ for some $\beta \in\Ext_\A^{r-4,185+r}(\widetilde{P}^{62})$
  with $2\le r\le4$.
  Recall
  \begin{enumerate}
    \item[\eqref{eq:P62t:Ext}~(8)\quad]
      $\Ext_\A^{0,189}(\widetilde{P}^{62}) =0$.
    \item[\eqref{eq:P62t:Ext}~(9)\quad]
      $\Ext_\A^{1,190}(\widetilde{P}^{62}) =0$.
    \item[\eqref{eq:P62t:Ext}~(10)\quad]
      $\Ext_\A^{2,191}(\widetilde{P}^{62})
      =\Z/2(\bar{\widetilde{e}_{31}h_5}h_7)
      \dsum \Z/2(\bar{e_{47}}h_4h_7)$.
  \end{enumerate}
  \begin{align*}
  \end{align*}
  Thus the only possibility is
  $d_2(\delta_1\bar{\widetilde{e}_{31}h_5}h_7 +\delta_2\bar{e_{47}}h_4h_7)
  =\widetilde{q}_\ast(\tau)$ where $(\delta_1,\delta_2)\ne(0,0)$.
  By \autoref{thm:P62t:Ext:2-191:d2},
  \begin{align*}
    d_2(\delta_1\bar{\widetilde{e}_{31}h_5}h_7 +\delta_2\bar{e_{47}}h_4h_7)
    &=\delta_1(\widetilde{q}_\ast(\widehat{D_3}h_7) +\bar{\widetilde{e}_{31}h_0h_4^2}h_7)
      +\delta_2(\widetilde{q}_\ast(\widehat{D_3}h_7) +\bar{\widetilde{e}_{31}h_0h_4^2}h_7)\\
    &=0\quad\mbox{or}\quad
      \widetilde{q}_\ast(\widehat{D_3}h_7) +\bar{\widetilde{e}_{31}h_0h_4^2}h_7.
  \end{align*}
  By \eqref{eq:P62t:Ext}~(6) and \eqref{eq:qtau} above,
  we see $\widetilde{q}_\ast(\tau) =
  \varepsilon_1\widetilde{q}_\ast(\widehat{h}_4h_0c_4)
  +\varepsilon_2\widetilde{q}_\ast(\widehat{D_3}h_7)
  +\varepsilon_3\widetilde{q}_\ast(\widehat{g}_4h_0)
  +\varepsilon_4\widetilde{q}_\ast(\widehat{h}_7D_3)
  \ne\widetilde{q}_\ast(\widehat{D_3}h_7) +\bar{\widetilde{e}_{31}h_0h_4^2}h_7$.
  This contradiction proves that $AF(\widehat{\theta}_5\theta_6)\ge5$.
  This concludes the proof of \autoref{thm:t5ht6:AF-ge-5}.
\end{proof}

We proceed to prove \autoref{thm:dr-h6}.
We want to prove for this proposition that in the ASS for spheres,
$d_r(h_6^3) =0$ for $r=2$ and $3$. That
$d_2(h_6^3) =0$ is quite easy as $d_2(h_6^3 =h_6h_6^2) =h_0h_5^2h_6^2 =0$.
We have $d_3(h_6^3) \in E_3^{6,194}$ which is a subquotient of $E_2^{6,194} =\Ext_\A^{6,194}$.
Recall
\begin{enumerate}
  \item[\eqref{eq:Ext}~(8)\quad]
    $\Ext_\A^{6,194}
    =\Z/2(h_0^2g_4)\dsum\Z/2(h_5^2D_3(1))$.
\end{enumerate}
So $d_3(h_6^3) =\varepsilon_1h_0^2g_4 +\varepsilon_2h_5^2D_3(1)$ with
$\varepsilon_i =0$ or $1$, $1\le i\le 2$. We have to show $\varepsilon_1 =\varepsilon_2 =0$.
For this purpose we first prove the following case that says that
$\varepsilon_1 =0$ and $\varepsilon_2 =1$ can not happen.
\begin{lemm}
  \label{thm:S0:d3h63-ne-h52D31}
  In the ASS for $\pi_\ast^S$,
  $d_3(h_6^3) \ne h_5^2D_3(1)$.
\end{lemm}
\begin{proof}
  Recall that the composite of stable maps
  $S^{63} \to{\theta_5}S^1 \to{2\iota}S^1$
  is zero. So for the following cofibration
  \begin{align}
    \label{eq:d2h63:cofibration}
    S^0 &\to{i_1}S^0\U_{2\iota}e^1  \to{q_1}S^1,
  \end{align}
  there is a coextension
  \begin{align}
    \label{eq:t5bar}
    S^{63} \to{\bar{\theta}_5}S^0\U_{2\iota}S^1.
  \end{align}
  We have
  $h_6^3 \in\Ext_\A^{3,192}$, and
  $(i_1)_\ast(h_6^3) \in\Ext_\A^{3,192}(S^0\U_{2\iota}e^1)$ for
  the $S^0 \to{i_1}S^0\U_{2\iota}e^1$ in \eqref{eq:d2h63:cofibration}.
  We claim the following.
  \begin{note}
    \label{eq:d2h63:claim}
    \begin{enumerate}
      \item $(i_1)_\ast(h_6^3) \ne0$ in $\Ext_\A^{3,192}(S^0\U_{2\iota}e^1)$.
      \item The composite
        $S^{189} \to{\theta_6}S^{63} \to{\bar{\theta}_5}S^0\U_{2\iota}e^1$
        is detected by $(i_1)_\ast(h_6^3)$ in the ASS for $\pi_\ast^S(S^0\U_{2\iota}e^1)$.
      \item $\Ext_\A^{4,193}(S^0\U_{2\iota}e^1) =0$.
    \end{enumerate}
  \end{note}
  \noindent
  Assuming \eqref{eq:d2h63:claim}, we show $d_3(h_6^3) \ne h_5^2D_3(1)$
  in the ASS for $\pi_\ast^S$ as follows.
  If $d_3(h_6^3) =h_5^2D_3(1)$ in the ASS for $\pi_\ast^S$, then since
  $(i_1)_\ast(h_6^3)\ne0$ is an infinite cycle
  in the ASS for $\pi_\ast^S(S^0\U_{2\iota}e^1)$ (by \eqref{eq:d2h63:claim}),
  we must have either of the following two cases.
  \begin{note}
    \label{eq:d2h63:cases}
    \begin{enumerate}
      \item $(i_1)_\ast(h_5^2D_3(1)) =0$ in $\Ext_\A^{6,194}(S^0\U_{2\iota}e^1)$, or
      \item $(i_1)_\ast(h_5^2D_3(1)) \ne0$ in $\Ext_\A^{6,194}(S^0\U_{2\iota}e^1)$ and
        $d_2(\alpha) =(i_1)_\ast(h_5^2D_3(1))$ for some $\alpha \in\Ext_\A^{4,193}(S^0\U_{2\iota}e^1)$
        in the ASS for $\pi_\ast^S(S^0\U_{2\iota}e^1)$.
    \end{enumerate}
  \end{note}
  \noindent
  By \eqref{eq:d2h63:claim} (3), \eqref{eq:d2h63:cases} (2) is impossible.
  To show that \eqref{eq:d2h63:cases} (1) is impossible,
  recall from \eqref{eq:M2:Ext:long-exact-sequence} that
  since \eqref{eq:d2h63:cofibration} induces a short exact sequence
  $0-> H_\ast(S^0) \to{(i_1)_\ast} H_\ast(S^0\U_{2\iota}e^1) \to{(q_1)_\ast} H_\ast(S^1)->0$,
  there is a long exact sequence
  \begin{align}
    \label{eq:d2h63:M2:Ext:long-exact-sequence}
    \cdots
    -> \Ext_\A^{s-1,t-1}
    \to{(2\iota)_\ast} \Ext_\A^{s,t}
    \to{(i_1)_\ast} \Ext_\A^{s,t}(S^0\U_{2\iota}e^1)
    \to{(q_1)_\ast} \Ext_\A^{s,t-1}
    \to{(2\iota)_\ast} \Ext_\A^{s+1,t}
    ->\cdots
  \end{align}
  So $(i_1)_\ast(h_5^2D_3(1)) =0$ in $\Ext_\A^{6,194}(S^0\U_{2\iota}e^1)$ for
  $(s,t)=(6,194)$ implies that
  there is $\beta\in\Ext_\A^{5,193}$ such that
  $(2\iota)_\ast(\beta) =h_0\beta =h_5^2D_3(1)$ in $\Ext_\A^{6,194}$.
  From \eqref{eq:Ext}~(2) and (8), there is no such $\beta \in\Ext_\A^{5,193}$.
  So \eqref{eq:d2h63:cases} (1) is impossible.
  This shows that $d_3(h_6^3) \ne h_5^2D_3(1)$.

  To complete the proof of \autoref{thm:S0:d3h63-ne-h52D31} it remains to show \eqref{eq:d2h63:claim}.
  \eqref{eq:d2h63:claim} (1), (3) can also be proved by using
  the long exact sequence in \eqref{eq:d2h63:M2:Ext:long-exact-sequence}.
  To show $(i_1)_\ast(h_6^3)\ne0$ as claimed in \eqref{eq:d2h63:claim} (1),
  note that $h_6^3 \in\Ext_\A^{3,192}$ for $(s,t) =(3,192)$, and
  by \eqref{eq:Ext}~(3), $\Ext_\A^{2,191} =0$. So, by \eqref{eq:d2h63:M2:Ext:long-exact-sequence},
  $(i_1)_\ast(h_6^3) \ne0$ in $\Ext_\A^{3,192}(S^0\U_{2\iota}e^1)$.
  This proves \eqref{eq:d2h63:claim} (1).
  To show $\Ext_\A^{4,193}(S^0\U_{2\iota}e^1) =0$ for $(s,t) =(4,193)$
  as claimed in \eqref{eq:d2h63:claim} (3), recall the following.
  \begin{enumerate}
    \item[\eqref{eq:Ext}~(1)\quad] $\Ext_\A^{4,192} =\Z/2(g_4)$.
    \item[\eqref{eq:Ext}~(2)\quad] $\Ext_\A^{5,193} =\Z/2(h_0g_4)$.
    \item[\eqref{eq:Ext}~(4)\quad] $\Ext_\A^{s,s+189} =\Z/2(h_0^{s-3}h_6^3)$, $3\le s\le4$.
  \end{enumerate}
  So
  $\Ext_\A^{4,192} =\Z/2(g_4)$ with
  $(2\iota)_\ast(g_4) =h_0g_4\ne0$ in $\Ext_\A^{5,193}$, and
  $\Ext_\A^{4,193} =\Z/2(h_0h_6^3)$ with $(2\iota)_\ast(h_6^3) =h_0h_6^3$.
  From \eqref{eq:d2h63:M2:Ext:long-exact-sequence} we thus see
  $\Ext_\A^{4,193}(S^0\U_{2\iota}e^1) =0$. This proves \eqref{eq:d2h63:claim} (3).
  Finally we show that $S^{188} \to{\theta_6}S^{62} \to{\bar{\theta}_5}S^0\U_{2\iota}e^1$
  is detected by $(i_1)_\ast(h_6^3)$
  as claimed in \eqref{eq:d2h63:claim} (2).
  $S^{62} \to{\bar{\theta}_5}S^0\U_{2\iota}e^1$ is detected by
  $(i_1)_\ast(h_6)$ by \cite{adams_non-existence_1960}.
  Then since $(i_1)_\ast(h_6)h_6^2 =(i_1)_\ast(h_6^3)\ne0$ as in \eqref{eq:d2h63:claim} and since
  $\Ext_\A^{s,190+s} =0$ for $0\le s\le1$ (\eqref{eq:Ext}~(6))
  we see that $(i_1)_\ast(h_6^3)$ is not a boundary in the spectral sequence and that 
  $\bar{\theta}_5\theta_6$ is detected by $(i_1)_\ast(h_6^3)$
  in the ASS for $\pi_\ast(S^0\U_{2\iota}e^1)$.
  This proves \autoref{thm:S0:d3h63-ne-h52D31}.
\end{proof}

Now we can prove \autoref{thm:dr-h6}.
\begin{proof}[Proof of \autoref{thm:dr-h6}]
  Recall the transfer map $P \to{t}S^0$ and
  the induced homomorphism $\Ext_\A^{s,t}(P) \to{t_\ast}\Ext_\A^{s+1,t+1}$
  in \eqref{eq:t_ast}.
  There is a class $\widehat{h}_5h_5h_7 \in\Ext_\A^{2,191}(P^{62})$ such that
  $t_\ast(\widehat{h}_5h_5h_7) =h_5^2h_7 =h_6^3$.
  Since $d_2(h_6^3) =0$ in the ASS for $\pi_\ast^S$, and
  $d_2(\widehat{h}_5h_5h_7)=\widehat{h}_5h_4^2h_0h_7=0$
  in the ASS for $\pi_\ast^S(P)$, we can
  consider the following diagram.
  \begin{align*}
    \xymatrix{
    \Ext_\A^{2,191}(P^{62})\ar[r]^{t_\ast}\ar[d]^{d_3}&
      \Ext_\A^{3,192}\ar[d]^{d_3}\\
    \Ext_\A^{5,193}(P^{62})\ar[r]^{t_\ast}&
      \Ext_\A^{6,194}
    }
  \end{align*}
  Suppose $d_3(h_6^3)=\alpha\ne0$ in $\Ext_\A^{6,194}$ in the ASS for $\pi_\ast^S$.
  We show this would lead to a contradiction.
  $d_3(h_6^3) =\alpha$ in the ASS for $\pi_\ast^S$ implies
  $d_3(\widehat{h}_5h_5h_7) =\alpha'\ne0$ in $\Ext_\A^{5,193}(P^{62})$
  with $t_\ast(\alpha')=\alpha$ in the ASS for $\pi_\ast^S(P^{62})$.
  We recall the following. 
  \begin{enumerate}
    \item[\eqref{eq:P62:Ext}~(i)~(9)\quad] 
      $\Ext_\A^{5,193}(P^{62})
        = \Z/2(\widehat{h}_5h_5D_3(1))
        \dsum \Z/2(\widehat{h}_4h_0^2c_4)
        \dsum \Z/2(\widehat{h}_4h_1f_3)
        \dsum \Z/2(\gamma_{61}h_7)
        \dsum \Z/2(\widehat{g}_4h_0^2)$.
  \end{enumerate}
  Therefore
  \begin{align}
    \label{eq:alphap}
    d_3(\widehat{h}_5h_5h_7)
      =\alpha'
      = \varepsilon_1\widehat{h}_5h_5D_3(1)
      +\varepsilon_2\widehat{h}_4h_0^2c_4
      +\varepsilon_3\widehat{h}_4h_1f_3
      +\varepsilon_4\gamma_{61}h_7
      &+\varepsilon_5\widehat{g}_4h_0^2,
      \quad\mbox{$\varepsilon_i=0$ or $1$, $1\le i\le 5$,}\\
      &\mbox{with $(\varepsilon_1,\varepsilon_2,\varepsilon_3,\varepsilon_4,\varepsilon_5)\ne(0,0,0,0,0)$.}\notag
  \end{align}
  Recall the complex $\widetilde{P}^{62}$ and the map
  $P^{62} \to{\widetilde{q}}\widetilde{P}^{62}$.
  Consider the following diagram.
  \begin{align*}
    \xymatrix{
    \Ext_\A^{2,191}(P^{62})\ar[r]^{\widetilde{q}_\ast}\ar[d]^{d_3}&
      \Ext_\A^{2,191}(\widetilde{P}^{62})\ar[d]^{d_3}\\
    \Ext_\A^{5,193}(P^{62})\ar[r]^{\widetilde{q}_\ast}&
      \Ext_\A^{5,193}(\widetilde{P}^{62})
    }
  \end{align*}
  We recall the following.
  \begin{enumerate}
    \item[\eqref{eq:P62t:Ext}~(7)~(a)\quad]
      $\widetilde{q}_\ast(\widehat{h}_5h_5D_3(1)) =0$ in $\Ext_\A^{5,193}(\widetilde{P}^{62})$.
    \item[\eqref{eq:P62t:Ext}~(7)~(b)\quad]
      $\widetilde{q}_\ast(\widehat{h}_4h_0^2c_4),
      \widetilde{q}_\ast(\widehat{h}_4h_1f_3),
      \widetilde{q}_\ast(\gamma_{61}h_7),
      \widetilde{q}_\ast(\widehat{g}_4h_0^2),
      \bar{\widetilde{e}_{31}h_4^2}h_0^2h_7$ are linearly independent in $\Ext_\A^{5,193}(\widetilde{P}^{62})$.
  \end{enumerate}
  Since $\widetilde{q}_\ast(\widehat{h}_5h_5h_7) =0$,
  we have $\widetilde{q}_\ast(\alpha') =\widetilde{q}_\ast(d_3(\widehat{h}_5h_5h_7)) =d_3(\widetilde{q}_\ast(\widehat{h}_5h_5h_7)) =0$
  in the ASS for $\pi_\ast^S(\widetilde{P}^{62})$.
  There are two cases for this to be possible.
  \begin{enumerate}
    \item If $\widetilde{q}_\ast(\alpha') =0$ in $\Ext_\A^{5,193}(\widetilde{P}^{62})$,
      from \eqref{eq:P62:Ext}~(i)~(9), \eqref{eq:alphap} and \eqref{eq:P62t:Ext}~(7)~(a), (b)
      we see that $\alpha' =\widehat{h}_5h_5D_3(1)$,
      that is, $\varepsilon_1 =1$, $\varepsilon_i=0$ for $2\le i\le5$.
      Then $d_3(h_6^3) =t_\ast(\alpha' =\widehat{h}_5h_5D_3(1)) =h_5^2D_3(1)$
      in the ASS for $\pi_\ast^S$.
      By \autoref{thm:S0:d3h63-ne-h52D31} this is impossible.
    \item If $\widetilde{q}_\ast(\alpha')\ne0$ in $\Ext_\A^{5,193}(\widetilde{P}^{62})$, then
      $\widetilde{q}_\ast(\alpha')$ must be a boundary
      in the ASS for $\pi_\ast^S(\widetilde{P}^{62})$.
      Let $\alpha'' =\widetilde{q}_\ast(\alpha')$,
      then there is
      $\beta \in\Ext_\A^{3,192}(\widetilde{P}^{62})$ with
      $d_2(\beta)=\alpha''$.
      By
      \begin{enumerate}
        \item[\eqref{eq:P62t:Ext}~(11)\quad]
          $\Ext_\A^{3,192}(\widetilde{P}^{62})
          =\Z/2(\bar{\widetilde{e}_{31}h_5}h_0h_7)
          \dsum \Z/2(\bar{e_{47}}h_4h_7h_0)
          \dsum\Z/2(\widetilde{q}_\ast(\bar{e_{62}h_1}h_6^2))$
      \end{enumerate}
      we have $d_2(\widetilde{q}_\ast(\bar{e_{62}h_1}h_6^2)) =0$
      ($\bar{e_{62}h_1}h_6^2 \in\Ext_\A^{3,192}(P^{62})$ is a pullback of
      $\widehat{h}_6h_0h_6^2 \in\Ext_\A^{3,192}(P)$, and
      $d_2(\widehat{h}_6) =\widehat{h}_5h_0h_5$ and $h_5h_6=0$ in the ASS for $\pi_\ast^S(P)$).
      And by \autoref{thm:P62t:Ext:2-191:d2},
      \begin{align*}
        d_2(\bar{\widetilde{e}_{31}h_5}h_7)
        &=\widetilde{q}_\ast(\widehat{D_3}h_7)
        +\bar{\widetilde{e}_{31}h_0h_4^2}h_7,\\
        d_2(\bar{e_{47}}h_4h_7)
        &=\widetilde{q}_\ast(\widehat{D_3}h_7)
        +\bar{\widetilde{e}_{31}h_0h_4^2}h_7.
      \end{align*}
      Therefore
      \begin{align*}
        d_2(\bar{\widetilde{e}_{31}h_5}h_7h_0)
        &=\widetilde{q}_\ast(\widehat{D_3}h_7)h_0
        +\bar{\widetilde{e}_{31}h_0h_4^2}h_0h_7
        = \bar{\widetilde{e}_{31}h_0h_4^2}h_0h_7,\\
        d_2(\bar{e_{47}}h_4h_7h_0)
        &=\widetilde{q}_\ast(\widehat{D_3}h_7)h_0
        +\bar{\widetilde{e}_{31}h_0h_4^2}h_0h_7
        = \bar{\widetilde{e}_{31}h_0h_4^2}h_0h_7.
      \end{align*}
      It follows that the only non-trivial boundary in the ASS for $\pi_\ast^S(\widetilde{P}^{62})$
      in $\Ext_\A^{5,193}(\widetilde{P}^{62})$ is $\bar{\widetilde{e}_{31}h_0h_4^2}h_0h_7$.
      Since $\widetilde{q}_\ast(\alpha')
      =\varepsilon_2\widetilde{q}_\ast(\widehat{h}_4h_0^2c_4)
      +\varepsilon_3\widetilde{q}_\ast(\widehat{h}_4h_1f_3)
      +\varepsilon_4\widetilde{q}_\ast(\gamma_{61}h_7)
      +\varepsilon_5\widetilde{q}_\ast(\widehat{g}_4h_0^2)
      \ne\bar{\widetilde{e}_{31}h_0h_4^2}h_0h_7$
      for any $\alpha' \in\Ext_\A^{5,193}(P^{62})$ (by \eqref{eq:P62t:Ext}~(7)~(b)),
      it follows that $\widetilde{q}_\ast(\alpha')$ is not a boundary
      in the ASS for $\pi_\ast^S(\widetilde{P}^{62})$.
  \end{enumerate}
  This completes the proof of \autoref{thm:dr-h6}.
\end{proof}

Finally we prove \autoref{thm:d2-h1g4}.
\begin{proof}[Proof of \autoref{thm:d2-h1g4}]
  In the ASS for $\pi_\ast^S$ we have $d_2(h_1)=0$,
  so $d_2(h_1g_4) =h_1d_2(g_4) \in\Ext_\A^{7,195}$.
  There are two cases to consider.
  \begin{enumerate}
    \item If $d_2(g_4) =0$ in $\Ext_\A^{6,193}$ then $d_2(h_1g_4) =0$
      and the proof is done.
    \item If $d_2(g_4)\ne0$ in $\Ext_\A^{6,193}$,
      then since $\Ext_\A^{6,193} =\Z/2(h_5V_0)$ (by \eqref{eq:Ext}~(7)) we see
      $d_2(g_4) =h_5V_0$.
      And from \eqref{eq:Ext}~(9),
      $h_1h_5V_0\ne h_0^3g_4$.
      So $d_2(h_1g_4) =h_1d_2(g_4) =h_1h_5V_0 \ne h_0^3g_4$ in $\Ext_\A^{7,195}$.
  \end{enumerate}
  This proves \autoref{thm:d2-h1g4}.
\end{proof}

The remainder of this section is devoted to proving \autoref{thm:P62t:Ext:2-191:d2}.
  First we prove
  \begin{enumerate}
    \item[\eqref{thm:P62t:Ext:2-191:d2}~(1)\quad]
      In the ASS for $\pi_\ast^S(\widetilde{P}^{62})$, $d_2(\bar{\widetilde{e}_{31}h_5}h_7) =\widetilde{q}_\ast(\widehat{D_3}h_7) +\bar{\widetilde{e}_{31}h_0h_4^2}h_7$.
  \end{enumerate}
  Since $d_2(h_7) =h_0h_6^2$ in the ASS for $\pi_\ast^S$ (by \eqref{eq:d2-hi}),
  we have the following differential
  in the ASS for $\pi_\ast^S(\widetilde{P}^{62})$.
  \begin{align}
    \label{eq:d2-e31h5h7:first}
    d_2(\bar{\widetilde{e}_{31}h_5}h_7)
    &= d_2(\bar{\widetilde{e}_{31}h_5})h_7 +\bar{\widetilde{e}_{31}h_5}d_2(h_7)\\
    &= d_2(\bar{\widetilde{e}_{31}h_5})h_7 +\bar{\widetilde{e}_{31}h_5}h_0h_6^2\notag\\
    &= d_2(\bar{\widetilde{e}_{31}h_5})h_7.\notag
  \end{align}
  Here $\bar{\widetilde{e}_{31}h_5}h_0h_6^2 =0$ in $\Ext_\A^{4,192}(\widetilde{P}^{62})$
  because,
  from \eqref{eq:P62t:Ext:new} (1) we have
  $\bar{\widetilde{e}_{31}h_5} =\{\widetilde{e}_{31}\lambda_{31}\}\in\Ext_\A^{1,63}(\widetilde{P}^{62})$,
  and we see that
  $\bar{\widetilde{e}_{31}h_5}h_0h_6^2
  = \bar{\widetilde{e}_{31}h_5}h_6h_0h_6
  = \{\widetilde{e}_{31}\lambda_{31}\lambda_{63}\lambda_0\lambda_{63}\}
  = 0$.
  Thus to compute $d_2(\bar{\widetilde{e}_{31}h_5}h_7)$ we need to
  compute $d_2(\bar{\widetilde{e}_{31}h_5})$
  in the ASS for $\pi_\ast^S(\widetilde{P}^{62})$.
  Consider the following cofibration.
  \begin{align}
    \label{eq:d2-e31h5h7:cofibration}
    \widetilde{P}^{62} &\to{q_{17}'}\widetilde{P}^{62}_{17} \to{\Delta}\Susp P^{16}.
  \end{align}
  Recall
  \begin{enumerate}
    \item[\eqref{eq:P62t:Ext}~(4)\quad]
      $\Ext_\A^{1,63}(\widetilde{P}^{62})
      =\Z/2(\bar{\widetilde{e}_{31}h_5})
      \dsum \Z/2(\bar{e_{47}}h_4)$.
    \item[\eqref{eq:P62-17t:Ext}~(1)\quad]
      $\Ext_\A^{1,63}(\widetilde{P}^{62}_{17})
      =\Z/2(\bar{\widetilde{e}_{31}}h_5) \dsum \Z/2(\bar{e_{47}}h_4)$.\\
      $\Ext_\A^{1,63}(\widetilde{P}^{62}) \to{(q_{17}')_\ast}\Ext_\A^{1,63}(\widetilde{P}^{62}_{17})$ is given by
      $(q_{17}')_\ast(\bar{e_{31}h_5}) =\bar{e_{31}}h_5$ and
      $(q_{17}')_\ast(\bar{e_{47}}h_4) =\bar{e_{47}}h_4$.
    \item[\eqref{eq:P62t:Ext}~(3)\quad]
      $\Ext_\A^{3,64}(\widetilde{P}^{62})
        = \Z/2(\bar{\widetilde{e}_{31}h_4^2}h_0)
        \dsum \Z/2(\widetilde{q}_\ast(\widehat{D_3}))$
    \item[\eqref{eq:P62-17t:Ext}~(3)\quad]
      $\Ext_\A^{3,64}(\widetilde{P}^{62}_{17})
      =\Z/2(\bar{\widetilde{e}_{31}}h_0h_4^2)
      \dsum \Z/2( (\widetilde{q}_{17})_\ast(\bar{e_{23}}h_0h_3h_5))$.\\
      $\Ext_\A^{3,64}(\widetilde{P}^{62}) \to{(q_{17}')_\ast}\Ext_\A^{3,64}(\widetilde{P}^{62}_{17})$ is given by
      $(q_{17}')_\ast(\bar{\widetilde{e}_{31}h_4^2}h_0) =\bar{\widetilde{e}_{31}}h_0h_4^2$, and
      $(q_{17}')_\ast(\widetilde{q}_\ast(\widehat{D_3})) =(\widetilde{q}_{17})_\ast(\bar{e_{23}}h_0h_3h_5)$.
  \end{enumerate}
  Note that
  $\Ext_\A^{1,63}(\widetilde{P}^{62}) \to{(q_{17}')_\ast}\Ext_\A^{1,63}(\widetilde{P}^{62}_{17})$ and
  $\Ext_\A^{3,64}(\widetilde{P}^{62}) \to{(q_{17}')_\ast}\Ext_\A^{3,64}(\widetilde{P}^{62}_{17})$ are both isomorphisms.
  By naturality of the differentials in the Adams spectral sequence,
  we have the following equation in the ASS for $\pi_\ast^S(\widetilde{P}^{62}_{17})$.
  \begin{align}
    \label{eq:d2-e31h5h7:q17:first}
    (q_{17}')_\ast(d_2(\bar{\widetilde{e}_{31}h_5}))
    &=d_2( (q_{17}')_\ast(\bar{\widetilde{e}_{31}h_5}))\\
    &=d_2( \bar{\widetilde{e}_{31}}h_5)
      &\mbox{(by \eqref{eq:P62-17t:Ext}~(1))}\notag\\
    &=d_2( \bar{\widetilde{e}_{31}})h_5 +\bar{\widetilde{e}_{31}}d_2(h_5)\notag\\
    &=d_2( \bar{\widetilde{e}_{31}})h_5 +\bar{\widetilde{e}_{31}}h_0h_4^2
      &\mbox{(by \eqref{eq:P62-17t:Ext}~(3))}\notag
  \end{align}
  Thus to compute $(q_{17}')_\ast(d_2(\bar{\widetilde{e}_{31}h_5}))$
  we need to compute $d_2(\bar{\widetilde{e}_{31}})$
  in the ASS for $\pi_\ast^S(\widetilde{P}^{62}_{17})$.
  Recall
  \begin{enumerate}
    \item[\eqref{eq:P62-17t:Ext}~(2)\quad]
      $\Ext_\A^{2,32}(\widetilde{P}^{62}_{17})=\Z/2( (\widetilde{q}_{17})_\ast(\bar{e_{23}}h_0h_3))$.
  \end{enumerate}
  We claim
  \begin{note}
    \label{eq:d2-e31h5h7:claim}
    $d_2(\bar{\widetilde{e}_{31}}) \ne0$ in the ASS for $\pi_\ast^S(\widetilde{P}^{62}_{17})$.
  \end{note}
  \noindent
  Assuming \eqref{eq:d2-e31h5h7:claim}, from \eqref{eq:P62-17t:Ext}~(2) we must have
  $d_2(\bar{\widetilde{e}_{31}}) =(\widetilde{q}_{17})_\ast(\bar{e_{23}}h_0h_3)$.
  Then we have
  \begin{align}
    \label{eq:d2-e31h5h7:q17}
    (q_{17}')_\ast(d_2(\bar{\widetilde{e}_{31}h_5}))
    &=d_2( \bar{\widetilde{e}_{31}})h_5 +\bar{\widetilde{e}_{31}}h_0h_4^2
      &\mbox{(by \eqref{eq:d2-e31h5h7:q17:first})}\\
    &=(\widetilde{q}_{17})_\ast(\bar{e_{23}}h_0h_3)h_5 +\bar{\widetilde{e}_{31}}h_0h_4^2
      &\mbox{(by $d_2(\bar{\widetilde{e}_{31}}) =(\widetilde{q}_{17})_\ast(\bar{e_{23}}h_0h_3)$)}\notag\\
    &=(\widetilde{q}_{17})_\ast(\bar{e_{23}}h_0h_3h_5) +\bar{\widetilde{e}_{31}}h_0h_4^2 \notag\\
    &=(q_{17}')_\ast(\widetilde{q}_\ast(\widehat{D_3}) +\bar{\widetilde{e}_{31}h_4^2}h_0)
      &\mbox{(by \eqref{eq:P62-17t:Ext}~(3))}\notag
  \end{align}
  And since
  $\Ext_\A^{3,64}(\widetilde{P}^{62}) \to{(q_{17}')_\ast}\Ext_\A^{3,64}(\widetilde{P}^{62}_{17})$ is an isomorphism,
  we have
  \begin{align}
    \label{eq:d2-e31h5h7:e31h5}
    d_2(\bar{\widetilde{e}_{31}h_5})
    &=\widetilde{q}_\ast(\widehat{D_3}) +\bar{\widetilde{e}_{31}h_4^2}h_0.
  \end{align}
  We then have
  the following \eqref{eq:d2-e31h5h7:conclusion}.
  \begin{align}
    \label{eq:d2-e31h5h7:conclusion}
    d_2(\bar{\widetilde{e}_{31}h_5}h_7)
    &=d_2(\bar{\widetilde{e}_{31}h_5})h_7
      &\mbox{(by \eqref{eq:d2-e31h5h7:first})}\\
    &=(\widetilde{q}_\ast(\widehat{D_3}) +\bar{\widetilde{e}_{31}h_4^2}h_0)h_7
      &\mbox{(by \eqref{eq:d2-e31h5h7:e31h5})}\notag\\
    &=\widetilde{q}_\ast(\widehat{D_3}h_7) +\bar{\widetilde{e}_{31}h_4^2}h_0h_7.\notag
  \end{align}
  This proves \eqref{thm:P62t:Ext:2-191:d2}~(1) modulo \eqref{eq:d2-e31h5h7:claim}.
  The proof of \eqref{eq:d2-e31h5h7:claim} is given as follows.

  We use the $\widetilde{P}^{62}_{17} \to{\Delta}\Susp P^{16}$ map
  in \eqref{eq:d2-e31h5h7:cofibration}.
  Recall from \eqref{eq:P62t:diff} that
  in the bigraded differential module
  $ H_\ast(\widetilde{P}^{62})\tensor\Lambda$ we have
  the differential
  $\delta(\widetilde{e}_{31})
  =e_{15}\lambda_{15} \in H_\ast(\widetilde{P}^{62})\tensor\Lambda$
  (put $\lambda_I =1$ in \eqref{eq:P62t:diff}).
  We have $\Ext_\A^{1,32}(\Susp P^{16}) \iso\Ext_\A^{1,31}(P^{16})$, and
  recall
  \begin{enumerate}
    \item[\eqref{eq:P62-k:Ext}~(6)\quad]
      $\Ext_\A^{1,31}(P^{16}) =\Z/2(\bar{e_{15}}h_4)$.
  \end{enumerate}
  So we have
  \begin{align}
    \label{eq:d2-e31h5h7:Delta}
    \Delta_\ast(\bar{\widetilde{e}_{31}} = \{\widetilde{e}_{31}\})
    = \{\delta(\widetilde{e}_{31})\}
    = \{e_{15}\lambda_{15}\}
    =\bar{e_{15}}h_4 \in\Ext_\A^{1,31}(P^{16}) \iso\Ext_\A^{1,32}(\Susp P^{16}).
  \end{align}
  Recall
  \begin{enumerate}
    \item[\eqref{eq:P62-k:Ext}~(7)\quad]
      $\Ext_\A^{2,16}(P^{16}) =\Z/2(\bar{e_7}h_0h_3)$.
    \item[\eqref{eq:P62-k:Ext}~(8)\quad]
      $\Ext_\A^{3,32}(P^{16})
        =\Z/2(\bar{e_{15}}h_0h_3^2)$.
  \end{enumerate}
  Since $\bar{e_{15}} \in\Ext_\A^{0,15}(P^{16})$ and $\bar{e_7} \in\Ext_\A^{0,7}(P^{16})$
  are pullbacks of $\widehat{h}_4 \in\Ext_\A^{0,15}(P)$ and $\widehat{h}_3 \in\Ext_\A^{0,7}(P)$,
  by naturality of the differentials in the Adams spectral sequence,
  we have the following equation in the ASS for $\pi_\ast^S(P^{16})$.
  \begin{align}
    \label{eq:P16:e31:d2:Delta}
    \Delta_\ast(d_2(\bar{\widetilde{e}_{31}}))
    &=d_2(\Delta_\ast(\bar{\widetilde{e}_{31}}))\\
    &=d_2(\bar{e_{15}}h_4)
      &\mbox{(by \eqref{eq:d2-e31h5h7:Delta})}\notag\\
    &=d_2(\bar{e_{15}})h_4 +\bar{e_{15}}d_2(h_4)\notag\\
    &=\bar{e_7}h_0h_3h_4 +\bar{e_{15}}h_0h_3^2
      &\mbox{(by \eqref{eq:P62-k:Ext}~(7) and \eqref{eq:d2-hi})}\notag\\
    &=\bar{e_{15}}h_0h_3^2 \ne0
      &\mbox{(by \eqref{eq:P62-k:Ext}~(8))}\notag.
  \end{align}
  It follows that we must have $d_2(\bar{\widetilde{e}_{31}})\ne0$
  in the ASS for $\pi_\ast^S(\widetilde{P}^{62}_{17})$.
  This proves \eqref{eq:d2-e31h5h7:claim}.

  This completes the proof of \autoref{thm:P62t:Ext:2-191:d2} (1).

  Next we prove
  \begin{enumerate}
    \item[\eqref{thm:P62t:Ext:2-191:d2}~(2)\quad]
      In the ASS for $\pi_\ast^S(\widetilde{P}^{62})$,
      $d_2(\bar{e_{47}}h_4h_7) =\widetilde{q}_\ast(\widehat{D_3}h_7) +\bar{\widetilde{e}_{31}h_0h_4^2}h_7$.
  \end{enumerate}
  Since $d_2(h_4) =h_0h_3^2$ and $d_2(h_7) =h_0h_6^2$  (by \eqref{eq:d2-hi}),
  we have the following differential
  in the ASS for $\pi_\ast^S(\widetilde{P}^{62})$.
  \begin{align}
    \label{eq:d2-e47h4h7:first}
    d_2(\bar{e_{47}}h_4h_7)
    &=d_2(\bar{e_{47}})h_4h_7 +\bar{e_{47}}d_2(h_4)h_7 +\bar{e_{47}}h_4d_2(h_7)\\
    &=d_2(\bar{e_{47}})h_4h_7 +\bar{e_{47}}h_0h_3^2h_7 +\bar{e_{47}}h_4h_0h_6^2 \notag\\
    &=d_2(\bar{e_{47}})h_4h_7. \notag
  \end{align}
  Here $\bar{e_{47}}h_4h_0h_6^2 =\bar{e_{47}}h_0h_4h_6^2 =0$
  in $\Ext_\A^{4,192}(\widetilde{P}^{62})$
  because $h_4h_6^2 =0$, and
  $\bar{e_{47}}h_0h_3^2h_7 =0$ in $\Ext_\A^{4,192}(\widetilde{P}^{62})$ because
  $\bar{e_{47}}h_0h_3^2h_7 =\bar{e_{47}}h_3^2h_0h_7 =\{e_{47}\lambda_7^2\lambda_0\lambda_{127}\}$
  and $\delta(e_{55}\lambda_7\lambda_0\lambda_{127}) =e_{47}\lambda_7^2\lambda_0\lambda_{127}$
  in $ H_\ast(\widetilde{P}^{62})\tensor\Lambda$.
  So to compute $d_2(\bar{e_{47}}h_4h_7)$ we need
  to compute $d_2(\bar{e_{47}})$ in the ASS for $\pi_\ast^S(\widetilde{P}^{62})$.
  Consider the following diagram where the
  horizontal maps are cofibrations.
  \begin{align}
    \label{eq:d2-e47h4h7:cofibration}
    \begin{split}
    \xymatrix{
    P^{32}\ar[r]^i\ar[d]&
    P^{62}\ar[r]^{q_{33}}\ar[d]^{\widetilde{q}}&  P^{62}_{33}\ar[r]^\Delta\ar@{=}[d]& \Susp P^{32}\\
    \widetilde{P}^{32}\ar[r]^{\widetilde{i}}&
      \widetilde{P}^{62}\ar[r]^{q_{33}'}&  {P}^{62}_{33}
    }
    \end{split}
  \end{align}
  We also recall the following.
  \begin{enumerate}
    \item[\eqref{eq:P62-k:Ext}~(4)\quad]
      $\Ext_\A^{0,47}(P^{62}_{33})=\Z/2(\bar{e_{47}})$.
    \item[\eqref{eq:P62-k:Ext}~(5)\quad]
      $\Ext_\A^{2,48}(P^{62}_{33}) =\Z/2(\bar{e_{39}}h_0h_3)$.
    \item[\eqref{eq:P62t:Ext}~(1)\quad]
      $\Ext_\A^{0,47}(\widetilde{P}^{62}) =\Z/2(\bar{e_{47}})$.
    \item[\eqref{eq:P62t:Ext}~(2)\quad]
      $\Ext_\A^{2,48}(\widetilde{P}^{62})
      = \Z/2(\bar{e_{39}h_3}h_0)
      \dsum \Z/2(\widetilde{q}_\ast(\widehat{h}_4h_0h_5))$.
    \item[\eqref{eq:P62t:Ext}~(3)\quad]
      $\bar{e_{39}h_3}h_0h_4 =\widetilde{q}_\ast(\widehat{D_3}) +\bar{\widetilde{e}_{31}h_0h_4^2}$ in $\Ext_\A^{3,64}(\widetilde{P}^{62})$.
    \item[\eqref{eq:P62-17t:Ext}~(4)\quad]
      $\Ext_\A^{2,48}(\widetilde{P}^{62}) \to{(q_{33}')_\ast}\Ext_\A^{2,48}(\widetilde{P}^{62}_{33} =P^{62}_{33})$
      is given by
      $(q_{33}')_\ast(\bar{e_{39}h_3}h_0) =\bar{e_{39}h_3}h_0$,
      $(q_{33}')_\ast(\widetilde{q}_\ast(\widehat{h}_4h_0h_5)) =0$.
    \item[\eqref{eq:P62-17t:Ext}~(5)\quad]
      $\Ext_\A^{0,47}(\widetilde{P}^{62}) \to{(q_{33}')_\ast}\Ext_\A^{0,47}(\widetilde{P}^{62}_{33} =P^{62}_{33})$
      is given by
      $(q_{33}')_\ast(\bar{e_{47}}) =\bar{e_{47}}$.
  \end{enumerate}
  By naturality of the differentials in the Adams spectral sequence,
  we have the following equations in the ASS for $\pi_\ast^S(P^{62}_{33})$.
  \begin{align}
    \label{eq:d2-e47h4h7:q33:first}
    (q_{33}')_\ast(d_2(\bar{e_{47}}))
    &= d_2( (q_{33}')_\ast(\bar{e_{47}}))\\
    &= d_2(\bar{e_{47}})
      &\mbox{(by \eqref{eq:P62-17t:Ext}~(5))}\notag
  \end{align}
  We claim
  \begin{note}
    \label{eq:d2-e47h4h7:claim}
    In the ASS for $\pi_\ast^S(P^{62}_{33})$, $d_2(\bar{e_{47}}) \ne0$.
  \end{note}
  \noindent
  Assuming \eqref{eq:d2-e47h4h7:claim}, by \eqref{eq:P62-k:Ext}~(5)
  we must have $d_2(\bar{e_{47}}) =\bar{e_{39}}h_0h_3$
  in the ASS for $\pi_\ast^S(P^{62}_{33})$.
  From this and \eqref{eq:d2-e47h4h7:q33:first}
  we have $(q_{33}')_\ast(d_2(\bar{e_{47}})) =\bar{e_{39}}h_0h_3$,
  and together with \eqref{eq:P62-17t:Ext}~(4) we see
  $d_2(\bar{e_{47}}) =\bar{e_{39}h_3}h_0 +\varepsilon\widetilde{q}_\ast(\widehat{h}_4h_0h_5)$
  in the ASS for $\pi_\ast^S(\widetilde{P}^{62})$,
  with $\varepsilon=0$ or $1$.
  So we have the following equation in the ASS for $\pi_\ast^S(\widetilde{P}^{62})$.
  \begin{align}
    \label{eq:P62t:e47h4h7:d2:result}
    d_2(\bar{e_{47}}h_4h_7)
    &= d_2(\bar{e_{47}})h_4h_7
      &\mbox{(by \eqref{eq:d2-e47h4h7:first})}\\
    &= (\bar{e_{39}h_3}h_0 +\varepsilon\widetilde{q}_\ast(\widehat{h}_4h_0h_5))h_4h_7
      &\mbox{(by $d_2(\bar{e_{47}}) =\bar{e_{39}h_3}h_0 +\varepsilon\widetilde{q}_\ast(\widehat{h}_4h_0h_5)$)}\notag\\
    &= \bar{e_{39}h_3}h_0h_4h_7\notag\\
    &= (\widetilde{q}_\ast(\widehat{D_3}) +\bar{\widetilde{e}_{31}h_0h_4^2})h_7
      &\mbox{(by \eqref{eq:P62t:Ext}~(3))}\notag\\
    &= \widetilde{q}_\ast(\widehat{D_3}h_7) +\bar{\widetilde{e}_{31}h_0h_4^2}h_7\notag
  \end{align}
  This proves \eqref{thm:P62t:Ext:2-191:d2}~(2) modulo \eqref{eq:d2-e47h4h7:claim}.
  The proof of \eqref{eq:d2-e47h4h7:claim} is given as follows.

  We use the map $P^{62}_{33} \to{\Delta}\Susp P^{32}$ in \eqref{eq:d2-e47h4h7:cofibration}.
  Recall from \eqref{eq:P:diff} that
  $\delta(e_{47}) =e_{31}\lambda_{15}$ in $ H_\ast(P)\tensor\Lambda$.
  We have $\Ext_\A^{1,46}(\Susp P^{32}) \iso\Ext_\A^{1,47}(P^{32})$.
  Recall
  \begin{enumerate}
    \item[\eqref{eq:P62-k:Ext}~(10)\quad]
      $\Ext_\A^{1,47}(P^{32})
        =\Z/2(\bar{e_{31}}h_4)$
  \end{enumerate}
  So we have
  \begin{align}
    \label{eq:d2-e47h4h7:Delta}
    \Delta_\ast(\bar{e_{47}} = \{e_{47}\})
    = \{\delta(e_{47})\}
    = \{e_{31}\lambda_{15}\}
    =\bar{e_{31}}h_4 \in\Ext_\A^{1,47}(P^{32}) \iso\Ext_\A^{1,46}(\Susp P^{32}).
  \end{align}
  Recall
  \begin{enumerate}
    \item[\eqref{eq:P62-k:Ext}~(11)\quad] 
      $\Ext_\A^{2,32}(P^{32}) =\Z/2(\bar{e_{15}}h_0h_4)$.
    \item[\eqref{eq:P62-k:Ext}~(12)\quad]
      $\Ext_\A^{3,48}(P^{32})
      =\Z/2(\bar{e_{31}}h_0h_3^2) 
      \dsum\Z/2(\bar{e_{15}}h_0h_4^2)$.
  \end{enumerate}
  Since $\bar{e_{31}} \in\Ext_\A^{0,31}(P^{32})$
  is a pullback of $\widehat{h}_5 \in\Ext_\A^{0,31}(P)$,
  by naturality of the differentials in the Adams spectral sequence,
  we have the following equations in the ASS for $\pi_\ast^S(P^{32})$.
  \begin{align}
    \label{eq:P32:e47:d2}
    \Delta_\ast(d_2(\bar{e_{47}}))
    &=d_2(\Delta_\ast(\bar{e_{47}}))\\
    &=d_2(\bar{e_{31}}h_4)\notag\\
    &=d_2(\bar{e_{31}})h_4 +\bar{e_{31}}d_2(h_4)\notag\\
    &=\bar{e_{15}}h_0h_4^2 +\bar{e_{31}}h_0h_3^2
      &\mbox{(by \eqref{eq:P62-k:Ext}~(11) and \eqref{eq:d2-hi})}\notag\\
    &\ne0
      &\mbox{(by \eqref{eq:P62-k:Ext}~(12))}\notag
  \end{align}
  It follows that $d_2(\bar{e_{47}}) \ne0$ in the ASS for $\pi_\ast^S(P^{62}_{33})$.
  This proves \eqref{eq:d2-e47h4h7:claim}.

  This completes the proof of \autoref{thm:P62t:Ext:2-191:d2} (2).

\newpage
\section{Modification of the proof without $\theta_6$}\label{se:t6-not-exist}
In this section we discuss how to modify our proof of \autoref{thm:main}
without using the knowledge of $\theta_6$.
Recall that the Kervaire invariant element
$\theta_6 \in\pi_{126}^S$ is a stable homotopy element
detected by $h_6^2 \in\Ext_\A^{2,128}$ in the ASS for $\pi_\ast^S$.
In \autoref{se:intro} through \autoref{se:pi}
we have assumed the existence of a $\theta_6 \in\pi_{126}^S$ when
proving the following theorems.
\begin{note}
  \label{eq:t6:used}
  \autoref{thm:intro:boundary},
  \autoref{thm:t5t6-detected},
  \autoref{thm:t5ht6-detected},
  \autoref{thm:t5ht6:AF-ge-5},
  \autoref{thm:qt5t6:AF-ge-5}, and
  \autoref{thm:S0:d3h63-ne-h52D31}.
\end{note}
\noindent
Let $X =S^{124}\U_\eta e^{126}$.
We will show later in \autoref{thm:t6bar} that
there is a stable map $X \to{\bar{\theta}_6}S^0$
that is detected by $h_6^2$ on the top cell $e^{126}$.
Then we will show that the methods that
we proved the results in \eqref{eq:t6:used} assuming the existence of a $\theta_6$
still works when $S^{126} \to{\theta_6}S^0$ is replaced by $X \to{\bar{\theta}_6}S^0$.

To define the stable map $X \to{\bar{\theta}_6}S^0$
we need to use the
quadratic construction $D_2(Y)$ of a stable complex $Y$ with base point
(\cite{milgram_unstable_1974,barratt_kervaire_1983}).
These are recalled as follows.
Given a stable complex with base point $(Y,\ast)$,
for each $n\ge0$ denote
the points in $S^n\times(Y\Smash Y)$ by $(\lambda,x\Smash y)$,
for $\lambda \in S^n$ and $x,y\in Y$.
Define an involution $S^n\times(Y\Smash Y) \to{T}S^n\times(Y\Smash Y)$ by
\begin{align*}
  T(\lambda,x\Smash y) &=(-\lambda,y\Smash x),
  \quad\mbox{where $-\lambda \in S^n$ is the antipodal point of $\lambda$.}
\end{align*}
Then the group $\{1,T\}$ defines a $\Z_2$-action on $S^n\times(Y\Smash Y)$.
Let $\bar{D_2^n}(Y) =S^n\times_{\Z_2}(Y\Smash Y) =S^n\times(Y\Smash Y)/\{1,T\}$.
Let $D_2^n(Y) =\frac{\bar{D_2^n}(Y)}{S^n\times_{\Z_2}\{\ast\}}$
with base point $\ast =S^n\times_{\Z_2}\{\ast\}$.
Denote the points in $D_2^n(Y)$ by $[\lambda,x\Smash y]$,
where $\lambda\in S^n$ and $x,y\in Y$.
It is not difficult to see that $D_2^n(Y) \subset D_2^{n+1}(Y)$ for each $n\ge0$.
\begin{note}
  \label{eq:D2Y}
  The quadratic construction of $Y$ is the complex $D_2(Y) =\U_{n\ge0} D_2^n(Y)$.
\end{note}
\begin{note}
  \label{eq:D2f}
  Given a based stable map $Y \to{f}X$ where $Y,X$ are stable complexes with base points,
  $D_2(Y) \to{D_2(f)}D_2(X)$ is the base point preserving map defined by
  \begin{align*}
    D_2(f)([\lambda,x\Smash y]) &=[\lambda,f(x)\Smash f(y)],
    \quad\mbox{where $\lambda\in S^n$, $x,y\in Y$.}
  \end{align*}
\end{note}
\noindent
Given a stable map $Y \to{f}S^0$, consider the induced map $D_2(Y) \to{D_2(f)}D_2(S^0)$.
It is not difficult to see that
$D_2(S^0) =P\U\{\ast\}$ where $P$ is the infinite real projective space.
Let $P\U\{\ast\} \to{g}S^0=\{\ast,1\}$ be the map defined by $g(\ast) =\ast$, $g(P) =1$.
\begin{note}
  \label{eq:gamma-f}
  Let $D_2(Y) \to{\gamma(f)}S^0$ be the following composite.
  \begin{align*}
    \gamma(f): D_2(Y) &\to{D_2(f)}D_2(S^0) \to{g}S^0.
  \end{align*}
\end{note}

To describe the stable map $X \to{\bar{\theta}_6}S^0$, we need to
recall from \cite{barratt_kervaire_1983} the following \autoref{thm:Sq2n4} and \autoref{thm:phi-nonzero}.
To state these results,
suppose there is a Kervaire invariant element $\theta_k \in\pi_n^S$ with $2\theta_k =0$, where
$k\ge1$ and $n=2^{k+1}-2$.
Let $Y' =S^n\U_{2\iota}e^{n+1}$,
and let $Y' =S^n\U_{2\iota}e^{n+1} \to{f}S^0$ be an extension of $S^n \to{\theta_k}S^0$.
From \cite{barratt_kervaire_1983},
$D_2(Y')$ has a cell structure as shown in the following diagram.
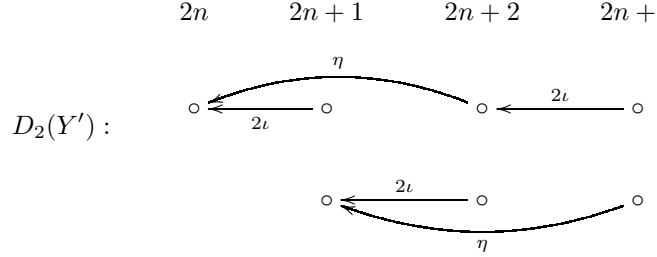
\begin{align*}
  D_2(Y'):\qquad
  \begin{split}
  \xymatrix{
  2n& 2n+1& 2n+2& 2n+3\\
  \circ&
    \circ\ar[l]^{2\iota}&
    \circ\ar@/_1pc/[ll]_\eta&
    \circ\ar[l]_{2\iota}\\
  & \circ&
    \circ\ar[l]_{2\iota}&
    \circ\ar@/^1pc/[ll]^\eta
  }
  \end{split}
\end{align*}
\begin{note}
  \label{eq:Ch}
  Let $X' =S^{2n} \U_\eta e^{2n+2} \U_{2\iota} e^{2n+3}\subset D_2(Y')$ and
  consider the map $X' \to{h =\gamma(f)|_{X'}}S^0$, where
  $D_2(Y') \to{\gamma(f)}S^0$ is as in \eqref{eq:gamma-f}.
  From the definition of $D_2(Y')\to{D_2(f)}D_2(S^0)$ and \eqref{eq:gamma-f}
  we see that the restriction of 
  $X' \to{h}S^0$
  to $S^{2n}$ is $S^{2n}\to{\theta_k^2}S^0$.
\end{note}
\noindent
Consider the mapping cone $C_h$. In \cite{barratt_kervaire_1983} it is shown that
\begin{thm}
  \label{thm:Sq2n4}
  $H^0(C_h) \to{Sq^{2n+4}}H^{2n+4}(C_h)$ is non-zero.
\end{thm}
Here $n=2^{k+1}-2$ so $Sq^{2n+4}=Sq^{2^{k+2}}$, $k\ge1$.
By Adams's decomposition of $Sq^{2^l}$ (\cite{adams_non-existence_1960})
we see the following.
\begin{coro}
  \label{thm:phi-nonzero}
  $H^0(C_h) \to{\Phi_{k+1,k+1}}H^{2n+3}(C_h)$ is non-zero.
\end{coro}

To describe the stable map $X \to{\bar{\theta}_6}S^0$
we will only use \autoref{thm:phi-nonzero} for
$k=5$, $n=62$.
(Recall that there is a $\theta_5 \in\pi_{62}^S$ with $2\theta_5 =0$.)
So $X' =S^{124} \U_\eta e^{126} \U_{2\iota} e^{127}$.
$X \to{\bar{\theta}_6}S^0$ is defined in the following proposition.
\begin{prop}
  \label{thm:t6bar}
  Let $X =S^{124}\U_\eta e^{126}\subset X'$, and
  let $X \to{\bar{\theta}_6 =h|_X}S^0$.
  Then $X \to{\bar{\theta}_6}S^0$ has the following properties.
  \begin{enumerate}
    \item $X \to{\bar{\theta}_6}S^0$ is detected on the top cell by $h_6^2$.
    \item
      Let $S^{126} \to{\bar{2\iota}}X =S^{124} \U_\eta e^{126}$ be a coextension
      of $S^{126} \to{2\iota}S^{126}$. Then the composite
      $S^{126} \to{\bar{2\iota}}X \to{\bar{\theta}_6}S^0$ is zero.
  \end{enumerate}
\end{prop}
It is not difficult to see that
\eqref{thm:t6bar}~(1) follows from \autoref{thm:phi-nonzero}.
It is also easy to see that we have \eqref{thm:t6bar}~(2) because
$X' \to{h}S^0$ is an extension of $X \to{\bar{\theta}_6}S^0$.

We proceed to describe the modifications for the theorems in \eqref{eq:t6:used},
starting from the following
\autoref{thm:t5t6-detected},
\autoref{thm:t5ht6-detected},
\autoref{thm:t5ht6:AF-ge-5}, and
\autoref{thm:qt5t6:AF-ge-5}.
\theoremstyle{nonumberplain}
\begin{thm:t5t6-detected*}
  $\theta_5\theta_6$ in $\pi_{188}^S$
  is detected by $h_0^2g_4 +\varepsilon h_5^2D_3(1)\ne0$
  in $\Ext_\A^{6,194}$,
  where $\varepsilon=0\mbox{ or }1$.
\end{thm:t5t6-detected*}
\begin{thm:t5ht6-detected*}
  Let $\widehat{\theta}_5\in\pi_{62}^S(P^{62})$ be as in \eqref{eq:theta_n-hat}.
  In $\pi_{188}^S(P^{62})$,
  $\widehat{\theta}_5\theta_6\ne0$ and is detected in the ASS for $\pi_\ast^S$
  by a class
  \begin{align*}
    \tau
    =\widehat{g}_4h_0^2,
    +\varepsilon\widehat{h}_5h_5D_3(1)
    +\bar{\tau}\ne0
    \quad\mbox{in $\Ext_\A^{5,193}(P^{62})$}
  \end{align*}
  where
  $\varepsilon=0$ or $1$.
  Here
  $\bar{\tau}\in\Ext_\A^{5,193}(P^{62})$ is a class such that
  $t_\ast(\bar{\tau})=0$ in $\Ext_\A^{6,194}$.
  Moreover,
  $t_\ast(\widehat{g}_4h_0^2)=h_0^2g_4$, and
  $t_\ast(\widehat{h}_5h_5D_3(1))= h_5^2D_3(1)$.
\end{thm:t5ht6-detected*}
\begin{thm:t5ht6:AF-ge-5*}
  For the composite $S^{188} \to{\theta_6}S^{62} \to{\widehat{\theta}_5}P^{62}$,
  $AF(\widehat{\theta}_5\theta_6)\ge5$.
\end{thm:t5ht6:AF-ge-5*}
\newtheorem{thm:qt5t6:AF-ge-5}{\autoref{thm:qt5t6:AF-ge-5}}
\begin{thm:qt5t6:AF-ge-5*}
  The composite
  $S^{188} \to{\theta_6}S^{62} \to{\widehat{\theta}_5} P^{62}\to{\widetilde{q}} \widetilde{P}^{62}$
  has $AF\ge5$
  in the ASS for $\pi_\ast^S(\widetilde{P}^{62})$.
\end{thm:qt5t6:AF-ge-5*}
Denote $\Susp^{62}(X =S^{124}\U_\eta e^{126}) \iso S^{186}\U_\eta e^{188}$ still by $X$,
and $X =S^{186}\U_\eta e^{188} \to{\Susp^{62}\bar{\theta}_6}S^{62}$ still by $\theta_6$.
In each of \autoref{thm:t5t6-detected},
\autoref{thm:t5ht6-detected},
\autoref{thm:t5ht6:AF-ge-5}, and
\autoref{thm:qt5t6:AF-ge-5},
replace $S^{188}$ with $X$, and $S^{188} \to{\theta_6}S^{62}$ with $X \to{\bar{\theta}_6}S^{62}$.
To see that the arguments in the original proofs still work after this modification,
consider the cofibration $S^{186} \to{i_1}X \to{q_1}S^{188}$.
From \eqref{eq:Ch} we see that the composite
$S^{186} \to{i_1}X \to{\bar{\theta}_6}S^{62}$ is $S^{186} \to{\theta_5^2}S^{62}$.
Since $\theta_5$ is detected by $h_5^2$ in the ASS for $\pi_\ast^S$ and since $h_5^4 =0$ in $\Ext_\A^{\ast,\ast}$,
we see $AF(\theta_5^2 =\bar{\theta}_6i_1)\ge5$.
If $AF(\theta_5^2) =5$ then $\theta_5^2$ is detected by a cohomology class $\alpha \in\Ext_\A^{5,129}$.
From \eqref{eq:Ext:known} we deduce
\begin{note}
  $\Ext_\A^{5,129} =0$.
\end{note}
\noindent
Therefore we have
\begin{note}
  \label{eq:t6i:AF-ge-6}
  $AF(\theta_5^2 =\bar{\theta}_6i_1)\ge6$.
\end{note}
\noindent
It follows that the composite
$S^{186} \to{i_1}X \to{\bar{\theta}_6}S^{62} \to{\widehat{\theta}_5}P^{62}$
has $AF \ge7$.
Then it is not difficult to see that all the arguments for proving
\autoref{thm:t5t6-detected},
\autoref{thm:t5ht6-detected},
\autoref{thm:t5ht6:AF-ge-5}, and
\autoref{thm:qt5t6:AF-ge-5}
carry through
if we replace $S^{188} \to{\theta_6}S^{62}$ by
$X \to{\bar{\theta}_6}S^{62}$.

We proceed to describe the modification for the following
\newtheorem{thm:S0:d3h63-ne-h52D31}{\autoref{thm:S0:d3h63-ne-h52D31}}
\begin{thm:S0:d3h63-ne-h52D31*}
  In the ASS for $\pi_\ast^S$,
  $d_3(h_6^3) \ne h_5^2D_3(1)$.
\end{thm:S0:d3h63-ne-h52D31*}
In the proof of \autoref{thm:S0:d3h63-ne-h52D31} we have assumed
a Kervaire invariant element $S^{189} \to{\theta_6}S^{63}$ in the following content.
\begin{enumerate}
  \item[\eqref{eq:d2h63:claim}~(2)] The composite
    $S^{189} \to{\theta_6}S^{63} \to{\bar{\theta}_5}S^0\U_{2\iota}e^1$
    is detected by $(i_1)_\ast(h_6^3)$ in the ASS for $\pi_\ast^S(S^0\U_{2\iota}e^1)$.
\end{enumerate}
We replace
$S^{189}$ and $S^{189} \to{\theta_6}S^{63}$
by $\Susp X=S^{187}\U_\eta e^{189}$ and $\Susp X \to{\Susp\bar{\theta}_6}S^{63}$ respectively
in the proof of \eqref{eq:d2h63:claim}~(2).
From \eqref{eq:t6i:AF-ge-6} it is not difficult to see that the modified proof of \eqref{eq:d2h63:claim} carries through.
This completes the modification for \autoref{thm:S0:d3h63-ne-h52D31}.

Finally we describe the modification for \autoref{thm:intro:boundary}.
\begin{thm:intro:boundary*}
  $h_0^3g_4$ is a boundary in the ASS for spheres.
\end{thm:intro:boundary*}
We have
assumed a Kervaire invariant element $\theta_6$
in the proof of \autoref{thm:intro:boundary} when
using \autoref{thm:t5t6-detected} and the property $2\theta_5\theta_6 =0$.
We replace \autoref{thm:t5t6-detected} with the modified version,
and replace each $S^{188} \to{\theta_6}S^{62}$ by $X \to{\bar{\theta}_6}S^{62}$.
We also have $2\theta_5\bar{\theta}_6 =0$ because $2\theta_5 =0$.
Then the modified proof for \autoref{thm:intro:boundary} carries through.

The completes the modifications to the proof of \autoref{thm:main} with
$S^{126} \to{\theta_6}S^0$ replaced by $X \to{\bar{\theta}_6}S^0$.


\end{document}